# FINITE NEUTROSOPHIC COMPLEX NUMBERS

W. B. Vasantha Kandasamy
Florentin Smarandache

2011

# CONTENTS

# PREFACE

In this book for the first time the authors introduce the notion of real neutrosophic complex numbers. Further the new notion of finite complex modulo integers is defined. For every $C(Z_n)$ the complex modulo integer $i_F$ is such that $i_F^2 = n - 1$. Several algebraic structures on $C(Z_n)$ are introduced and studied.

Further the notion of complex neutrosophic modulo integers is introduced. Vector spaces and linear algebras are constructed using these neutrosophic complex modulo integers.

This book is organized into 5 chapters. The first chapter introduces real neutrosophic complex numbers. Chapter two introduces the notion of finite complex numbers; algebraic structures like groups, rings etc are defined using them. Matrices and polynomials are constructed using these finite complex numbers.

Chapter three introduces the notion of neutrosophic complex modulo integers. Algebraic structures using neutrosophic complex modulo integers are built and around 90 examples are given. Some probable applications are suggested in chapter four and chapter five suggests around 160 problems some of which are at research level.

We thank Dr. K.Kandasamy for proof reading and being extremely supportive.

W.B.VASANTHA KANDASAMY
FLORENTIN SMARANDACHE

**Chapter One**

# REAL NEUTROSOPHIC COMPLEX NUMBERS

In this chapter we for the first time we define the notion of integer neutrosophic complex numbers, rational neutrosophic complex numbers and real neutrosophic complex numbers and derive interesting properties related with them.

Throughout this chapter Z denotes the set of integers, Q the rationals and R the reals. I denotes the indeterminacy and $I^2 = I$. Further i is the complex number and $i^2 = -1$ or $i = \sqrt{-1}$. Also $\langle Z \cup I \rangle = \{a + bI \mid a, b \in Z\} = N(Z)$ and $ZI = \{aI \mid a \in Z\}$. Similarly $\langle Q \cup I \rangle = \{a + bI \mid a, b \in Q\} = N(Q)$ and $QI = \{aI \mid a \in Q\}$. $\langle R \cup I \rangle = \{a + bI \mid a, b \in R\} = N(R)$ and $RI = \{aI \mid a \in R\}$. For more about neutrosophy and the neutrosophic or indeterminate I refer [9-11, 13-4].

Let $C(\langle Z \cup I \rangle) = \{a + bI + ci + dIi \mid a, b, c, d \in Z\}$ denote the integer complex neutrosophic numbers or integer neutrosophic complex numbers. If in C $\langle Z \cup I \rangle = a = c = d = 0$ then we get pure neutrosophic numbers $ZI = \{aI \mid a \in Z\}$. If $c = d = 0$ we get the neutrosophic integers $\{a + bI \mid a, b \in Z\} = N(Z)$. If $b = d = 0$ then we get $\{a + ci\}$ the collection of complex integers J. Likewise $P = \{dIi \mid d \in Z\}$ give the collection of pure neutrosophic complex integers. However dI that is pure neutrosophic collection ZI is a subset of P. Some of the subcollection will have a nice algebraic structure. Thus neutrosophic complex integer is a 4-tuple {a +bI + ci + dIi | a, b, c, d ∈ Z}. We give operations on them. Let x = a + bI + ci + dIi and y = m + nI + si + tIi be in C(⟨Z ∪ I⟩). Now x + y = (a + bI + ci + dIi) + (m + nI + si + tIi) = (a + m) + (b + n)I + (c + s)i + (d+t)Ii.

(We can denote dIi by idI or iId or Iid or diI) we see x + y is again in C(⟨Z ∪ I⟩). We see $0_i^I$ = 0 + 0I + 0i + 0Ii acts as the additive identity.

Thus $0_i^I$ + x = x + $0_i^I$ = x for every x ∈ C(⟨Z ∪ I⟩).

In view of this we have the following theorem.

**THEOREM 1.1**: *C(⟨Z ∪ I⟩) = {a + bI+ ci + idI | a, b, c, d ∈ Z} is integer complex neutrosophic group under addition.*

Proof is direct and hence left as an exercise to the reader.

Let x = a + bI + ci + idI and y = m + nI + ti + isI be in C(⟨Z ∪ I⟩). To find the product xy = (a + bI + ci + idI) (m + nI + ti + isI) = am + mbI + mci + imdI + anI + bnI + nciI + indI + ati + ibtI + cti$^2$ + i$^2$tdI + iasI + ibis + i$^2$csI + i$^2$dsI (using the fact I$^2$ = I and i$^2$ = –1).

= am + mbI + mci + imdI + anI + bnI + incI + indI + iat + ibtI – ct – tdI + iasI + ibis – csI – dsI

= (am – ct) + (mb + an + bn – td – cs – ds)I + i (mc + at) + i(md + nc + nd + bt + as + bs)I.

Clearly xy ∈ C (⟨Z ∪ I⟩).

Now $1_i^I$ = 1 + 0I + 0i + i0I acts as the multiplicative identity. For x $1_i^I$ = $1_i^I$ x = x for every x ∈ C (⟨Z ∪ I⟩). Thus C (⟨Z ∪ I⟩) is a monoid under multiplication. No element in C (⟨Z ∪ I⟩) has inverse with respect to multiplication.

Hence without loss of generality we can denote $1_i^I$ by 1 and $0_i^I$ by 0.

**THEOREM 1.2**: *$C(\langle Z \cup I \rangle)$ is a integer complex neutrosophic monoid or integer neutrosophic complex monoid commutative monoid under multiplication.*

Proof is simple and hence left as an exercise to the reader.

Also it is easily verified that for x, y, z ∈ C ($\langle Z \cup I \rangle$)

x.(y+z) = x.y+x.z and (x+y)z = x.z + y.z. Thus product distributes over the addition.

In view of theorems 1.2 and 1.1 we have the following theorem.

**THEOREM 1.3:** *Let $S = (C (\langle Z \cup I \rangle), +, \times); = \{a + bI + ci + idI$ where $a, b, c, d \in Z\}$; under addition + and multiplication $\times$ is a integer neutrosophic complex commutative ring with unit of infinite order.*

Thus $S = (C (\langle Z \cup I \rangle), +, \times)$ is a ring. S has subrings and ideals. Further C ($\langle Z \cup I \rangle$) has subrings which are not ideals. This is evident from the following theorem the proof of which is left to the reader.

**THEOREM 1.4:** *Let C ($\langle Z \cup I \rangle$) be the integer complex neutrosophic ring.*
i) *$nZI \subseteq C(\langle Z \cup I \rangle)$ is a integer neutrosophic subring of $C(\langle Z \cup I \rangle)$ and is not an ideal of $C(\langle Z \cup I \rangle)$ (n = 1, 2, …)*
ii) *$nZ \subseteq C (\langle Z \cup I \rangle)$ is an integer subring of $C (\langle Z \cup I \rangle)$ which is not an ideal of $C (\langle Z \cup I \rangle)$, (n = 1, 2, …)*
iii) *Let $C (Z) = \{a + ib \mid a, b \in Z\} \subseteq C(\langle Z \cup I \rangle)$, C(Z) is again a complex integer subring which is not an ideal of $C(\langle Z \cup I \rangle)$.*
iv) *Let $S = \{a + bI + ic + idI \mid a, b, c, d \in nZ; 2 < n < \infty\} \subseteq C (\langle Z \cup I \rangle)$.*

*S is a integer complex neutrosophic subring and S is also an ideal of C ($\langle Z \cup I \rangle$).*

The proof of all these results is simple and hence is left as an exercise to the reader.

We can define ideals and also quotient rings as in case of integers Z.

Consider J = {a + bi + cI + idI | a, b, c, d ∈ 2Z} ⊆ C (⟨Z ∪ I⟩) be the ideal of C (⟨Z ∪ I⟩).

Consider $\frac{C(\langle Z \cup I \rangle)}{J}$ = {J, 1 + J, i + J, iI +J, I + J, 1 + i + J, 1+I+J, 1+iI+J, 1+i+I+J, 1+i+iI+J, 1+I+iI+J, i+I+iI+J, 1+i+I+iI+J, i+I+J, i+iI+J, I+iI+J}. Clearly order of P = $\frac{C(\langle Z \cup I \rangle)}{J}$ is $2^4$.

We see P is not an integral domain P has zero divisors. Likewise if we consider S = {a+bi+cI+idI | a, b, c, d, ∈ nZ} ⊆ C (⟨Z ∪ I⟩) S is an ideal of C (⟨Z ∪ I⟩).

Clearly $\frac{C(\langle Z \cup I \rangle)}{S}$ is again not an integral domain and the number of elements in $\frac{C(\langle Z \cup I \rangle)}{S}$ is $n^4$.

We can define different types of integer neutrosophic complex rings using C(⟨Z ∪ I⟩).

**DEFINITION 1.1:** *Let C (⟨Z ∪ I⟩) be the integer complex neutrosophic ring. Consider S = {($x_1$, …, $x_n$) | $x_i$ ∈ C (⟨Z ∪ I⟩); 1 ≤ i ≤ n}; S is again a integer complex neutrosophic 1 × n row matrix ring. The operation in S is taken component wise and we get a ring. This ring is not an integral domains has zero divisors.*

We will give some examples of them.

***Example 1.1:*** Let S = {($x_1$, $x_2$, $x_3$) | $x_t$ = $a_t$ + $b_t$i + $c_t$I + $id_t$I | $a_t$, $b_t$, $c_t$, $d_t$ ∈ Z; 1 ≤ t ≤ 3} be the integer neutrosophic complex ring. S has zero divisors, subrings and ideals.

***Example 1.2:*** Let V = {($x_1$, $x_2$,…, $x_{12}$) | $x_t$ = $a_t$ + $ib_t$ + $c_t$I + $id_t$I where $a_t$, $b_t$, $c_t$, $d_t$ ∈ Z; 1 ≤ t ≤ 12} be the integer neutrosophic complex 1 × 12 matrix ring. V has zero divisors, no units, no idempotents, subrings and ideals. V is of infinite order.

**DEFINITION 1.2**: *Let $S = \{(a_{ij}) \mid a_{ij} \in C(\langle Z \cup I \rangle); 1 \leq i, j \leq n\}$ be a collection of $n \times n$ complex neutrosophic integer matrices. S is a ring of $n \times n$ integer complex neutrosophic ring of infinite order and is non commutative. S has zero divisors, units, idempotents, subrings and ideals.*

We give examples of them.

*Example 1.3*: Let
$$M = \left\{ \begin{pmatrix} a_1 & a_2 \\ a_3 & a_4 \end{pmatrix} \middle| a_i \in C(\langle Z \cup I \rangle); 1 \leq i \leq 4 \right\}$$
be a $2 \times 2$ complex neutrosophic integer ring. M has subrings which are not ideals.

For we see
$$N = \left\{ \begin{pmatrix} a_1 & 0 \\ a_2 & 0 \end{pmatrix} \middle| a_1, b_1 \in C(\langle Z \cup I \rangle) \right\} \subseteq M$$
is a integer complex neutrosophic subring of M which is only a left ideal of M. Clearly N is not a right ideal for
$$\begin{pmatrix} x & 0 \\ y & 0 \end{pmatrix} \begin{pmatrix} a & b \\ c & d \end{pmatrix} = \begin{pmatrix} xa & xb \\ yc & yd \end{pmatrix} \notin N.$$

However
$$\begin{pmatrix} a & b \\ c & d \end{pmatrix} \begin{pmatrix} x & 0 \\ y & 0 \end{pmatrix} = \begin{pmatrix} ax+by & 0 \\ cx+dy & 0 \end{pmatrix}$$

is in N. Hence N is a left ideal and not a right ideal.
Consider
$$T = \left\{ \begin{pmatrix} x & y \\ 0 & 0 \end{pmatrix} \middle| x, y \in C(\langle Z \cup I \rangle) \right\} \subseteq M$$

is a subring but is only a right ideal as
$$\begin{pmatrix} x & y \\ 0 & 0 \end{pmatrix} \begin{pmatrix} a & b \\ c & d \end{pmatrix} = \begin{pmatrix} xa+yc & xb+yd \\ 0 & 0 \end{pmatrix} \in T$$

but

$$\begin{pmatrix} a & b \\ c & d \end{pmatrix} \begin{pmatrix} x & y \\ 0 & 0 \end{pmatrix} = \begin{pmatrix} ax & ay \\ cx & dy \end{pmatrix}$$

is not in T hence is only a right ideal of M.

*Example 1.3:* Let

$$M = \left\{ \begin{pmatrix} a_1 & a_2 & a_3 & a_4 \\ a_5 & a_6 & a_7 & a_8 \\ a_9 & a_{10} & a_{11} & a_{12} \\ a_{13} & a_{14} & a_{15} & a_{16} \end{pmatrix} \middle| a_i \in C(\langle Z \cup I \rangle); 1 \le i \le 16 \right\}$$

be a integer complex neutrosophic ring of $4 \times 4$ matrices.

M is not commutative. M has zero divisors.

$$I_{4 \times 4} = \begin{pmatrix} 1 & 0 & 0 & 0 \\ 0 & 1 & 0 & 0 \\ 0 & 0 & 1 & 0 \\ 0 & 0 & 0 & 1 \end{pmatrix}$$

in M is such that $I_{4 \times 4}$ is the multiplicative identity in M.

Consider

$$P = \left\{ \begin{pmatrix} b_1 & b_2 & b_3 & b_4 \\ 0 & b_5 & b_6 & b_7 \\ 0 & 0 & b_8 & b_9 \\ 0 & 0 & 0 & b_{10} \end{pmatrix} \middle| b_i \in C(\langle Z \cup I \rangle); 1 \le i \le 10 \right\} \subseteq M$$

is an integer complex neutrosophic subring which is not a left ideal or right ideal of M.

Consider

$$T = \left\{ \begin{pmatrix} a_1 & 0 & 0 & 0 \\ a_2 & a_3 & 0 & 0 \\ a_4 & a_5 & a_6 & 0 \\ a_7 & a_8 & a_9 & a_{10} \end{pmatrix} \middle| a_i \in C(\langle Z \cup I \rangle); 1 \le i \le 10 \right\} \subseteq M$$

is an integer complex neutrosophic subring only and not a left ideal or a right ideal of M.

Now we can proceed onto define polynomial integer complex neutrosophic ring.

**DEFINITION 1.4:** *Let*

$$V = \left\{ \sum_{i=0}^{\infty} a_i x^i \,\middle|\, a_i \in C(\langle Z \cup I \rangle) \right\}$$

*be the collection of all polynomials in the variable x with coefficients from the integer complex neutrosophic integral ring $C(\langle Z \cup I \rangle)$ with the following type of addition and multiplication.*

*If $p(x) = a_0 + a_1 x + \ldots + a_n x^n$ and $q(x) = b_0 + b_1 x + \ldots + b_n x^n$ where $a_i, b_i \in C(\langle Z \cup I \rangle)$; $1 \leq i \leq n$; are in V then*
$p(x) + q(x) = (a_0 + a_1 x + \ldots + a_n x^n) + (b_0 + b_1 x + \ldots + b_n x^n)$
$= (a_0 + b_0) + (a_1 + b_1)x + \ldots + (a_n + b_n)x^n \in V$.
*The $0 = 0 + 0x + \ldots + 0x^n$ is the zero integer complex neutrosophic polynomial in V.*
*Now*
$p(x) \cdot q(x) = a_0 b_0 + (a_0 b_1 + a_1 b_0) x + \ldots + a_n b_n x^{2n}$
*is in $C(\langle Z \cup I \rangle)$. $1 = 1 + 0x + \ldots + 0x^n$ in V is such that $p(x) \cdot 1 = 1 \cdot p(x) = p(x)$.*
*(V, +, .) is defined as the integer complex neutrosophic polynomial ring.*

We just enumerate some of the properties enjoyed by V.
   (i) V is a commutative ring with unit.
   (ii) V is an infinite ring.

We can define irreducible polynomials in V as in case of usual polynomials.
   $p(x) = x^2 - 2 \in V$ we see $p(x)$ is irreducible in V.
   $q(x) = x^3 - 3 \in V$ is also irreducible in V.
   $q(x) = x^{12} - 5$ is also irreducible in V.
   $p(x) = x^2 + 4$ is irreducible in V.

Thus in V we can define reducibility and irreducibility in polynomials in V.

Let $p(x) \in V$ if $p(x) = (x - a_1) \ldots (x - a_t)$ where $a_i \in C(\langle Z \cup I \rangle)$; $1 \leq i \leq t$ then $p(x)$ is reducible linearly in a unique manner except for the order in which $a_i$'s occur.

It is infact not an easy task to define relatively prime or greatest common divisor of two polynomials with coefficients from $C(\langle Z \cup I \rangle)$. For it is still difficult to define g c d of two elements in $C(\langle Z \cup I \rangle)$.

For if $a = 5 + 3I$ and $b = (7 + 5I)$ we can say g c d $(a, b) = 1$, if $a = 3 + i - 4I$ and $b = 4i - 2I$ then also g c d $(a, b) = 1$.

If $a = 3 + 3i + 6I + 9iI$ and $b = 12I + 18i + 24$ then g c d $(a, b) = 3$ and so on.

So it is by looking or working with a, b in $C(\langle Z \cup I \rangle)$ we can find g c d.

Now having seen the problem we can not put any order on $C(\langle Z \cup I \rangle)$. For consider i and I we cannot order them for i is the complex number and I is an indeterminate so no relation can be obtained between them. Likewise $1 + i$ and I and so on.

Concept of reducibility and irreducibility is an easy task but other concepts to be obtained in case of neutrosophic complex integer polynomials is a difficult task.

Thus

$$C(\langle Z \cup I \rangle)[x] = \left\{ \sum_{i=0}^{\infty} a_i x^i \,\middle|\, a_i \in C(\langle Q \cup I \rangle) \right\}$$

is a commutative integral domain. Let $p(x)$ and $q(x) \in C(\langle Z \cup I \rangle)[x]$, we can define degree of $p(x)$ as the highest power of x in $p(x)$ with non zero coefficients from $C(\langle Z \cup I \rangle)$.

So if $\deg(p(x)) = n$ and $\deg q(x) = m$ and if $n < m$ then we can divide $q(x)$ by $p(x)$ and find

$$\frac{q(x)}{p(x)} = r(x) + \frac{s(x)}{p(x)}$$

where $\deg s(x) < n$.

This division also is carried out as in case of usual polynomials but due to the presence of four tuples the division process is not simple.

$Z[x] \subseteq C (\langle Z \cup I \rangle) [x] = S$. Thus $Z[x]$ is only a subring and $Z[x]$ is an integral domain in $C (\langle Z \cup I \rangle) [x] = S$. Likewise $ZI[x] \subseteq C (\langle Z \cup I \rangle) [x]$ is again an integral domain which is a subring.

Consider $P = \{a + bI \mid a, b \in Z\} \subseteq C(\langle Z \cup I \rangle) [x]$, $P$ is again a subring which is not an ideal of $C(\langle Z \cup I \rangle) [x]$.

Also $C(Z) = \{a + ib \mid a, b \in Z\}$ is again a subring of $C(\langle Z \cup I \rangle)[x]$.

Further

$$C(Z)[x] = \left\{ \sum_{i=0}^{\infty} a_i x^i \,\middle|\, a_i \in C(Z) \right\} \subseteq C (\langle Z \cup I \rangle) [x]$$

is only a subring of $C (\langle Z \cup I \rangle) [x]$ which is not an ideal.

Likewise

$$C(ZI)[x] = \left\{ \sum_{i=0}^{\infty} a_i x^i \,\middle|\, a_i \in C(ZI) ; a_i = a + ib\ a, b \in ZI \right\} \subseteq S$$

is only a subring of $S$ and is also an ideal of $S$.

Several such properties enjoyed by $C (\langle Z \cup I \rangle) [x]$ can be derived without any difficulty.

We can also define the notion of prime ideal as in case of $C (\langle Z \cup I \rangle) [x]$.

Now we can also define semigroups using $C (\langle Z \cup I \rangle)$.

Consider

$$\left\{ \begin{bmatrix} a_1 & a_2 \\ \vdots & \vdots \\ a_{n-1} & a_n \end{bmatrix} \,\middle|\, a_i \in C (\langle Z \cup I \rangle); 1 \leq i \leq n \right\} = H;$$

$H$ is a group under addition but multiplication cannot be defined on $H$. So $H$ is not a ring but only a semigroup. Thus any collection of $m \times n$ matrices with entries from $C (\langle Z \cup I \rangle)$ ($m \neq n$) is only an abelian group under addition and is not a ring as multiplication cannot be defined on that collection.

We can replace Z by Q then we get $C(\langle Z \cup I \rangle)$ to be rational complex neutrosophic numbers. Also $C(\langle Z \cup I \rangle) \subseteq C(\langle Q \cup I \rangle)$.

$C(\langle Q \cup I \rangle) = \{a + bi + cI + idI \mid a, b, c, d \in Q\}$ is a ring. Infact $C(\langle Q \cup I \rangle)$ has no zero divisors.

For if we take
$$x = 6 + 2i - 3I + 4iI$$
and
$$y = a + bi + cI + diI$$
in $C(\langle Q \cup I \rangle)$.

$$\begin{aligned} xy &= (6 + 2i - 3I + 4iI)(a + bi + cI + diI) \\ &= 6a + 2ai - 3aI + 4aiI \\ &\quad 6bi - 2b - 3biI - 4bI \\ &\quad 6cI + 2ciI - 3cI + 4ciI \\ &\quad 6diI - 2dI - 3diI - 4dI \\ &= (6a - 2b) + I(2a + 6b) + (-3a - 4b + 6c - 3c - 2d - 4d)I + (4a - 3b + 2c + 4c + 6d - 3d) \\ &= 0 \end{aligned}$$

$6a = 2b \quad 2a + 6b = 0$
$b = 3a \quad a = -3b$ this is possible only when $a = b = 0$
$3c - 6d = 0 \quad c = 2d$
$6c + 3d = 0 \quad d = -2c$
So $c = d = 0$.

Thus $a = b = c = d = 0$. But in general $C(\langle Q \cup I \rangle)$ is not a field.

This field contains subfields like Q, $S = \{a + ib \mid a, b \in Q\}$ $C(\langle Q \cup I \rangle)$ contains also subrings.

We can build algebraic structures using $C(\langle Q \cup I \rangle)$. We call $C(\langle Q \cup I \rangle)$ as the rational complex neutrosophic like field. $C(\langle Q \cup I \rangle)$ is of characteristic zero $C(\langle Q \cup I \rangle)$ is not a prime like field for it has subfields of characteristic zero.

$C(\langle Q \cup I \rangle)[x]$ is defined as the neutrosophic complex rational polynomial;

$$C(\langle Q \cup I \rangle)[x] = \left\{ \sum_{i=0}^{\infty} a_i x^i \,\middle|\, a_i \in C(\langle Q \cup I \rangle) \right\}.$$

$C(\langle Q \cup I \rangle)$ is not a field, it is only an integral domain we can derive the polynomial properties related with rational complex neutrosophic polynomials in $C(\langle Q \cup I \rangle)[x]$.

$C(\langle Q \cup I \rangle)$ [x] can have ideals. Now consider $T = \{(x_1, x_2, \ldots, x_n) \mid x_i \in C(\langle Q \cup I \rangle); 1 \leq i \leq n\}$; T is a rational complex neutrosophic $1 \times n$ matrix. T is only a ring for T contains zero divisors but T has no idempotents; T has ideals for take $P = \{(x_1, x_2, x_3, 0, \ldots, 0) \mid x_i \in C(\langle Q \cup I \rangle); 1 \leq i \leq 3\} \subseteq T$ is a subring as well as an ideal of T.

T has several ideals, T also has subrings which are not ideals. For take $S = \{(x_1, x_2, x_3, \ldots, x_n) \mid x_i \in C(\langle Z \cup I \rangle); 1 \leq i \leq n\} \subseteq T$; S is a subring of T and is not an ideal of T. We can have several subrings of T which are not ideals of T.

Now we can define $M = \{A = (a_{ij}) \mid A$ is a $n \times n$ rational complex neutrosophic matrix with $a_{ij} \in C(\langle Q \cup I \rangle); 1 \leq i, j \leq n\}$ to be the $n \times n$ rational complex neutrosophic matrix ring. M also has zero divisors, units, ideals, and subrings. For consider $N = \{$collection of all upper triangular $n \times n$ matrices with elements from $C(\langle Q \cup I \rangle)\} \subseteq M$; N is a subring of M and N is not an ideal of M. We have $T = \{$all diagonal $n \times n$ matrices with entries from $C(\langle Q \cup I \rangle)\} \subseteq M$; T is an ideal of M.

All usual properties can be derived with appropriate modifications. Now if we replace Q by R we get $C(\langle R \cup I \rangle)$ to be the real complex neutrosophic ring. $C(\langle R \cup I \rangle)$ is not a field called the like field of real complex neutrosophic numbers. $C(\langle R \cup I \rangle)$ is not a prime field. It has subfields and subrings which are not subfields.

All properties enjoyed by $C(\langle Q \cup I \rangle)$ can also be derived for $C(\langle R \cup I \rangle)$. We see $C(\langle R \cup I \rangle) \supset C(\langle Q \cup I \rangle) \supset C(\langle Z \cup I \rangle)$.

We construct polynomial ring with real complex neutrosophic coefficients and $1 \times n$ matrix ring with real complex neutrosophic matrices. Likewise the $n \times n$ real complex neutrosophic matrix ring can also be constructed. The latter two will have zero divisors and units where as the first ring has no zero divisors it is an integral domain.

Now we have seen real complex neutrosophic like field, C($\langle R \cup I \rangle$).

We proceed onto define vector spaces, set vector spaces, group vector spaces of complex neutrosophic numbers.

**DEFINITION 1.4:** *Let V be a additive abelian group of complex neutrosophic numbers. Q be the field. If V is vector space over Q then define V to be a ordinary complex neutrosophic vector space over the field Q.*

We will give examples of them.

*Example 1.5:* Let
$$V = \left\{ \begin{bmatrix} a_1 & a_2 \\ a_3 & a_4 \\ a_5 & a_6 \\ a_7 & a_8 \\ a_9 & a_{10} \end{bmatrix} \middle| a_i \in C(\langle Q \cup I \rangle); 1 \le i \le 10 \right\}$$
be an ordinary complex neutrosophic vector space over Q.

*Example 1.6:* Let
$$M = \left\{ \begin{bmatrix} a_1 & a_2 \\ a_3 & a_4 \end{bmatrix} \middle| a_i \in C(\langle Q \cup I \rangle); 1 \le i \le 4 \right\}$$
be the ordinary complex neutrosophic vector space over Q.

*Example 1.7:* Let
$$P = \left\{ \begin{pmatrix} a_1 & a_2 & a_3 \\ a_4 & a_5 & a_6 \end{pmatrix} \middle| a_i \in C(\langle Q \cup I \rangle); 1 \le i \le 6 \right\}$$
be the ordinary complex neutrosophic vector space over Q.

We can as in case of usual vector spaces define subspaces and basis of P over Q.

However it is pertinent to mention here that we can have other types of vector spaces defined depending on the field we choose.

**DEFINITION 1.5:** *Let V be the complex neutrosophic additive abelian group. Take $F = \{a + bi \mid a, b \in Q; i^2 = -1\}$; if V is a vector space over the complex field F; then we call V to be complex - complex neutrosophic vector space.*

We will give examples of them.

*Example 1.8:* Let

$$V = \left\{ \begin{bmatrix} a_1 & a_2 & a_3 & a_4 \\ a_5 & a_6 & a_7 & a_8 \end{bmatrix} \middle| a_i \in C(\langle Q \cup I \rangle); 1 \le i \le 8 \right\}$$

be a complex - complex neutrosophic vector space over the field $F = \{a + b_i \mid a, b \in Q\}$.

Take

$$W = \left\{ \begin{pmatrix} a_1 & a_2 & a_3 & a_4 \\ 0 & 0 & 0 & 0 \end{pmatrix} \middle| a_i \in C(\langle Q \cup I \rangle); 1 \le i \le 4 \right\}$$

is a complex - complex neutrosophic vector subspace of V over F. Infact V has several such subspaces.

Take

$$W_1 = \left\{ \begin{pmatrix} a_1 & 0 & 0 & a_2 \\ 0 & a_3 & 0 & 0 \end{pmatrix} \middle| a_i \in C(\langle Q \cup I \rangle); 1 \le i \le 3 \right\} \subseteq V,$$

a complex - complex neutrosophic vector subspace of V over the rational complex field F.

Suppose

$$W_2 = \left\{ \begin{pmatrix} 0 & 0 & 0 & 0 \\ a_1 & 0 & a_2 & 0 \end{pmatrix} \middle| a_1, a_2 \in C(\langle Q \cup I \rangle) \right\} \subseteq V$$

is a complex - complex neutrosophic vector subspace of V over F.

Consider

$$W_3 = \left\{ \begin{pmatrix} 0 & a_1 & a_2 & 0 \\ 0 & 0 & 0 & a_3 \end{pmatrix} \middle| a_1, a_2, a_3 \in C(\langle Q \cup I \rangle) \right\} \subseteq V;$$

$W_3$ is a complex - complex neutrosophic vector subspace of V over F.

Clearly
$$V = W_1 \cup W_2 \cup W_3$$
$$= W_1 + W_2 + W_3.$$

Thus V is a direct sum of subspaces.

*Example 1.9:* Let

$$V = \left\{ \begin{bmatrix} a_1 & a_2 & a_3 & a_4 \\ a_5 & a_6 & a_7 & a_8 \\ a_9 & a_{10} & a_{11} & a_{12} \\ a_{13} & a_{14} & a_{15} & a_{16} \\ a_{17} & a_{18} & a_{19} & a_{20} \end{bmatrix} \middle| a_i \in C(\langle Q \cup I \rangle); 1 \leq i \leq 20 \right\}$$

be a complex - complex neutrosophic vector space over the complex field $F = \{a + bi \mid a, b \in Q\}$. V has subspaces and V can be written as a direct sum / union of subspaces of V over F.

Now we can define complex neutrosophic - neutrosophic like vector space or neutrosophic - neutrosophic vector space over the neutrosophic like field $\langle Q \cup I \rangle$ or $\langle R \cup I \rangle$.

**DEFINITION 1.6:** *Let V be an additive abelian group of complex neutrosophic numbers. Let $F = \langle Q \cup I \rangle$ be the neutrosophic like field of rationals. If V is a like vector space over F then we define V to be a neutrosophic - neutrosophic complex like vector space over the field F (complex neutrosophic - neutrosophic vector space over the field F or neutrosophic complex neutrosophic vector space over the field F).*

We will give examples of this situation.

*Example 1.10:* Let

$$V = \left\{ \begin{bmatrix} a_1 & a_2 \\ a_3 & a_4 \\ a_5 & a_6 \\ \vdots & \vdots \\ a_{15} & a_{16} \end{bmatrix} \middle| a_i \in C(\langle Q \cup I \rangle); 1 \le i \le 16 \right\}$$

be a neutrosophic complex neutrosophic like vector space over the neutrosophic like field $F = \langle Q \cup I \rangle$.

*Example 1.11:* Let

$$M = \left\{ \begin{bmatrix} a_1 & a_2 & a_3 & a_4 \\ a_5 & a_6 & a_7 & a_8 \\ a_9 & a_{10} & a_{11} & a_{12} \\ a_{13} & a_{14} & a_{15} & a_{16} \end{bmatrix} \middle| a_i \in C(\langle Q \cup I \rangle); 1 \le i \le 16 \right\}$$

be a neutrosophic - neutrosophic complex like vector space over the neutrosophic like field $\langle Q \cup I \rangle = F$.

Take

$$P_1 = \left\{ \begin{bmatrix} a_1 & a_2 & 0 & 0 \\ a_3 & a_4 & 0 & 0 \\ 0 & 0 & 0 & 0 \\ 0 & 0 & 0 & 0 \end{bmatrix} \middle| a_i \in C(\langle Q \cup I \rangle); 1 \le i \le 4 \right\} \subseteq V,$$

$$P_2 = \left\{ \begin{bmatrix} 0 & 0 & a_1 & a_2 \\ 0 & 0 & a_3 & a_4 \\ 0 & 0 & 0 & 0 \\ 0 & 0 & 0 & 0 \end{bmatrix} \middle| a_i \in C(\langle Q \cup I \rangle); 1 \le i \le 4 \right\} \subseteq V,$$

$$P_3 = \left\{ \begin{bmatrix} 0 & 0 & 0 & 0 \\ 0 & 0 & 0 & 0 \\ a_1 & a_2 & 0 & 0 \\ a_4 & a_5 & 0 & 0 \end{bmatrix} \middle| a_i \in C(\langle Q \cup I \rangle); 1 \le i \le 4 \right\} \subseteq V$$

and

$$P_4 = \left\{ \begin{bmatrix} 0 & 0 & 0 & 0 \\ 0 & 0 & 0 & 0 \\ 0 & 0 & a_1 & a_2 \\ 0 & 0 & a_3 & a_4 \end{bmatrix} \middle| a_i \in C(\langle Q \cup I \rangle); 1 \le i \le 4 \right\} \subseteq V$$

be subspace of V.

Clearly $V = P_1 + P_2 + P_3 + P_4$ and $P_i \cap P_j = (0)$ if $i \ne j$, so V is a direct sum of subspaces of V over $F = \langle Q \cup I \rangle$.

Consider

$$V_1 = \left\{ \begin{bmatrix} a_1 & 0 & 0 & a_2 \\ 0 & 0 & 0 & 0 \\ 0 & a_3 & a_4 & 0 \\ 0 & 0 & a_5 & a_6 \end{bmatrix} \middle| a_i \in C(\langle Q \cup I \rangle); 1 \le i \le 6 \right\} \subseteq V,$$

$$V_2 = \left\{ \begin{bmatrix} a_1 & a_2 & a_3 & 0 \\ a_4 & a_5 & 0 & 0 \\ 0 & 0 & 0 & 0 \\ 0 & 0 & 0 & a_6 \end{bmatrix} \middle| a_i \in C(\langle Q \cup I \rangle); 1 \le i \le 6 \right\} \subseteq V,$$

$$V_3 = \left\{ \begin{bmatrix} a_1 & 0 & 0 & 0 \\ 0 & a_2 & a_3 & a_4 \\ a_5 & 0 & 0 & 0 \\ 0 & 0 & 0 & 0 \end{bmatrix} \middle| a_i \in C(\langle Q \cup I \rangle); 1 \le i \le 6 \right\} \subseteq V,$$

$$V_4 = \left\{ \begin{bmatrix} a_1 & 0 & 0 & 0 \\ 0 & 0 & a_2 & 0 \\ a_3 & a_4 & a_5 & a_6 \\ 0 & 0 & 0 & a_7 \end{bmatrix} \middle| a_i \in C(\langle Q \cup I \rangle); 1 \le i \le 7 \right\} \subseteq V,$$

and

$$V_5 = \left\{ \begin{bmatrix} a_1 & 0 & 0 & a_2 \\ 0 & 0 & a_3 & 0 \\ 0 & 0 & 0 & a_4 \\ a_5 & a_6 & a_7 & 0 \end{bmatrix} \middle| a_i \in C(\langle Q \cup I \rangle); 1 \le i \le 7 \right\} \subseteq V$$

is such that $V = V_1 \cup V_2 \cup V_3 \cup V_4 \cup V_5$ but $V_i \cap V_j \ne (0)$ if $i \ne j$ so V is only a pseudo direct sum of the subspaces of V over F. Now we have defined neutrosophic complex - neutrosophic like vector spaces over the neutrosophic like field.

We now proceed onto define special complex neutrosophic like vector space over the complex neutrosophic like field.

**DEFINITION 1.7:** *Let V be an abelian group under addition of complex neutrosophic numbers. Let $F = C(\langle Q \cup I \rangle)$ be the complex - neutrosophic like field of rationals; if V is a vector space over F then we define V to be a special complex neutrosophic like vector space over the complex neutrosophic rational like field $F = C(\langle Q \cup I \rangle)$.*

We will give examples of them.

*Example 1.12:* Let

$$V = \left\{ \begin{bmatrix} a_1 \\ a_2 \\ \vdots \\ a_{15} \end{bmatrix} \right\}$$

where $a_i \in C(\langle Q \cup I \rangle); 1 \le i \le 15\}$ be the special complex neutrosophic like vector space over the complex neutrosophic like field of rationals $C(\langle Q \cup I \rangle) = F$.

It is easily verified V has subspaces. The dimension of V over F is 15.
For take

$$B = \left\{ \begin{bmatrix} 1 \\ 0 \\ \vdots \\ 0 \end{bmatrix}, \begin{bmatrix} 0 \\ 1 \\ 0 \\ \vdots \\ 0 \end{bmatrix}, \begin{bmatrix} 0 \\ 0 \\ 1 \\ 0 \\ \vdots \\ 0 \end{bmatrix}, \dots, \begin{bmatrix} 0 \\ \vdots \\ 0 \\ 0 \\ 1 \\ 0 \end{bmatrix}, \begin{bmatrix} 0 \\ \vdots \\ 0 \\ 0 \\ 0 \\ 1 \end{bmatrix} \right\} \subseteq V$$

is B is a basis of V over F.

*Example 1.13:* Let

$$V = \left\{ \begin{bmatrix} a_1 & a_2 & a_3 & a_4 \\ a_5 & a_6 & a_7 & a_8 \\ a_9 & a_{10} & a_{11} & a_{12} \end{bmatrix} \middle| a_i \in C(\langle Q \cup I \rangle); 1 \leq i \leq 12 \right\}$$

be a special complex neutrosophic like vector space over the field $F = C(\langle Q \cup I \rangle)$.

Now we can define like subfield subspace of a vector space.

**DEFINITION 1.8:** *Let V be an additive abelian group of complex neutrosophic numbers. V be a special complex neutrosophic vector space over the complex neutrosophic like field*
 $F = C(\langle Q \cup I \rangle) = \{a + bi + cI + idI \mid a, b, c, d \in Q\}$.
 *Let $W \subseteq V$, W also a proper subgroup of V and $K \subseteq F$. K the neutrosophic like field $\langle Q \cup I \rangle \subseteq C(\langle Q \cup I \rangle) = F$.*
 *If W is a vector space over K then we define W to be a neutrosophic subfield complex neutrosophic vector subspace of V over the neutrosophic like subfield K of F.*

We will give examples of this situation.

*Example 1.14:* Let

$$V = \left\{ \begin{bmatrix} a_1 & a_2 & a_3 \\ a_4 & a_5 & a_6 \\ a_7 & a_8 & a_9 \\ a_{10} & a_{11} & a_{12} \\ a_{13} & a_{14} & a_{15} \end{bmatrix} \middle| a_i \in C(\langle Q \cup I \rangle); 1 \le i \le 15 \right\}$$

be a special complex neutrosophic vector space over the complex neutrosophic like field $F = C(\langle Q \cup I \rangle)$.

Nothing is lost if we say neutrosophic field also for we have defined so but an indeterminate or neutrosophic field need not have a real field structure like a neutrosophic group is not a group yet we call it a group.

Take

$$W = \left\{ \begin{bmatrix} a_1 & a_2 & a_3 \\ 0 & 0 & 0 \\ 0 & 0 & 0 \\ 0 & 0 & 0 \\ a_4 & a_5 & a_6 \end{bmatrix} \middle| a_i \in C(\langle Q \cup I \rangle); 1 \le i \le 6 \right\} \subseteq V;$$

take $K = \langle Q \cup I \rangle \subseteq F = C(\langle Q \cup I \rangle)$; W is a neutrosophic special complex neutrosophic vector subspace of V over the neutrosophic like subfield $\langle Q \cup I \rangle$ of $F = C(\langle Q \cup I \rangle)$.

Consider

$$M = \left\{ \begin{bmatrix} a_1 & 0 & 0 \\ 0 & a_2 & 0 \\ 0 & 0 & a_3 \\ a_4 & 0 & 0 \\ 0 & a_5 & 0 \end{bmatrix} \middle| a_i \in C(\langle Q \cup I \rangle); 1 \le i \le 5 \right\} \subseteq V$$

take $K = \{a + ib \mid a, b \in Q\} \subseteq F = C(\langle Q \cup I \rangle)$ a complex subfield of $C(\langle Q \cup I \rangle)$. W is a special complex neutrosophic complex subvector space of V over the rational complex subfield K of F.

Suppose

$$S = \left\{ \begin{bmatrix} a_1 & 0 & a_2 \\ 0 & a_3 & 0 \\ a_4 & 0 & a_5 \\ 0 & a_6 & 0 \\ a_7 & 0 & a_8 \end{bmatrix} \middle| a_i \in C(\langle Q \cup I \rangle); 1 \leq i \leq 8 \right\} \subseteq V.$$

Take $Q \subseteq C(\langle Q \cup I \rangle) = F$ the field of rationals as a subfield of F. S is a ordinary special neutrosophic complex vector subspace of V over the rational subfield Q of $F = C(\langle Q \cup I \rangle)$.
All properties can be derived for vector space over $C(\langle Q \cup I \rangle)$.

It is interesting to see that only these special neutrosophic complex vector spaces has many properties which in general is not true over the usual field Q.

**THEOREM 1.5:** *Let V be a special neutrosophic complex vector space over $C(\langle Q \cup I \rangle)$. V has only 2 subfields over which vector subspaces can be defined.*

When we say this it is evident from the example 1.14, hence left as an exercise to the reader.

**THEOREM 1.6:** *Let V be an ordinary neutrosophic complex vector space over the field Q. V has no subfield vector subspace.*

Proof easily follows from the fact Q is a prime field.

**THEOREM 1.7:** *Let V be a neutrosophic complex neutrosophic vector space over the neutrosophic field $\langle Q \cup I \rangle = F$. V has only one subfield over which vector subspaces can be defined.*

**THEOREM 1.8:** *Let V be a complex neutrosophic complex vector space over the rational complex field $F = \{a + ib \mid a, b \in Q\}$. V has only one subfield over which subvector spaces can be defined.*

Proof follows from the fact $Q \subseteq F$.
We see
$$Q \subseteq \langle Q \cup I \rangle \subseteq C(\langle Q \cup I \rangle),$$
$$Q \subseteq C(Q) = \{a + ib \mid a, b \in Q\} \subseteq C(\langle Q \cup I \rangle).$$

These spaces behave in a very different way which is evident from the following example.

*Example 1.15:* Let

$$V = \left\{ \begin{bmatrix} a_1 & a_2 & a_3 \\ a_4 & a_5 & a_6 \\ a_7 & a_8 & a_9 \end{bmatrix} \middle| a_i \in C(\langle Q \cup I \rangle); 1 \leq i \leq 9 \right\}$$

be a special complex neutrosophic vector space over the complex neutrosophic like field $F = C(\langle Q \cup I \rangle)$.

It is easily verified V is of dimension 9 over $F = C(\langle Q \cup I \rangle)$. However V has special complex neutrosophic subspaces of dimensions 1, 2, 3, 4, 5, 6, 7 and 8.

Now consider

$$W = \left\{ \begin{bmatrix} a_1 & a_2 & 0 \\ 0 & a_4 & 0 \\ 0 & 0 & a_4 \end{bmatrix} \middle| a_i \in C(\langle Q \cup I \rangle); 1 \leq i \leq 4 \right\} \subseteq V.$$

W is a special subfield neutrosophic complex neutrosophic vector subspace of V over the neutrosophic subfield $\langle Q \cup I \rangle$ of F. Clearly dimension of W over $\langle Q \cup I \rangle$ is not finite. Suppose W is considered as a special neutrosophic complex vector subspace of V over $F = C(\langle Q \cup I \rangle)$ then dimension of W over F is four.

W as a special complex neutrosophic complex vector subspace over the rational complex field $K = \{a + ib \mid a, b \in Q\} \subseteq C(\langle Q \cup I \rangle)$ which is also of infinite dimension over K. W as a special neutrosophic complex ordinary vector subspace over the rational field Q is also of infinite dimension over the rational subfield Q of F.

Consider

$$T = \left\{ \begin{bmatrix} a_1 & 0 & a_2 \\ 0 & a_3 & 0 \\ a_4 & 0 & a_5 \end{bmatrix} \middle| a_i \in Q; 1 \le i \le 5 \right\} \subseteq V.$$

T is a special neutrosophic complex ordinary vector subspace of V over Q dimension 5.

Clearly T is not defined over the subfield, $\langle Q \cup I \rangle$ of C(Q) or $C(\langle Q \cup I \rangle)$.

Consider

$$S = \left\{ \begin{bmatrix} 0 & a_1 & 0 \\ a_2 & 0 & a_3 \\ 0 & a_4 & 0 \end{bmatrix} \middle| a_i \in \langle Q \cup I \rangle; 1 \le i \le 4 \right\} \subseteq V$$

be a special neutrosophic-neutrosophic complex vector subspace of V over the neutrosophic rational subfield $\langle Q \cup I \rangle$. S is of dimension four over $\langle Q \cup I \rangle$. S is also a special neutrosophic complex ordinary vector subspace of V over the rational field Q and S is of infinite dimension over Q. Clearly S is not a vector subspace over C(Q) = {a + ib | a, b ∈ Q} or $C(\langle Q \cup I \rangle)$.

Consider

$$P = \left\{ \begin{bmatrix} a_1 & a_2 & a_3 \\ 0 & a_4 & a_5 \\ 0 & 0 & a_6 \end{bmatrix} \middle| a_i \in C(Q) \right.$$

= {a + ib | a, b ∈ Q}, $1 \le i \le 6\} \subseteq V$, P is a special complex neutrosophic complex vector subspace of V over the rational complex field C(Q) of dimension 6. P is also a special complex neutrosophic ordinary vector subspace over the field of rationals Q of infinite dimension.

Clearly P is not a special complex neutrosophic complex vector subspace of V over $\langle Q \cup I \rangle$ or $C(\langle Q \cup I \rangle)$. As it is not defined over the two fields.

Thus we have seen several of the properties about special neutrosophic complex vector spaces defined over C ($\langle Q \cup I \rangle$). We now proceed onto define the linear algebra structures.

**DEFINITION 1.9:** *Let V be a ordinary complex neutrosophic vector space over the rationals Q. If on V a product '.' can be defined and V is compatible with respect to the product '.', then we call V to be a ordinary complex neutrosophic linear algebra over the rational field Q.*

We provide examples of them.

***Example 1.16:*** Let
$$V = \left\{ \begin{bmatrix} a_1 & a_2 \\ a_3 & a_4 \end{bmatrix} \middle| a_i \in C(\langle Q \cup I \rangle); 1 \leq i \leq 4 \right\}$$

be an ordinary complex neutrosophic linear algebra over Q. V is of infinite dimension over Q. V has linear subalgebras. It is interesting to notice V is also a vector space but all ordinary complex neutrosophic vector spaces over Q need not in general be linear algebras.

In view of this we have the following theorem.

**THEOREM 1.9:** *Let V be an ordinary neutrosophic complex linear algebra over the rationals Q, then V is an ordinary neutrosophic complex vector space. If V is an ordinary neutrosophic complex vector space over Q then V in general is not an ordinary neutrosophic complex linear algebra over Q.*

The proof is straight forward hence left as an exercise for the reader. On similar lines we can define neutrosophic - complex neutrosophic linear algebra over $\langle Q \cup I \rangle$, neutrosophic complex - complex linear algebra over $C(Q) = \{a + ib \mid a, b \in Q\}$ and special neutrosophic complex linear algebra over $C(\langle Q \cup I \rangle)$.

We give only examples of them as the definition is a matter of routine.

*Example 1.17:* Let
$$V = \left\{ \sum_{i=0}^{\infty} a_i x^i \,\middle|\, a_i \in C(\langle Q \cup I \rangle) \right\}$$
be a neutrosophic complex neutrosophic linear algebra over the neutrosophic field $F = \langle Q \cup I \rangle$.

*Example 1.18:* Let $M = \{$All $5 \times 5$ upper triangular matrices with entries from $C(\langle Q \cup I \rangle)\}$ be a neutrosophic complex neutrosophic linear algebra over the neutrosophic field $\langle Q \cup I \rangle = F$

*Example 1.19:* Let $P = \{$all $10 \times 10$ matrices with complex neutrosophic entries from $C(\langle Q \cup I \rangle)\}$ be a neutrosophic complex neutrosophic linear algebra over the neutrosophic field $F = \langle Q \cup I \rangle$. $M = \{$all $10 \times 10$ upper triangular matrices with complex neutrosophic entries from $C(\langle Q \cup I \rangle)\} \subseteq P$ is the neutrosophic complex neutrosophic linear subalgebra of P over the field $F = \langle Q \cup I \rangle$.

This space is of infinite dimension over $\langle Q \cup I \rangle = F$.

*Example 1.20:* Let
$$P = \left\{ \sum_{i=0}^{\infty} a_i x^i \,\middle|\, a_i \in C(\langle Q \cup I \rangle) \right\}$$
be a complex neutrosophic complex linear algebra over the complex rational field $C(Q) = \{a + bi \mid a, b \in Q\}$.

*Example 1.21:* Let
$$V = \left\{ \begin{pmatrix} a_1 & a_2 & a_3 & a_4 \\ 0 & a_5 & a_6 & a_7 \\ 0 & 0 & a_8 & a_9 \\ 0 & 0 & 0 & a_{10} \end{pmatrix} \,\middle|\, a_i \in C(\langle Z \cup I \rangle); 1 \le i \le 10 \right\}$$
be a complex neutrosophic complex linear algebra over the rational complex field $C(Q) = \{a + bi \mid a, b \in Q\}$.

***Example 1.22:*** Let M = {all 8 × 8 lower triangular complex neutrosophic matrices with entries from C (⟨Q ∪ I⟩)} be a complex neutrosophic complex linear algebra over the complex field C(Q) = {a + bi | a, b ∈ Q}.

***Example 1.23:*** Let

$$T = \left\{ \begin{bmatrix} a_1 & a_2 & a_3 \\ a_4 & a_5 & a_6 \\ a_7 & a_8 & a_9 \end{bmatrix} \middle| a_i \in C(\langle Q \cup I \rangle); 1 \le i \le 9 \right\}$$

be an ordinary neutrosophic complex linear algebra over the rational field Q.

***Example 1.24:*** Let A = {all 10 × 10 upper triangular matrices with complex neutrosophic entries from the complex neutrosophic field} be an ordinary neutrosophic complex linear algebra over the field of rationals Q.

***Example 1.25:*** Let

$$B = \left\{ \sum_{i=0}^{\infty} a_i x^i \middle| a_i \in C(\langle Q \cup I \rangle) \right\}$$

be a special neutrosophic complex linear algebra over the neutrosophic complex field C (⟨Q ∪ I⟩).

***Example 1.26:*** Let

$$C = \left\{ \begin{bmatrix} a_1 & a_2 & a_3 \\ 0 & a_4 & a_5 \\ 0 & 0 & a_6 \end{bmatrix} \middle| a_i \in C(\langle Q \cup I \rangle); 1 \le i \le 6 \right\}$$

be a special neutrosophic complex linear algebra over the neutrosophic complex field C (⟨Q ∪ I⟩).

Clearly C is of dimension 6 over C (⟨Q ∪ I⟩). We can derive almost all properties of vector spaces in case of linear algebras with simple appropriate modifications.

We can also define the notion of linear transformation and linear operators.

We can define linear transformation T of neutrosophic complex vector spaces from V to W only if V and W are defined over the same field. Also T (I) = I is a basic criteria for the transformation for the indeterminacy cannot be mapped on to any other element.

Likewise if W is replaced by V this linear transformation T of V to V becomes the linear operator. All properties associated with linear transformation and linear operators of neutrosophic complex vector spaces can be easily derived in case of these operators.

We can now define the characteristic values and characteristic vectors as in case of complex neutrosophic vector spaces / linear algebras defined over the fields.

**DEFINITION 1.10:** *Let V be a special complex neutrosophic vector space over the complex neutrosophic field C (⟨Q ∪ I⟩). Let T be a special linear operator on V. A complex neutrosophic value of T is a scalar c in C (⟨Q ∪ I⟩) so that there is a non zero neutrosophic complex vector $\alpha$ in V with $T\alpha = c\alpha$.*

*If c is the special characteristic value of T then*
  *i) any $\alpha$ such that $T\alpha = c\alpha$ is called the characteristic neutrosophic complex vector of T associated with c.*
  *ii) The collection of all $\alpha$ such that $T\alpha = c\alpha$ is called the neutrosophic complex characteristic space associated with c.*

*If the complex neutrosophic field C(⟨Q ∪ I⟩) is replaced by ⟨Q ∪ I⟩ or C(Q) = {a + bi | a, b ∈ Q} or Q we get the characteristic value c as neutrosophic number a + bI or c + di or a respectively and the associated characteristic space of them would be a neutrosophic complex neutrosophic subspace or complex neutrosophic complex subspace or ordinary complex neutrosophic subspace respectively.*

The following theorem is left an exercise for the reader to prove.

**THEOREM 1.10:** *Let T be a complex neutrosophic linear operator on a finite dimensional special (ordinary or complex or neutrosophic) complex neutrosophic vector space V and let c be a scalar in C(⟨Q ∪ I⟩) (or Q or C(Q) or ⟨Q ∪ I⟩ respectively),*
*The following are equivalent*
  *i)  c is the characteristic value of T.*
  *ii) The operator (T – cI) is singular.*
  *iii) det (T – cI) = 0.*

We can define diagonlizable linear operator as in case of other linear operator. Also we will show by an example how to find the characteristic polynomial in case of special (ordinary or complex or neutrosophic) matrices.

Let
$$M = \begin{bmatrix} a_1 + b_1 i + c_1 I + d_1 iI & a_2 + b_2 i + c_2 I + d_2 iI \\ a_3 + b_3 i + c_3 I + d_3 iI & a_4 + b_4 i + c_4 I + d_4 iI \end{bmatrix},$$

be a 2 × 2 matrix with entries over the field C (⟨Q ∪ I⟩) such that the matrix (M – $cI_{2\times 2}$) is non invertible.

All results related with linear operators of special (ordinary or complex or neutrosophic) vector spaces can be derived as in case of usual vector spaces with simple appropriate modifications.

Now we proceed onto define the concept of linear functionals. We can define four types of linear functionals.
  Let V be a special neutrosophic complex vector space over the neutrosophic complex field F = C(⟨Q ∪ I⟩). The special linear functional on V is a map (a linear transformation) f : V → F such that
$$F (c \alpha + \beta ) = cf (\alpha) + f (\beta)$$
$\alpha, \beta \in$ V and C ∈ C (⟨Q ∪ I⟩).

We will first illustrate this by an example.

*Example 1.27:* Let

$$V = \left\{ \begin{pmatrix} a_1 & a_2 & a_3 \\ a_4 & a_5 & a_6 \end{pmatrix} \middle| a_i \in C(\langle Q \cup I \rangle); 1 \leq i \leq 6 \right\}$$

be a special neutrosophic complex vector space over the neutrosophic complex field $F = C(\langle Q \cup I \rangle)$.

Define $f : V \to F$ by

$$f\left( \begin{pmatrix} a_1 & a_2 & a_3 \\ a_4 & a_5 & a_6 \end{pmatrix} \right) = \sum_{i=0}^{6} a_i.$$

Clearly f is a special linear functional on V.

If V is a neutrosophic complex neutrosophic vector space over the neutrosophic field $K = \langle Q \cup I \rangle$. We define $f : V \to K$ and $f(v) \in \langle Q \cup I \rangle$ that is only neutrosophic number.

*Example 1.28:* Let

$$V = \left\{ \begin{bmatrix} a_1 \\ a_2 \\ a_3 \\ a_4 \\ a_5 \\ a_6 \end{bmatrix} \middle| a_i \in C(\langle Q \cup I \rangle); 1 \leq i \leq 6 \right\}$$

be a neutrosophic complex neutrosophic vector space over the neutrosophic field $\langle Q \cup I \rangle = K$.
Define $f : V \to K$ by

$$f\left( \begin{bmatrix} a_1 \\ a_2 \\ \vdots \\ a_6 \end{bmatrix} \right) = (a_{11} + \ldots + a_{16}) + (c_{11} + \ldots + c_{16})I$$

where $a_1 = a_{11} + b_{11}i + c_{11}I + id_{11}I$
$a_2 = a_{12} + b_{12}i + c_{12}I + id_{12}I$
$\vdots$
and $a_6 = a_{16} + b_{16}i + c_{16}I + id_{16}I$.

Clearly

$$f\left(\begin{bmatrix} a_1 \\ a_2 \\ \vdots \\ a_6 \end{bmatrix}\right) = (a_{11} + \ldots + a_{16}) + (c_{11} + \ldots + c_{16})I$$

is in $\langle Q \cup I \rangle$ is a neutrosophic linear functional on V.

*Example 1.29:* Let $V = \{(a_1, a_2, a_3) \mid a_i \in C (\langle Q \cup I \rangle); 1 \leq i \leq 3\}$ be a complex neutrosophic complex vector space over the rational complex field $K = C(Q) = \{a + ib \mid a, b \in Q\}$. Let $f : V \to K$ defined by $f(a_1, a_2, a_3) = (a_{11} + a_{12} + a_{13}) + i(b_{11} + b_{12} + b_{13})$ where
$a_1 = a_{11} + b_{11}i + c_{11}I + id_{11}I$
$a_2 = a_{12} + b_{12}i + c_{12}I + id_{12}I$
and $a_3 = a_{13} + b_{13}i + c_{13}I + id_{13}I$

where $a_{ij}$, $b_{ij}$, $c_{ij}$ and $d_{ij}$ are in Q; $1 \leq j \leq 3$.

Clearly f is a complex linear functional on V. Now we proceed onto give an example of a ordinary linear functional on V.

*Example 1.30:* Let

$$V = \left\{ \begin{pmatrix} a_1 & a_2 \\ a_3 & a_4 \end{pmatrix} \middle| a_j = a_{1i} + ib_{1i} + c_{1i}I + id_{1i}I \right.$$

where $a_{1i}$, $b_{1i}$, $c_{1i}$, $d_{1i} \in Q$; $1 \leq i \leq 3$ and $1 \leq j \leq 4\}$ be an ordinary neutrosophic complex vector space over the rational field Q.
Define $f : V \to Q$ by

$$f\left(\begin{pmatrix} a_1 & a_2 \\ a_3 & a_4 \end{pmatrix}\right) = a_{11} + a_{12} + a_{13};$$

f is an ordinary linear functional on V.

Now having seen the definitions of linear functionals interested reader can derive all properties related with linear functionals with appropriate changes. We can now define set neutrosophic complex vector spaces, semigroup neutrosophic complex vector spaces and group neutrosophic complex vector spaces. Also the corresponding linear algebras.

**DEFINITION 1.11:** *Let $V \subseteq C(\langle Q \cup I \rangle)$ be a proper subset of complex neutrosophic rationals. $S \subseteq Q$ be a subset of S. We define V to a set complex neutrosophic vector space over the set $S \subseteq Q$ if for all $v \in V$ and $s \in S$, vs and sv $\in V$.*

We give examples of this situation.

*Example 1.31:* Let

$$V = \left\{ (a_1, a_2, a_3), \begin{bmatrix} a_1 & a_2 \\ a_3 & a_4 \\ a_5 & a_6 \end{bmatrix}, \sum_{i=0}^{20} a_i x^i \middle| a_i \in C(\langle Q \cup I \rangle); 0 \leq i \leq 20 \right\}$$

be set vector space of neutrosophic complex rationals over $Z \subseteq Q$ or set complex neutrosophic rational vector space over the set Z.

*Example 1.32:* Let

$$V = \left\{ \begin{bmatrix} a_1 & a_2 \\ 0 & a_3 \end{bmatrix}, \begin{bmatrix} 0 & a_1 & a_2 & a_3 \\ a_4 & 0 & a_5 & 0 \end{bmatrix} \middle| a_i \in C(\langle Q \cup I \rangle); \right.$$

$1 \leq i \leq 5$} be a set complex neutrosophic vector space over the set $S = \{0, 1, 2, 3, 4, 12, 17, -5, -9, -23\} \subseteq Q$.

*Example 1.33:* Let

$$V = \left\{ \sum_{i=0}^{8} a_i x^i, (a_1, a_2, a_3, \ldots, a_{20}), \begin{bmatrix} a_1 \\ a_2 \\ \vdots \\ a_{12} \end{bmatrix} \;\middle|\; a_i \in C(\langle Q \cup I \rangle); \right.$$

$0 \leq i \leq 20\}$ be a set neutrosophic complex vector space over the set $S = \{0, -1, 1\}$.

Now having seen examples of set neutrosophic complex vector spaces we now proceed onto define set neutrosophic complex vector subspaces of V over the set S.

Let V be a set complex neutrosophic vector space over the set S. Suppose $W \subseteq V$ and if W itself is a set complex neutrosophic vector space over the set S then we define W to be a set complex neutrosophic vector subspace of V over S.

We will illustrate this situation by some examples.

*Example 1.34:* Let

$$V = \left\{ (a_1, a_2, a_3), \begin{bmatrix} a_1 \\ a_2 \\ a_3 \\ a_4 \\ a_5 \end{bmatrix}, \sum_{i=0}^{3} a_i x^i \;\middle|\; a_i \in C(\langle Q \cup I \rangle); 0 \leq i \leq 5 \right\}$$

be set vector space of neutrosophic complex rationals over the set $S = \{0, 1\}$.

Consider

$$W = \left\{ (0, a_1, 0), \begin{bmatrix} a_1 \\ 0 \\ a_2 \\ 0 \\ a_3 \end{bmatrix} \;\middle|\; a_1, a_2, a_3 \in C(\langle Q \cup I \rangle) \right\} \subseteq V,$$

W is a set vector subspace of neutrosophic complex rationals over the set S ={0, 1}.

Take

$$M = \left\{ (a_1, a_2, a_3), \begin{bmatrix} a_1 \\ a_2 \\ a_3 \\ a_4 \\ a_5 \end{bmatrix}, \sum_{i=0}^{3} a_i x^i \;\middle|\; a_i \in C(\langle Z \cup I \rangle) \right.$$

$\subseteq C(\langle Q \cup I \rangle); 0 \le i \le 5\} \subseteq V$ is a set complex neutrosophic vector subspace of V over the set S = {0, 1}.

*Example 1.35:* Let

$$V = \left\{ \sum_{i=0}^{25} a_i x^i, \begin{bmatrix} a_1 & a_2 & a_3 \\ a_4 & a_5 & a_6 \\ \vdots & \vdots & \vdots \\ a_{28} & a_{29} & a_{30} \end{bmatrix}, \begin{bmatrix} a_1 & a_2 & \dots & a_{10} \\ a_{11} & a_{12} & \dots & a_{20} \\ a_{21} & a_{22} & \dots & a_{30} \\ a_{31} & a_{32} & \dots & a_{40} \end{bmatrix} \right.$$

$a_i \in C(\langle Q \cup I \rangle); 0 \le i \le 40\}$ be a set neutrosophic complex vector space over the set $\{-5, 4, 2, 1, 3, 8, 10, 0\} \subseteq S$.

Consider

$$W = \left\{ \sum_{i=0}^{12} a_i x^i, \begin{bmatrix} a_1 & a_2 & a_3 \\ 0 & 0 & 0 \\ 0 & 0 & 0 \\ \vdots & \vdots & \vdots \\ 0 & 0 & 0 \\ a_4 & a_5 & a_6 \end{bmatrix}, \begin{bmatrix} a_1 & a_2 & \dots & a_{10} \\ 0 & 0 & \dots & 0 \\ \vdots & \vdots & & \vdots \\ a_{21} & a_{22} & \dots & a_{30} \\ 0 & 0 & \dots & 0 \end{bmatrix} \right.$$

$a_j, a_i \in C(\langle Q \cup I \rangle); 0 \le a_i \le 12, a_j = 1, 2, 3, 4, 5, 6, 7, 8, 9, 10, 21, 22, \dots, 30\} \subseteq V$ is a set neutrosophic complex vector subspace of V over the set S.

Now as in case of usual set vector spaces we can derive all the related properties we can also define the notion of subset neutrosophic complex vector subspace over a subset, which is simple and left as an exercise to the reader.

We give examples of this structure.

*Example 1.36:* Let

$$V = \left\{ \begin{bmatrix} a_1 & a_2 & a_3 & a_4 \\ a_5 & a_6 & a_7 & a_8 \end{bmatrix}, \begin{bmatrix} a_1 & a_2 & a_3 \\ a_4 & a_5 & a_6 \\ a_7 & a_8 & a_9 \end{bmatrix}, \begin{bmatrix} a_1 & a_2 & a_3 & a_4 \\ a_5 & a_6 & a_7 & a_8 \\ a_9 & a_{10} & a_{11} & a_{12} \\ a_{13} & a_{14} & a_{15} & a_{16} \end{bmatrix} \right.$$

$a_i \in C (\langle Q \cup I \rangle); 1 \le i \le 16\}$ be a set neutrosophic complex vector space over the set $S = 3Z \cup 5Z \cup 7Z \cup 13Z$.
    Consider

$$P = \left\{ \begin{bmatrix} a_1 & a_2 & 0 & 0 \\ 0 & 0 & a_3 & a_4 \end{bmatrix}, \begin{bmatrix} a_1 & 0 & a_2 \\ 0 & a_3 & 0 \\ a_4 & 0 & a_5 \end{bmatrix}, \begin{bmatrix} a_1 & 0 & a_2 & 0 \\ 0 & a_4 & 0 & a_5 \\ a_3 & 0 & a_6 & 0 \\ 0 & a_7 & 0 & a_8 \end{bmatrix} \right.$$

$a_i \in C (\langle Q \cup I \rangle); 1 \le i \le 8\} \subseteq V$; P is a subset neutrosophic complex vector subspace of V over the subset $T = 3Z \cup 13Z \subset S$.
    Also

$$M = \left\{ \begin{bmatrix} 0 & 0 & 0 & 0 \\ a_1 & a_2 & a_3 & a_4 \end{bmatrix}, \begin{bmatrix} a_1 & a_2 & a_3 \\ 0 & 0 & 0 \\ 0 & 0 & a_4 \end{bmatrix}, \begin{bmatrix} a_1 & 0 & 0 & 0 \\ 0 & a_2 & 0 & 0 \\ 0 & 0 & a_3 & 0 \\ 0 & 0 & 0 & a_4 \end{bmatrix} \right.$$

$a_i \in C (\langle Q \cup I \rangle); 1 \le i \le 4\} \subseteq V$, M is a subset neutrosophic complex vector subspace of V over the set $T = 3Z \cup 5Z \subseteq S$.

*Example 1.37:* Let

$$V = \left\{ \sum_{i=0}^{29} a_i x^i, \begin{bmatrix} a_1 & a_2 \\ a_3 & a_4 \\ \vdots & \vdots \\ a_{11} & a_{12} \end{bmatrix}, \begin{bmatrix} a_1 & a_2 \\ a_3 & a_4 \end{bmatrix}, (a_1, a_2, a_3, a_4) \right.$$

$a_i \in C (\langle Q \cup I \rangle); 0 \leq i \leq 29\}$ be a set neutrosophic complex vector space over the set $S = 3Z^+ \cup 5Z \cup 7Z^+$.

Consider

$$M = \left\{ \sum_{i=0}^{20} a_i x^i, \begin{bmatrix} 0 & a_1 \\ a_2 & 0 \end{bmatrix}, \begin{bmatrix} a_1 & 0 \\ a_2 & 0 \\ 0 & 0 \\ 0 & 0 \\ 0 & 0 \\ a_3 & a_4 \end{bmatrix}, (a_1, 0, a_2, 0) \right.$$

$a_i \in C (\langle Q \cup I \rangle); 0 \leq i \leq 20\} \subseteq V$; is a subset neutrosophic complex vector subspace of V over the subset $T = \{3Z^+ \cup 5Z\} \subseteq S$. Take

$$P = \left\{ \sum_{i=0}^{10} a_i x^i, \begin{bmatrix} 0 & 0 \\ a_1 & 0 \end{bmatrix}, \begin{bmatrix} a_1 & 0 \\ a_2 & 0 \\ a_3 & 0 \\ a_4 & 0 \\ a_5 & 0 \\ a_6 & 0 \end{bmatrix}, (0, 0, a_1, a_2) \right.$$

$a_i \in C (\langle Q \cup I \rangle); 0 \leq i \leq 10\} \subseteq V$, is a subset vector complex neutrosophic subspace of V over the subset $T = 7Z^+ \subseteq S$.

We can define set linear transformation, set linear operator and set basis of a set neutrosophic complex vector space over the set S, which is left as an exercise as it can be carried out as a matter of routine. We can also define special set neutrosophic complex vector space over the complex neutrosophic subset of $C(\langle Q \cup I \rangle)$.

*Example 1.38:* Let

$$V = \left\{ \begin{bmatrix} a_1 & a_2 \\ a_3 & a_4 \end{bmatrix}, \begin{bmatrix} a_1 \\ a_2 \\ \vdots \\ a_8 \end{bmatrix}, \begin{pmatrix} a_1 & a_2 & a_3 \\ a_4 & a_5 & a_6 \\ a_7 & a_8 & a_9 \end{pmatrix} \middle| a_i \in C(\langle Q \cup I \rangle); 1 \leq i \leq 9 \right\}$$

be a set neutrosophic of complex neutrosophic vector over a set $S = \{1, i, 1 + I, 3i - 1, 2I+1+5i + 18I\} \subseteq C(\langle Q \cup I \rangle)$. We observe $S \not\subseteq Q$ so we call such vector spaces as special set vector spaces.

**DEFINITION 1.12:** *Let V be a set vector space of complex neutrosophic rationals (set complex neutrosophic rationals vector space) over a set $S \subseteq C(\langle Q \cup I \rangle)$ and $S \not\subseteq Q$ then we define V to be a special set vector space of complex neutrosophic rationals over the set S if vs and sv $\in$ V for all s $\in$ S and v $\in$ V.*

We will illustrate this by some examples.

*Example 1.39:* Let

$$V = \left\{ \begin{bmatrix} a_1 & a_2 \\ a_3 & a_4 \end{bmatrix}, \begin{bmatrix} a_1 & a_2 & \ldots & a_7 \\ a_8 & a_9 & \ldots & a_{14} \end{bmatrix}, (a_1, a_2, \ldots, a_{10}) \middle| \right.$$

$a_i \in C(\langle Q \cup I \rangle); 1 \leq i \leq 14\}$ be a special set complex neutrosophic rational vector space over the set $S = C(\langle Z \cup I \rangle) \subseteq C(\langle Q \cup I \rangle)$.

*Example 1.40:* Let

$$M = \left\{ \begin{bmatrix} a_1 & a_2 & a_3 \\ a_4 & a_5 & a_6 \\ a_7 & a_8 & a_9 \end{bmatrix}, \begin{bmatrix} a_1 & a_2 & \ldots & a_9 \\ a_{10} & a_{11} & \ldots & a_{18} \end{bmatrix}, \begin{bmatrix} a_1 \\ a_2 \\ \vdots \\ a_7 \end{bmatrix} \right.$$

$a_i \in C(\langle Q \cup I \rangle)$; $1 \leq i \leq 18$} be a special set neutrosophic complex vector space over the set $S = C(\langle 3Z \cup I \rangle) \subseteq C(\langle Q \cup I \rangle)$.

*Example 1.41:* Let

$$V = \left\{ (a_1, a_2, \ldots, a_9), \begin{bmatrix} a_1 & a_2 & a_3 & a_4 \\ a_5 & a_6 & a_7 & a_8 \\ a_9 & a_{10} & a_{11} & a_{12} \end{bmatrix}, \sum_{i=0}^{12} a_i x^i \right\}$$

$a_i \in C(\langle Q \cup I \rangle)$; $0 \leq i \leq 12$} be a special set complex neutrosophic vector space over the set $C(\langle 3Z \cup I \rangle)$.

*Example 1.42:* Let

$$M = \left\{ \begin{bmatrix} a_1 & a_2 & a_3 \\ a_4 & a_5 & a_6 \end{bmatrix}, \begin{bmatrix} a_1 & a_2 \\ a_3 & a_4 \\ a_5 & a_6 \\ a_7 & a_8 \\ a_9 & a_{10} \end{bmatrix}, \sum_{i=0}^{25} a_i x^i, \begin{bmatrix} a_1 & a_2 \\ a_3 & a_4 \\ a_5 & a_6 \end{bmatrix} \right\}$$

$a_i \in C(\langle Q \cup I \rangle)$; $0 \leq i \leq 25$} be a special set complex neutrosophic vector space over the set $S = \{C(\langle 5Z \cup I \rangle) \cup C(\langle 13Z \cup I \rangle)\}$. Take

$$P = \left\{ \sum_{i=0}^{10} a_i x^i, \begin{bmatrix} a_1 & 0 & a_2 \\ 0 & a_3 & 0 \end{bmatrix}, \middle| a_i \in C(<Q \cup I>); 0 \leq i \leq 10 \right\} \subseteq V,$$

P is a special set neutrosophic complex vector subspace of V over S. Consider

$$W = \left\{ \begin{bmatrix} a_1 & 0 \\ 0 & a_2 \\ a_3 & 0 \end{bmatrix}, \sum_{i=0}^{20} a_i x^i, \begin{bmatrix} 0 & a_1 \\ a_2 & 0 \\ 0 & a_3 \\ a_4 & 0 \\ 0 & a_5 \end{bmatrix} \middle| a_i \in C(\langle Q \cup I \rangle); 0 \leq i \leq 20 \right\}$$

⊆ V, W is a special set complex neutrosophic vector subspace of V over the set S.

Consider

$$B = \left\{ \begin{bmatrix} a_1 & a_2 \\ 0 & 0 \\ a_3 & a_4 \end{bmatrix}, \begin{bmatrix} 0 & 0 \\ a_1 & a_2 \\ 0 & 0 \\ a_3 & a_4 \\ 0 & 0 \end{bmatrix}, \sum_{i=0}^{5} a_i x^i \;\middle|\; a_i \in C(\langle Q \cup I \rangle); \right.$$

$0 \le i \le 5\} \subseteq V$, take $T = \{C(\langle 25Z \cup I \rangle) \cup C(\langle 39Z \cup I \rangle)\} \subseteq S$; we see B is a special subset complex neutrosophic vector subspace of V over the subset T of S.

Now having seen examples of special set subspaces and special subset vector subspaces of complex neutrosophic rationals we proceed onto define the notion of semigroup complex neutrosophic vector space and special semigroup complex neutrosophic vector space.

Just we mention in case of special set neutrosophic complex vector space also one can define special set linear transformations provided both the special set vector spaces are defined over the same set of complex neutrosophic numbers. Further the basis, direct sum of subspace and other properties can be easily derived as a matter of routine. All these work is left as exercises to the reader.

**DEFINITION 1.13:** *Let V be any subset of complex neutrosophic numbers and S be any additive semigroup with zero. We call V to be a semigroup neutrosophic complex vector space over S if the following conditions hold good.*

  i)  *$vs = sv \in V$ for all $s \in S$ and $v \in V$.*
  ii) *$0.v = 0 \in V$ for all $v \in V$ and $0 \in S$; $0 \in V$ is a zero vector.*
  iii) *$(s_1 + s_2) v = s_1 v + s_2 v$ for all $s_1, s_2 \in S$ and $v \in V$.*

We will first illustrate this situation by some examples.

*Example 1.43:* Let

$$V = \left\{ \begin{bmatrix} a_1 & a_2 \\ a_3 & a_4 \\ \vdots & \\ a_9 & a_{10} \end{bmatrix}, \begin{bmatrix} a_1 & a_2 \\ a_3 & a_4 \end{bmatrix}, \begin{bmatrix} a_1 & a_2 & \ldots & a_{10} \\ a_{11} & a_{12} & \ldots & a_{20} \\ a_{21} & a_{22} & \ldots & a_{30} \end{bmatrix} \right\}$$

$a_i \in C(\langle Q \cup I \rangle)$; $1 \leq i \leq 30\}$ be a semigroup neutrosophic complex vector space over the semigroup $S = Z^+ \cup \{0\}$.

*Example 1.44:* Let

$$M = \left\{ \sum_{i=0}^{40} a_i x^i, \begin{bmatrix} a_1 & a_2 & a_3 & a_4 \\ a_5 & a_6 & a_7 & a_8 \\ a_9 & a_{10} & a_{11} & a_{12} \\ a_{13} & a_{14} & a_{15} & a_{16} \end{bmatrix}, \begin{bmatrix} a_1 & a_2 \\ a_3 & a_4 \\ \vdots & \vdots \\ a_{41} & a_{42} \end{bmatrix}, \begin{pmatrix} a_1 & a_2 & \ldots & a_{20} \\ a_{21} & a_{22} & \ldots & a_{40} \end{pmatrix} \right\}$$

$a_i \in C(\langle Q \cup I \rangle)$; $0 \leq i \leq 40\}$ be a semigroup neutrosophic complex vector space over the additive semigroup $S = 5Z^+ \cup \{0\}$.

*Example 1.45:* Let

$$P = \left\{ \begin{bmatrix} a_1 & a_2 \\ a_3 & a_4 \end{bmatrix}, \begin{bmatrix} a_1 \\ a_2 \\ \vdots \\ a_{10} \end{bmatrix}, \begin{pmatrix} a_1 & a_2 & a_3 & a_4 & a_5 \\ a_6 & a_7 & a_8 & a_9 & a_{10} \end{pmatrix} \right\}$$

$a_i \in C(\langle Q \cup I \rangle)$; $1 \leq i \leq 10\}$ be a semigroup neutrosophic complex vector space over the semigroup $S = 2Z^+ \cup \{0\}$.

Consider

$$V = \left\{ \begin{bmatrix} a_1 & a_2 \\ 0 & 0 \end{bmatrix}, \begin{bmatrix} a_1 \\ 0 \\ 0 \\ 0 \\ \vdots \\ a_9 \\ a_{10} \end{bmatrix}, \begin{pmatrix} 0 & a_1 & 0 & 0 & a_2 \\ a_3 & 0 & 0 & 0 & 0 \end{pmatrix} \middle| a_1, a_2, a_3, a_9, a_{10} \in \right.$$

$C(\langle Q \cup I \rangle) \subseteq P$; V is a semigroup neutrosophic complex vector subspace of P over the semigroup S.

Also

$$B = \left\{ \begin{bmatrix} 0 & a_2 \\ a_1 & 0 \end{bmatrix}, \begin{pmatrix} a_1 & a_1 & a_1 & 0 & 0 \\ 0 & 0 & 0 & a_1 & a_1 \end{pmatrix} \middle| \right.$$

$a_1, a_2 \in C(\langle Q \cup I \rangle) \subseteq B$ is a subsemigroup neutrosophic complex vector subspace of P over the subsemigroup $A = 8Z^+ \cup \{0\}$ of the semigroup S.

Here also all properties of semigroup linear transformation of vector spaces can be obtained provided they are defined over the same, semigroup, semigroup linear operator and basis can be defined as in case of usual semigroup vector spaces.

*Example 1.46:* Let

$$V = \left\{ \begin{bmatrix} a_1 & a_2 \\ a_3 & a_4 \end{bmatrix}, \begin{bmatrix} a_1 & a_4 & a_7 & a_{10} & a_{13} \\ a_2 & a_5 & a_8 & a_{11} & a_{14} \\ a_3 & a_6 & a_9 & a_{12} & a_{15} \end{bmatrix}, \begin{bmatrix} a_1 & a_2 \\ a_3 & a_4 \\ \vdots & \vdots \\ a_9 & a_{10} \end{bmatrix} \middle| \right.$$

$a_i \in C(\langle Q \cup I \rangle); \ 1 \le i \le 15\}$ be a semigroup complex neutrosophic vector space over the semigroup $S = 3Z^+ \cup \{0\}$.

Consider

$$W_1 = \left\{ \begin{bmatrix} a_1 & a_2 \\ a_3 & a_4 \end{bmatrix} \middle| \; a_i \in C(\langle Q \cup I \rangle); \; 1 \leq i \leq 4 \right\} \subseteq V,$$

$$W_2 = \left\{ \begin{bmatrix} a_1 & a_4 & a_7 & a_{10} & a_{13} \\ a_2 & a_5 & a_8 & a_{11} & a_{14} \\ a_3 & a_6 & a_9 & a_{12} & a_{15} \end{bmatrix} \middle| \; a_i \in C(\langle Q \cup I \rangle); 1 \leq i \leq 15 \right\} \subseteq V$$

and

$$W_3 = \left\{ \begin{bmatrix} a_1 & a_2 \\ a_3 & a_4 \\ \vdots & \vdots \\ a_9 & a_{10} \end{bmatrix} \middle| \; a_i \in C(\langle Q \cup I \rangle); \; 1 \leq i \leq 10 \right\} \subseteq V$$

be semigroup complex neutrosophic vector subspaces of V over the semigroup $S = 3Z^+ \cup \{0\}$. It is easily seen $V = W_1 + W_2 + W_3$ and $W_i \cap W_j = \phi$; $1 \leq i, j \leq 3$. Thus V is a direct sum of semigroup complex neutrosophic vector subspaces of V over S.

Let

$$W_1 = \left\{ \begin{bmatrix} a_1 & a_2 \\ a_3 & a_4 \end{bmatrix}, \begin{bmatrix} a_1 & 0 & 0 & 0 & a_3 \\ 0 & 0 & a_2 & 0 & 0 \\ 0 & 0 & 0 & a_4 & 0 \end{bmatrix} \right\}$$

$a_i \in C(\langle Q \cup I \rangle); \; 1 \leq i \leq 4\} \subseteq V,$

$$W_2 = \left\{ \begin{bmatrix} a_1 & 0 \\ 0 & a_2 \end{bmatrix}, \begin{bmatrix} a_1 & a_2 & a_3 & a_4 & a_5 \\ a_6 & a_7 & a_8 & a_9 & a_{10} \\ a_{11} & a_{12} & a_{13} & a_{14} & a_{15} \end{bmatrix}, \begin{bmatrix} 0 & 0 \\ 0 & 0 \\ 0 & 0 \\ 0 & 0 \\ a_1 & a_2 \end{bmatrix} \right\}$$

$a_i \in C(\langle Q \cup I \rangle); \; 1 \leq i \leq 15\} \subseteq V$ and

$$W_3 = \left\{ \begin{bmatrix} a_1 & 0 & 0 & 0 & a_3 \\ 0 & 0 & a_2 & 0 & 0 \\ 0 & 0 & 0 & a_4 & 0 \end{bmatrix}, \begin{bmatrix} a_1 & a_2 \\ a_3 & a_4 \\ a_5 & a_6 \\ a_7 & a_8 \\ a_9 & a_{10} \end{bmatrix} \right.$$

$a_i \in C(\langle Q \cup I \rangle);\ 1 \le i \le 10\} \subseteq V$; clearly $W_1 + W_2 + W_3 = V$ and $W_i \cap W_j \ne \phi$ if $i \ne j$; $1 \le i, j \le 3$. Thus V is only a pseudo direct union of subspaces.

We define special semigroup complex vector as follows.

**DEFINITION 1.14:** *Let V be a semigroup neutrosophic complex vector space over the complex neutrosophic additive semigroup S. Then we define V to be a special semigroup neutrosophic complex vector space over the semigroup S.*

We will illustrate this situation by some examples.

*Example 1.47:* Let

$$V = \left\{ \begin{bmatrix} a_1 & a_2 \\ a_3 & a_4 \end{bmatrix}, \begin{bmatrix} a_1 & a_2 & a_3 \\ a_4 & a_5 & a_6 \\ \vdots & \vdots & \vdots \\ a_{28} & a_{29} & a_{30} \end{bmatrix}, \begin{pmatrix} a_1 & a_2 & \ldots & a_{20} \\ a_{21} & a_{22} & \ldots & a_{40} \end{pmatrix} \right.$$

$a_i \in C(\langle Q \cup I \rangle);\ 1 \le i \le 40\}$ be a special semigroup neutrosophic complex vector space over the complex neutrosophic additive semigroup $S = C(\langle 3Z \cup I \rangle)$.

*Example 1.48:* Let

$$M = \left\{ \begin{bmatrix} a_1 & a_2 & a_3 \\ a_4 & a_5 & a_6 \\ \vdots & \vdots & \vdots \\ a_{43} & a_{44} & a_{45} \end{bmatrix}, \begin{bmatrix} a_1 & a_2 & \ldots & a_{10} \\ a_{11} & a_{12} & \ldots & a_{20} \end{bmatrix}, \sum_{i=0}^{28} a_i x^i \right.$$

$a_i \in C(\langle Q \cup I \rangle); \ 0 \le i \le 45\}$ be a special semigroup complex neutrosophic vector space over the semigroup under addition $S = C(\langle 5Z \cup I \rangle)$.

*Example 1.49:* Let

$$M = \left\{ \sum_{i=0}^{12} a_i x^i, \begin{bmatrix} a_1 & a_2 & a_3 \\ a_4 & a_5 & a_6 \\ a_7 & a_8 & a_9 \end{bmatrix}, \begin{bmatrix} a_1 & a_2 \\ a_3 & a_4 \\ \vdots & \vdots \\ a_{31} & a_{32} \end{bmatrix} \middle| a_i \in C(\langle Q \cup I \rangle); 0 \le i \le 32 \right\}$$

be a special semigroup complex neutrosophic vector space over the neutrosophic complex semigroup $S = C(\langle 10Z \cup I \rangle)$ under addition.

Take

$$V = \left\{ \sum_{i=0}^{6} a_i x^i, \begin{bmatrix} a_1 & a_2 & a_3 \\ 0 & a_4 & a_5 \\ 0 & 0 & a_6 \end{bmatrix} \middle| a_i \in C(\langle Q \cup I \rangle); 0 \le i \le 6 \right\} \subseteq M$$

be a special semigroup complex neutrosophic vector subspace of M over the neutrosophic complex semigroup S.

Take

$$P = \left\{ \sum_{i=0}^{6} a_i x^i, \begin{bmatrix} a_1 & 0 \\ a_2 & 0 \\ \vdots & \vdots \\ a_{16} & 0 \end{bmatrix} \middle| a_i \in C(\langle Q \cup I \rangle); 0 \le i \le 16 \right\} \subseteq M$$

be a special subsemigroup complex neutrosophic subvector space of M over the subsemigroup $T = C(\langle 40Z \cup I \rangle) \subseteq S$ under addition. We can also write M as a direct sum of subspaces as well as a pseudo direct sum of special semigroup vector subspaces.

Now we proceed onto define the notion of set linear algebra, special set linear algebra, semigroup linear algebra and special semigroup linear algebra using complex neutrosophic rationals.

**DEFINITION 1.15:** *Let V be a set neutrosophic complex vector space over the set S. We define V to be a set complex neutrosophic linear algebra over the set S if s (a + b) = sa + sb for all a, b ∈ V and s ∈ S.*

We give examples of them.

*Example 1.50:* Let

$$V = \left\{ \begin{bmatrix} a_1 & a_2 \\ a_3 & a_4 \\ \vdots & \vdots \\ a_{13} & a_{14} \end{bmatrix} \middle| a_i \in C(\langle Q \cup I \rangle); 1 \le i \le 14 \right\}$$

be a set complex neutrosophic linear algebra over the set S = $(5Z \cup 7Z)$.

*Example 1.51:* Let

$$M = \left\{ \begin{bmatrix} a_1 & a_2 \\ a_3 & a_4 \end{bmatrix} \middle| a_i \in C(\langle Q \cup I \rangle); 1 \le i \le 4, \times \right\}$$

be a set complex neutrosophic linear algebra over the set S = $3Z \cup 7Z \cup 5Z$.

*Example 1.52:* Let

$$M = \left\{ \begin{bmatrix} a_1 & a_2 & \ldots & a_{12} \\ a_{13} & a_{14} & \ldots & a_{24} \end{bmatrix} \middle| a_i \in C(\langle Q \cup I \rangle); 1 \le i \le 24 \right\}$$

be a set neutrosophic complex linear algebra over the set S = $5Z \cup 2Z \cup 17Z$.

$$W = \left\{ \begin{bmatrix} a_1 & 0 & a_3 & \ldots & a_{11} & 0 \\ 0 & a_{14} & 0 & \ldots & 0 & a_{24} \end{bmatrix} \middle| a_1, a_3, \ldots, a_{24} \in C(\langle Q \cup I \rangle) \right\}$$

⊆ M is a set neutrosophic complex linear subalgebra of M over the set S. Take

$$R = \left\{ \begin{bmatrix} a_1 & a_1 & \ldots & a_1 \\ b_1 & b_1 & \ldots & b_1 \end{bmatrix} \middle| a_1, b_1 \in C(\langle Q \cup I \rangle) \right\} \subseteq M;$$

R is a set neutrosophic complex linear algebra over the set $S = 5Z \cup 2Z \cup 17Z$.

Consider

$$A = \left\{ \begin{bmatrix} a_1 & 0 & a_3 & 0 & a_5 & 0 & ... & a_{11} & 0 \\ a_2 & 0 & a_4 & 0 & a_6 & 0 & ... & a_{12} & 0 \end{bmatrix} \middle| \begin{array}{l} a_i \in C(\langle Q \cup I \rangle); \\ 1 \le i \le 12 \end{array} \right\}$$

$\subseteq M$ be a subset complex neutrosophic linear subalgebra of M over the subset $T = 5Z \cup 2Z \subseteq S$.

Suppose

$$B = \left\{ \begin{bmatrix} a_1 & a_2 & ... & a_{12} \\ 0 & 0 & ... & 0 \end{bmatrix} \middle| a_i \in C(\langle Q \cup I \rangle); 1 \le i \le 12 \right\} \subseteq M;$$

B is a subset neutrosophic complex linear subalgebra of M over the subset $R = 17Z \subseteq S$.

*Example 1.53:* Let

$$M = \left\{ \begin{bmatrix} a_1 & a_2 & a_3 \\ a_4 & a_5 & a_6 \\ a_7 & a_8 & a_9 \\ a_{10} & a_{11} & a_{12} \\ a_{13} & a_{14} & a_{15} \\ a_{16} & a_{17} & a_{18} \end{bmatrix} \middle| a_i \in C(\langle Q \cup I \rangle); 1 \le i \le 18 \right\}$$

be a set neutrosophic complex linear algebra over the set $S = 3Z^+ \cup \{0\}$.

Consider

$$W_1 = \left\{ \begin{bmatrix} a_1 & a_2 & a_3 \\ 0 & 0 & 0 \\ a_4 & a_5 & a_6 \\ 0 & 0 & 0 \\ 0 & 0 & 0 \\ 0 & 0 & 0 \end{bmatrix} \middle| a_i \in C(\langle Q \cup I \rangle); 1 \le i \le 6 \right\}$$

a proper subset of M,

$$W_2 = \left\{ \begin{bmatrix} 0 & 0 & 0 \\ a_1 & a_2 & a_3 \\ 0 & 0 & 0 \\ a_4 & a_5 & a_6 \\ 0 & 0 & 0 \\ 0 & 0 & 0 \end{bmatrix} \middle| a_i \in C(\langle Q \cup I \rangle); 1 \le i \le 6 \right\} \subseteq M,$$

$$W_3 = \left\{ \begin{bmatrix} 0 & 0 & 0 \\ 0 & 0 & 0 \\ 0 & 0 & 0 \\ 0 & 0 & 0 \\ a_1 & a_2 & a_3 \\ 0 & 0 & 0 \end{bmatrix} \middle| a_1, a_2, a_3 \in C(\langle Q \cup I \rangle) \right\} \subseteq M$$

and

$$W_4 = \left\{ \begin{bmatrix} 0 & 0 & 0 \\ 0 & 0 & 0 \\ 0 & 0 & 0 \\ 0 & 0 & 0 \\ 0 & 0 & 0 \\ a_1 & a_2 & a_3 \end{bmatrix} \middle| a_i \in C(\langle Q \cup I \rangle); 1 \le i \le 3 \right\} \subseteq M$$

be set neutrosophic complex linear subalgebras of M over the set S.

Clearly $M = W_1 + W_2 + W_3 + W_4$ and $W_i \cap W_j = (0)$ if $i \ne j$; $1 \le i, j \le 4$. Thus M is the direct sum of sublinear algebras of M. Now we can also define special set linear algebra of complex neutrosophic numbers.

We will give only examples of them and their substructures as it is a matter of routine to define them.

*Example 1.54:* Let

$$V = \left\{ \begin{bmatrix} a_1 & a_{10} \\ a_2 & a_{11} \\ \vdots & \vdots \\ a_9 & a_{18} \end{bmatrix} \middle| a_i \in C(\langle Q \cup I \rangle); 1 \leq i \leq 18 \right\}$$

be a special set linear algebra of complex neutrosophic rationals over the set $S = C(\langle Z \cup I \rangle)$.

*Example 1.55:* Let

$$M = \left\{ \begin{bmatrix} a_1 & a_2 & \ldots & a_{10} \\ a_{11} & a_{12} & \ldots & a_{20} \\ a_{21} & a_{22} & \ldots & a_{30} \end{bmatrix} \middle| a_i \in C(\langle Q \cup I \rangle); 1 \leq i \leq 30 \right\}$$

be a special set neutrosophic complex linear algebra over the set $S = C(\langle 3Z \cup I \rangle) \cup C(\langle 5Z \cup I \rangle)$.

*Example 1.56:* Let

$$V = \left\{ \begin{bmatrix} a_1 & a_2 & \ldots & a_{40} \\ a_{11} & a_{12} & \ldots & a_{80} \end{bmatrix} \middle| a_i \in C(\langle Q \cup I \rangle); 1 \leq i \leq 80 \right\}$$

be a special set complex neutrosophic linear algebra of V over the set $S = C(\langle 3Z \cup I \rangle) \cup C(\langle 5Z \cup I \rangle)$.

Take

$$M = \left\{ \begin{bmatrix} a_1 & a_2 & 0 & \ldots & 0 & a_3 \\ a_4 & a_5 & 0 & \ldots & 0 & a_6 \end{bmatrix} \middle| a_i \in C(\langle Q \cup I \rangle); 1 \leq i \leq 6 \right\}$$

$\subseteq$ V be a special set complex neutrosophic linear subalgebra of V over the set S.

Consider

$$T = \left\{ \begin{bmatrix} a_1 & a_2 & ... & a_{40} \\ 0 & 0 & ... & 0 \end{bmatrix} \middle| a_i \in C(\langle Q \cup I \rangle); 1 \leq i \leq 40 \right\} \subseteq V,$$

be a special set complex neutrosophic linear subalgebra of V over S.

Take

$$A = \left\{ \begin{bmatrix} 0 & 0 & ... & a_{40} \\ a_1 & a_1 & ... & 0 \end{bmatrix} \middle| a_i \in C(\langle Q \cup I \rangle); 1 \leq i \leq 40 \right\} \subseteq V$$

and subset $T = C(\langle 3Z \cup I \rangle) \subseteq S$.
A is a subset special neutrosophic complex sublinear algebra of V over the subset $T = C(\langle 3Z \cup I \rangle) \subseteq S$.

Having seen examples of substructures in case of set special neutrosophic complex linear algebra over the set S, we now proceed onto give examples of semigroup complex neutrosophic linear algebras and their substructures.

*Example 1.57:* Let

$$V = \left\{ \begin{bmatrix} a_1 & a_2 & a_3 \\ a_4 & a_5 & a_6 \\ a_7 & a_8 & a_9 \end{bmatrix} \middle| a_i \in C(\langle Q \cup I \rangle); 1 \leq i \leq 9 \right\}$$

be a semigroup linear algebra of complex neutrosophic numbers over the semigroup $S = Z^+ \cup \{0\}$.

*Example 1.58:* Let P = {all 10 × 10 neutrosophic complex numbers from $C(\langle Q \cup I \rangle)$} be a semigroup neutrosophic complex linear algebra of complex numbers over the semigroup $S = 5Z$.

Clearly V = {all 10 × 10 upper triangular neutrosophic complex numbers with entries from $C(\langle Q \cup I \rangle)$} $\subseteq$ P is a semigroup neutrosophic complex linear subalgebra of P over S.

Take W = {all 10 × 10 diagonal neutrosophic complex matrices with entries from $C(\langle Q \cup I \rangle)$} $\subseteq$ P; W is a subsemigroup complex neutrosophic linear subalgebra of P over the subsemigroup $T = 15Z$, a subsemigroup of $5Z = S$.

*Example 1.58:* Let

$$V = \left\{ \begin{bmatrix} a_1 & a_2 & a_3 & a_4 & a_5 \\ a_6 & a_7 & a_8 & a_9 & a_{10} \\ a_{11} & a_{12} & a_{13} & a_{14} & a_{15} \end{bmatrix} \middle| a_i \in C(\langle Q \cup I \rangle); 1 \le i \le 15 \right\}$$

be a semigroup complex neutrosophic linear algebra over the semigroup $S = 2Z$.

Let

$$M_1 = \left\{ \begin{bmatrix} a_1 & a_2 & a_3 & a_4 & a_5 \\ 0 & 0 & 0 & 0 & 0 \\ a_6 & a_7 & a_8 & a_9 & a_{10} \end{bmatrix} \middle| a_i \in C(\langle Q \cup I \rangle); 1 \le i \le 10 \right\}$$

$\subseteq V$ be a semigroup complex neutrosophic linear subalgebra of V over $S = 2Z$.

$$M_2 = \left\{ \begin{bmatrix} 0 & 0 & 0 & 0 & 0 \\ a_1 & a_2 & 0 & 0 & 0 \\ 0 & 0 & 0 & 0 & 0 \end{bmatrix} \middle| a_1, a_2 \in C(\langle Q \cup I \rangle) \right\} \subseteq V,$$

is a semigroup complex neutrosophic linear subalgebra of V over $S = 2Z$.

$$M_3 = \left\{ \begin{bmatrix} 0 & 0 & 0 & 0 & 0 \\ 0 & 0 & a_1 & a_2 & a_2 \\ 0 & 0 & 0 & 0 & 0 \end{bmatrix} \middle| a_1, a_2 \in C(\langle Q \cup I \rangle) \right\} \subseteq V$$

is a semigroup complex neutrosophic linear subalgebra of V over $S = 2Z$.

We see $V = W_1 + W_2 + W_3$ where $W_i \cap W_j = (0)$, $1 \le i, j \le 3$. Thus V is a direct sum of semigroup linear subalgebras of complex neutrosophic numbers over $S = 2Z$.

Let

$$M_1 = \left\{ \begin{bmatrix} a_1 & a_2 & 0 & 0 & 0 \\ 0 & 0 & 0 & 0 & 0 \\ 0 & 0 & 0 & a_3 & a_4 \end{bmatrix} \middle| a_i \in C(\langle Q \cup I \rangle); 1 \leq i \leq 4 \right\} \subseteq V,$$

$$M_2 = \left\{ \begin{bmatrix} a_1 & 0 & 0 & a_2 & a_3 \\ 0 & a_4 & 0 & 0 & 0 \\ a_5 & 0 & 0 & 0 & a_6 \end{bmatrix} \middle| a_i \in C(\langle Q \cup I \rangle); 1 \leq i \leq 6 \right\} \subseteq V,$$

$$M_3 = \left\{ \begin{bmatrix} a_1 & 0 & 0 & a_2 & a_3 \\ a_4 & a_5 & 0 & 0 & 0 \\ a_6 & a_7 & 0 & 0 & a_7 \end{bmatrix} \middle| a_i \in C(\langle Q \cup I \rangle); 1 \leq i \leq 7 \right\} \subseteq V,$$

$$M_4 = \left\{ \begin{bmatrix} 0 & a_1 & a_2 & a_3 & 0 \\ a_4 & a_5 & 0 & a_6 & a_7 \\ 0 & 0 & a_8 & 0 & 0 \end{bmatrix} \middle| a_i \in C(\langle Q \cup I \rangle); 1 \leq i \leq 8 \right\} \subseteq V,$$

$$M_5 = \left\{ \begin{bmatrix} a_1 & 0 & 0 & 0 & a_2 \\ 0 & 0 & a_3 & 0 & a_4 \\ a_4 & 0 & a_5 & 0 & a_6 \end{bmatrix} \middle| a_i \in C(\langle Q \cup I \rangle); 1 \leq i \leq 6 \right\} \subseteq V$$

be semigroup complex neutrosophic linear subalgebras of V over the semigroup S.

We see

$$V = \bigcup_{i=1}^{5} W_i$$

but $W_i \cap W_j \neq (0)$; if $i \neq j$, $1 \leq i, j \leq 5$. Thus V is the pseudo direct sum of semigroup complex neutrosophic linear subalgebras of V over S.

*Example 1.60:* Let

$$A = \left\{ \begin{bmatrix} a_1 & a_2 & a_3 & a_4 & a_5 \\ a_6 & a_7 & a_8 & a_9 & a_{10} \\ a_{11} & a_{12} & a_{13} & a_{14} & a_{15} \\ a_{16} & a_{17} & a_{18} & a_{19} & a_{20} \\ a_{21} & a_{22} & a_{23} & a_{24} & a_{25} \end{bmatrix} \middle| a_i \in C(\langle Q \cup I \rangle); 1 \le i \le 25 \right\}$$

be a semigroup neutrosophic complex linear algebra over the semigroup $S = 10Z$.

Take

$$M = \left\{ \begin{bmatrix} a_1 & 0 & a_2 & 0 & a_3 \\ 0 & a_4 & 0 & a_5 & 0 \\ a_6 & 0 & a_7 & 0 & a_8 \\ 0 & a_9 & 0 & a_{10} & 0 \\ a_{11} & 0 & a_{12} & 0 & a_{13} \end{bmatrix} \middle| a_i \in C(\langle Q \cup I \rangle); 1 \le i \le 13 \right\} \subseteq A;$$

is a subsemigroup neutrosophic complex linear subalgebra of A over the subsemigroup $40Z = T$ of $S = 10Z$.

We can have several such subsemigroup complex neutrosophic linear subalgebras of A over $T \subseteq S$.

Now we give examples of special semigroup neutrosophic complex linear algebras and their substructures.

*Example 1.61:* Let

$$V = \left\{ \begin{bmatrix} a_1 & a_2 & a_3 & a_4 \\ a_5 & a_6 & a_7 & a_8 \\ \vdots & \vdots & \vdots & \vdots \\ a_{61} & a_{62} & a_{63} & a_{64} \end{bmatrix} \middle| a_i \in C(\langle Q \cup I \rangle); 1 \le i \le 64 \right\}$$

be a special semigroup complex neutrosophic linear algebra over the semigroup $C(\langle Z \cup I \rangle)$.

*Example 1.62:* Let

$$M = \left\{ \begin{bmatrix} a_1 & a_2 & \cdots & a_{10} \\ a_5 & a_6 & \cdots & a_{20} \\ \vdots & \vdots & \vdots & \vdots \\ a_{51} & a_{52} & \cdots & a_{60} \end{bmatrix} \middle| a_i \in C(\langle Q \cup I \rangle); 1 \le i \le 60 \right\}$$

be a special semigroup neutrosophic complex linear algebra over the semigroup C ($\langle Z \cup I \rangle$).

*Example 1.63:* Let

$$V = \left\{ \begin{bmatrix} a_1 & a_2 & a_3 & a_4 \\ a_5 & a_6 & a_7 & a_8 \\ a_9 & a_{10} & a_{11} & a_{12} \\ a_{13} & a_{14} & a_{15} & a_{16} \end{bmatrix} \middle| a_i \in C(\langle Q \cup I \rangle); 1 \le i \le 16, + \right\}$$

be a special semigroup neutrosophic complex linear algebra over the semigroup S = C ($\langle Q \cup I \rangle$).

We see V is of finite dimension and dimension of V is 16 over S.

Take M = {set of all 4 × 4 upper triangular matrices with entries from C ($\langle Q \cup I \rangle$)} $\subseteq$ V; V is a special semigroup complex neutrosophic linear subalgebra of V over S.

Consider

$$H = \left\{ \begin{bmatrix} a_1 & 0 & 0 & 0 \\ 0 & a_2 & 0 & 0 \\ 0 & 0 & a_3 & 0 \\ 0 & 0 & 0 & a_4 \end{bmatrix} \middle| a_i \in C(\langle Q \cup I \rangle); 1 \le i \le 4 \right\} \subseteq V,$$

H is a special subsemigroup complex neutrosophic linear subalgebra of V over the subsemigroup T = C ($\langle Z \cup I \rangle$) $\subseteq$ C ($\langle Q \cup I \rangle$).

Take

$$W_1 = \left\{ \begin{bmatrix} a_1 & a_2 & 0 & 0 \\ a_3 & a_4 & 0 & 0 \\ 0 & 0 & 0 & 0 \\ 0 & 0 & 0 & 0 \end{bmatrix} \middle| a_1, a_2, a_3, a_4 \in C(\langle Q \cup I \rangle) \right\} \subseteq V,$$

$W_1$ is a special semigroup of complex neutrosophic linear subalgebra of V over C ($\langle Q \cup I \rangle$).

$$W_2 = \left\{ \begin{pmatrix} 0 & 0 & a_1 & a_2 \\ 0 & 0 & 0 & 0 \\ 0 & 0 & 0 & 0 \\ 0 & 0 & 0 & 0 \end{pmatrix} \middle| a_i \in C(\langle Q \cup I \rangle), 1 \leq i \leq 2 \right\} \subseteq V$$

is also a sublinear algebra of V.

$$W_3 = \left\{ \begin{pmatrix} 0 & 0 & 0 & 0 \\ 0 & 0 & 0 & 0 \\ a_1 & a_2 & 0 & 0 \\ a_3 & a_4 & 0 & 0 \end{pmatrix} \middle| a_i \in C(\langle Q \cup I \rangle), 1 \leq i \leq 4 \right\} \subseteq V$$

is also a special semigroup complex neutrosophic linear subalgebra of V over C ($\langle Q \cup I \rangle$).

Let

$$W_4 = \left\{ \begin{pmatrix} 0 & 0 & 0 & 0 \\ 0 & 0 & a_1 & a_2 \\ 0 & 0 & 0 & 0 \\ 0 & 0 & 0 & 0 \end{pmatrix} \middle| a_1, a_2 \in C(\langle Q \cup I \rangle) \right\} \subseteq V;$$

be a special semigroup neutrosophic complex linear subalgebra of V over the semigroup S.

$$W_5 = \left\{ \begin{pmatrix} 0 & 0 & 0 & 0 \\ 0 & 0 & 0 & 0 \\ 0 & 0 & a_1 & a_2 \\ 0 & 0 & 0 & 0 \end{pmatrix} \middle| a_1, a_2 \in C(\langle Q \cup I \rangle) \right\} \subseteq V$$

be a special semigroup neutrosophic complex linear subalgebra of V over the semigroup S.

Clearly $W_1 + W_2 + W_3 + W_4 + W_5 \neq V$ with $W_i \cap W_j = (0)$; $i \neq j$; $1 \leq i, j \leq 5$, but still V is not a direct sum of subspaces. On the other hand suppose we add the special semigroup complex neutrosophic linear subalgebra.

$$W_6 = \left\{ \begin{pmatrix} 0 & 0 & 0 & 0 \\ 0 & 0 & 0 & 0 \\ 0 & 0 & 0 & 0 \\ 0 & 0 & a_1 & a_2 \end{pmatrix} \middle| a_1, a_2 \in C(\langle Q \cup I \rangle) \right\} \subseteq V,$$

then $V = W_1 + W_2 + W_3 + W_4 + W_5 + W_6$ and $W_i \cap W_j = (0)$ if $i \neq j$, $1 \leq i, j \leq 6$.

Thus V is a direct sum of subspaces.

*Example 1.64:* Let

$$V = \left\{ \begin{bmatrix} a_1 & a_2 & a_3 \\ a_4 & a_5 & a_6 \\ a_7 & a_8 & a_9 \end{bmatrix} \middle| a_i \in C(\langle Q \cup I \rangle); 1 \leq i \leq 9 \right\}$$

be a special semigroup neutrosophic complex linear algebra over the neutrosophic complex semigroup $S = C(\langle Q \cup I \rangle)$.

$$W_1 = \left\{ \begin{bmatrix} a_1 & 0 & a_3 \\ 0 & a_2 & 0 \\ 0 & 0 & 0 \end{bmatrix} \middle| a_1, a_2, a_3 \in C(\langle Q \cup I \rangle) \right\} \subseteq V;$$

$$W_2 = \left\{ \begin{bmatrix} a_1 & a_2 & 0 \\ 0 & a_3 & 0 \\ 0 & 0 & a_4 \end{bmatrix} \middle| a_i \in C(\langle Q \cup I \rangle); 1 \leq i \leq 4 \right\} \subseteq V,$$

$$W_3 = \left\{ \begin{bmatrix} a_1 & 0 & a_3 \\ a_4 & 0 & a_5 \\ 0 & 0 & a_2 \end{bmatrix} \middle| a_i \in C(\langle Q \cup I \rangle); 1 \leq i \leq 5 \right\} \subseteq V$$

and $\quad W_4 = \left\{ \begin{bmatrix} a_1 & 0 & 0 \\ 0 & a_3 & a_2 \\ a_4 & a_5 & a_6 \end{bmatrix} \middle| a_i \in C(\langle Q \cup I \rangle); 1 \leq i \leq 6 \right\} \subseteq V$

be the collection of special semigroup complex neutrosophic linear subalgebra of V. We see

$$V = \bigcup_{i=1}^{4} W_i$$

and $W_i \cap W_j \neq (0)$; if $i \neq j$, $1 \leq i, j \leq 4$, thus V is only a pseudo direct sum of subspaces of V.

We can as in case of semigroup vector spaces define the notion of linear transformations of special semigroup complex neutrosophic linear algebras only if both these linear algebras are defined over the same complex neutrosophic semigroup. We can also define linear operator of special semigroup complex neutrosophic linear algebras over the semigroup of complex neutrosophic numbers.

Further the notion of basis and dimension can also be defined. Now we proceed onto define the notion of group neutrosophic complex vector space and other related concepts.

**DEFINITION 1.16:** *Let V be a set of complex neutrosophic numbers with zero, which is non empty. Let G be a group under addition.*

*We call V to be a group neutrosophic complex number vector space over G if the following conditions are true.*

i. For every $v \in V$ and $g \in V$ $gv$ and $vg$ are in $V$.
ii. $0.v = 0$ for every $v \in V$, $0$ is the additive identity of $G$.

We give examples of them.

*Example 1.65:* Let

$$V = \left\{ \begin{bmatrix} a & b \\ c & d \end{bmatrix}, \begin{bmatrix} a \\ b \\ c \\ d \end{bmatrix}, (a,b,c,d) \,\middle|\, a,b,c,d \in C(\langle Q \cup I \rangle) \right\}$$

be a group complex neutrosophic vector space over the group $G = Z$.

*Example 1.66:* Let

$$V = \left\{ \sum_{i=0}^{21} a_i x^i, \begin{bmatrix} a_1 \\ a_2 \\ a_3 \\ a_4 \end{bmatrix}, \begin{pmatrix} a_1 & a_2 & \dots & a_8 \\ a_9 & a_{10} & \dots & a_{16} \end{pmatrix} \,\middle|\, \begin{array}{l} a_i \in C(\langle Q \cup I \rangle), \\ 0 \leq i \leq 21 \end{array} \right\}$$

be a group complex neutrosophic vector space over the group $G = 3Z$.

*Example 1.67:* Let

$$V = \left\{ \begin{bmatrix} a_1 & a_2 & a_3 & a_4 \\ a_5 & a_6 & a_7 & a_8 \\ \vdots & \vdots & \vdots & \vdots \\ a_{45} & a_{46} & a_{47} & a_{48} \end{bmatrix} \,\middle|\, a_i \in C(\langle Q \cup I \rangle), 1 \leq i \leq 48 \right\}$$

be a group complex neutrosophic vector space over the group $G = Q$.

*Example 1.68:* Let

$$V = \left\{ \sum_{i=0}^{20} a_i x^i, \begin{bmatrix} a_1 \\ a_2 \\ \vdots \\ a_{12} \end{bmatrix}, (a_1, a_2, \ldots, a_{48}) \,\Big|\, a_i \in C(\langle Q \cup I \rangle), 0 \leq i \leq 20 \right\}$$

be a group complex neutrosophic vector space over the group G = Q.

Consider

$$H = \left\{ \sum_{i=0}^{15} a_i x^i, \begin{bmatrix} a_1 \\ a_2 \\ a_3 \\ 0 \\ \vdots \\ 0 \end{bmatrix}, (0, a_1, 0, a_2, 00\ldots 00) \,\Big|\, \begin{matrix} a_i \in C(\langle Q \cup I \rangle), \\ 0 \leq i \leq 15 \end{matrix} \right\}$$

$\subseteq$ V; H is a group complex neutrosophic vector subspace over the group G = Q.

Take

$$P = \left\{ \sum_{i=0}^{12} a_i x^i, (a_1, a_2, a_3, a_4, 0, \ldots, 0) \,\Big|\, a_i \in C(\langle Q \cup I \rangle), 0 \leq i \leq 12 \right\}$$

$\subseteq$ V; P is a group complex neutrosophic vector subspace of V over G.

*Example 1.69:* Let

$$V = \left\{ \begin{bmatrix} a_1 & a_2 & \ldots & a_{10} \\ a_{11} & a_{12} & \ldots & a_{20} \end{bmatrix}, \begin{bmatrix} a_1 & a_2 & a_3 \\ a_4 & a_5 & a_6 \\ a_7 & a_8 & a_9 \end{bmatrix}, \sum_{i=0}^{9} a_i x^i \right.$$

$a_i \in C(\langle Q \cup I \rangle); 0 \leq i \leq 20\}$ be a group complex neutrosophic vector space over the group G = Z.

Consider

$$W_1 = \left\{ \begin{bmatrix} a_1 & a_2 & \cdots & a_{10} \\ a_{11} & a_{12} & \cdots & a_{20} \end{bmatrix} \middle| a_i \in C(\langle Q \cup I \rangle), 1 \le i \le 20 \right\} \subseteq V,$$

$$W_2 = \left\{ \begin{bmatrix} a_1 & a_2 & a_3 \\ a_4 & a_5 & a_6 \\ a_7 & a_8 & a_9 \end{bmatrix} \middle| a_i \in C(\langle Q \cup I \rangle), 1 \le i \le 9 \right\} \subseteq V$$

and

$$W_3 = \left\{ \sum_{i=0}^{9} a_i x^i \middle| a_i \in C(\langle Q \cup I \rangle), 0 \le i \le 9 \right\} \subseteq V$$

be group complex neutrosophic vector subspaces of V over G.

Further $V = W_1 + W_2 \cup W_3$; $W_i \cap W_j \ne (0)$; if $i \ne j$, $1 \le i, j \le 3$. Let

$$B_1 = \left\{ \begin{bmatrix} a_1 & a_2 & \cdots & a_{10} \\ 0 & 0 & \cdots & 0 \end{bmatrix}, \sum_{i=0}^{9} a_i x^i \middle| a_i \in C(\langle Q \cup I \rangle), 0 \le i \le 10 \right\} \subseteq V$$

$$B_2 = \left\{ \begin{bmatrix} a_1 & a_2 & \cdots & a_{10} \\ a_{11} & a_{12} & \cdots & a_{20} \end{bmatrix} \middle| a_i \in C(\langle Q \cup I \rangle), 1 \le i \le 20 \right\} \subseteq V$$

and

$$B_3 = \left\{ \begin{bmatrix} a_1 & a_3 & 0 & \cdots & 0 & a_9 \\ a_4 & a_2 & 0 & \cdots & 0 & a_{10} \end{bmatrix}, \begin{bmatrix} a_1 & a_2 & a_3 \\ a_4 & a_5 & a_6 \\ a_7 & a_8 & a_9 \end{bmatrix}, \sum_{i=0}^{9} a_i x^i \middle| \right.$$

$a_i \in C(\langle Q \cup I \rangle)$; $1 \le i \le 10 \} \subseteq V$ be group complex neutrosophic vector subspaces of V over the group G.

Clearly $V = B_1 \cup B_2 \cup B_3$; $B_i \cap B_j \ne \phi$; $1 \le i, j \le 3$, hence V is a pseudo direct union of subspaces.

We say a group complex neutrosophic vector space V over a group G is said to be a group complex neutrosophic linear algebra over the group G if V is a group under addition. We give examples of group complex neutrosophic linear algebra over group G.

*Example 1.70:* Let
$$V = \left\{ \sum_{i=0}^{28} a_i x^i \,\middle|\, a_i \in C(\langle Q \cup I \rangle), 0 \leq i \leq 28 \right\}$$

be a group complex neutrosophic linear algebra over the group $G = Z$.

*Example 1.71:* Let
$$V = \left\{ \begin{bmatrix} a_1 & a_2 & a_3 \\ a_4 & a_5 & a_6 \end{bmatrix} \,\middle|\, a_i \in C(\langle Q \cup I \rangle), 1 \leq i \leq 6 \right\}$$

be a group neutrosophic complex linear algebra over the group $G = 3Z$.

*Example 1.72:* Let
$$V = \left\{ \begin{bmatrix} a_1 & a_2 & a_3 & a_4 \\ a_5 & a_6 & a_7 & a_8 \\ a_9 & a_{10} & a_{11} & a_{12} \\ a_{13} & a_{14} & a_{15} & a_{16} \end{bmatrix} \,\middle|\, a_i \in C(\langle Q \cup I \rangle), 1 \leq i \leq 16 \right\}$$

be a group neutrosophic complex linear algebra over the group $G = Z$.

$$W = \left\{ \begin{bmatrix} a_1 & 0 & a_2 & 0 \\ 0 & a_3 & 0 & a_4 \\ a_5 & 0 & a_6 & 0 \\ 0 & a_7 & 0 & a_8 \end{bmatrix} \,\middle|\, a_i \in C(\langle Q \cup I \rangle), 1 \leq i \leq 8 \right\} \subseteq V;$$

be a group neutrosophic complex linear subalgebra of V over the group G = Z.

$$M = \left\{ \begin{bmatrix} a_1 & 0 & 0 & 0 \\ a_2 & a_3 & 0 & 0 \\ a_4 & a_5 & a_6 & 0 \\ a_7 & a_8 & a_9 & a_{10} \end{bmatrix} \middle| a_i \in C(\langle Q \cup I \rangle), 1 \leq i \leq 10 \right\} \subseteq V$$

be a group neutrosophic complex neutrosophic complex linear subalgebra of V over the group G = Z.

*Example 1.73:* Let

$$V = \left\{ \begin{bmatrix} a_1 & a_2 & a_3 \\ a_4 & a_5 & a_6 \\ \vdots & \vdots & \vdots \\ a_{25} & a_{26} & a_{27} \end{bmatrix} \middle| a_i \in C(\langle Q \cup I \rangle), 1 \leq i \leq 27 \right\}$$

be a group complex neutrosophic linear algebra over the group G = Z.

Take

$$P_1 = \left\{ \begin{bmatrix} a_1 & a_2 & a_3 \\ 0 & 0 & 0 \\ 0 & 0 & 0 \\ \vdots & \vdots & \vdots \\ 0 & 0 & 0 \end{bmatrix} \middle| a_i \in C(\langle Q \cup I \rangle), 1 \leq i \leq 3 \right\} \subseteq V,$$

$$P_2 = \left\{ \begin{bmatrix} 0 & 0 & 0 \\ a_1 & 0 & 0 \\ \vdots & \vdots & \vdots \\ a_8 & 0 & 0 \end{bmatrix} \middle| a_i \in C(<Q \cup I>), 1 \leq i \leq 8 \right\} \subseteq V,$$

$$P_3 = \left\{ \begin{bmatrix} 0 & 0 & 0 \\ 0 & a_1 & a_2 \\ 0 & a_3 & a_4 \\ 0 & a_5 & a_6 \\ 0 & 0 & 0 \\ \vdots & \vdots & \vdots \\ 0 & 0 & 0 \end{bmatrix} \middle| a_i \in C(\langle Q \cup I \rangle), 1 \le i \le 6 \right\} \subseteq V,$$

$$P_4 = \left\{ \begin{bmatrix} 0 & 0 & 0 \\ 0 & 0 & 0 \\ 0 & 0 & 0 \\ 0 & 0 & 0 \\ 0 & a_1 & 0 \\ \vdots & \vdots & \vdots \\ 0 & a_5 & 0 \end{bmatrix} \middle| a_i \in C(\langle Q \cup I \rangle), 1 \le i \le 5 \right\} \subseteq V$$

and

$$P_5 = \left\{ \begin{bmatrix} 0 & 0 & 0 \\ 0 & 0 & 0 \\ 0 & 0 & 0 \\ 0 & 0 & 0 \\ 0 & 0 & a_1 \\ \vdots & \vdots & \vdots \\ 0 & 0 & a_5 \end{bmatrix} \middle| a_i \in C(\langle Q \cup I \rangle), 1 \le i \le 5 \right\} \subseteq V$$

be group complex neutrosophic linear subalgebras of V over the group G. Now we have $V = P_1 + P_2 + P_3 + P_4 + P_5$ and $P_i \cap P_j = (0)$ if $i \ne j$, $1 \le i, j \le 5$. Thus V is the direct sum of group complex neutrosophic vector subspaces $P_1, P_2, \ldots, P_5$ of V over G.

*Example 1.74:* Let

$$V = \left\{ \begin{bmatrix} a_1 & a_2 & a_3 \\ a_4 & a_5 & a_6 \\ a_7 & a_8 & a_9 \\ a_{10} & a_{11} & a_{12} \end{bmatrix} \middle| a_i \in C(\langle Q \cup I \rangle), 1 \leq i \leq 12 \right\}$$

be a special semigroup linear algebra of neutrosophic complex numbers over the complex neutrosophic semigroup $S = C(\langle Z \cup I \rangle)$. Take

$$P = \left\{ \begin{bmatrix} a_1 & a_2 & a_3 \\ a_4 & a_5 & a_6 \\ a_7 & a_8 & a_9 \\ a_{10} & a_{11} & a_{12} \end{bmatrix} \middle| a_i \in C(\langle Z \cup I \rangle), 1 \leq i \leq 12 \right\} \subseteq V,$$

P is a special semigroup linear subalgebra of complex neutrosophic numbers over the neutrosophic complex semigroup $S = C(\langle Z \cup I \rangle)$. We see P cannot be used to find a direct sum of sublinear algebras however, P can be used in the pseudo direct union of sublinear algebras.

Further we define special subsemigroup pseudo complex neutrosophic linear subalgebra of V.

We say a proper subset T of V is a pseudo special semigroup complex neutrosophic linear subalgebra of V over the pseudo subsemigroup B of S if B is just a semigroup of reals and T is also only reals.

We proceed onto give examples of this situation.

*Example 1.75:* Let

$$V = \left\{ \begin{bmatrix} a_1 & a_2 & a_3 \\ a_4 & a_5 & a_6 \end{bmatrix} \middle| a_i \in C(\langle Z \cup I \rangle), 1 \leq i \leq 6 \right\}$$

be a special semigroup complex neutrosophic linear algebra over the complex neutrosophic semigroup $S = C(\langle Z \cup I \rangle)$.

Let

$$P = \left\{ \begin{bmatrix} a_1 & a_2 & a_3 \\ a_4 & a_5 & a_6 \end{bmatrix} \middle| a_i \in Q, 1 \leq i \leq 6 \right\} \subseteq V,$$

P is not a special semigroup linear subalgebra over the semigroup $S = C(\langle Z \cup I \rangle)$, but P is a semigroup linear algebra over the semigroup $Z = T$. Then P is defined / called as the pseudo special semigroup of complex neutrosophic linear subalgebra over the pseudo subsemigroup $T = Z$.

Infact P is also a pseudo special semigroup complex neutrosophic linear subalgebra over the pseudo special subsemigroup $M = 3Z^+ \cup \{0\}$.
Thus we have infinite number of pseudo special semigroup complex neutrosophic linear subalgebras over the pseudo special subsemigroups N of S.

Now we can also have pseudo special neutrosophic semigroup linear subalgebras over the pseudo neutrosophic subsemigroup B of S.
We will illustrate this situation also by examples.

*Example 1.76:* Let

$$V = \left\{ \begin{bmatrix} a_1 & a_2 \\ a_3 & a_4 \end{bmatrix} \middle| a_i \in C(\langle Q \cup I \rangle), 1 \leq i \leq 4 \right\}$$

be a special complex neutrosophic linear algebra over the complex neutrosophic semigroup $S = C(\langle Z \cup I \rangle)\}$.
Take

$$W = \left\{ \begin{bmatrix} a_1 & a_2 \\ a_3 & a_4 \end{bmatrix} \middle| a_i \in (\langle Q \cup I \rangle), 1 \leq i \leq 4 \right\} \subseteq V$$

be the pseudo neutrosophic subsemigroup special complex neutrosophic linear subalgebra of V over the pseudo neutrosophic subsemigroup $B = \langle Z \cup I \rangle$ of S.

Now consider

$$A = \left\{ \begin{bmatrix} a_1 & a_2 \\ a_3 & a_4 \end{bmatrix} \middle| a_i \in C(\langle Z \cup I \rangle), 1 \le i \le 4 \right\}$$

be a pseudo special subsemigroup complex neutrosophic linear subalgebra of V over the pseudo neutrosophic subsemigroup $\langle 3Z \cup I \rangle$ of the neutrosophic complex semigroup S.

Consider

$$N = \left\{ \begin{bmatrix} a_1 & a_2 \\ a_3 & a_4 \end{bmatrix} \middle| a_i \in (\langle Z \cup I \rangle), 1 \le i \le 4 \right\} \subseteq V,$$

N is a pseudo special subsemigroup linear subalgebra of neutrosophic complex numbers of V over the neutrosophic subsemigroup $T = \{3Z \cup I\} \subseteq S$.

We can have several such examples. The following theorem is sufficient to prove these.

**THEOREM 1.11:** *Let V be a special semigroup neutrosophic complex linear algebra over the complex neutrosophic semigroup S.*
1. *V has pseudo special neutrosophic subsemigroup complex neutrosophic linear subalgebras over the pseudo neutrosophic subsemigroup of S.*
2. *V has pseudo special ordinary subsemigroup complex neutrosophic linear subalgebra over the pseudo real subsemigroup of S.*

*Example 1.77:* Let

$$V = \left\{ \begin{bmatrix} a_1 & a_2 \\ a_3 & a_4 \\ \vdots & \vdots \\ a_{15} & a_{16} \end{bmatrix} \middle| a_i \in C(\langle Q \cup I \rangle), 1 \le i \le 16 \right\}$$

be a special semigroup complex neutrosophic linear algebra over the semigroup $S = C(\langle Q \cup I \rangle)$. Clearly dimension of V

over S is 16; however if S is replaced by T = C ($\langle Z \cup I \rangle$) $\subseteq$ S then the dimension of V over T is infinite.

So we can have special subsemigroup linear algebras to be of infinite dimension even when the dimension of the special semigroup linear algebra is finite dimension over the semigroup S but of infinite dimension over the subsemigroup of S.

*Example 1.78:* Let

$$V = \left\{ \begin{bmatrix} a & a \\ a & a \\ a & a \\ a & a \\ a & a \\ a & a \end{bmatrix} \middle| a \in C(\langle Q \cup I \rangle) \right\}$$

be a special semigroup linear algebra of complex neutrosophic number over C ($\langle Q \cup I \rangle$). Dimension of V is one.

Clearly V has no special semigroup linear subalgebras but V has pseudo ordinary special subsemigroup linear subalgebras and pseudo neutrosophic special subsemigroup linear subalgebras.

For

$$M = \left\{ \begin{bmatrix} a & a \\ a & a \\ a & a \\ a & a \\ a & a \\ a & a \end{bmatrix} \middle| a \in C(\langle Z \cup I \rangle) \right\} \subseteq V,$$

is a pseudo neutrosophic subsemigroup linear subalgebra over the pseudo neutrosophic subsemigroup T = $\langle Z \cup I \rangle$ $\subseteq$ C($\langle Q \cup I \rangle$).

Likewise

$$N = \left\{ \begin{bmatrix} a & a \\ a & a \\ a & a \\ a & a \\ a & a \\ a & a \end{bmatrix} \middle| a \in Z \right\}$$

is a pseudo ordinary subsemigroup neutrosophic complex linear subalgebra over the pseudo ordinary subsemigroup $B = Z \subseteq S$.

Thus we can say, V is a special semigroup neutrosophic complex linear algebra over S to be simple if it has no proper special semigroup neutrosophic complex linear subalgebra over S.

*Example 1.79:* Let

$$V = \left\{ \begin{bmatrix} a & a & a \\ a & a & a \\ a & a & a \end{bmatrix} \middle| a \in C(\langle Q \cup I \rangle) \right\}$$

be a special semigroup neutrosophic complex linear algebra over the semigroup $S = C(\langle Q \cup I \rangle)$. V is simple. However V has special subsemigroup complex neutrosophic linear subalgebras M over the subsemigroup $T = C(\langle Z \cup I \rangle) \subseteq C(\langle Q \cup I \rangle) = S$ where

$$M = \left\{ \begin{bmatrix} a & a & a \\ a & a & a \\ a & a & a \end{bmatrix} \middle| a \in C(\langle Z \cup I \rangle) \right\} \subseteq V.$$

Also we consider

$$N = \left\{ \begin{bmatrix} b & b & b \\ b & b & b \\ b & b & b \end{bmatrix} \middle| a \in C(\langle 3Z \cup I \rangle) \right\} \subseteq V;$$

N is a special subsemigroup neutrosophic complex linear subalgebra of V over the subsemigroup $B = \{C(\langle 2Z \cup I \rangle)\} \subseteq S$.

**DEFINITION 1.17:** *Let V be a ordinary (neutrosophic or special) semigroup complex neutrosophic linear algebra over the ordinary semigroup S (or neutrosophic semigroup or complex neutrosophic semigroup). If on V we can define a product and V is compatible with product (V is a semigroup with respect to another operation apart from addition) then we define V to be a ordinary (neutrosophic or special) semigroup double linear algebra over the ordinary semigroup S (or neutrosophic semigroup or neutrosophic complex semigroup).*

We will illustrate this situation by some simple examples.

*Example 1.80:* Let

$$V = \left\{ \begin{bmatrix} a_1 & a_2 & a_3 \\ a_4 & a_5 & a_6 \\ a_7 & a_8 & a_9 \end{bmatrix} \middle| a_i \in C(\langle Q \cup I \rangle), 1 \leq i \leq 9 \right\}$$

be a ordinary semigroup complex neutrosophic linear algebra over S = Z. V is clearly a ordinary semigroup complex neutrosophic double linear algebra over S where S = Z.

*Example 1.81:* Let

$$P = \left\{ \sum_{i=0}^{\infty} a_i x^i \middle| a_i \in C(\langle Q \cup I \rangle) \right\}$$

be a ordinary semigroup complex neutrosophic double linear algebra over the semigroup $S = 3Z^+ \cup \{0\}$.

We can define substructure, basis, linear operator and linear transformation which is just a matter of routine.

We have the following interesting result.

**THEOREM 1.12:** *Every ordinary semigroup complex neutrosophic double linear algebra over the semigroup S is a ordinary semigroup complex neutrosophic linear algebra over S but however in general a ordinary semigroup complex neutrosophic linear algebra is not a double linear algebra.*

For the latter part of the proof we give an example.

***Example 1.82:*** Let

$$M = \left\{ \begin{bmatrix} a_1 & a_2 \\ a_3 & a_4 \\ a_5 & a_6 \\ a_7 & a_8 \\ a_9 & a_{10} \end{bmatrix} \middle| a_i \in C(\langle Q \cup I \rangle); 1 \le i \le 10 \right\}$$

be an ordinary semigroup complex neutrosophic linear algebra over the semigroup $S = Z^+ \cup \{0\}$. Clearly M is not an ordinary semigroup complex neutrosophic double linear algebra over S.

***Example 1.83:*** Let

$$V = \left\{ \sum_{i=0}^{\infty} a_i x^i \middle| a_i \in C(\langle Q \cup I \rangle) \right\}$$

be an ordinary semigroup complex neutrosophic double linear algebra over the semigroup $S = 3Z^+ \cup \{0\}$.

Consider

$$M = \left\{ \sum_{i=0}^{\infty} a_i x^i \middle| a_i \in C(\langle Z \cup I \rangle) \right\} \subseteq V;$$

M is also an ordinary semigroup complex neutrosophic double linear subalgebra of V over the semigroup S.

Take

$$P = \left\{ \sum_{i=0}^{20} a_i x^i \middle| a_i \in C(\langle Q \cup I \rangle) \right\} \subseteq V;$$

P is only a pseudo ordinary semigroup complex neutrosophic double linear subalgebra of V over S or ordinary semigroup complex neutrosophic pseudo double linear subalgebra of V over S for on P product cannot be defined.

This same concept of double linear algebra can be easily extended to the case of special semigroup complex neutrosophic linear algebras, complex semigroup complex neutrosophic linear algebras and neutrosophic semigroup complex neutrosophic linear algebras. The definition is a matter of routine so we give only examples of them.

*Example 1.84:* Let

$$M = \left\{ \sum_{i=0}^{\infty} a_i x^i \,\middle|\, a_i \in C(\langle Q \cup I \rangle) \right\}$$

be a neutrosophic semigroup neutrosophic complex double linear algebra over the neutrosophic semigroup $S = \langle Z \cup I \rangle$.
Take

$$P = \left\{ \sum_{i=0}^{100} a_i x^i \,\middle|\, a_i \in C(\langle Z \cup I \rangle) \right\} \subseteq M;$$

P is only a neutrosophic semigroup neutrosophic complex linear subalgebra which is not a double linear subalgebra of M. We call T a neutrosophic semigroup neutrosophic complex pseudo double linear subalgebra of M. Thus apart from double linear subalgebras we can also have pseudo double linear subalgebras of M.

*Example 1.85:* Let

$$V = \left\{ \begin{bmatrix} a_1 & a_{11} \\ a_2 & a_{12} \\ \vdots & \vdots \\ a_{10} & a_{20} \end{bmatrix} \,\middle|\, a_i \in C(\langle Q \cup I \rangle), 1 \le i \le 20 \right\}$$

be a complex semigroup neutrosophic complex semigroup neutrosophic complex linear algebra over the complex semigroup $C(Z) = \{a + bi \mid a, b \in Z\} = S$.
Clearly V is not a double linear algebra.

*Example 1.86:* Let

$$S = \left\{ \begin{bmatrix} a_1 & a_2 \\ a_3 & a_4 \end{bmatrix} \,\middle|\, a_i \in C(\langle Q \cup I \rangle), 1 \le i \le 4 \right\}$$

be a complex semigroup complex neutrosophic double linear algebra over the complex semigroup $T = C(Q) = \{a + ib \mid a, b \in Q\}$.

We see
$$X = \left\{ \begin{bmatrix} 0 & a \\ b & 0 \end{bmatrix} \middle| a, b \in C(\langle Q \cup I \rangle) \right\} \subseteq S$$

is only a complex semigroup complex neutrosophic pseudo double linear subalgebra of S over T.

For we see if
$$A = \begin{bmatrix} 0 & a \\ b & 0 \end{bmatrix} \text{ and } B = \begin{bmatrix} 0 & c \\ d & 0 \end{bmatrix}$$

in X then
$$AB = \begin{bmatrix} 0 & a \\ b & 0 \end{bmatrix} \begin{bmatrix} 0 & c \\ d & 0 \end{bmatrix} = \begin{bmatrix} ad & 0 \\ 0 & bc \end{bmatrix} \notin X.$$

Hence X is only a pseudo linear subalgebra of S over T.

Let
$$V = \left\{ \begin{bmatrix} a & b \\ 0 & c \end{bmatrix} \middle| a, b, c \in C(\langle Q \cup I \rangle) \right\} \subseteq S,$$

be again a complex semigroup neutrosophic complex double linear subalgebra of S over the semigroup T.

Let us consider
$$L = \left\{ \begin{bmatrix} a & 0 \\ 0 & b \end{bmatrix} \middle| a, b \in C(\langle Q \cup I \rangle) \right\} \subseteq S,$$

L is again a complex subsemigroup neutrosophic complex double linear subalgebra over the complex subsemigroup $C(Z) = \{a + bi \mid a, b \in Z\} \subseteq C(Q) = T$.

Now having seen neutrosophic and complex double linear algebra we now proceed onto give examples of special semigroup neutrosophic complex double linear algebras.

*Example 1.87:* Let

$$V = \left\{ \begin{bmatrix} a_1 & a_2 & a_3 \\ a_4 & a_5 & a_6 \\ a_7 & a_8 & a_9 \end{bmatrix} \middle| a_i \in C(\langle Q \cup I \rangle), 1 \leq i \leq 9 \right\}$$

be a special complex neutrosophic double linear algebra over the neutrosophic complex semigroup $S = \{C(\langle Q \cup I \rangle)\}$.

Consider

$$A = \left\{ \begin{bmatrix} a_1 & a_2 & a_3 \\ 0 & a_4 & a_5 \\ 0 & 0 & a_6 \end{bmatrix} \middle| a_i \in C(\langle Q \cup I \rangle), 1 \leq i \leq 6 \right\} \subseteq V,$$

A is a special complex neutrosophic double linear subalgebra of V over the complex semigroup $S = \{C(\langle Q \cup I \rangle)\}$.

Take

$$B = \left\{ \begin{bmatrix} 0 & 0 & a_1 \\ 0 & a_2 & 0 \\ a_3 & 0 & 0 \end{bmatrix} \middle| a_i \in C(\langle Q \cup I \rangle), 1 \leq i \leq 3 \right\} \subseteq V;$$

B is a special neutrosophic complex pseudo linear subalgebra of V over S.

For if
$$x = \begin{bmatrix} 0 & 0 & a_1 \\ 0 & a_2 & 0 \\ a_3 & 0 & 0 \end{bmatrix} \text{ and } y = \begin{bmatrix} 0 & 0 & b_1 \\ 0 & b_2 & 0 \\ b_3 & 0 & 0 \end{bmatrix} \in B;$$

$$xy = \begin{bmatrix} 0 & 0 & a_1 \\ 0 & a_2 & 0 \\ a_3 & 0 & 0 \end{bmatrix} \begin{bmatrix} 0 & 0 & b_1 \\ 0 & b_2 & 0 \\ b_3 & 0 & 0 \end{bmatrix}$$

$$= \begin{bmatrix} a_1b_3 & 0 & 0 \\ 0 & a_2b_2 & 0 \\ 0 & 0 & a_1b_3 \end{bmatrix} \notin B$$

so B is only a pseudo linear subalgebra of V over S.

Now consider

$$L = \left\{ \begin{bmatrix} a_1 & 0 & 0 \\ a_2 & a_3 & 0 \\ a_4 & a_5 & a_6 \end{bmatrix} \middle| a_i \in C(\langle Q \cup I \rangle), 1 \le i \le 6 \right\} \subseteq V,$$

be a special subsemigroup complex neutrosophic double linear subalgebra of V over the special subsemigroup $T = C(\langle Z \cup I \rangle) \subseteq S$; of S.

Let

$$P = \left\{ \begin{bmatrix} a_1 & a_2 & a_3 \\ a_4 & a_5 & a_6 \\ a_7 & a_8 & a_9 \end{bmatrix} \middle| a_i \in C(\langle Z \cup I \rangle), 1 \le i \le 9 \right\} \subseteq V$$

be a pseudo special complex subsemigroup complex neutrosophic double linear subalgebra of V over the complex subsemigroup $N = \{(a + ib) \mid a, b \in Z\} \subseteq S$, N a pseudo special subsemigroup of S.

Let

$$C = \left\{ \begin{bmatrix} a_1 & a_2 & a_3 \\ a_4 & a_5 & a_6 \\ a_7 & a_8 & a_9 \end{bmatrix} \middle| a_i \in (\langle Z \cup I \rangle), 1 \le i \le 9 \right\} \subseteq V$$

be a pseudo neutrosophic special subsemigroup of double linear subalgebra of V over the neutrosophic subsemigroup $E = \{\langle Z \cup I \rangle\} \subseteq S$, E is a pseudo neutrosophic special subsemigroup of V over S.

We can define basis, special double linear transformation, double operator and so on as in case of usual semigroup linear algebras in case of semigroup neutrosophic complex linear algebras.

This task is also left as an exercise to the reader.

Now we will proceed onto continue to define and work to group neutrosophic complex vector spaces / linear algebras.

*Example 1.88:* Let

$$V = \left\{ \begin{bmatrix} a_1 & a_2 & a_3 & a_4 \\ a_5 & a_6 & a_7 & a_8 \\ a_9 & a_{10} & a_{11} & a_{12} \end{bmatrix} \middle| a_i \in C(\langle Q \cup I \rangle), 1 \leq i \leq 12 \right\}$$

be a group complex neutrosophic linear algebra over the group $S = Z$.

Let

$$M = \left\{ \begin{bmatrix} a_1 & a_2 & a_3 & a_4 \\ a_5 & a_6 & a_7 & a_8 \\ a_9 & a_{10} & a_{11} & a_{12} \end{bmatrix} \middle| a_i \in C(\langle Z \cup I \rangle), 1 \leq i \leq 12 \right\} \subseteq V$$

be a pseudo semigroup complex neutrosophic linear subalgebra of V over the semigroup $T = Z^+ \cup \{0\} \subseteq Z$.

That is T is a Smarandache special definite group as it has a proper set which is a semigroup. We define a group neutrosophic complex linear algebra V to be a group neutrosophic complex double linear algebra if V is endowed with another operation product.

We will illustrate this situation by some examples.

*Example 1.89:* Let

$$M = \left\{ \begin{bmatrix} a_1 & a_2 & a_3 \\ a_4 & a_5 & a_6 \\ a_7 & a_8 & a_9 \end{bmatrix} \middle| a_i \in C(\langle Q \cup I \rangle), 1 \leq i \leq 9 \right\}$$

be a group neutrosophic complex double linear algebra over the group $G = Q$. We can find substructures of M. Clearly M is also compatible with respect to matrix multiplication.

*Example 1.90:* Let

$$V = \left\{ \begin{bmatrix} a_1 & a_2 & a_3 & a_4 \\ a_5 & a_6 & a_7 & a_8 \\ a_9 & a_{10} & a_{11} & a_{12} \\ a_{13} & a_{14} & a_{15} & a_{16} \end{bmatrix} \middle| a_i \in C(\langle Q \cup I \rangle), 1 \le i \le 16 \right\}$$

be a group neutrosophic complex double linear algebra over the group $G = Q$.

Consider

$$M = \left\{ \begin{bmatrix} a_1 & a_2 & a_3 & a_4 \\ a_5 & a_6 & a_7 & a_8 \\ a_9 & a_{10} & a_{11} & a_{12} \\ a_{13} & a_{14} & a_{15} & a_{16} \end{bmatrix} \middle| a_i \in C(\langle Z \cup I \rangle), 1 \le i \le 16 \right\}$$

$\subseteq$ V, M is a subgroup neutrosophic complex double linear subalgebra of V over the subgroup $H = Z \subseteq Q = G$.

Take

$$N = \left\{ \begin{bmatrix} a_1 & 0 & 0 & 0 \\ 0 & a_2 & 0 & 0 \\ 0 & 0 & a_3 & 0 \\ 0 & 0 & 0 & a_4 \end{bmatrix} \middle| a_i \in C(\langle Q \cup I \rangle), 1 \le i \le 4 \right\} \subseteq V,$$

be a group neutrosophic complex double linear subalgebra of V over the group $G = Q$.

Let

$$L = \left\{ \begin{bmatrix} a_1 & a_2 & 0 & 0 \\ a_3 & a_4 & 0 & 0 \\ 0 & 0 & 0 & 0 \\ 0 & 0 & 0 & 0 \end{bmatrix} \middle| a_i \in C(\langle Q \cup I \rangle), 1 \le i \le 4 \right\} \subseteq V,$$

L is also a group neutrosophic complex double linear subalgebra of V over the group $G = Q$.

Let

$$A = \left\{ \begin{bmatrix} 0 & 0 & a_1 & a_2 \\ 0 & 0 & a_3 & a_4 \\ 0 & 0 & 0 & 0 \\ 0 & 0 & 0 & 0 \end{bmatrix} \middle| a_i \in C(\langle Q \cup I \rangle), 1 \le i \le 4 \right\} \subseteq V$$

be a group complex neutrosophic linear subalgebra of V, infact a double linear subalgebra of V over G.

For if

$$x = \begin{bmatrix} 0 & 0 & a & b \\ 0 & 0 & c & d \\ 0 & 0 & 0 & 0 \\ 0 & 0 & 0 & 0 \end{bmatrix} \text{ and } y = \begin{bmatrix} 0 & 0 & e & f \\ 0 & 0 & g & h \\ 0 & 0 & 0 & 0 \\ 0 & 0 & 0 & 0 \end{bmatrix} \in A$$

then

$$xy = \begin{bmatrix} 0 & 0 & a & b \\ 0 & 0 & c & d \\ 0 & 0 & 0 & 0 \\ 0 & 0 & 0 & 0 \end{bmatrix} \times \begin{bmatrix} 0 & 0 & e & f \\ 0 & 0 & g & h \\ 0 & 0 & 0 & 0 \\ 0 & 0 & 0 & 0 \end{bmatrix} = \begin{bmatrix} 0 & 0 & 0 & 0 \\ 0 & 0 & 0 & 0 \\ 0 & 0 & 0 & 0 \\ 0 & 0 & 0 & 0 \end{bmatrix}.$$

Thus A is a group complex neutrosophic double linear subalgebra of V over G.

Let

$$C = \left\{ \begin{bmatrix} 0 & 0 & 0 & a_1 \\ 0 & 0 & a_2 & 0 \\ 0 & a_3 & 0 & 0 \\ a_4 & 0 & 0 & 0 \end{bmatrix} \middle| a_i \in C(\langle Q \cup I \rangle), 1 \leq i \leq 4 \right\} \subseteq V$$

be group neutrosophic complex pseudo double linear subalgebra of V over the group G. For if

$$x = \begin{bmatrix} 0 & 0 & 0 & a_1 \\ 0 & 0 & a_2 & 0 \\ 0 & a_3 & 0 & 0 \\ a_4 & 0 & 0 & 0 \end{bmatrix} \text{ and } y = \begin{bmatrix} 0 & 0 & 0 & b_1 \\ 0 & 0 & b_2 & 0 \\ 0 & b_3 & 0 & 0 \\ b_4 & 0 & 0 & 0 \end{bmatrix}$$

are in C.
Now consider the product

$$x.y = \begin{bmatrix} 0 & 0 & 0 & a_1 \\ 0 & 0 & a_2 & 0 \\ 0 & a_3 & 0 & 0 \\ a_4 & 0 & 0 & 0 \end{bmatrix} \begin{bmatrix} 0 & 0 & 0 & b_1 \\ 0 & 0 & b_2 & 0 \\ 0 & b_3 & 0 & 0 \\ b_4 & 0 & 0 & 0 \end{bmatrix}$$

$$= \begin{bmatrix} a_1 b_4 & 0 & 0 & 0 \\ 0 & a_2 b_3 & 0 & 0 \\ 0 & 0 & a_3 b_2 & 0 \\ 0 & 0 & 0 & a_4 b_1 \end{bmatrix} \notin C.$$

So $C \subseteq V$ is not a double linear subalgebra, hence C is only a group complex neutrosophic pseudo double linear subalgebra of V over G.

Now having seen examples of them we can proceed on to define basis, linear operator, linear transformation, direct sum, pseudo direct sum as in case of semigroup neutrosophic complex double linear algebras. This task is left as an exercise to the reader.

We can define fuzzy neutrosophic complex groups, semigroups etc in two ways. For if we want to induct complex fuzzy number we see $-1 = i^2$ so if product is to be defined we need to induct $i^2 = -1$ but if we wish to work only with min max then we do in the following way.

Let $C(\langle[0,1] \cup [0,I]\rangle) = \{a + bi + cI + idI$ where $a, b, c, d \in [0, 1]\}$, we define $C(\langle[0,1] \cup [0,I]\rangle)$ to be the fuzzy complex neutrosophic numbers.

We define max or min operation on $C(\langle[0,1] \cup [0,I]\rangle)$ as follows:

If $x = a + bi + cI + idI$ and $y = m + ni + tI + isI$ are in $C(\langle[0,1] \cup [0,I]\rangle)$;
then
$$\min(x, y) = \min(a+bi+cI+idI, m+ni+tI+isI)$$
$$= \min(a, m) + \min(bi, ni) + \min(cI, tI) + \min(idI, isI).$$

It is easily verified min $(x, y)$ is again in $C(\langle[0,1] \cup [0,I]\rangle)$. Consider $x = 0.7 + 0.61i + 0.23I + i(0.08)I$ and $y = 0.9 + 0.23i + 0.193I + i(0.7)I$ in $C(\langle[0,1] \cup [0,I]\rangle)$.

Now min $\{x, y\} = \min\{0.7 + 0.61i + 0.23I + i(0.08)I, 0.9 + 0.23i + 0.193I + I(0.7)I\}$
$= \min\{0.7, 0.9\} + \min\{0.61i, 0.23i\} + \min\{0.23I, 0.193I\} + \min\{i(0.08)I, I(0.7)I\}$
$= 0.7 + 0.23i + 0.23I + i(0.08)I$.

Thus $\{C(\langle[0,1] \cup [0,I]\rangle), \min\}$ is a semigroup. Likewise we can define the operation of max on $C(\langle[0,1] \cup [0,I]\rangle)$ and $\{C(\langle[0,1] \cup [0,I]\rangle) \max\}$ is also a semigroup. These semigroups will be known as fuzzy neutrosophic complex semigroup.

Now we define special fuzzy complex neutrosophic group, semigroup and ring as follows:

**DEFINITION 1.18:** *Let V be a ordinary semigroup neutrosophic complex vector space over the semigroup S. Let $\eta$ be a map from V into $C(\langle[0,1] \cup [0,I]\rangle)$ such that $(V, \eta)$ is a fuzzy semigroup neutrosophic complex vector space or semigroup neutrosophic complex fuzzy vector space [ ].*

Likewise one can define for neutrosophic semigroup neutrosophic complex vector space, complex semigroup neutrosophic complex vector space, special semigroup neutrosophic complex vector space over the neutrosophic semigroup, complex semigroup and complex neutrosophic semigroup respectively the fuzzy analogue. Further we can define for set complex neutrosophic vector space and group complex neutrosophic vector space over a set and group respectively also the fuzzy analogue.

Thus defining these fuzzy notions is a matter of routine and hence left as an exercise to the reader.

We give examples of these situations.

*Example 1.91:* Let

$$V = \left\{ \begin{bmatrix} a_1 & a_2 & a_3 \\ a_4 & a_5 & a_6 \\ a_7 & a_8 & a_9 \end{bmatrix} \middle| a_i \in C(\langle Z \cup I \rangle), 1 \le i \le 9 \right\}$$

be a semigroup complex neutrosophic linear algebra over the semigroup $S = Z$.

Let $\eta : V \to C(\langle [0,1] \cup [0,I] \rangle)$

$$\eta\left(\begin{bmatrix} a_1 & a_2 & a_3 \\ a_4 & a_5 & a_6 \\ a_7 & a_8 & a_9 \end{bmatrix}\right) = \begin{cases} \dfrac{1}{a_i} & a_i \ne 0 \text{ if } a_i \in Z \\ \dfrac{i}{a_i} & \text{if } a_i \in C(Z) \\ & a_i \ne 0 \\ \dfrac{I}{a_i} & \text{if } a_i \in \langle Z \cup I \rangle \\ \dfrac{iI}{a_i} & \text{if } a_i \in C(\langle Z \cup I \rangle) \\ & a_i \ne 0 \end{cases}$$

If $a_i = 0$ then 1. If every $a_i = 0$ then

$$\eta \left[ \begin{pmatrix} 0 & 0 & 0 \\ 0 & 0 & 0 \\ 0 & 0 & 0 \end{pmatrix} \right] = 1.$$

$$\eta(A+B) \geq \min(\eta(x), \eta(y))$$

where

$$A = \begin{pmatrix} a_1 & a_2 & a_3 \\ a_4 & a_5 & a_6 \\ a_7 & a_8 & a_9 \end{pmatrix} \text{ and } B = \begin{pmatrix} b_1 & b_2 & b_3 \\ b_4 & b_5 & b_6 \\ b_7 & b_8 & b_9 \end{pmatrix}$$

for every $A, B \in V$ and $\eta(rx) \geq \eta(x)$, $r \in Z$.

*Example 1.92:* Let

$$V = \left\{ \sum_{i=0}^{9} a_i x^i, \begin{bmatrix} a_1 \\ a_2 \\ a_3 \\ a_4 \end{bmatrix}, (a_1, a_2, \ldots, a_{10}) \middle| a_i \in C(\langle Q \cup I \rangle), 0 \leq i \leq 10 \right\}$$

be a set complex neutrosophic vector space over the set $S = 3Z \cup 2Z \cup 5Z$.

Define $\eta : V \to C(\langle [0,1] \cup [0,I] \rangle)$

$$\eta(a_i) = \begin{cases} \dfrac{1}{a_i} & \text{if } a_i \neq 0 \\ 1 & \text{if } a_i = 0 \\ \dfrac{i}{a_i} & \text{if } ia_i \neq 0 \text{ is complex} \\ \dfrac{I}{a_i} & \text{if } Ia_i \neq 0 \text{ is neutrosophic} \\ \dfrac{iI}{a_i} & \text{if } iIa_i. \end{cases}$$

with this stipulation $(V, \eta)$ is a set fuzzy complex neutrosophic vector space over Z.

Our main criteria is to built the map $\eta$ such that $\eta$ (i) = i, $\eta$ (I) = I and $\eta$ (iI) = iI with this ordinary mapping $\eta : V \to C (\langle[0,1] \cup [0,I]\rangle)$ we get the fuzzy complex neutrosophic structures.

Now other way of defining fuzzy structures of complex neutrosophic elements is as follows.

**DEFINITION 1.19 :** *Let $V = \{(x_1, \ldots, x_n) \mid x_i \in C (\langle[0, 1] \cup [0, I]\rangle), 1 \leq i \leq n\}$; V is a fuzzy complex neutrosophic semigroup with min (or max) function.*
*'or' used in the mutually exclusive sense.*

*Example 1.93:* Let $V = \{(x_1, x_2, x_3) \mid x_i \in C (\langle[0,1] \cup [0,I]\rangle); 1 \leq i \leq 3\}$, V is a fuzzy complex neutrosophic semigroup under max.

*Example 1.94:* Let

$$V = \left\{ \begin{bmatrix} x_1 \\ x_2 \\ \vdots \\ x_{10} \end{bmatrix} \middle| \; x_i \in C (\langle[0,1] \cup [0,I]\rangle); 1 \leq i \leq 10 \right\}$$

be a fuzzy complex neutrosophic semigroup under min function.

*Example 1.95:* Let M = {all $10 \times 3$ fuzzy neutrosophic complex matrices with entries from $C (\langle[0,1] \cup [0,I]\rangle)$}. M under max (or min) is a semigroup.
For if

$$X = \begin{bmatrix} a_1 & a_2 & a_3 \\ a_4 & a_5 & a_6 \\ \vdots & \vdots & \vdots \\ a_{28} & a_{29} & a_{30} \end{bmatrix} \text{ and } Y = \begin{bmatrix} b_1 & b_2 & b_3 \\ b_4 & b_5 & b_6 \\ \vdots & \vdots & \vdots \\ b_{28} & b_{29} & b_{30} \end{bmatrix}$$

be in M,

$$\min(X, Y) = \min\left(\begin{bmatrix} a_1 & a_2 & a_3 \\ a_4 & a_5 & a_6 \\ \vdots & \vdots & \vdots \\ a_{28} & a_{29} & a_{30} \end{bmatrix}, \begin{bmatrix} b_1 & b_2 & b_3 \\ b_4 & b_5 & b_6 \\ \vdots & \vdots & \vdots \\ b_{28} & b_{29} & b_{30} \end{bmatrix}\right)$$

$$= \begin{bmatrix} \min\{a_1,b_1\} & \min\{a_2,b_2\} & \min\{a_3,b_3\} \\ \min\{a_4,b_4\} & \min\{a_5,b_5\} & \min\{a_6,b_6\} \\ \vdots & \vdots & \vdots \\ \min\{a_{28},b_{28}\} & \min\{a_{29},b_{29}\} & \min\{a_{30},b_{30}\} \end{bmatrix}.$$

(If min is replaced by max still we get a semigroup).

*Example 1.96:* Let

$$V = \left\{ \begin{bmatrix} a_1 & a_2 & a_3 \\ a_4 & a_5 & a_6 \\ a_7 & a_8 & a_9 \end{bmatrix} \middle| a_i \in C(\langle[0,1]\cup[0,I]\rangle); 1\le i \le 9 \right\}.$$

V is a semigroup under max (or min).

Now we can define only set semivector space over the set S.

**DEFINITION 1.20:** *Let V be a set of fuzzy complex neutrosophic elements. S = {0, 1} be a set; V is a set semivector space over the set S if*

  i) *For $v \in V$ and $s \in S$ we have $vs = sv \in V$.*
  ii) *$0.v = 0 \in V$.*
     *$1.v = v \in V$,*

*we define V to be a fuzzy neutrosophic complex semivector space over the set S.*

We will give examples of set fuzzy neutrosophic complex semivector space over the set S.

**Example 1.97:** Let $V = \{(x_1, x_2, x_3, x_4),$

$$\begin{bmatrix} x_1 & x_2 & x_3 \\ x_4 & x_5 & x_6 \\ x_7 & x_8 & x_9 \end{bmatrix}, \begin{bmatrix} x_1 \\ x_2 \\ x_3 \\ x_4 \\ x_5 \end{bmatrix}, \sum_{i=0}^{6} a_i x^i \mid x_i, a_j \in C\ (\langle[0,1] \cup [0,I]\rangle),$$

$1 \leq i \leq 9;\ 0 \leq j \leq 6\}$ be a set of fuzzy neutrosophic complex elements. V is a set fuzzy neutrosophic complex semivector space over the set $S = \{0,1\}$.

**Example 1.98:** Let $V = \{(x_1, x_2, \ldots, x_9),$

$$\begin{bmatrix} x_1 & x_2 \\ x_3 & x_4 \\ \vdots & \vdots \\ x_{19} & x_{20} \end{bmatrix}, \begin{bmatrix} a_1 & a_2 & a_3 & a_4 \\ a_5 & a_6 & a_7 & a_8 \\ a_9 & a_{10} & a_{11} & a_{12} \\ a_{13} & a_{14} & a_{15} & a_{16} \end{bmatrix} \mid x_i, a_j \in C\ (\langle[0,1] \cup [0,I]\rangle),$$

$1 \leq i \leq 20;\ 1 \leq j \leq 16\}$ be a set fuzzy complex neutrosophic semivector space over the set $S = \{0, 1\}$.

**Example 1.99:** Let

$$V = \left\{ \begin{bmatrix} a_1 & a_2 & a_3 & a_4 & a_5 \\ a_6 & a_7 & a_8 & a_9 & a_{10} \\ \vdots & \vdots & \vdots & \vdots & \vdots \\ a_{56} & a_{57} & a_{58} & a_{59} & a_{60} \end{bmatrix} \right.$$

where $a_i \in C\ (\langle[0,1] \cup [0,I]\rangle),\ 1 \leq i \leq 60\}$ be semigroup under max operation V is a set fuzzy complex neutrosophic semilinear algebra over the set $S = \{0, 1\}$.

*Example 1.100:* Let

$$M = \left\{ \begin{bmatrix} a_1 & a_2 & a_3 & a_4 \\ a_5 & a_6 & a_7 & a_8 \\ \vdots & \vdots & \vdots & \vdots \\ a_{77} & a_{78} & a_{79} & a_{80} \end{bmatrix} \right.$$

| $a_i \in C(\langle[0,1] \cup [0,I]\rangle)$, $1 \leq i \leq 80$} be a set fuzzy neutrosophic complex semi linear algebra over the set $X = \{0, 1\}$.

Now we can define substructures, basis, linear transformation and linear operator for these set fuzzy complex neutrosophic semivector spaces / semilinear algebras over the set S.
 This task is left as an exercise to the reader.

**DEFINITION 1.21:** *Let V be a set with elements from C ($\langle[0,1] \cup [0,I]\rangle$) and S be a semigroup with min or max or product. We define V to be a semigroup fuzzy neutrosophic complex semivector space over S if the following conditions hold good;*
 *i)   $s.v \in V$ for all $v \in V$ and $s \in S$.*
 *ii)  $0.v = 0 \in V$ for all $v \in V$ and $0 \in S$.*

*Example 1.101:* Let

$$V = \left\{ \begin{bmatrix} a_1 & a_2 \\ a_3 & a_4 \end{bmatrix}, (a_1, a_2, \ldots, a_9), \begin{bmatrix} a_1 \\ a_2 \\ \vdots \\ a_{20} \end{bmatrix} \middle| a_i \in C(\langle[0,1] \cup [0,I]\rangle); 1 \leq i \leq 20 \right\}$$

be a semigroup fuzzy complex neutrosophic semivector space over the semigroup $S = \{0, 1\}$ under multiplication.
 Results in this direction can be got without any difficulty. Thus one can define structures on C ($\langle[0,1] \cup [0,I]\rangle$) but the scope is limited. However if we consider a map $\eta : V \to$ C ($\langle[0,1] \cup [0,I]\rangle$) we can have almost define all algebraic fuzzy neutrosophic complex structures without any difficulty.

**Chapter Two**

# FINITE COMPLEX NUMBERS

We know i is a complex number where $i^2 = -1$ or $\sqrt{-1} = i$ where $-1$ is from the set of reals. However here we define finite complex numbers in modulo integers which are in a finite set up. Through out this book $Z_n$ will denote the set of modulo integers $\{0, 1, 2, \ldots, n-1\}$ and $-1$ is $(n-1)$, $-2 = (n-2)$ and so on. Thus we use the fact $-1 = n-1$ to define the finite complex modulo numbers.

**DEFINITION 2.1:** *Let $C(Z_n) = \{a + bi_F \mid a, b \in Z_n, i_F$ is the finite complex modulo number such that $i_F^2 = n-1, n < \infty\}$ we define $i_F$ as the finite complex modulo number. $C(Z_n)$ is the finite complex modulo integer numbers.*

It is interesting to note that since finite values in $Z_n$ are dependent on n so also the finite complex number is also dependent on $Z_n$ for every n.

We give examples of them.

*Example 2.1:* Let $C(Z_2) = \{i_F, 1, 0, 1 + i_F\}$. We see $i_F^2 = -1 = 2-1 = 1$. Also $(i_F + 1)^2 = 1 + i_F^2 + 2i_F = 1+1 = 0$. $C(Z_2)$ is a ring.

***Example 2.2:*** Let C ($Z_3$) = {a + b$i_F$ | a, b ∈ $Z_3$} = {0, 1, $i_F$, 2, 2$i_F$, 1 + $i_F$, 2+$i_F$, 2 + 2$i_F$ | $i_F^2$ = –1 = 3 – 1 = 2} $(2i_F)^2$ = –4 = –1 = 2

$$(i_F + 1)^2 = 1 + i_F^2 + 2i_F$$
$$= 1 + 2 + 2i_F = 2i_F \pmod 3.$$
$$(1 + 2i_F)^2 = 1 + (2i_F)^2 + 22i_F \pmod 3$$
$$= 1 + 2 + i_F \pmod 3$$
$$= i_F \pmod 3.$$
$$(2 + i_F)^2 = 4 + (i_F)^2 + 4i_F \pmod 3$$
$$= 4 + 2 + i_F$$
$$= i_F \pmod 3.$$
$$(2i_F + 2)^2 = 4 + 8i_F + (2i_F)^2 \pmod 3$$
$$= 1 + 2i_F + 2 \pmod 3$$
$$= 2i_F \pmod 3.$$

We give the tables associated with C($Z_2$) and C($Z_3$).

Multiplication table C($Z_2$).

| X | 0 | 1 | $i_F$ | 1 + $i_F$ |
|---|---|---|---|---|
| 0 | 0 | 0 | 0 | 0 |
| 1 | 0 | 1 | $i_F$ | 1 + $i_F$ |
| $i_F$ | 0 | $i_F$ | 1 | $i_F$ + 1 |
| 1 + $i_F$ | 0 | $i_F$ + 1 | $i_F$ + 1 | 0 |

Table for C($Z_3$)

| X | 0 | 1 | 2 | $i_F$ | 2$i_F$ | 1+ $i_F$ | 2+ $i_F$ | 1+2$i_F$ | 2+2$i_F$ |
|---|---|---|---|---|---|---|---|---|---|
| 0 | 0 | 0 | 0 | 0 | 0 | 0 | 0 | 0 | 0 |
| 1 | 0 | 1 | 2 | $i_F$ | 2$i_F$ | $i_F$+1 | 2+$i_F$ | 1+2$i_F$ | 2+2$i_F$ |
| 2 | 0 | 2 | 1 | 2$i_F$ | $i_F$ | 2+2$i_F$ | 1+2$i_F$ | 2+$i_F$ | $i_F$+1 |
| $i_F$ | 0 | $i_F$ | 2$i_F$ | 2 | 1 | $i_F$+2 | 2$i_F$+2 | $i_F$+1 | 2$i_F$+1 |
| 2$i_F$ | 0 | 2$i_F$ | $i_F$ | 1 | 2 | 2$i_F$+1 | $i_F$+1 | 2$i_F$+2 | $i_F$+2 |
| 1+$i_F$ | 0 | 1+$i_F$ | 2+2$i_F$ | $i_F$+2 | 2$i_F$+1 | 2$i_F$ | 1 | 2 | $i_F$ |
| 2+$i_F$ | 0 | 2+$i_F$ | 1+2$i_F$ | 2$i_F$+2 | $i_F$+1 | 1 | $i_F$ | 2$i_F$ | 2 |
| 1+2$i_F$ | 0 | 1+2$i_F$ | 2+$i_F$ | $i_F$+1 | 2$i_F$+2 | 2 | 2$i_F$ | $i_F$ | 1 |
| 2+2$i_F$ | 0 | 2+2$i_F$ | 1+$i_F$ | 2$i_F$+1 | $i_F$+2 | $i_F$ | 2 | 1 | 2$i_F$ |

We see C($Z_2$) is a ring with zero divisor of characteristic two where as C($Z_3$) is a field of characteristic three and o(C($Z_3$)) = 9.

*Example 2.3:* Let $C(Z_4) = \{a + bi_F \mid a, b \in Z_4\} = \{0, 1, 2, 3, i_F, 2i_F, 3i_F, 1+ i_F, 2+ i_F, 3+i_F, 1+2i_F, 1+3i_F, 2i_F+2, 2i_F+3, 3i_F+2, 3i_F+3\}$ is a complex ring of order 16.

*Example 2.4:* $C(Z_5) = \{a + bi_F \mid a, b \in Z_5\} = \{0, 1, 2, 3, 4, i_F, 2i_F, 3i_F, 4i_F, \ldots, 4 + 4i_F\}$ is only a finite complex ring for $(1 + 2i_F)(2 + i_F) = 0$ is a zero divisor in $C(Z_5)$.

*Example 2.5:* Consider $C(Z_7) = \{a + bi_F \mid a, b \in Z_7\}$ be the finite complex ring $C(Z_7)$ is a field. For take $(a + bi_F)(c + di_F) = 0$ where $a + bi_F, c + di_F \in C(Z_7)$ with $a, b, c, d \in Z_7 \setminus \{0\}$.
$(a + bi_F)(c + di_F) = 0$ implies

$$ac + 6bd = 0 \qquad (1)$$
$$ad + bc = 0 \qquad (2)$$

$(1) \times a + (2) \times b$ gives
$$a^2 c + 6bda = 0$$
$$b^2 c + bda = 0$$
$$c^2 (a^2 + b^2) = 0 \; c \neq 0$$
this forces $a^2 + b^2 = 0$ in $Z_7$. But for no a and b in $Z_7 \setminus (0)$ we have $a^2 + b^2 = 0$ so $C(Z_7)$ is a finite complex field of characteristic seven.

We see $C(Z_{11})$ is again a finite complex field of characteristic eleven $C(Z_{13})$ is not finite complex field only a ring for $9 + 4i_F$ and $4 + 9i_F$ in $C(Z_{13})$ is such that $(9 + 4i_F)(3 + 9i_F) = 0$.

In view of all these we have the following theorem which gurantees when $C(Z_p)$ is not a field.

**THEOREM 2.1:** *Let $C(Z_p)$ be the finite complex number ring. $C(Z_p)$ is not a field if and only if there exists $a, b \in Z_p \setminus \{0\}$ with $a^2 + b^2 \equiv p$ or $a^2 + b^2 \equiv 0 \pmod{p}$ where p is a prime.*

*Proof:* Let $C(Z_p) = \{a + bi_F \mid a, b \in Z_p, i_F^2 = p - 1\}$. To show $C(Z_p)$ is not a field it is enough if we show $C(Z_p)$ has zero divisors. Suppose $C(Z_p)$ has zero divisors say $a + bi_F$ and $c + di_F$ in $C(Z_p)$ ($a, b, c, d \in Z_p \setminus \{0\}$) is such that $(a + bi_F)(c+di_F) = 0$, then $ac + bdi_F^2 + (ad + bc) i_F = 0$ that is $ac + (p–1) bd + (ad + bc) i_F = 0$ this forces

$$ac + (p - 1)bd = 0 \qquad (i)$$
and $$bc + ad = 0 \qquad (ii)$$

(i) × a + (ii) b gives

$$a^2c + (p–1) bad = 0 +$$
$$b^2c + bad = 0$$

(since a, b, c, d are in $Z_p \setminus \{0\}$ and p is a prime, no item is zero) gives $a^2c + b^2c = 0$ that is

$$c(a^2 + b^2) = 0 \text{ as } c \in Z_p \setminus \{0\}$$

and $c^{-1}$ exists as p is a prime.

We see $c(a^2 + b^2) = 0$ is possible only if $a^2 + b^2 \equiv 0$.

Conversely if $a^2 + b^2 = 0$ than we have
$$(a + bi_F)(b + ai_F) = 0$$
for consider,
$$(a + bi_F)(b + ai_F) = ab + (p–1)ab + (a^2+b^2) i_F$$
$$pab + (a^2+b^2)i_F \equiv 0 \pmod{p}$$
as $a, b \in Z_p$ and given $a^2 + b^2 = 0$. Thus $C(Z_p)$ has zero divisors hence $C(Z_p)$ is a ring and not a field.

**THEOREM 2.2:** *Let $C(Z_p)$ be the commutative finite complex ring, p a prime, $C(Z_p)$ is a field if and only if $Z_p$ has no two distinct elements a and b different from zero such that $a^2 + b^2 \equiv 0 \pmod{p}$.*

*Proof:* Follows from the fact that if $C(Z_p)$ has no zero divisors it is a commutative integral domain which is finite hence is a field by [ ].

Now we will derive other properties related with these finite complex rings / fields. Let $C(Z) = \{a + bi \mid a, b \in Z\}$ be the collection of complex numbers.

Clearly $C(Z)$ is a ring.

Consider the ideal generated by $2 + 2i$, denote it by I.

$$\frac{C(Z)}{I} = \{I, 1 + I, i + I, 1 + i + I\}$$

where $(1 + i)^2 + I = 2 + 2i + I = I$ as here $i^2 \equiv 1 \pmod{2}$ and $1 + 1 = 0 \bmod 2$.

Now we denote it by $i_F$ and $\dfrac{C(Z)}{I} \cong C(Z_2)$. Consider $3 + 3I \in C(Z)$. Let J be the ideal generated by $3 + 3I$.

$$\dfrac{C(Z)}{J = \langle 3+3I \rangle} = \{J, 1 + J, 2 + J, i + J, 2i + J, i + 1 + I, 2 + i + J, \\ 2i + 1 + J, 2 + 2i + J\}.$$

Here also $2 + 1 \equiv 0 \pmod 3$ and $2i + 1 \equiv 0 \pmod 3$ further $i^2 = -1 = 2$. We use i for $i_F$ for one can understand from the fact $i_F^2 = n - 1$. So we can say
$$\dfrac{C(Z)}{J} \cong C(Z_3)$$
where $i_F^2 = 2$. Thus we get a relation between
$$C(Z_n) \text{ and } \dfrac{C(Z)}{J = \langle n + ni \rangle}.$$

We have in $C(Z_n)$ zero divisors, idempotents, nilpotent everything depending on n.

We call $C(Z_p)$ when $C(Z_p)$ is a field as the complex Galois field to honour Galois. However for every prime p, $C(Z_p)$ need not in general be a complex Galois field. Here we give some properties about $C(Z_p)$, p a prime or otherwise.

Consider
$$x = a + bi_F \in C(Z_p)$$
then $\quad \overline{x} = a + (n - 1)i_F b \in C(Z_p)$
is defined as the conjugate of x and vice versa. We see

$$\begin{aligned} x \cdot \overline{x} &= (a + i_F b)(a + (n - 1)i_F b) \\ &= a^2 + i_F ab + (n - 1)i_F ab + (n - 1)b^2 i_F^2 \\ &= a^2 + 0 + (n - 1)^2 b^2 \\ &= a^2 + b^2. \end{aligned}$$

Recall if $a^2 + b^2 \equiv 0 \pmod n$ then $C(Z_n)$ has a nontrivial zero divisor. We can add and multiply finite complex number using $i_F^2 = (n - 1)$ if elements are from $C(Z_n)$.

Thus if $x = a + i_F b$ and $y = c + i_F d$ are in $C(Z_n)$ then
$$x + y = (a + i_F b) + (c + i_F d)$$
$$= (a + c) \bmod n + i_F (b + d) \bmod n.$$
Likewise
$$x \cdot y = (a + i_F b)(c + i_F d)$$
$$= ac + i_F bc + i_F da + i_F^2 bd$$
$$= ac + i_F bc + i_F da + (n-1)bd$$
$$= (ac + (n-1)bd) \pmod n + i_F (bc + da) \pmod n.$$
The operation $+$ and $\times$ are commutative and associative.

It is important to note they are modulo n integers and they are not orderable ($n \geq 2$) we can have for $a + ib = x$, $x^{-1}$ exist or need not exist in $C(Z_n)$. For in $C(Z_2)$ we see $x = 1 + i_F \in C(Z_2)$ and $x^{-1}$ does not exist as $(1 + i_F)^2 = 0 \pmod 2$

Consider $x = 3 + i_F 4$ in $C(Z_7)$ we have $y = 6 + i_F 6$ in $C(Z_7)$ such that $xy = 1$.

Consider $xy = (3 + i_F 4)(6 + i_F 6)$
$$= 18 + i_F^2 24 + i_F 18 + 24 i_F$$
$$= 4 + 3 i_F + 4 i_F + 6 \times 3$$
$$= 4 + 7 i_F + 4$$
$$= 8 + 7 i_F$$
$$= 1 \pmod 7 \text{ as } 7 i_F = 0 \bmod 7 \text{ and } 8 \equiv 1 \bmod 7.$$

Other properties of usual complex numbers are not true in case of finite complex modulo numbers.

We can give a graphical representation of complex modulo integers in 3 layers. The inner most layer consists of real modulo integers makes as in the figure.

The outer layer consists of complex modulo numbers, where as the outer most layer is the mixed complex modulo integer we represent the graph or diagram for $C(Z_2)$

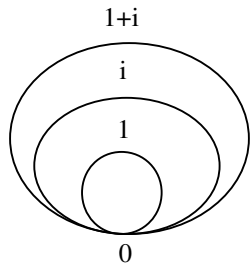

The diagram for C ($Z_3$)

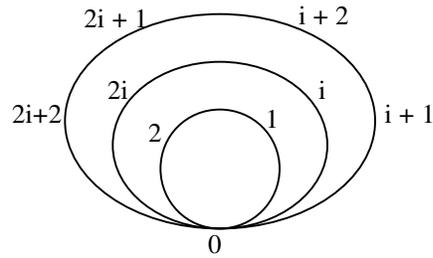

The diagram for C ($Z_4$)

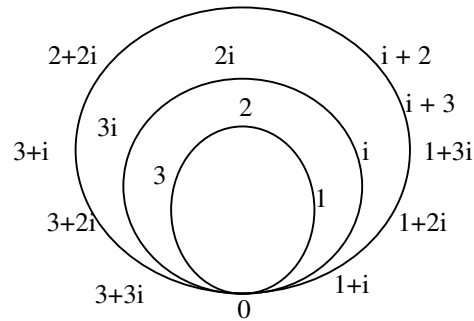

o(C($Z_4$)) = 16 and is a ring we can draw for any C($Z_n$); n ≥ 2.

We cannot do any plane geometry using these finite complex modulo integers but we can get several algebraic results on C($Z_n$); n ≥ 2.

Our main motivation to introduce these concepts is to introduce finite neutrosophic complex modulo number / integers.

Just before we proceed to define these concepts, we define substructures in C($Z_n$); n ≥ 2.

**DEFINITION 2.2:** *Let C ($Z_n$) be the ring of finite complex modulo integer / numbers. Let H $\subseteq$ C($Z_n$) (n ≥ 2) if H itself is a ring under the operations of C($Z_n$). We call H a subring of finite complex modulo integers C($Z_n$).*

We give some examples of them.

**Example 2.6:** Let $C(Z_2)$ be the ring of finite complex modulo integers $C(Z_2) = \{a + bi_F \mid a, b \in Z_2\}$ $P = \{0, 1+i_F \mid i_F^2 = 1\} \subseteq C(Z_2)$ is a subring of $C(Z_2)$.

**Example 2.7:** Let $C(Z_3) = \{a + bi_F \mid a, b \in Z_3\}$, be the ring. $P = \{1, 2, 0\} \subseteq C(Z_3)$ is a subring of $C(Z_3; P = Z_3)$ is the ring of modulo integers three. We see $i_F^2 = 2$, so we cannot find a subring in $C(Z_3)$.

It is pertinent to mention here that for $2 \in Z_2$ is zero $2 \in Z_3$ is such that $2^2 = 1$, $2 \in Z_4$ is such that $2^2 = 0$, $2 \in Z_5$ is such that $2^3 = 1 \pmod{7}$ and so on.

Likewise $i_F$ is finite, $i_F^2 = 1 \pmod 2$, $(i_F)^2 = 2 \pmod 3$, $(i_F)^2 = 3 \pmod 4$ and so on $(i_F)^2 = (n-1) \pmod n$.

**Example 2.8:** Let $C(Z_6) = \{a + bi_F \mid a, b \in Z_6, i_F^2 = 5\}$, $3 + 3i_F \in C(Z_6)$;

$$\begin{aligned}(3 + 3i_F)^2 &= 9 + 9i_F^2 + 2.3^2 i_F \\ &= 3 + 9 \times 5 \\ &= 3 + 3 \times 5 = 3 + 3 \\ &= 0 \pmod 6.\end{aligned}$$

For $1 + i_F \in C(Z_6)$ we have
$$(1+i_F)^2 = 1 + 2i_F + i_F^2 = 1 + 2i_F + 5 = 2i_F.$$

We have an interesting theorem.

**THEOREM 2.3:** *Let $C(Z_{2p})$ be the finite complex ring, p a prime $p > 2$. $C(Z_{2p})$ is a ring and $(1+i_F)^2 = 2i_F$.*

*Proof:* Consider $(1+i_F)^2 = 1 + 2i_F + i_F^2 = 1 + 2i_F + 2p - 1 = 2i_F$ as $2p = 0$.

Now to show $C(Z_{2p})$ is a ring it is enough if we prove the existence of a zero divisor. Take $p + pi_F$ in $C(Z_{2p})$,
$$\begin{aligned}(p + i_F p)^2 &= p^2 (1+i_F)^2 \\ &= p^2 (2i_F) = 2p^2 \cdot i_F \\ &= 2p (pi_F) \equiv 0 \pmod{2p}.\end{aligned}$$
Thus $C(Z_{2p})$ has zero divisors, hence $C(Z_{2p})$ is a ring.

**COROLLARY 2.1:** *Let $C(Z_n)$ be any complex finite ring / field. If $x = 1+i_F \in C(Z_n)$ then $x^2 = 2i_F$.*

*Proof:* Consider $1 + i_F = x \in C(Z_n)$, now
$$(1 + i_F)^2 = 2 + i_F^2 + 2i_F = 1 + n - 1 + 2i_F$$
$$= 2i_F \pmod{n}.$$
Hence the claim.

It is infact a difficult task to find subrings and ideals in $C(Z_n)$ where $o(C(Z_n)) = n^2$.

*Example 2.9:* Let $C(Z_4)$ be the finite complex ring. $I = \{0, 2 + 2i_F\} \subseteq C(Z_4)$ is an ideal of $C(Z_4)$. $P = \{0, 1, 2, 3\} \subseteq C(Z_4)$ is only a subring of $C(Z_4)$ and not an ideal. We see $S = \{0, i_F, 2i_F, 3i_F\} \subseteq C(Z_4)$ S is not even a subring only a set under multiplication as $i_F^2 = 3$.

So in general we cannot get any nice algebraic structure using $\{0, i_F, 2i_F, 3i_F\} = S$; S can only be a group under addition and not even a semigroup under multiplication.

**THEOREM 2.4:** *Let $C(Z_n)$ be the ring, $Z_n \subseteq C(Z_n)$ is a subring of $C(Z_n)$ and not an ideal of $C(Z_n)$.*

Proof is direct and hence left as an exercise to the reader.

**THEOREM 2.5:** *Let $C(Z_{2p})$ be the ring. $P = \{0, p + pi_F\} \subseteq C(Z_{2p})$ is the ideal of $C(Z_{2p})$, p a prime.*

*Proof:* Let $C(Z_{2p})$ be the given ring, p a prime. Consider $S = \{0, p + pi_F\} \subseteq C(Z_{2p})$. Clearly S under addition is an abelian group.
Now
$$(p + pi_F)(a + b\, i_F) = ap + ap\, i_F + pb\, i_F + pb\, (i_F)^2$$
$$= ap + (ap + pb) i_F + (p - 1) pb$$
$$= p(a + (p - 1)b) + p(a + b) i_F.$$
$$= p + p\, i_F \pmod{2p}$$
(using simple number theoretic techniques).

*Example 2.10:* Let $C(Z_{26}) = R$ be a complex ring of modulo integers. $C(Z_{26}) = \{a + i_F b \mid a, b \in Z_{26}\}$. Consider $P = \{0, 13 +$ $13i_F\} \subseteq R$; to show P is an ideal of R. (P, +) is an additive abelian group. To show (P, ×) is a semigroup $(13 + 13\ i_F) \cdot 0 = 0$, $(13 + 13\ i_F)(13 + 13i_F) = 13^2 + 2.13.13i + 13^2 \cdot i_F^2$ ($i_F^2 = 25$)

$$= 13^2 + 13^2 \cdot (25)$$
$$= 13^2 (1+25)$$
$$= 13^2 \cdot 26 \pmod{26} = 0.$$

Thus (P, +, ×) is a ring. Consider $(a + bi_F)$ P; $a, b \in Z_{26}$; if $a = 3$ and $b = 5$ both a and b odd in $Z_{26}$, then

$(3 + 5i_F)(13+13i_F) = 3 \times 13 + 5.13i_F + 3.13i + 5.13 i_F^2 \pmod{26}$

$$= 39 + 104\ i_F + 1625 \pmod{26}$$
$$= 1664 + 104\ i_F \pmod{26}$$
$$= 0 + 0i_F \pmod{26} \in P.$$

Now take $a = 5$ odd and $b = 8$ even. $a + i_F b = 5 + 8i_F$ in R

$(5 + 8i_F)(13 + 13i_F)$

$$= 65 + 8 \times 13i_F + 13\ i_F.5 + 8.13 \times 25 \pmod{26}$$
$$= 65 + 65\ i_F \pmod{26}$$
$$= 13 + 13i_F \pmod{26}$$

Hence $13 + 13i_F$ is in P. Hence the claim.

***Example 2.11:*** Let $R = C(Z_{11})$ be a complex modulo integer ring with $i_F^2 = 10$. Consider $5 + 5i_F$ in R.

$(5+5i_F)^2 = 25 + 2.5^2 \cdot i_F + 5^2 \cdot i_F^2 \pmod{11}$

$$= 3 + 6i_F + 8 = 6i_F.$$

It is interesting to see $(a + ai_F)^2$ is always a $bi_F$ only an imaginary or a complex modulo number further $(6 + 6i_F)^2 = 6i_F$.

We have a nice interesting number theoretic result.

**THEOREM 2.6:** *Let $R = C(Z_p)$, p a prime be the ring of complex modulo integers $i_F^2 = p - 1$, then*

*(i) $(a + ai_F)^2 = bi_F$, $b, a \in Z_p$ and*

*(ii) $\left( \dfrac{p+1}{2} + \dfrac{p+1}{2} i_F \right)^2 = \left( \dfrac{(p+1)}{2} \right) i_F = \dfrac{(p+1)i_F}{2}.$*

The proof uses only simple number theoretic techniques hence left as an exercise to the reader.

***Example 2.12:*** Let $S = C(Z_{20}) = \{a + i_F b \mid a, b \in Z_{20}\}$ be a complex ring of modulo integers.

Consider
$$\begin{aligned}
x &= 8 + 12i_F \in S, x^2 = (8 + 12i_F)^2 \\
&= 8^2 + 144 \times 19 + 2 \times 8 \times 12i_F \pmod{20} \\
&= 4 + 16 + 12i_F \pmod{20} \\
&= 12i_F \pmod{20}.
\end{aligned}$$
Thus if $x = a + (20 - a) i_F \in C(Z_{20})$; then $x^2 = (a + (20 - a)i_F)^2 = a^2 + (20 - a)^2 \times 19 + 2a(20 - a)i_F \pmod{20}$
$$\begin{aligned}
&= a^2 + a^2 \times 19 + (2 \times 20a - 2a^2) i_F \pmod{20} \\
&= a^2(1 + 19) + ci_F \pmod{20} \\
&= 0 + ci_F \pmod{20}; c \in Z_{20}.
\end{aligned}$$

Inview of this example we have an interesting result.

**THEOREM 2.7:** *Let $S = C(Z_n)$ be a complex ring of modulo integers $x = a + (n-a)i_F \in C(Z_n)$ then $x^2 = bi_F$ for some $b \in Z_n$, $i_F^2 = n - 1$.*

Simple number theoretic method yields the solution.

***Example 2.13:*** Let $C(Z_{12}) = S$ be a modulo integer ring. Consider $x = (6 + 6i_F)$ in $S$, $x^2 = 0$ is a zero divisor in $S$. Take $y = 3i_F + 9 \in S$, $y^2 = 9 \times 11 + 81 + 2.3.9i_F = 99 + 81 + 54i_F \pmod{12} = 6i_F$. $a = 4+8i_F$ in $S$.
$a^2 = 16 + 64 \times 11 + 2.4.8\, i_F = 4i_F$ in $S$.

If we consider $b = 2+6i_F$ in $S$, $b^2 = (2+6i_F)^2 = 4 + 36 \times 11 + 2.2.6i_F = 4$ for $2.6 \in Z_{12}$ is such that $2.6 = 0$. If $n = 3+8i_F$ then $n^2 = (3 + 8i_F)^2 = 9 + 64 \times 11 + 2.3.8\, i_F = 5$.

Thus we see every $x = a + bi_F$ in $S$ for which $a.b = 0 \bmod 12$ is such that $x^2 \in Z_{12}$ is a real value.

**THEOREM 2.8:** *Let $R = C(Z_n)$ be a ring of complex modulo integers (n not a prime). Every $x = a + bi_F$ in which $a.b = 0 \pmod{n}$ gives $x^2$ to be a real value.*

Proof follows from simple number theoretic techniques. We see this ring $S = C(Z_n)$ has zero divisors, units and idempotents when n is not a prime. Now we can use these rings and build polynomial rings and matrix rings. We will give examples of them. Let $S = C(Z_n)$ be a complex modulo integer ring.

$S[x] = [C(Z_n)][x]$ is a polynomial ring with complex modulo integer coefficients.

For we see if

$$S[x] = C(Z_3)[x] = \left\{\sum_{i=0}^{\infty} a_i x^i\right.$$

where $a_i = m_i + n_i\, i_F$; $m_i, n_i \in Z_3\}$ we see $S[x]$ is a commutative ring of infinite order but of characteristic three.

If $\quad p(x) = (2 + i_F) + (1 + 2i_F)\, x + (2 + 2\, i_F)x^7$
and $\quad q(x) = i_F + (2 + i_F)x^3 + 2x^6$
are in $S[x]$,

$$\begin{aligned}
p(x)\, q(x) &= [(2 + i_F) + (1 + 2i_F)x + (2 + 2\, i_F)\, x^7]\, [i_F + (2 + i_F)x^3 + 2x^6] \\
&= (2 + i_F)i_F + (2i_F + 1)i_F x + (2 + 2i_F)i_F\, x^7 + (2 + i_F)^2\, x^3 + (1+2i_F)(2+i_F)x^4 + (2+2\, i_F)(2+ i_F)\, x^{10} + (2+2i_F)2x^{13} + (2+i_F)2x^6 + (1+2i_F)\, 2x^7 \\
&= (2i_F+2) + (1 + i_F)x + (1 + 2i_F)x^7 + (1 + 2 + 4i_F)x^3 + (2 + i_F + 1 + i_F)x^4 + (1 + i_F + 2i_F + 1)x^{10} + (1 + i_F)x^{13} + (1 + 2i_F)x^6 + (2 + i_F)x^7 \\
&= (2 + 2i_F) + (1+ i_F)x + i_F\, x^3 + 2\, i_F\, x^4 + 2x^{10} + (1 + i_F)x^{13} + (1 + 2i_F)x^6 + (2 + i_F)x^7.
\end{aligned}$$

This is the way the product of two polynomials is performed.

***Example 2.14:*** Let $S[x] = C(Z_8)[x] = \left\{\sum_{i=0}^{90} a_i x^i \,\middle|\, a_i \in C(Z_8) = \right.$ $\{x + i_F\, y \mid x, y \in Z_8\}$; $0 \le i \le 90\}$ be the set of polynomials in the variable x with coefficients from $C(Z_8)$. $S[x]$ is an additive abelian group of complex modulo integer polynomials of finite order. Clearly $S[x]$ is not closed with respect to product.

Clearly $S[x]$ contains $Z_8[x] = \left\{\sum_{i=0}^{90} a_i x^i \,\middle|\, a_i \in Z_8\right\} \subseteq S[x]$, the usual polynomial group under addition of polynomials. $Z_8 i_F[x] = \left\{\sum_{i=0}^{90} a_i x^i \text{ such that } a_i = x_i\, i_F \text{ with } x_i \in Z_8 \text{ and } i_F^2 = 7\right\}$ $\subseteq$ S[x] is a pure complex modulo integer group of finite order contained in S[x] under addition. Thus S[x] = $Z_8$[x] + $Z_8$ $i_F$[x] with $Z_8$[x] $\cap$ $Z_8$ $i_F$ [x] = {0} is the direct sum of subgroups. Other than these subgroups S[x] has several subgroups.

*Example 2.15:* Let

$$S[x] = C[Z_2][x] = \left\{\sum_{i=0}^{3} a_i x^i \mid a_i = x + y i_F; x, y \in Z_2\right\}$$

= {0, 1, 2, $i_F$, $2i_F$, x, 2x, $i_F$x, $2i_F$x, $x^3$, $2x^3$, $i_F$ $x^3$ and so on} be a complex abelian group under addition of modulo complex integers of finite order.

*Example 2.16:* Let

$$S[x] = \left\{\sum_{i=0}^{\infty} a_i x^i \;\middle|\; a_i \in C(Z_{80}) = \{a + b i_F \mid a, b \in Z_{80}, i_F^2 = 79\}\right\}$$

be the set of polynomials in the variable x with complex modulo integer coefficients. S[x] is a semigroup under multiplication.

We have the following interesting theorems; the proof of which is simple and hence is left as exercise to the reader.

**THEOREM 2.9:** *Let S [x] = $\left\{\sum_{i=0}^{n} a_i x^i \;\middle|\; n < \infty, a_i \in C(Z_n) = \{a + i_F b \mid i_F^2 = n-1 \text{ and } a, b \in Z_n\}\right.$. S [x] is group of complex modulo integer polynomials under addition.*

**THEOREM 2.10:** *Let S [x] = $\left\{\sum_{i=0}^{\infty} a_i x^i \;\middle|\; a_i \in C(Z_n) = \{a + b i_F \mid a, b \in Z_n \text{ and } i_F^2 = n-1\}\right.$.*
  *(1) S [x] is a group under addition of infinite order.*
  *(2) S [x] is a commutative monoid under multiplication.*

*Example 2.17:* Let S[x] = $C(Z_{18})[x]$ = $\left\{\sum_{i=0}^{\infty} a_i x^i \;\middle|\; a_i \in C(Z_{18}) = \{x + y i_F \mid x, y \in Z_{18} \text{ and } i_F^2 = 17\}\right.$ be the collection of all polynomials. S[x] is a ring. S[x] has zero divisors, units and idempotents.

It is interesting to note that these polynomial rings with complex modulo integers can be studied as in case of usual polynomials like the concept of reducibility, irreducibility finding roots etc. This can be treated as a matter of routine and can be done with simple and appropriate changes.

Now we proceed onto study and introduce the notion of matrix ring with complex modulo integer entries. We will illustrate this by some examples.

***Example 2.18:*** Let $V = \{(x_1, x_2, \ldots, x_9) \mid x_i \in C(Z_7) = \{a + bi_F \mid a, b \in Z_7 \text{ and } i_F^2 = 6\}; 1 \le i \le 9\}$ be the $1 \times 9$ row matrix (vector). Clearly V is a group under addition of finite order and V is a semigroup under product of finite order. V as a semigroup has zero divisors. This semigroup $(V, \times)$ is a Smarandache semigroup.

***Example 2.19:*** Let $P = \{(x_1, x_2, x_3, x_4) \mid x_i \in C(Z_{12}) = \{a + bi_F \mid a, b \in Z_{12} \text{ and } i_F^2 = 11, 1 \le i \le 4\}$ be the collection of $1 \times 4$ row matrices with complex modulo integer entries. P is a group under addition of finite order and $(P, \times)$ is a semigroup of finite order having ideals, subsemigroups, zero divisors, and units.

Now we show how addition is carried out. Let
$$x = (7 + 3i_F, 4 + 2i_F, i_F, 7)$$
and $\qquad y = (2, 5i_F, 9 + 2i_F, 10 + i_F)$
be in $(P, +)$.

$x + y = (7 + 3i_F, 4 + 2i_F, i_F, 7) + (2, 5i_F, 9 + 2i_F, 10 + i_F)$
$\quad = (7 + 3i_F + 2, 4 + 2i_F + 5i_F, i_F + 9 + 2i_F + 7, 10 + i_F)$ (mod 12)
$\quad = (9 + 3i_F, 4 + 7i_F, 3i_F + 9, 5 + i_F) \in P$.

$(0, 0, 0, 0)$ acts as the additive identity of P.

Now $(P, \times)$ is a commutative semigroup of finite order; $(1, 1, 1, 1)$ is the multiplicative identity of $(P, \times)$.

Suppose $x = (3 i_F, 2 + i_F, 8, 5+7 i_F)$ and $y = (3+ i_F, 8, 9+ i_F, 8i_F)$ are in P to find
$x.y = (3 i_F, 2+ i_F, 7, 5+7 i_F) \times (3 + i_F, 8, 9 + i_F, 8i_F)$
$\quad = (3i_F \times 3+i_F, 2+i_F \times 8, 8\times(9+ i_F), 5+7i_F \times 8i_F)$ (mod 12)
$\quad = (9i_F + 3\times 11, 16 + 8i_F, 72 + 8i_F, 40i_F + 56 \times 11)$ (mod 12)
$\quad$ (using $i_F^2 = 11$)

= $(9i_F + 9, 4+8i_F, 8i_F, 4+4i_F)$ is in P. Clearly (P, ×) has zero divisors, units, idempotents, ideals and subsemigroups.

Interested reader can study the associated properties of the group of row matrices with complex modulo number entries and semigroup of row matrices of complex modulo integer entries.

Now let

$$Y = \begin{bmatrix} x_1 \\ x_2 \\ \vdots \\ x_m \end{bmatrix}$$

where $x_i$'s are in C $(Z_n)$; $1 \leq i \leq m$, that is $x_i = a_i + b_i i_F$ with $i_F^2 =$ n–1. Y is called / defined as the column matrix or column vector of complex modulo integers.

We see if we get a collection of m × 1 column vectors, that collection is group under addition. Infact multiplication is not defined on that collection.

We will illustrate this situation by some examples.

*Example 2.20:* Let

$$P = \left\{ \begin{bmatrix} x_1 \\ x_2 \\ x_3 \\ x_4 \\ x_5 \end{bmatrix} \right\}$$

$x_i \in$ C $(Z_5) = \{a + bi_F \mid a, b \in Z_5 ; i_F^2 = 4, 1 \leq i \leq 5\}$. P is a group under addition of finite order.

For if

$$x = \begin{bmatrix} 2+4i_F \\ 3i_F \\ 2i_F + 1 \\ 3 \\ 4i_F + 3 \end{bmatrix} \text{ and } y = \begin{bmatrix} 0 \\ 2+4i_F \\ 4 \\ 3+2i_F \\ 1+i_F \end{bmatrix}$$

are in P to find x + y.

$$x + y = \begin{bmatrix} 2+4i_F \\ 3i_F \\ 2i_F+1 \\ 3 \\ 4i_F+3 \end{bmatrix} + \begin{bmatrix} 0 \\ 2+4i_F \\ 4 \\ 3+2i_F \\ 1+i_F \end{bmatrix}$$

$$= \begin{bmatrix} 2+4i_F+0 \\ 3i_F+2+4i_F \\ 1+2i_F+4 \\ 3+3+2i_F \\ 3+4i_F+1+i_F \end{bmatrix} \pmod 5 = \begin{bmatrix} 2+4i_F \\ 2i_F+2 \\ 2i_F \\ 1+2i_F \\ 4 \end{bmatrix}$$

is in P.

Now $\begin{bmatrix} 0 \\ 0 \\ 0 \\ 0 \\ 0 \\ 0 \end{bmatrix}$ acts as the additive identity. Multiplication or product cannot be defined on P.

*Example 2.21:* Let

$$M = \left\{ \begin{bmatrix} x_1 \\ x_2 \\ \vdots \\ x_{20} \end{bmatrix} \right\}$$

$x_i \in C(Z_{24}) = \{a + bi_F \mid a, b \in Z_{24}\}$; $i_F^2 = 23$; $1 \le i \le 20\}$ be a group of $20 \times 1$ complex modulo integer column matrix. M is of finite order and is commutative. Now we proceed to define m × n (m ≠ n) complex modulo integer matrix and work with them.

*Example 2.22:* Let

$$P = \left\{ \begin{bmatrix} x_1 & x_2 & x_3 \\ x_4 & x_5 & x_6 \\ x_7 & x_8 & x_9 \\ x_{10} & x_{11} & x_{12} \\ x_{13} & x_{14} & x_{15} \end{bmatrix} \middle| x_i \in C(Z_{20}) = \{a + bi_F \mid a, b \in Z_{20}\}; \right.$$

$i_F^2 = 19$; $1 \le i \le 15\}$ be a $5 \times 3$ complex modulo integer matrix. P is a group under addition and in P we cannot define product P is defined as a $5 \times 3$ matrix group of complex modulo integers.

*Example 2.23:* Let

$$W = \left\{ \begin{bmatrix} a_1 & a_2 & \ldots & a_{11} \\ a_{12} & a_{13} & \ldots & a_{22} \end{bmatrix} \middle| a_i \in C(Z_{19}) = \{a + bi_F \mid a, b \in Z_{19}\}; \right.$$

$i_F^2 = 18$; $1 \le i \le 22\}$ be a $2 \times 11$ complex modulo integer group of finite order under addition.

$$\begin{bmatrix} 0 & 0 & \ldots & 0 \\ 0 & 0 & \ldots & 0 \end{bmatrix} \in W$$

acts as the additive identity.

*Example 2.24:* Let

$$S = \left\{ \begin{bmatrix} a_1 & a_2 & a_3 & a_4 & a_{17} \\ a_5 & a_6 & a_7 & a_8 & a_{18} \\ a_9 & a_{10} & a_{11} & a_{12} & a_{19} \\ a_{13} & a_{14} & a_{15} & a_{16} & a_{20} \end{bmatrix} \middle| a_i \in C(Z_3) = \{a + bi_F \mid \right.$$

$a, b \in Z_3\}$; $i_F^2 = 2$; $1 \le i \le 20\}$ be a group of $4 \times 5$ complex modulo integer matrices under addition.

Consider

$$P = \left\{ \begin{bmatrix} 0 & 0 & 0 & 0 & 0 \\ a_1 & a_2 & a_3 & a_4 & a_5 \\ 0 & 0 & 0 & 0 & 0 \\ a_6 & a_7 & a_8 & a_9 & a_{10} \end{bmatrix} \middle| a_i \in C(Z_3); 1 \leq i \leq 10 \right\} \subseteq S$$

is a subgroup of S of finite order.

Also

$$T = \left\{ \begin{bmatrix} 0 & a_3 & 0 & a_7 & 0 \\ a_1 & 0 & a_5 & 0 & a_9 \\ a_3 & 0 & a_6 & 0 & a_{10} \\ 0 & a_4 & 0 & a_8 & 0 \end{bmatrix} \middle| a_i \in C(Z_3); 1 \leq i \leq 10 \right\} \subseteq S$$

is a subgroup of S of finite order.

Now

$$P \cap T = \left\{ \begin{bmatrix} 0 & 0 & 0 & 0 & 0 \\ a_1 & 0 & a_2 & 0 & a_3 \\ 0 & 0 & 0 & 0 & 0 \\ 0 & a_4 & 0 & a_5 & 0 \end{bmatrix} \middle| a_i \in C(Z_3); 1 \leq i \leq 5 \right\}$$

is a subgroup of S.

Thus almost all results true in case of general groups can be easily derived in case of group of complex modulo integer matrix groups.

*Example 2.25:* Let

$$P = \left\{ \begin{bmatrix} a_1 & a_2 \\ a_3 & a_4 \end{bmatrix} \middle| a_i \in C(Z_{50}) \right\}$$

$= \{a + bi_F \mid a, b \in Z_{50}\}; i_F^2 = 49; 1 \leq i \leq 4\}$ be the group of complex modulo integers under addition.

$$S = \left\{ \begin{bmatrix} a_1 & a_2 \\ 0 & 0 \end{bmatrix} \middle| a_i \in C(Z_{50}); 1 \leq i \leq 2 \right\} \subseteq P$$

is a subgroup of P of finite order.

*Example 2.26:* Let S = {all 15 × 15 matrices with entries from C ($Z_{42}$)} be the group of complex modulo integers under addition which is commutative and is of finite order.

*Example 2.27:* Let

$$M = \left\{ A = \begin{bmatrix} a_1 & a_2 & a_3 \\ a_4 & a_5 & a_6 \\ a_7 & a_8 & a_9 \end{bmatrix} \middle| a_i \in C(Z_7) = \{a + bi_F \mid a, b \in Z_7\}; \right.$$

$i_F^2 = 6; 1 \leq i \leq 9\}$ be the set of 3 × 3 matrices. If $|A| \neq 0$, M is a group under matrix multiplication. M is a finite non commutative group. If we relax the condition $|A| \neq 0$ then M is only finite non commutative semigroup of complex modulo integers.

M has subgroups also. Also properties of finite groups can be easily extended to the case of complex modulo integers finite groups without any difficulty.

Infact M = {collection of all m × m square matrices with entries from C ($Z_n$) (n < ∞)} is a non commutative complex modulo integer finite ring of characteristic n. Thus we have a nice ring structure on M. So for finite non commutative rings one can easily make use of them. We can have subrings, ideals, zero divisors, units, idempotents and so on as in case of usual rings. This work is also a matter of routine and hence left as an exercise to the reader.

Now having seen examples of rings, semigroups and rings using matrices of complex modulo integers we proceed onto define the notion of complex modulo integers vector spaces and linear algebra.

**DEFINITION 2.3:** *Let V be a complex modulo integer group under addition with entries from $C(Z_p)$; p a prime. If V is a vector space over $Z_p$ then we define V to be a complex modulo integer vector space over $Z_p$.*

We will first illustrate this situation by some examples.

***Example 2.28:*** Let $V = \left\{ \sum_{i=0}^{20} a_i x^i \mid a_i \in C(Z_5) = \{a + bi_F \mid a, b \in Z_5\}; i_F^2 = 4; 0 \le i \le 20 \right\}$ be a complex modulo integer vector space over the field $Z_5$.

***Example 2.29:*** Let

$$V = \left\{ \begin{bmatrix} x_1 & x_2 \\ x_3 & x_4 \\ \vdots & \vdots \\ x_{25} & x_{26} \end{bmatrix} \middle| a_i \in C(Z_{43}) = \{a + bi_F \mid a, b \in Z_{43}\}; \right.$$

$i_F^2 = 42; 0 \le i \le 26\}$ be a complex modulo integer vector space over the field $Z_{43}$.

***Example 2.30:*** Let

$$V = \left\{ \begin{bmatrix} a_1 & a_2 & a_3 & a_4 \\ a_5 & a_6 & a_7 & a_8 \\ a_9 & a_{10} & a_{11} & a_{12} \\ a_{13} & a_{14} & a_{15} & a_{16} \\ a_{17} & a_{18} & a_{19} & a_{20} \end{bmatrix} \right.$$

$a_i \in C(Z_{41}) = \{a + bi_F \mid a, b \in Z_{41}\}; i_F^2 = 40; 1 \le i \le 20\}$ be a complex modulo integer vector space over the field $F = Z_{41}$.

We can define subspaces, basis, direct sum, pseudo direct sum, linear transformation and linear operators, which is a matter of routine so left as an exercise to the reader.
However we proceed onto give some examples.

***Example 2.31:*** Let

$$P = \left\{ \begin{bmatrix} a_1 & a_2 & a_3 \\ a_4 & a_5 & a_6 \\ a_7 & a_8 & a_9 \\ a_{10} & a_{11} & a_{12} \\ a_{13} & a_{14} & a_{15} \end{bmatrix} \right.$$

$a_i \in C(Z_{23}) = \{a + bi_F \mid a, b \in Z_{23}\}$; $i_F^2 = 22$; $1 \leq i \leq 15\}$ be a complex modulo integer vector space over the field $Z_{23} = F$. Consider

$$M = \left\{ \begin{bmatrix} a_1 & a_2 & a_3 \\ 0 & 0 & 0 \\ a_4 & a_5 & a_6 \\ 0 & 0 & 0 \\ a_7 & a_8 & a_9 \end{bmatrix} \middle| a_i \in C(Z_{23}); 1 \leq i \leq 9 \right\} \subseteq P;$$

M is a complex modulo integer vector subspace of P over $F = Z_{23}$.

$$S = \left\{ \begin{bmatrix} a_1 & 0 & a_6 \\ a_2 & 0 & a_7 \\ a_3 & 0 & a_8 \\ a_4 & 0 & a_9 \\ a_5 & 0 & a_{10} \end{bmatrix} \middle| a_i \in C(Z_{23}); 1 \leq i \leq 10 \right\} \subseteq P,$$

is again a complex modulo integer vector subspace of P over $Z_{23} = F$. We see

$$M \cap S = \left\{ \begin{bmatrix} a_1 & 0 & a_2 \\ 0 & 0 & 0 \\ a_4 & 0 & a_3 \\ 0 & 0 & 0 \\ a_5 & 0 & a_6 \end{bmatrix} \middle| a_i \in C(Z_{23}); 1 \leq i \leq 6 \right\} \subseteq P$$

is again a vector subspace of complex modulo integers over the field $Z_{23} = F$.

In view of this we have a nice theorem the proof of which is left as an exercise to the reader.

**THEOREM 2.11:** *Let V be a complex modulo integer vector space over the field $Z_p = F$; p a prime. If $W_1, W_2, \ldots, W_t$ ($t < \infty$) are complex modulo integer vector subspaces of V over $F = Z_p$ then $W = \bigcap_{i=1}^{t} W_i$ is a complex modulo integer vector subspace of V over F.*

*Example 2.32:* Let

$$V = \left\{ \begin{bmatrix} a_1 & a_2 & a_3 \\ a_4 & a_5 & a_6 \\ a_7 & a_8 & a_9 \\ a_{10} & a_{11} & a_{12} \end{bmatrix} \middle| a_i \in C(Z_{11}) = \{a + bi_F \mid a, b \in Z_{11}\}; \right.$$

$i_F^2 = 10; 1 \le i \le 12\}$ be a complex modulo integer vector space over the field $F = Z_{11}$.

Consider

$$W_1 = \left\{ \begin{bmatrix} a_1 & a_2 & 0 \\ 0 & 0 & 0 \\ 0 & 0 & 0 \\ 0 & a_3 & a_4 \end{bmatrix} \middle| a_i \in C(Z_{11}); 1 \le i \le 4 \right\} \subseteq V,$$

$$W_2 = \left\{ \begin{bmatrix} 0 & 0 & a_1 \\ 0 & 0 & 0 \\ 0 & 0 & 0 \\ a_2 & 0 & 0 \end{bmatrix} \middle| a_i \in C(Z_{11}); 1 \le i \le 2 \right\} \subseteq V,$$

$$W_3 = \left\{ \begin{bmatrix} 0 & 0 & 0 \\ a_1 & 0 & a_2 \\ 0 & a_3 & 0 \\ 0 & 0 & 0 \end{bmatrix} \middle| a_i \in C(Z_{11}); 1 \le i \le 3 \right\} \subseteq V$$

and

$$W_4 = \left\{ \begin{bmatrix} 0 & 0 & 0 \\ 0 & a_1 & 0 \\ a_2 & 0 & a_3 \\ 0 & 0 & 0 \end{bmatrix} \middle| a_i \in C(Z_{11}); 1 \le i \le 3 \right\} \subseteq V,$$

be complex modulo integer vector subspaces of V over the field $Z_{11}$. Clearly $V = \bigcup_{i=1}^{4} W_i$ and $W_j \cap W_i = (0)$ if $i \neq j$, $1 \leq i, j \leq 4$.

Thus V is a direct sum or direct union of complex modulo integer subspaces over $Z_{11}$.

Consider

$$P_1 = \left\{ \begin{bmatrix} a_1 & 0 & a_2 \\ 0 & a_3 & 0 \\ 0 & 0 & 0 \\ 0 & 0 & 0 \end{bmatrix} \middle| a_i \in C(Z_{11}); 1 \leq i \leq 3 \right\} \subseteq V,$$

$$P_2 = \left\{ \begin{bmatrix} a_1 & 0 & 0 \\ a_2 & 0 & a_3 \\ 0 & 0 & 0 \\ a_4 & 0 & 0 \end{bmatrix} \middle| a_i \in C(Z_{11}); 1 \leq i \leq 4 \right\} \subseteq V,$$

$$P_3 = \left\{ \begin{bmatrix} a_1 & a_2 & 0 \\ 0 & a_3 & a_4 \\ 0 & 0 & 0 \\ a_5 & 0 & 0 \end{bmatrix} \middle| a_i \in C(Z_{11}); 1 \leq i \leq 5 \right\} \subseteq V,$$

$$P_4 = \left\{ \begin{bmatrix} a_1 & 0 & 0 \\ 0 & 0 & 0 \\ 0 & a_2 & a_3 \\ a_4 & 0 & 0 \end{bmatrix} \middle| a_i \in C(Z_{11}); 1 \leq i \leq 5 \right\} \subseteq V,$$

$$P_5 = \left\{ \begin{bmatrix} a_1 & 0 & 0 \\ 0 & 0 & a_3 \\ a_4 & 0 & 0 \\ 0 & a_2 & 0 \end{bmatrix} \middle| a_i \in C(Z_{11}); 1 \leq i \leq 4 \right\} \subseteq V$$

and

$$P_6 = \left\{ \begin{bmatrix} a_1 & 0 & 0 \\ 0 & 0 & 0 \\ a_4 & 0 & 0 \\ a_2 & a_3 & a_5 \end{bmatrix} \middle| a_i \in C(Z_{11}); 1 \le i \le 4 \right\} \subseteq V$$

be complex modulo integer vector subspaces of V over the complex field $Z_{11}$. Clearly $V = \bigcup_{i=1}^{6} P_i$, but $P_i \cap P_j \ne (0)$ for $i \ne j$; $1 \le i, j \le 6$. Thus V is only a pseudo direct sum of complex modulo integer vector subspaces of V over $Z_{11}$.

***Example 2.33:*** Let $V = \left\{ \begin{bmatrix} a_1 & a_2 & a_3 \\ a_4 & a_5 & a_6 \end{bmatrix} \middle| a_i \in C(Z_5) = \{a + bi_F \mid a, b \in Z_5\}; i_F^2 = 4 \right\}$ be a complex modulo integer vector space over the field $Z_5$. Now consider the set

$$B = \left\{ \begin{bmatrix} 1 & 0 & 0 \\ 0 & 0 & 0 \end{bmatrix}, \begin{bmatrix} 0 & 1 & 0 \\ 0 & 0 & 0 \end{bmatrix}, \begin{bmatrix} 0 & 0 & 1 \\ 0 & 0 & 0 \end{bmatrix}, \right.$$

$$\begin{bmatrix} 0 & 0 & 0 \\ 1 & 0 & 0 \end{bmatrix}, \begin{bmatrix} 0 & 0 & 0 \\ 0 & 1 & 0 \end{bmatrix}, \begin{bmatrix} 0 & 0 & 0 \\ 0 & 0 & 1 \end{bmatrix}, \begin{bmatrix} i_F & 0 & 0 \\ 0 & 0 & 0 \end{bmatrix}, \begin{bmatrix} 0 & 0 & 0 \\ 0 & i_F & 0 \end{bmatrix},$$

$$\left. \begin{bmatrix} 0 & i_F & 0 \\ 0 & 0 & 0 \end{bmatrix}, \begin{bmatrix} 0 & 0 & i_F \\ 0 & 0 & 0 \end{bmatrix}, \begin{bmatrix} 0 & 0 & 0 \\ i_F & 0 & 0 \end{bmatrix}, \begin{bmatrix} 0 & 0 & 0 \\ 0 & 0 & i_F \end{bmatrix} \right\} \subseteq V.$$

B is a basis of V over $Z_5$. Clearly complex dimension of V over $Z_5$ is 12. Likewise we can find basis for V and determine the dimension of V over F.

***Example 2.34:*** Let

$$V = \left\{ \begin{bmatrix} a_1 & a_2 & a_3 & a_4 \\ a_5 & a_6 & a_7 & a_8 \end{bmatrix} \middle| a_i \in C(Z_{13}) = \{a + bi_F \mid a, b \in Z_{11}\}; \right.$$

$i_F^2 = 12\}$ be a complex modulo integer vector space over the field $Z_{13}$
and

$$W = \left\{ \begin{bmatrix} a_1 & a_2 & a_3 \\ a_4 & a_5 & a_6 \\ a_7 & a_8 & a_9 \end{bmatrix} \middle| a_i \in C(Z_{13}) = \{a + bi_F \mid a, b \in Z_{13}\}; \right.$$

$i_F^2 = 12\}$ be a complex modulo integer vector space over the same field $Z_{13}$. Only now we can define the notion of complex linear transformation from V to W.

Let $T : V \to W$ be a map such that

$$T\left( \begin{bmatrix} a_1 & a_2 & a_3 & a_4 \\ a_5 & a_6 & a_7 & a_8 \end{bmatrix} \right) = \begin{bmatrix} a_1 & a_2 & a_3 \\ a_4 & a_5 & a_6 \\ 0 & a_8 & a_7 \end{bmatrix}.$$

T is a complex linear transformation of V to W.

It is pertinent to mention here that most important factor is if T is a complex linear transformation we map i to i. Of course we can have other ways of mapping but those maps should give us the expected complex linear transformation.

We can define

$$\ker T = \left\{ \begin{bmatrix} a_1 & a_2 & a_3 & a_4 \\ a_5 & a_6 & a_7 & a_8 \end{bmatrix} \middle| T\left( \begin{bmatrix} a_1 & a_2 & a_3 & a_4 \\ a_5 & a_6 & a_7 & a_8 \end{bmatrix} \right) = \begin{pmatrix} 0 & 0 & 0 & 0 \\ 0 & 0 & 0 & 0 \end{pmatrix} \right\}$$

$\subseteq V$; ker T is again a subspace of V. We also can easily derive all results about complex linear transformations as in case of usual transformation only with some simple appropriate operations on them.

We can also define as in case of usual vector spaces the notion of invariant subspace by a complex linear operator and so on.

*Example 2.35:* Let

$$V = \left\{ \begin{bmatrix} a_1 & a_2 & a_3 \\ a_4 & a_5 & a_6 \\ a_7 & a_8 & a_9 \\ a_{10} & a_{11} & a_{12} \\ a_{13} & a_{14} & a_{15} \end{bmatrix} \middle| a_i \in C(Z_{29}) = \{a + bi_F \mid a, b \in Z_{29}\}; \right.$$

$i_F^2 = 28; 1 \leq i \leq 15\}$ be a complex modulo integer vector space over the field $Z_{29} = F$.

Let $T : V \to V$ be a map such that

$$T\left( \begin{bmatrix} a_1 & a_2 & a_3 \\ a_4 & a_5 & a_6 \\ a_7 & a_8 & a_9 \\ a_{10} & a_{11} & a_{12} \\ a_{13} & a_{14} & a_{15} \end{bmatrix} \right) = \begin{bmatrix} a_1 & 0 & 0 \\ 0 & a_2 & 0 \\ 0 & 0 & a_3 \\ a_4 & 0 & 0 \\ 0 & a_5 & 0 \end{bmatrix};$$

T is a linear operator on V.

Now let

$$W = \left\{ \begin{bmatrix} a_1 & 0 & 0 \\ 0 & a_2 & 0 \\ 0 & 0 & a_3 \\ a_4 & 0 & 0 \\ 0 & a_5 & 0 \end{bmatrix} \middle| a_i \in C(Z_{29}); 1 \leq i \leq 5 \right\} \subseteq V;$$

W is a complex modulo integer vector subspace of V over $Z_{29}$. Now we see $T(W) \subseteq W$ so W is invariant under the complex linear operator T.

Consider

$$P = \left\{ \begin{bmatrix} a_1 & a_2 & a_3 \\ 0 & 0 & 0 \\ a_4 & a_5 & a_6 \\ 0 & 0 & 0 \\ a_7 & a_8 & a_9 \end{bmatrix} \middle| a_i \in C(Z_{29}); 1 \leq i \leq 9 \right\} \subseteq V,$$

P is a complex modulo integer vector subspace of V over $Z_{29}$ but $T(P) \subseteq P$, so P is not a complex vector subspace of V which is invariant under the linear operator T of V. Interested reader can study the notion of $\text{Hom}_{Z_p}(V,V)$ and $\text{Hom}_{Z_p}(W,V)$; p a prime. Such study is also interesting.

*Example 2.36:* Let $V = \left\{ \begin{bmatrix} a_1 & a_2 \\ a_3 & a_4 \end{bmatrix} \middle| a_i \in C(Z_{17}) = \{a + i_F b \mid a, b \in Z_{179}\}; i_F^2 = 16 \right\}$ be a complex modulo integer vector space over the field $Z_{17}$.

If

$$W_1 = \left\{ \begin{bmatrix} a_1 & 0 \\ 0 & 0 \end{bmatrix} \middle| a_i \in C(Z_{17}) \right\} \subseteq V,$$

$$W_2 = \left\{ \begin{bmatrix} 0 & a_2 \\ 0 & 0 \end{bmatrix} \middle| a_i \in C(Z_{17}) \right\} \subseteq V,$$

$$W_3 = \left\{ \begin{bmatrix} 0 & 0 \\ a_3 & 0 \end{bmatrix} \middle| a_i \in C(Z_{17}) \right\} \subseteq V,$$

and

$$W_4 = \left\{ \begin{bmatrix} 0 & 0 \\ 0 & a_4 \end{bmatrix} \middle| a_i \in C(Z_{17}) \right\} \subseteq V$$

be subspaces of V. We say V is a direct sum of $W_1$, $W_2$, $W_3$ and $W_4$ and V is spanned by $W_1$, $W_2$, $W_3$ and $W_4$.

We can have several properties about linear operators like projections and the related results in case of complex modulo integers.

We can define modulo integer linear algebras if the modulo integer vector space is endowed with a product. Infact all complex modulo integer linear algebras are vector spaces but vector spaces in general are not linear algebras. Several of the examples given prove this. We now give examples of complex modulo integer linear algebras.

***Example 2.37:*** Let $V = \left\{ \sum_{i=0}^{\infty} a_i x^i \mid a_i \in C(Z_{11}) = \{a + bi_F \mid a, b \in Z_{11}; i_F^2 = 10\}\right\}$ be a complex modulo integer linear algebra over $Z_{11}$. Clearly V is of infinite dimension over $Z_{11}$.

***Example 2.38:*** Let $M = \{(x_1, x_2, x_3, x_4, x_5, x_6) \mid x_i \in C(Z_{17}) = \{a + bi_F \mid a, b \in Z_{17}; i_F^2 = 16\}; 1 \le i \le 6\}$ be a complex modulo integer linear algebra over $Z_{17}$. Clearly M is finite dimensional and M is also of finite order.

***Example 2.39:*** Let

$$V = \left\{ \begin{bmatrix} a_1 & a_2 & a_3 \\ a_4 & a_5 & a_6 \\ a_7 & a_8 & a_9 \end{bmatrix} \mid a_i \in C(Z_{53}) = \{a + bi_F \mid a, b \in Z_{53}\};\right.$$

$\left. i_F^2 = 52; 1 \le i \le 9 \right\}$ be a complex modulo linear algebra over $Z_{53}$. V is finite dimensional and has only finite number of elements in it. However V has complex modulo integer linear subalgebras.

***Example 2.40:*** Let

$$V = \left\{ \begin{bmatrix} a_1 & a_2 \\ a_3 & a_4 \end{bmatrix} \mid a_i \in C(Z_{61}) = \{a + bi_F \mid a, b \in Z_{61}; i_F^2 = 60\};\right.$$

$1 \le i \le 4\}$ be a complex modulo integer linear algebra over the field $F = Z_{61}$.

Define $T : V \to V$ by

$$T\left(\begin{bmatrix} a_1 & a_2 \\ a_3 & a_4 \end{bmatrix}\right) = \begin{bmatrix} a_1 & 0 \\ 0 & a_2 \end{bmatrix};$$

T is a linear operator on V.

Further if

$$W = \left\{ \begin{bmatrix} a_1 & 0 \\ a_2 & a_3 \end{bmatrix} \mid a_i \in C(Z_{61}); 1 \le i \le 3 \right\} \subseteq V;$$

be a complex modulo integer linear subalgebra of V over $Z_{61} = F$.

$T(W) \not\subseteq W$ so T does not keep W as a invariant subspace of V over $Z_{61}$.

*Example 2.41:* Let

$$V = \left\{ \begin{bmatrix} a_1 & a_2 \\ a_3 & a_4 \end{bmatrix} \middle| a_i \in C(Z_{13}); 1 \le i \le 4 \right\}$$

and

$$W = \left\{ \begin{bmatrix} a_1 & a_2 & a_3 & a_4 \\ a_5 & a_6 & a_7 & a_8 \\ a_9 & a_{10} & a_{11} & a_{12} \\ a_{13} & a_{14} & a_{15} & a_{16} \end{bmatrix} \middle| a_i \in C(Z_{13}); 1 \le i \le 16 \right\}$$

be any two complex modulo integer linear algebras over the field $F = Z_{13}$.

Define $T : V \to W$ by

$$T(\begin{bmatrix} a_1 & a_2 \\ a_3 & a_4 \end{bmatrix}) = \begin{bmatrix} a_1 & 0 & 0 & 0 \\ 0 & a_2 & 0 & 0 \\ 0 & 0 & a_3 & 0 \\ 0 & 0 & 0 & a_4 \end{bmatrix};$$

T is a linear transformation from V to W.

*Example 2.42:* Let

$$V = \left\{ \begin{bmatrix} a_1 & a_2 & a_3 & a_4 \\ a_5 & a_6 & a_7 & a_8 \\ a_9 & a_{10} & a_{11} & a_{12} \\ a_{13} & a_{14} & a_{15} & a_{16} \end{bmatrix} \middle| a_i \in C(Z_{37}) \right.$$

$= \{a + bi_F \mid a, b \in Z_{37}; i_F^2 = 36\}; 1 \le i \le 16\}$ be a complex modulo integer linear algebra over the field $Z_{37} = F$.

Define T : V → V by

$$T\left(\begin{bmatrix} a_1 & a_2 & a_3 & a_4 \\ a_5 & a_6 & a_7 & a_8 \\ a_9 & a_{10} & a_{11} & a_{12} \\ a_{13} & a_{14} & a_{15} & a_{16} \end{bmatrix}\right) = \begin{bmatrix} a_1 & 0 & a_2 & 0 \\ 0 & a_3 & 0 & a_4 \\ a_5 & 0 & a_6 & 0 \\ 0 & a_7 & 0 & a_8 \end{bmatrix};$$

T is a linear operator on V.

All properties related with usual operators can be easily derived in case of complex modulo integer linear operators.

We can define in case of linear operators of vector space / linear algebra of complex modulo integers we can define the related special characteristic values and characteristic vectors.

Let V be a complex modulo integer vector space over the field $Z_p$ and let T be a linear operator on V. A special characteristic value of T is a scalar c in $Z_p$ such that there is a non zero vector $\alpha$ in V with $T\alpha = c\alpha$. If c of characteristic value of T, then

(a) any $\alpha$ such that $T\alpha = c\alpha$ is called a characteristic vector of T associated with value c.
(b) collection of all $\alpha$ such that $T\alpha = c\alpha$ is called the special characteristic space associated with c.

The fact to be remembered is that in case of usual vector spaces the associated matrix has its entries from the field F over which the space is defined but in case of complex modulo integers this may not be possible in general. We say in case of special characteristic equation the roots may not in general be in the field $Z_p$ but in C ($Z_p$) where the vector space is defined having its entries. Thus this is the major difference between usual vector space and the complex modulo integer vector space.

Now we define the notion of special Smarandache complex modulo integer vector space / linear algebra over the S-ring.

**DEFINITION 2.4:** *Let V be an additive abelian group with entries from {C ($Z_n$) = {a + bi$_F$ | a, b ∈ $Z_n$}; $i_F^2 = n–1$} where $Z_n$ is a S-ring. If V is a vector space over the S-ring then we define*

*V* to be a *Smarandache complex modulo integer vector space over the S-ring* $Z_n$.

We give examples of them.

*Example 2.43:* Let
$$V = \left\{ \begin{bmatrix} a_1 & a_2 & a_3 \\ a_4 & a_5 & a_6 \end{bmatrix} \middle| a_i \in C(Z_{10}) \right.$$
$= \{a + bi_F \mid a, b \in Z_{10}; i_F^2 = 9\}; 1 \le i \le 6\}$ be a Smarandache complex modulo integer vector space over the S-ring $R = Z_{10}$.

*Example 2.44:* Let
$$V = \left\{ \begin{bmatrix} a_1 & a_2 & a_3 & a_4 \\ a_5 & a_6 & a_7 & a_8 \\ \vdots & \vdots & \vdots & \vdots \\ a_{37} & a_{38} & a_{39} & a_{40} \end{bmatrix} \middle| a_i \in C(Z_{14}) \right.$$
$= \{a + bi_F \mid a, b \in Z_{14}\}; i_F^2 = 13\}$ be a S-complex modulo integer vector space over the ring $S = Z_{14}$.

*Example 2.45:* Let
$$V = \left\{ \begin{bmatrix} a_1 & a_2 & a_3 & a_4 \\ a_5 & a_6 & a_7 & a_8 \\ a_9 & a_{10} & a_{11} & a_{12} \\ a_{13} & a_{14} & a_{15} & a_{16} \end{bmatrix} \middle| a_i \in C(Z_{46}) \right.$$
$= \{a + bi_F \mid a, b \in Z_{46}; i_F^2 = 45\}; 1 \le i \le 6\}$ be a S-complex modulo integer linear algebra over the S-ring $S = Z_{46}$.

*Example 2.46:* Let
$$M = \left\{ \begin{bmatrix} a_1 & a_2 & a_3 \\ a_4 & a_5 & a_6 \\ a_7 & a_8 & a_9 \end{bmatrix} \middle| a_i \in C(Z_{38}) \right.$$
$= \{a + bi_F \mid a, b \in Z_{38}; i_F^2 = 37\}; 1 \le i \le 9\}$ be a complex modulo integers linear algebra over $Z_{38}$.

Take
$$W = \left\{ \begin{bmatrix} a_1 & a_2 & a_3 \\ 0 & a_4 & a_5 \\ 0 & 0 & a_6 \end{bmatrix} \middle| a_i \in C(Z_{38}); 1 \le i \le 6 \right\} \subseteq M;$$

W is a Smarandache complex modulo integer linear subalgebra of M over $Z_{38}$. Infact M has many S-complex modulo integer linear subalgebras.

*Example 2.47:* Let $V = \{(x_1, x_2, x_3, x_4, x_5, x_6, x_7, x_8) \mid x_i \in C(Z_{34}) = \{a+bi_F \mid a, b \in Z_{34}; i_F^2 = 33\}; 1 \le i \le 8\}$ be a Smarandache complex modulo integer linear algebra over the S-ring $Z_{34} = S$.

Take
$P_1 = \{(x_1, x_2, 0, 0, 0, 0, 0, 0) \mid x_1, x_2 \in C(Z_{34})\} \subseteq V$,
$P_2 = \{(0, 0, x_3, 0, 0, 0, 0, 0) \mid x_3 \in C(Z_{34})\} \subseteq V$,
$P_3 = \{(0, 0, 0, x_4, x_5, 0, 0, 0) \mid x_5, x_4 \in C(Z_{34})\} \subseteq V$,
$P_4 = \{(0\ 0\ 0\ 0\ 0\ x_6\ 0\ 0) \mid x_6 \in C(Z_{34})\} \subseteq V$,
$P_5 = \{(0\ 0\ 0\ 0\ 0\ 0\ x_1\ 0) \mid x_7 \in C(Z_{34})\} \subseteq V$ and
$P_6 = \{(0\ 0\ \ldots\ 0, x_8) \mid x_8 \in C(Z_{34})\} \subseteq V$

be smarandache complex modulo integer linear subalgebras of V over the S-ring $Z_{34}$. Clearly $V = \bigcup_{i=1}^{6} P_i$ and $P_i \cap P_j = (0\ 0\ 0\ 0\ 0\ 0\ 0\ 0)$ if $i \ne j$; $1 \le i, j \le 6$. Thus V is a direct sum of subspaces.

Take
$M_1 = \{(x_1, x_2, 0, x_3\ 0\ 0\ 0\ 0) \mid x_i \in C(Z_{34}); 1 \le i \le 3\} \subseteq V$,
$M_2 = \{(x_1, x_2, 0, 0, x_3, x_4\ 0\ 0) \mid x_i \in C(Z_{34}); 1 \le i \le 4\} \subseteq V$,
$M_3 = \{(x_1, 0\ 0\ x_2\ 0\ x_3\ 0\ 0) \mid x_i \in C(Z_{34}); 1 \le i \le 3\} \subseteq V$,
$M_4 = \{(x_1, 0\ x_2, 0, x_3\ 0\ x_4) \mid x_i \in C(Z_{34}); 1 \le i \le 4\} \subseteq V$,
$M_5 = \{(x_1, 0\ 0\ 0\ 0\ 0\ x_3\ x_2) \mid x_i \in C(Z_{34}); 1 \le i \le 3\} \subseteq V$ and
$M_6 = \{(x_1, 0, x_2, 0, x_3\ 0\ x_4\ 0) \mid x_i \in C(Z_{34}); 1 \le i \le 4\} \subseteq V$,

be S-complex modulo integer linear subalgebras of V over the S-ring $S = Z_{34}$.
$$V = \bigcup_{i=1}^{6} M_i\ ;\ M_i \cap M_j \ne (0) \text{ if } i \ne j;\ 1 \le i, j \le 6.$$

Thus V is only a pseudo direct union of S-complex modulo integer sublinear algebras of V over $Z_{34} = S$.

*Example 2.48:* Let

$$V = \left\{ \begin{bmatrix} a_1 & a_2 \\ a_3 & a_4 \end{bmatrix} \middle| a_i \in C(Z_6) = \{a + i_F b \mid a, b \in Z_6; i_F^2 = 5\}; \right.$$

$1 \leq i \leq 4\}$ be a Smarandache complex modulo integer linear algebra over the S-ring $Z_6$.

Let

$$W = \left\{ \begin{pmatrix} a_1 & a_2 \\ 0 & a_3 \end{pmatrix} \middle| a_i \in C(Z_6); 1 \leq i \leq 3 \right\} \subseteq V,$$

W is a Smarandache complex modulo integer linear subalgebra of V over the S-ring $Z_6$.

Now let

$$P = \left\{ \begin{pmatrix} a_1 & 0 \\ 0 & a_2 \end{pmatrix} \middle| a_i \in C(Z_6); 1 \leq i \leq 2 \right\} \subseteq V,$$

P is a complex modulo integer linear algebra over the field F = $\{0, 3\} \subseteq Z_6$. Thus P is defined as a Smarandache pseudo complex modulo integer linear algebra over the field F of V.

Infact P is also a Smarandache pseudo complex modulo integer linear subalgebra of V over the field $\{0, 2, 4\} \subseteq Z_6$.

Also

$$T = \left\{ \begin{pmatrix} a & a \\ a & a \end{pmatrix} \middle| a \in C(Z_6) \right\} \subseteq V$$

is a pseudo Smarandache complex modulo integer linear subalgebra of V over the field F = $\{0, 2, 4\} \subseteq Z_6$.

*Example 2.49:* Let

$$V = \left\{ \begin{pmatrix} a & a & a & a \\ a & a & a & a \\ a & a & a & a \\ a & a & a & a \end{pmatrix} \middle| a \in C(Z_7) \right.$$

= $\{a + i_F b \mid a, b \in Z_7; i_F^2 = 6\}\}$ be a complex modulo integer vector space over the field $Z_7$. Clearly V has no subvector space with entries of the form $a = x + i_F y$; $x \neq 0$ and $y \neq 0$.

Inview of this we have the following theorem.

**THEOREM 2.12:** *Let*

$$V = \left\{ \begin{pmatrix} a & a & a & \cdots & a \\ a & a & a & \cdots & a \\ \vdots & \vdots & \vdots & & \vdots \\ a & a & a & \cdots & a \end{pmatrix} \middle| a \in C(Z_p); \right.$$

*p a prime} be the set of all $n \times n$ matrices complex modulo integer linear algebra over $Z_p$.*

*i) V has no complex modulo integer linear subalgebra with entries of the matrix $x = a + bi_F$, $a \neq 0$ and $b \neq 0$ $a, b \in Z_p$.*

*ii)* $M = \left\{ \begin{pmatrix} x & x & \cdots & x \\ x & x & \cdots & x \\ \vdots & \vdots & & \vdots \\ x & x & \cdots & x \end{pmatrix} \middle| x \in Z_p \right\} \subseteq V$ *is a pseudo complex modulo integer linear subalgebra of V over $Z_p$.*

The proof is straight forward and hence left as an exercise to the reader.

*Example 2.50:* Let

$$V = \left\{ \begin{pmatrix} a & a \\ a & a \end{pmatrix} \middle| a \in C(Z_5) \right\}$$

be a complex modulo integer linear algebra over the field $Z_5$.
Consider

$$X = \begin{pmatrix} 3i_F & 3i_F \\ 3i_F & 3i_F \end{pmatrix} \text{ and } Y = \begin{pmatrix} 4i_F & 4i_F \\ 4i_F & 4i_F \end{pmatrix}$$

in V. Now

$$XY = \begin{pmatrix} 3i_F & 3i_F \\ 3i_F & 3i_F \end{pmatrix} \begin{pmatrix} 4i_F & 4i_F \\ 4i_F & 4i_F \end{pmatrix}$$

$$= \begin{pmatrix} 2 \times 12i_F^2 & 2 \times 12i_F^2 \\ 2 \times 12i_F^2 & 2 \times 12i_F^2 \end{pmatrix} = \begin{pmatrix} 24 \times 4 & 24 \times 4 \\ 24 \times 4 & 24 \times 4 \end{pmatrix} = \begin{pmatrix} 1 & 1 \\ 1 & 1 \end{pmatrix}$$

(using $i_F^2 = 4$ and modulo adding 5).

Thus if
$$P = \left\{ \begin{pmatrix} a & a \\ a & a \end{pmatrix} \middle| a \in Z_5 i_F \right\} = \{ai_F \mid a \in Z_5\} \subseteq V$$
then P is not a linear subalgebra and product is defined on V but it is not in P. For take
$$x = \begin{pmatrix} i_F & i_F \\ i_F & i_F \end{pmatrix} \text{ and } y = \begin{pmatrix} 3i_F & 3i_F \\ 3i_F & 3i_F \end{pmatrix}$$
in P. Now
$$x.y = \begin{pmatrix} i_F & i_F \\ i_F & i_F \end{pmatrix} \begin{pmatrix} 3i_F & 3i_F \\ 3i_F & 3i_F \end{pmatrix}$$

$$= \begin{pmatrix} 3i_F^2 + 3i_F^2 & 3i_F^2 + 3i_F^2 \\ 3i_F^2 + 3i_F^2 & 3i_F^2 + 3i_F^2 \end{pmatrix} \pmod 5 \; (i_5^2 = 4) = \begin{pmatrix} 4 & 4 \\ 4 & 4 \end{pmatrix} \notin P.$$

Thus P cannot be a linear subalgebra but P can be complex modulo integer vector subspace so P is a pseudo complex modulo integer vector subspace of V over $Z_5$.

*Example 2.51:* Let
$$P = \left\{ \begin{pmatrix} a & a & a & a \\ a & a & a & a \\ a & a & a & a \\ a & a & a & a \end{pmatrix} \middle| a \in C(Z_{12}) = \right.$$

$\{a_1 + i_F b_1 \mid a_1, b_1 \in Z_{12}; i_F^2 = 11\}\}$ be a S-complex modulo integer linear algebra over the S-ring $Z_{12}$. P has pseudo S-complex modulo integer linear subalgebra given by

$$S = \left\{ \begin{bmatrix} a & a & a & a \\ a & a & a & a \\ a & a & a & a \\ a & a & a & a \end{bmatrix} \middle| a \in Z_{12} i_F \right\} \subseteq P$$

is not a subalgebra. Thus if $C(Z_n) = \{a + i_F b \mid a, b \in Z_n; i_F^2 = n-1\}$ then $Z_n \subseteq C(Z_n)$ and $Z_n$ is a subring called the pseudo complex modulo integer subring. However $Z_n i_F = \{a i_F \mid a \in Z_n, i_F^2 = n - 1\}$ is not a ring or a semigroup under product but only a group under addition for $a i_F \cdot b i_F$ ($a, b \in Z_n$) is $ab\, i_F^2 = ab(n-1)$ and $ab(n-1) \notin Z_n i_F$. Thus product is not defined on $Z_n i_F$.

Now all properties associated with usual vector spaces can be derived in case of complex modulo integer vector spaces. We just indicate how linear functionals on complex modulo integer of linear algebras / vector spaces defined over $Z_p$ is described in the following.

Let
$$V = \left\{ \begin{pmatrix} a_1 & a_2 \\ a_3 & a_4 \\ \vdots & \vdots \\ a_9 & a_{10} \end{pmatrix} \middle| a \in C(Z_{23}) \right.$$
$= \{a + b i_F \mid a, b \in Z_{23}; i_F^2 = 22\}; 1 \le i \le 10\}$ be a complex modulo integer vector space over the field $Z_{23}$.

We define $f : V \to Z_{23}$ as follows:
$$f\left( \begin{bmatrix} a_1 + b_1 i_F & a_2 + b_2 i_F \\ a_3 + b_3 i_F & a_4 + b_4 i_F \\ \vdots & \vdots \\ a_9 + b_9 i_F & a_{10} + b_{10} i_F \end{bmatrix} \right)$$
$= (a_1 + b_1 + a_2 + b_2 + \ldots + a_{10} + b_{10}) \bmod 23$.

Thus $f$ is a linear functional on $V$ over $Z_{23}$.

However we have to be careful to make appropriate changes while defining dual spaces and basis properties associated with linear functionals. Further we can have $f(v) = 0$ even if $v \ne 0$; and $v \in V$. Interested reader can develop it as a matter of routine.

If $A = (a_i + b_i i_F)$ is a $n \times n$ complex modulo integer matrix with entries from $Z_p$, $p$ a prime we can find characteristic values as in case of usual matrices.

**Chapter Three**

# NEUTROSOPHIC COMPLEX MODULO INTEGERS

In this chapter for the first time the authors introduce the notion of neutrosophic complex modulo integers. We work with these newly defined neutrosophic complex modulo integers to built algebraic structures.

Let $C(\langle Z_n \cup I \rangle) = \{a + bi_F + cI + i_F dI \mid a, b, c, d \in Z_n, i_F^2 = n - 1, I^2 = I, (iI)^2 = (n-1)I\}$ be the collection of neutrosophic complex modulo integers. Thus neutrosophic complex modulo integer is a 4-tuple (a, b, c, d) where first coordinate is a real value second coordinate signifies the complex number, the third coordinate the neutrosophic integer and the forth coordinate represents the complex neutrosophic integer.

***Example 3.1:*** Let $C(\langle Z_5 \cup I \rangle) = \{(a + i_F b + cI + i_F dI) \mid a, b, c, d \in Z_5, i_F^2 = 4, I^2 = I, (i_F I)^2 = 4I\}$ be the neutrosophic complex modulo integer. The order of $C(\langle Z_5 \cup I \rangle)$ is $5^4$.

***Example 3.2:*** Let $C(\langle Z_4 \cup I \rangle) = \{(a + i_F b + cI + i_F dI) \mid a, b, c, d \in Z_4, i_F^2 = 3, I^2 = I, (i_F I)^2 = 3I\}$ be the neutrosophic complex modulo integers.

***Example 3.3:*** Let $C(\langle Z_{12} \cup I\rangle) = \{(a + bi_F + cI + di_FI) \mid a, b, c, d \in Z_{12},\ i_F^2 = 11,\ I^2 = I,\ (i_FI)^2 = 11I\}$ be the neutrosophic complex modulo integers.

We give now algebraic structure using $C(\langle Z_n \cup I\rangle)$.

**DEFINITION 3.1:** *Let $C(\langle Z_n \cup I\rangle) = \{a + i_Fb + cI + i_FdI \mid a, b, c, d \in Z_n,\ i_F^2 = n-1,\ I^2 = I,\ (i_FI)^2 = (n-1)I\}$, suppose define addition '+' on $C(\langle Z_n \cup I\rangle)$ as follows if $x = \{a + i_Fb + cI + i_FdI\}$ and $y = \{m + i_Fp + sI + i_FrI\}$ are in $C(\langle Z_n \cup I\rangle)$ then*
$$\begin{aligned}x + y &= (a + i_Fb + cI + i_FdI) + (m + i_Fp + sI + i_FrI)\\ &= ((a + m)\ (mod\ n) + i_F\ (b + p)\ (mod\ n) +\\ &\quad (c + s)I\ (mod\ n) + i_FI\ (d + r)\ (mod\ n))\end{aligned}$$
*$(C(\langle Z_n \cup I\rangle), +)$ is a group for $0 + 0 + 0 + 0 = 0$ acts as the additive identity. For if $x = \{a + bi_F + cI + di_FI\} \in C(\langle Z_n \cup I\rangle)$ we have $-x = (n-a) + (n-b)i_F + (n-c)I + (n-d)i_FI \in C(\langle Z_n \cup I\rangle)$ such that $x + (-x) = 0$.*

*$(C(\langle Z_n \cup I\rangle), +)$ is defined as the group of neutrosophic complex modulo integers or neutrosophic complex modulo integers group.*

All these groups are of finite order and are commutative.

We will illustrate this situation by some examples.

***Example 3.4:*** Let $G = C(\langle Z_{11} \cup I\rangle) = \{(a + i_Fb + cI + i_FdI) \mid a, b, c, d \in Z_{11},\ I^2 = I,\ (i_F)^2 = 10,\ (i_FI)^2 = 10I\}$ be a neutrosophic complex modulo integer group under addition. $o(G) = 11^4$. G can have only subgroups of order 11 or $121 = 11^2$ or $11^3$ only.

***Example 3.5:*** Let $G = \{C(\langle Z_2 \cup I\rangle), +\} = \{0, 1, i_F, I, i_FI, 1+ i_F, 1 + I, I + i_F, 1 + i_F I, + i_F + I, i_F + i_FI, 1+ i_F + I, 1 + i_F + i_FI, 1 + I + i_F I, I + i_F + i_FI, 1+ i_F+I + i_FI\}$ be the group neutrosophic complex modulo integers of order $2^4 = 16$.

We have $H_1 = \{0, 1\}$, $H_2 = \{0, I\}$ and $H_3 = \{0, 1, I, 1 + I\}$ are subgroups of order two and four respectively.

*Example 3.6:* Let

$$P = \left\{ \begin{bmatrix} a_1 & a_2 & a_3 \\ a_4 & a_5 & a_6 \end{bmatrix} \,\middle|\, a_i \in C(\langle Z_3 \cup I \rangle) \right.$$

$= \{a + bi_F + cI + i_F dI \mid a, b, c, d \in Z_3,\ i_F^2 = 2,\ I^2 = I,\ (i_F I)^2 = 2I,\ 1 \le i \le 6\}$ be a group of neutrosophic complex modulo number $2 \times 3$ matrix under addition.

*Example 3.7:* Let

$$M = \left\{ \sum_{i=0}^{27} a_i x^i \,\middle|\, a_i \in C(\langle Z_{42} \cup I \rangle) \right.$$

$= \{a + bi_F + cI + i_F dI \mid a, b, c, d \in Z_{42},\ i_F^2 = 41,\ I^2 = I,\ (i_F I)^2 = 41I\};\ 0 \le i \le 27\}$ be a group of polynomials under addition of complex neutrosophic modulo integers.

Clearly M is of finite order.

*Example 3.8:* Let

$$P = \left\{ \begin{bmatrix} x_1 & x_2 \\ x_3 & x_4 \\ \vdots & \vdots \\ x_{43} & x_{44} \end{bmatrix} \,\middle|\, a_i \in C(\langle Z_{50} \cup I \rangle) \right.$$

$= \{a + bi_F + cI + i_F dI \mid a, b, c, d \in Z_{50},\ i_F^2 = 49,\ I^2 = I,\ (i_F I)^2 = 49I\}\}$ be the $22 \times 2$ group matrix under addition with entries from $C(\langle Z_{50} \cup I \rangle)$.

*Example 3.9:* Let $M = \{(x_1, x_2, \ldots, x_{11}) \mid x_i \in C(\langle Z_{15} \cup I \rangle) = \{a + bi_F + cI + i_F dI \mid a, b, c, d \in Z_{15},\ i_F^2 = 14,\ I^2 = I,\ (i_F I)^2 = 14I\};\ 1 \le i \le 11\}$ be a group of row matrix of neutrosophic complex modulo integers under addition.

*Example 3.10:* Let $M = \{(a, a, a, a, a, a) \mid a \in C(\langle Z_{11} \cup I \rangle) = \{a + bi_F + cI + di_F I \mid a, b, c, d \in Z_{11},\ i_F^2 = 10,\ I^2 = I,\ (i_F I)^2 = 10I\}\}$ be a row matrix of neutrosophic complex modulo integers. M is a group under product. (M does not contain $(0, 0, 0, 0, 0, 0)$).

*Example 3.11:* Let

$$M = \left\{ \begin{bmatrix} a \\ a \\ a \\ a \\ a \\ a \\ a \\ a \\ a \\ a \end{bmatrix} \;\middle|\; a_i \in C(\langle Z_{13} \cup I \rangle) \right\}$$

= {a + $i_F$ b + cI + d$i_F$I | a, b, c, d ∈ $Z_{13}$, $i_F^2$ = 12, $I^2$ = I, $(i_F I)^2$ = 12I}} be a group of neutrosophic complex modular integers under addition of finite order.

*Example 3.12:* Let

$$M = \left\{ \begin{bmatrix} a & b \\ c & d \end{bmatrix} \;\middle|\; a, b, c, d \in C(\langle Z_{15} \cup I \rangle) \right\}$$

= {m + n$i_F$ + rI + $i_F$sI | m, n, r, s ∈ $Z_{15}$, $i_F^2$ = 14, $I^2$ = I, $(i_F I)^2$ = 14I}} be a group of complex neutrosophic modulo integers under addition. Let

$$P = \left\{ \begin{bmatrix} a & b \\ 0 & 0 \end{bmatrix} \;\middle|\; a, b \in C(\langle Z_{15} \cup I \rangle) \right\} \subseteq M;$$

P is a subgroup under addition of M.

Take

$$W = \left\{ \begin{bmatrix} a & b \\ c & 0 \end{bmatrix} \;\middle|\; a, b, c \in C(\langle Z_{15} \cup I \rangle) \right\} \subseteq M$$

is also a group under addition.

Consider

B = $\left\{ \begin{bmatrix} a & b \\ c & d \end{bmatrix} \;\middle|\; a, b, c, d \in \right.$ {m + n$i_F$ + rI + $i_F$sI | m, n, r, s ∈ {0, 5, 10} ⊆ $Z_{15}$}} ⊆ M is a subgroup of M.

Now we have seen the concept of subgroups of a group.

***Example 3.13:*** Let $M = \{(a_1, \ldots, a_{10})$ where $a_i \in C(\langle Z_{40} \cup I \rangle) = \{a + i_F b + cI + i_F dI \mid a, b, c, d \in Z_{40}, i_F^2 = 39, I^2 = I, (i_F I)^2 = 39I\}$; $1 \leq i \leq 10\}$ be closed under product but is not a group under product. But M is a commutative semigroup of neutrosophic complex modulo integers under product. M has ideals and subsemigroups.

***Example 3.14:*** Let

$$S = \left\{ \begin{bmatrix} a_1 & a_2 & a_3 \\ a_4 & a_5 & a_6 \\ \vdots & \vdots & \vdots \\ a_{28} & a_{29} & a_{30} \end{bmatrix} \middle| a_i \in C(\langle Z_{43} \cup I \rangle) \right\}$$

$= \{a + bi_F + cI + i_F dI \mid a, b, c, d \in Z_{43}, i_F^2 = 42, I^2 = I, (i_F I)^2 = 42I\}$; $1 \leq i \leq 30\}$ be a complex neutrosophic modulo integer semigroup under addition.
      a) Find order of S.
      b) Find ideals if any in S.
      c) Can S have zero divisors?
      d) Find subsemigroups which are not ideals in S.

Study (a) to (d) to understand about neutrosophic complex modulo integer semigroup.

***Example 3.15:*** Let

$$M = \left\{ \begin{bmatrix} a & b \\ c & d \end{bmatrix} \middle| a, b, c, d \in C(\langle Z_{13} \cup I \rangle) \right\}$$

$= \{m + r i_F + nI + I s i_F \mid i_F^2 = 12, I^2 = I, \text{ and } (I i_F)^2 = 12I\}$ be a neutrosophic complex modulo integer semigroup under product of matrices.
i.   Find order of M.
ii.  Find subsemigroups which are not ideals of M.
iii. Is $P = \left\{ \begin{bmatrix} 0 & a \\ b & 0 \end{bmatrix} \middle| a, b \in C(\langle Z_{13} \cup I \rangle) \right\} \subseteq M$ an ideal of M?

iv. Find right ideals in M, which are not left ideals of M.
v. Prove M is a non-commutative semigroup.
vi. Is M a S-semigroup?

Answer such questions to understand this concept. However the last chapter of this book provides many problems for the reader. Further we are not interested in the study of complex neutrosophic modulo integer semigroups or groups.

What we are interested is the construction and study of complex neutrosophic modulo integer vector spaces / linear algebra and set complex neutrosophic modulo integer vector space / linear algebra and their particularized forms. Now we have seen examples of semigroups and groups built using the complex neutrosophic modulo integers.

We will first define vector spaces of complex neutrosophic modulo integers.

**DEFINITION 3.2:** *Let V be an additive abelian group of complex neutrosophic modulo integers using $C(\langle Z_p \cup I \rangle)$, p a prime. $Z_p$ be the field. If V is a vector space over $Z_p$ then V is defined as the complex neutrosophic modulo integer vector space over $Z_p$.*

We will illustrate this situation by an example.

*Example 3.16:* Let
$$V = \left\{ \sum_{i=0}^{25} a_i x^i \mid a_i \in C(\langle Z_7 \cup I \rangle) \right.$$
$= \{a + bi_F + cI + i_F dI \mid a, b, c, d \in Z_7\}, i_F^2 = 6, I^2 = I\}; 0 \leq i \leq 25\}$
be a complex neutrosophic modulo integer vector space over the field $Z_7$.

*Example 3.17:* Let
$$M = \left\{ \begin{bmatrix} a_1 & a_2 & a_3 \\ a_4 & a_5 & a_6 \\ a_7 & a_8 & a_9 \end{bmatrix} \middle| a_i \in C(\langle Z_{11} \cup I \rangle) \right.$$

= {a + bi$_F$+ cI + i$_F$dI | a, b, c, d ∈ Z$_{11}$}, i$_F^2$ = 10, I$^2$ = I; (i$_F$I)$^2$ = 10I}; 1 ≤ i ≤ 9} be a neutrosophic complex modulo integer vector space over the field Z$_{11}$.

*Example 3.18:* Let

$$P = \left\{ \begin{bmatrix} a_1 \\ a_2 \\ \vdots \\ a_{12} \end{bmatrix} \middle| a_i \in C(\langle Z_{43} \cup I \rangle) \right\}$$

= {a + bi$_F$+ cI + i$_F$dI | a, b, c, d ∈ Z$_{43}$}, i$_F^2$ = 42, I$^2$ = I; (i$_F$I)$^2$ = 42I}; 1 ≤ i ≤ 12} be a complex neutrosophic modulo integer vector space over the field F = Z$_{43}$.

*Example 3.19:* Let M = {(a$_1$, a$_2$, …, a$_{10}$) | a$_i$ ∈ C (⟨Z$_{47}$ ∪ I⟩) = {a + bi$_F$+ cI + i$_F$dI | a, b, c, d ∈ Z$_{47}$}, i$_F^2$ = 46, I$^2$ = I; (i$_F$I)$^2$ = 46I}; 1 ≤ i ≤ 10} be a neutrosophic complex modulo integer vector space over the field F = Z$_{47}$.

*Example 3.20:* Let

$$W = \left\{ \begin{pmatrix} a_1 & a_2 & \cdots & a_{12} \\ a_{13} & a_{14} & \cdots & a_{24} \end{pmatrix} \middle| a_i \in C(\langle Z_{23} \cup I \rangle) \right\}$$

= {a + bi$_F$+ cI + i$_F$dI | a, b, c, d ∈ Z$_{23}$}, i$_F^2$ = 22, I$^2$ = I; (i$_F$I)$^2$ = 22I}; 1 ≤ i ≤ 24} be a neutrosophic complex modulo integer vector space over the field Z$_{23}$ = F.

*Example 3.21:* Let

$$W = \left\{ \begin{bmatrix} p_1 & p_2 & p_3 \\ p_4 & p_5 & p_6 \\ \vdots & \vdots & \vdots \\ p_{43} & p_{44} & p_{45} \end{bmatrix} \middle| p_i \in C(\langle Z_5 \cup I \rangle) \right\}$$

= {a + i$_F$b + cI + i$_F$dI | a, b, c, d ∈ Z$_5$}, i$_F^2$ = 4, I$^2$ = I; (i$_F$I)$^2$ = 4I}; 1 ≤ i ≤ 45} be a neutrosophic complex modulo integer vector space over the field F = Z$_5$.

*Example 3.22:* Let

$$V = \left\{ \begin{bmatrix} a_1 & a_2 & a_3 \\ a_4 & a_5 & a_6 \\ a_7 & a_8 & a_9 \\ a_{10} & a_{11} & a_{12} \end{bmatrix} \middle| a_i \in C(\langle Z_{29} \cup I \rangle) \right.$$

$= \{a + bi_F + cI + i_F dI \mid a, b, c, d \in Z_{29}\}$, $i_F^2 = 28$, $I^2 = I$; $(i_F I)^2 = 28I\}$; $1 \leq i \leq 12\}$ be the neutrosophic complex modulo integer vector space over the field $Z_{29} = F$.

The definition of subvector space is a matter of routine we give a few examples of them.

*Example 3.23:* Let

$$V = \left\{ \begin{bmatrix} a_1 & a_2 \\ a_3 & a_4 \\ \vdots & \vdots \\ a_{21} & a_{22} \end{bmatrix} \middle| a_i \in C(\langle Z_{13} \cup I \rangle) \right.$$

$= \{a + bi_F + cI + i_F dI \mid a, b, c, d \in Z_{13}\}$, $i_F^2 = 12$, $I^2 = I$; $(i_F I)^2 = 12I\}$; $1 \leq i \leq 22\}$ be the neutrosophic complex modulo integer vector space over $Z_{13}$.

Consider

$$T = \left\{ \begin{bmatrix} a_1 & 0 \\ a_2 & 0 \\ \vdots & \vdots \\ a_{11} & 0 \end{bmatrix} \middle| a_i \in C(\langle Z_{13} \cup I \rangle); 1 \leq i \leq 11 \right\} \subseteq V,$$

T is a complex neutrosophic modulo integer vector subspace of V over $Z_{13}$.

$$M = \left\{ \begin{bmatrix} a_1 & a_2 \\ 0 & 0 \\ a_3 & a_4 \\ 0 & 0 \\ a_5 & a_6 \\ 0 & 0 \\ a_7 & a_8 \\ 0 & 0 \\ a_9 & a_{10} \\ 0 & 0 \\ a_{11} & a_{12} \end{bmatrix} \middle| a_i \in C(<Z_{13} \cup I>); 1 \le i \le 12 \right\} \subseteq V$$

be a neutrosophic complex modulo integer vector subspace of V over $Z_{13}$.

Now

$$T \cap M = \left\{ \begin{bmatrix} a_1 & 0 \\ 0 & 0 \\ a_2 & 0 \\ 0 & 0 \\ a_3 & 0 \\ 0 & 0 \\ a_4 & 0 \\ 0 & 0 \\ a_5 & 0 \\ 0 & 0 \\ a_6 & 0 \end{bmatrix} \middle| a_i \in C(<Z_{13} \cup I>); 1 \le i \le 6 \right\} \subseteq V$$

is a subspace of V. We have the following theorem the proof of which is left as an exercise to the reader.

**THEOREM 3.1:** *Let V be a neutrosophic complex modulo integer vector space over the field $F = Z_p$. Let $W_1, \ldots, W_t$ ($t < \infty$) be the collection of vector subspaces of V.*

$W = \cap W_i$ is a vector subspace of V. W can be the zero subspace.

*Example 3.24:* Let

$$V = \left\{ \begin{bmatrix} a_1 & a_2 & a_3 & a_4 \\ a_5 & a_6 & a_7 & a_8 \\ a_9 & a_{10} & a_{11} & a_{12} \end{bmatrix} \middle| a_i \in C(\langle Z_7 \cup I \rangle) \right.$$

$= \{a + i_Fb + cI + i_FdI \mid a, b, c, d \in Z_7 ; i_F^2 = 6, I^2 = I; (i_FI)^2 = 6I\}$; $1 \leq i \leq 12\}$ be a neutrosophic complex modulo integer vector space over $Z_7$.

Consider

$$W_1 = \left\{ \begin{bmatrix} a_1 & 0 & 0 & a_3 \\ 0 & 0 & 0 & 0 \\ a_2 & 0 & 0 & 0 \end{bmatrix} \middle| a_i \in C(<Z_7 \cup I>); 1 \leq i \leq 3 \right\} \subseteq V,$$

$$W_2 = \left\{ \begin{bmatrix} 0 & 0 & 0 & 0 \\ a_1 & 0 & 0 & a_2 \\ 0 & 0 & 0 & a_3 \end{bmatrix} \middle| a_i \in C(<Z_7 \cup I>); 1 \leq i \leq 3 \right\} \subseteq V,$$

$$W_3 = \left\{ \begin{bmatrix} 0 & a_1 & 0 & 0 \\ 0 & a_2 & 0 & 0 \\ 0 & 0 & 0 & 0 \end{bmatrix} \middle| a_1, a_2 \in C(<Z_7 \cup I>) \right\} \subseteq V,$$

$$W_4 = \left\{ \begin{bmatrix} 0 & 0 & a_2 & 0 \\ 0 & 0 & 0 & 0 \\ 0 & a_1 & 0 & 0 \end{bmatrix} \middle| a_1, a_2 \in C(<Z_7 \cup I>) \right\} \subseteq V \text{ and}$$

$$W_5 = \left\{ \begin{bmatrix} 0 & 0 & 0 & 0 \\ 0 & 0 & a_1 & 0 \\ 0 & 0 & a_2 & 0 \end{bmatrix} \middle| a_1, a_2 \in C(<Z_7 \cup I>) \right\} \subseteq V$$

be neutrosophic complex modulo integer vector subspaces of V over $Z_7$. Clearly $V = \bigcup_{i=1}^{5} W_i$ if $i \neq j$, $1 \leq i, j \leq 5$. Thus V is the direct sum of $W_1, W_2, \ldots, W_5$.

*Example 3.25:* Let

$$V = \left\{ \begin{bmatrix} a_1 & a_2 & a_3 & a_4 \\ a_5 & a_6 & a_7 & a_8 \\ a_9 & a_{10} & a_{11} & a_{12} \\ a_{13} & a_{14} & a_{15} & a_{16} \\ a_{17} & a_{18} & a_{19} & a_{20} \\ a_{21} & a_{22} & a_{23} & a_{24} \\ a_{25} & a_{26} & a_{27} & a_{28} \end{bmatrix} \middle| a_i \in C(\langle Z_{53} \cup I \rangle) \right\}$$

$= \{a + bi_F + cI + i_F dI \mid a, b, c, d \in Z_{53} ; i_F^2 = 52, I^2 = I; (i_F I)^2 = 52I\}$; $1 \leq i \leq 28\}$ be a complex neutrosophic modulo integer vector space over the field $F = Z_{53}$.

*Example 3.26:* Let

$$M = \left\{ \begin{bmatrix} a_1 & a_2 & a_3 & a_4 \\ a_5 & a_6 & a_7 & a_8 \\ \vdots & \vdots & \vdots & \vdots \\ a_{25} & a_{26} & a_{27} & a_{28} \end{bmatrix} \middle| a_i \in C(Z_{53}) \right\}$$

$= \{a + bi_F \mid i_F^2 = 52\}$ $1 \leq i \leq 8\} \subseteq V$, (V mentioned in Example 3.25) M is called as the complex pseudo neutrosophic complex modulo integer vector subspace of V over $Z_{53}$. Take

$$P = \left\{ \begin{bmatrix} a_1 & a_2 & a_3 & a_4 \\ a_5 & a_6 & a_7 & a_8 \\ \vdots & \vdots & \vdots & \vdots \\ a_{25} & a_{26} & a_{27} & a_{28} \end{bmatrix} \middle| a_i \in \langle Z_{53} \cup I \rangle \right\}$$

$= \{a + bI \mid I^2 = I, a, b \in Z_{53}\} \subseteq V$. P is a neutrosophic pseudo neutrosophic complex modulo integer vector subspace of V over $Z_{53}$.

Now if

$$T = \left\{ \begin{bmatrix} a_1 & a_2 & a_3 & a_4 \\ a_5 & a_6 & a_7 & a_8 \\ \vdots & \vdots & \vdots & \vdots \\ a_{25} & a_{26} & a_{27} & a_{28} \end{bmatrix} \middle| a_i \in Z_{53}; 1 \le i \le 28 \right\} \subseteq V;$$

then T is a real pseudo neutrosophic modulo integer vector subspace of V over $Z_{53}$.

It is interesting to note that V is not a linear algebra.

We can also have

$$B = \left\{ \begin{bmatrix} a_1 & a_2 & a_3 & a_4 \\ a_5 & a_6 & a_7 & a_8 \\ \vdots & \vdots & \vdots & \vdots \\ a_{25} & a_{26} & a_{27} & a_{28} \end{bmatrix} \middle| a_i \in Z_{53}I = \{aI \mid a \in Z_{53}\}; 1 \le i \le 28 \right\}$$

$\subseteq$ V is a pure neutrosophic pseudo neutrosophic modulo integer vector subspace of V over $Z_{53}$.

Finally if we take

$$C = \left\{ \begin{bmatrix} a_1 & a_2 & a_3 & a_4 \\ a_5 & a_6 & a_7 & a_8 \\ \vdots & \vdots & \vdots & \vdots \\ a_{25} & a_{26} & a_{27} & a_{28} \end{bmatrix} \middle| a_j \in Z_{53}i_F; \{ai_F \mid a \in Z_{53}\}; 1 \le j \le 28 \right\}$$

$\subseteq$ V is again pure complex pseudo neutrosophic complex modulo integer vector subspace.

However it is pertinent to mention here that if V is a linear algebra then V cannot contain pure complex pseudo neutrosophic complex modulo integer vector subspace. This is the one of the marked difference between a vector space and linear algebra. We have not so far defined the notion of neutrosophic complex modulo integer linear algebra. We just say a complex neutrosophic modulo integer vector space V is a complex neutrosophic modulo integer linear algebra of V is closed under a product. Inview of this we have the following theorem.

**THEOREM 3.2:** *Let V be a complex neutrosophic modulo integer linear algebra over the field $Z_p$, then V is a neutrosophic complex modulo integer vector space over the field $Z_p$. If V is a complex neutrosophic modulo integer vector space then V need not in general be a linear algebra.*

The first part follows from the very definition and the later part is proved by a counter example.

*Example 3.27:* Let

$$M = \left\{ \begin{bmatrix} a_1 & a_2 & a_3 \\ a_4 & a_5 & a_6 \\ a_7 & a_8 & a_9 \\ a_{10} & a_{11} & a_{12} \end{bmatrix} \middle| a_i \in C(\langle Z_{17} \cup I \rangle) \right.$$

$= \{a + bi_F + cI + i_F dI \mid a, b, c, d \in Z_{17}; i_F^2 = 16, I^2 = I; (i_F I)^2 = 16I\}; 1 \leq i \leq 12\}$ be a complex neutrosophic modulo vector space over the field $Z_{17}$. Clearly a product cannot be defined on M so M is only a vector space and not a linear algebra.

*Example 3.28:* Let

$$M = \left\{ \begin{bmatrix} a_1 & a_2 & a_3 & a_4 \\ a_5 & a_6 & a_7 & a_8 \\ a_9 & a_{10} & a_{11} & a_{12} \\ a_{13} & a_{14} & a_{15} & a_{16} \end{bmatrix} \middle| a_i \in C(\langle Z_{23} \cup I \rangle) \right.$$

$= \{a + bi_F + cI + i_F dI \mid a, b, c, d \in Z_{23}; i_F^2 = 22, I^2 = I; (i_F I)^2 = 22I\}; 1 \leq i \leq 16\}$ be a complex neutrosophic modulo integer linear algebra over the field $Z_{23}$.

We can define subalgebras and direct sum.

*Example 3.29:* Let $M = \left\{ \sum_{i=0}^{\infty} a_i x^i \middle| a_i \in C(\langle Z_{11} \cup I \rangle) = \{a + bi_F + cI + i_F dI \mid a, b, c, d \in Z_{11}\}, i_F^2 = i_F, I^2 = I \text{ and } (Ii_F)^2 = 10I\}$ be a

complex neutrosophic modulo integer linear algebra of infinite dimension over $Z_{11}$.

***Example 3.30:*** Let $V = \{(a_1, a_2, a_3, a_4, a_5, a_6, a_7, a_8) \mid a_i \in C(\langle Z_{17} \cup I\rangle) = \{a + bi_F + cI + i_F dI \mid a, b, c, d \in Z_{17}\}$, $i_F^2 = 16$, $I^2 = I$ and $(i_F I)^2 = 16I\}$; $1 \le i \le 8\}$ be a neutrosophic complex modulo integer linear algebra over the field $F = Z_{17}$. Clearly V is finite dimensional over $Z_{17}$.

Consider
$P_1 = \{(x_1, 0, x_2, 0, 0, 0, 0, 0) \mid x_1, x_2 \in C(\langle Z_{17} \cup I\rangle)\} \subseteq V$,
$P_2 = \{(0, x_1, 0, x_2, 0, 0, 0, 0) \mid x_1, x_2 \in C(\langle Z_{17} \cup I\rangle)\} \subseteq V$,
$P_3 = \{(0, 0, 0, 0, x_1, x_2, 0, 0) \mid x_1, x_2 \in C(\langle Z_{17} \cup I\rangle)\} \subseteq V$ and
$P_4 = \{(0\ 0\ 0\ 0\ 0\ 0\ x_1\ x_2) \mid x_1, x_2 \in C(\langle Z_{12} \cup I\rangle)\} \subseteq V$
be a collection of complex neutrosophic modulo integer linear subalgebras of V over $Z_{17}$. Clearly $V = \bigcup_{i=1}^{4} P_i$ ; $P_i \cap P_j = (0\ 0\ 0\ 0\ 0\ 0\ 0\ 0)$ if $i \ne j$, $1 \le i, j \le 4$. Thus V is a direct sum of linear subalgebras.

Take $N = \{(x_1, x_2, x_3, x_4, x_5, x_6, x_7, x_8) \mid x_i \in C(Z_{17})$ and $x_i = a_i i_F$; $a_i \in Z_{17}$; $i_F^2 = 16\} \subseteq V$; clearly N is a pseudo complex neutrosophic complex modulo integer pseudo vector subspace of V over $Z_{17}$. Clearly N is not a pseudo complex neutrosophic complex modulo linear subalgebra of V over $Z_{17}$.

Take $B = \{(x_1, x_2, x_3, x_4, x_5, x_6, x_7, x_8) \mid x_i \in C(Z_{17})$ where $x_i = a + bi_F$ with $a, b \in Z_{17}$; $i_F^2 = 16\} \subseteq B$; clearly B is a pseudo complex neutrosophic complex modulo integer linear subalgebra of V over $Z_{17}$.

For if $x = (2i_F, 0, i_F, 8i_F, 3i_F, 0, i_F, 8i_F)$ is in N we see x.y = $(2i_F, 0, i_F, 8i_F, 3i_F, 0, i_F, 8i_F).(3i_F, 9i_F, i_F, i_F, 0, 8i_F, i_F, 9i_F) = (6.16, 0, i_F\ i_F, 8\ i_F\ i_F, 3i_F\ 0, 0, 8i_F, i_F. i_F, 8i_F, 9i_F) = (11, 0, 16, 0, 0, 16, 4) \notin N$. Thus N is not closed under product so N is not a pure complex modulo integer linear algebra over $Z_{17}$.

Let $A = \{(x_1, x_2, \ldots, x_8) \mid x_i \in \langle Z_{17} \cup I\rangle = \{a + bI \mid a, b \in Z_{17}, I^2 = I; 1 \le i \le 8\} \subseteq V$; V is a pseudo neutrosophic complex neutrosophic modulo integer sublinear algebra of V over $Z_{17}$.

Consider $B = \{(x_1, x_2, \ldots, x_8) \mid x_i \in Z_{17} I = \{aI \mid a \in Z_{17}, I^2 = I; 1 \leq i \leq 8\} \subseteq V$ be a pure neutrosophic pseudo neutrosophic complex modulo integer linear subalgebra of V over $Z_{17}$.

Now we can proceed onto define other properties which can be treated as a matter of routine we would provide more illustrate examples.

*Example 3.31:* Let

$$M = \left\{ \begin{bmatrix} a_1 & a_2 & a_3 \\ a_4 & a_5 & a_6 \\ a_7 & a_8 & a_9 \\ a_{10} & a_{11} & a_{12} \end{bmatrix} \middle| a_i \in C(\langle Z_{11} \cup I \rangle) \right\}$$

$= \{a + bi_F + cI + i_F dI \mid a, b, c, d \in Z_{11}; i_F^2 = 10, I^2 = I; (i_F I)^2 = 10I\}; 1 \leq i \leq 12\}$ and

$$N = \left\{ \begin{bmatrix} a_1 & a_2 \\ a_3 & a_4 \\ a_5 & a_6 \\ a_7 & a_8 \\ a_9 & a_{10} \\ a_{11} & a_{12} \end{bmatrix} \middle| a_i \in C(\langle Z_{11} \cup I \rangle) \right\}$$

$= \{a + bi_F + cI + i_F dI \mid a, b, c, d \in Z_{11}; i_F^2 = 10, I^2 = I; (i_F I)^2 = 10I\}; 1 \leq i \leq 12\}$ be two neutrosophic complex modulo integer vector spaces over $Z_{11}$.

Consider $T: M \to N$ defined by

$$T\left( \begin{bmatrix} a_1 & a_2 & a_3 \\ a_4 & a_5 & a_6 \\ a_7 & a_8 & a_9 \\ a_{10} & a_{11} & a_{12} \end{bmatrix} \right) = \begin{bmatrix} a_1 & a_2 \\ a_3 & a_4 \\ a_5 & a_6 \\ a_7 & a_8 \\ a_9 & a_{10} \\ a_{11} & a_{12} \end{bmatrix}$$

for every

$$\begin{bmatrix} a_1 & a_2 & a_3 \\ a_4 & a_5 & a_6 \\ a_7 & a_8 & a_9 \\ a_{10} & a_{11} & a_{12} \end{bmatrix} \in M.$$

T is a linear transformation of M to N.

We can define several other linear transformations from M to N.

*Example 3.32:* Let

$$M = \left\{ \begin{bmatrix} a_1 & a_2 & a_3 \\ a_4 & a_5 & a_6 \\ a_7 & a_8 & a_9 \\ a_{10} & a_{11} & a_{12} \end{bmatrix} \middle| a_i \in C(\langle Z_{19} \cup I \rangle) \right.$$

$= \{a + bi_F + cI + i_F dI \mid a, b, c, d \in Z_{19}; i_F^2 = 18, I^2 = I; (i_F I)^2 = 18I\}; 1 \leq i \leq 12\}$ be a neutrosophic complex modulo integer vector space over the field $Z_{19}$. Define $T : M \to M$ by

$$T\left( \begin{bmatrix} a_1 & a_2 & a_3 \\ a_4 & a_5 & a_6 \\ a_7 & a_8 & a_9 \\ a_{10} & a_{11} & a_{12} \end{bmatrix} \right) = \begin{bmatrix} a_1 & 0 & a_3 \\ 0 & a_5 & 0 \\ a_7 & 0 & a_9 \\ 0 & a_{11} & 0 \end{bmatrix};$$

T is a linear operator on M. Several such linear operators can be defined.

Suppose

$$W = \left\{ \begin{bmatrix} a_1 & 0 & a_2 \\ 0 & a_3 & 0 \\ a_4 & 0 & a_5 \\ 0 & a_6 & 0 \end{bmatrix} \middle| a_i \in (<Z_{19} \cup I>); 1 \leq i \leq 6 \right\} \subseteq M$$

be a complex neutrosophic modulo integer vector subspace of M over $Z_{19}$.

We see $T(W) \subseteq W$; thus W is a invariant subspace of M over $Z_{19}$.

Consider

$$P = \left\{ \begin{bmatrix} a_1 & a_2 & a_3 \\ 0 & 0 & 0 \\ a_4 & a_5 & a_6 \\ 0 & 0 & 0 \end{bmatrix} \middle| a_i \in C(<Z_{19} \cup I>); 1 \leq i \leq 6 \right\} \subseteq M;$$

P is a complex neutrosophic modulo integer vector subspace of M over $Z_{19}$.

We see

$$T \left( \begin{bmatrix} a_1 & a_2 & a_3 \\ 0 & 0 & 0 \\ a_4 & a_5 & a_6 \\ 0 & 0 & 0 \end{bmatrix} \right) = \begin{bmatrix} a_1 & 0 & a_3 \\ 0 & 0 & 0 \\ a_4 & 0 & a_9 \\ 0 & 0 & 0 \end{bmatrix}.$$

Thus P is also invariant under T.

Consider $\eta : M \to M$ given by

$$\eta \left( \begin{bmatrix} a_1 & a_2 & a_3 \\ a_4 & a_5 & a_6 \\ a_7 & a_8 & a_9 \\ a_{10} & a_{11} & a_{12} \end{bmatrix} \right) = \begin{bmatrix} a_1 & a_2 & a_3 \\ a_4 & 0 & a_6 \\ a_7 & 0 & a_9 \\ a_{10} & 0 & a_{12} \end{bmatrix},$$

$\eta$ is a linear operator on M. But $\eta(W) \not\subseteq W$ so W is not invariant under $\eta$. Further $\eta(P) \not\subseteq P$ so P is also not invariant under $\eta$, we can derive several properties like relating nullity and rank T with dimension of V, V a complex neutrosophic modulo integer vector space and T a linear transformation from V to W; W also a complex neutrosophic modulo integer vector space over the same field as that of V. We can combine two linear operators and find $\text{Hom}_{Z_p}(V, V)$ and so on.

*Example 3.33:* Let

$$V = \left\{ \begin{bmatrix} a_1 & \ldots & a_{10} \\ a_{11} & \ldots & a_{20} \end{bmatrix} \middle| a_i \in C(\langle Z_{13} \cup I \rangle); 1 \leq i \leq 20 \right\}$$

be a complex neutrosophic modulo integer vector space over the field $Z_{13}$.

Define $f : V \to Z_{13}$ by
$$f\left(\begin{bmatrix} a_1 & a_2 & ... & a_{10} \\ a_{11} & a_{12} & ... & a_{20} \end{bmatrix}\right) = \sum_{i=1}^{20} \text{Re part of } a_i \ (\text{mod } 13);$$
f is a linear functional on V.

***Example 3.34:*** Let $V = \{(a_1, a_2) \mid a_i \in C (\langle Z_7 \cup I \rangle) = \{a + bi_F + cI + i_F dI \mid a, b, c, d \in Z_7; i_F^2 = 6, I^2 = I; (i_F I)^2 = 6I\}; 1 \le i \le 2\}$ be a neutrosophic complex modulo integer vector space over the field $Z_7$. We define $f : V \to Z_7$ by
$$f((a_1, a_2)) = \sum_{i=1}^{2} \text{Re part of } a_i \ (\text{mod } 7).$$
For instance $f((3 + 4i_F + 6I + 4i_F I), (6+2i_F + 0.I + 3i_F I)) = 3 + 6$ (mod 7) = 2 (mod 7).

Thus f is a linear functional on V. Suppose instead of defining neutrosophic complex modulo integer vector space V over a field $Z_p$ if we define V over $Z_n$, $Z_n$ a S-ring then we call V as a Smarandache neutrosophic complex modulo integer vector space over the S-ring $Z_n$.

We will give examples of this situation.

***Example 3.35:*** Let
$$V = \left\{ \begin{bmatrix} a_1 & a_2 & ... & a_8 \\ a_9 & a_{10} & ... & a_{16} \\ a_{17} & a_{18} & ... & a_{24} \end{bmatrix} \middle| a_i \in C (\langle Z_{78} \cup I \rangle) \right.$$
$= \{a + bi_F + cI + i_F dI \mid a, b, c, d \in Z_{78}\}, i_F^2 = 77, I^2 = I$ and $(i_F I)^2 = 77I\}, 1 \le i \le 24\}$ be a Smarandache complex neutrosophic modulo integer vector space over the S-ring $Z_{78}$.

***Example 3.36:*** Let
$$V = \left\{ \begin{bmatrix} a_1 & a_2 & a_3 \\ a_4 & a_5 & a_6 \\ \vdots & \vdots & \vdots \\ a_{28} & a_{29} & a_{30} \end{bmatrix} \middle| a_i \in C (\langle Z_{94} \cup I \rangle) \right.$$

= {a + bi$_F$ + cI + i$_F$dI | a, b, c, d ∈ Z$_{94}$}, i$_F^2$ = 93, I$^2$ = I and (i$_F$I)$^2$ = 93I}, 1 ≤ i ≤ 30} be the S-neutrosophic complex modulo integer vector space over the S-ring Z$_{94}$.

*Example 3.37:* Let P = {(a$_1$, a$_2$, a$_3$, a$_4$, …, a$_{20}$) | a$_i$ ∈ C(⟨Z$_{10}$ ∪ I⟩) = {a + bi$_F$ + cI + i$_F$dI | a, b, c, d ∈ Z$_{10}$}, i$_F^2$ = 9, I$^2$ = I and (i$_F$I)$^2$ = 9I}, 1 ≤ a$_i$ ≤ 20} be a S-neutrosophic complex modulo integer linear algebra over the S-ring Z$_{10}$.

*Example 3.38:* Let

$$M = \left\{ \begin{bmatrix} a_1 & a_2 & \ldots & a_7 \\ a_8 & a_9 & \ldots & a_{14} \end{bmatrix} \middle| a_i \in C(\langle Z_{34} \cup I \rangle) \right.$$

= {a + bi$_F$ + cI + i$_F$dI | a, b, c, d ∈ Z$_{34}$}, i$_F^2$ = 33, I$^2$ = I and (i$_F$I)$^2$ = 33I}, 1 ≤ i ≤ 14} be a S-neutrosophic complex modulo integer vector space over the S-ring Z$_{34}$.

It is important and interesting to note that every S-neutrosophic complex modulo integer vector space over the S-ring is in general not a S-neutrosophic complex modulo integer linear algebra but always a S-neutrosophic complex modulo integer linear algebra over a S-ring is a vector space over a S-ring.

*Example 3.39:* Let

$$V = \left\{ \begin{bmatrix} a_1 & a_2 & a_3 \\ a_4 & a_5 & a_6 \\ \vdots & \vdots & \vdots \\ a_{16} & a_{17} & a_{18} \end{bmatrix} \middle| a_i \in C(\langle Z_{46} \cup I \rangle) \right.$$

= {a + bi$_F$ + cI + i$_F$dI | a, b, c, d ∈ Z$_{46}$}, i$_F^2$ = 45, I$^2$ = I and (i$_F$I)$^2$ = 45I}, 1 ≤ i ≤ 18} be a S-neutrosophic complex modulo integer vector space over the S-ring Z$_{46}$.

Clearly V is not a S-neutrosophic complex modulo integer linear algebra over Z$_{46}$.

*Example 3.40:* Let

$$V = \left\{ \begin{bmatrix} a_1 & a_2 & \cdots & a_9 \\ a_{10} & a_{11} & \cdots & a_{18} \\ a_{19} & a_{20} & \cdots & a_{27} \end{bmatrix} \middle| a_i \in C(\langle Z_{14} \cup I \rangle) \right.$$

$= \{a + bi_F + cI + i_F dI \mid a, b, c, d \in Z_{14}\}$, $i_F^2 = 13$, $I^2 = I$ and $(i_F I)^2 = 13I\}$, $1 \le i \le 27\}$ be a Smarandache complex neutrosophic modulo integer vector space over the S-ring $Z_{14}$.

Consider

$$P = \left\{ \begin{bmatrix} a_1 & 0 & \cdots & 0 \\ 0 & 0 & \cdots & 0 \\ a_2 & 0 & \cdots & 0 \end{bmatrix} \middle| a_1, a_2 \in (\langle Z_{14} \cup I \rangle) \right\} \subseteq V$$

is a Smarandache complex neutrosophic modulo integer vector subspace of V over $Z_{14}$.

Take

$$M = \left\{ \begin{bmatrix} 0 & 0 & \cdots & 0 & a_1 & a_4 \\ 0 & 0 & \cdots & 0 & a_2 & a_5 \\ 0 & 0 & \cdots & 0 & a_3 & a_6 \end{bmatrix} \middle| a_i \in C(\langle Z_{14} \cup I \rangle); 1 \le i \le 6 \right\} \subseteq V$$

is a Smarandache complex neutrosophic modulo integer vector subspace of V over $Z_{14}$.

$$P \cap M = \left( \begin{bmatrix} 0 & 0 & \cdots & 0 \\ 0 & 0 & \cdots & 0 \\ 0 & 0 & \cdots & 0 \end{bmatrix} \right) \subseteq V$$

is the zero subspace of V over $Z_{14}$.

*Example 3.41:* Let

$$P = \left\{ \begin{bmatrix} a_1 & a_2 & a_3 \\ a_4 & a_5 & a_6 \\ a_7 & a_8 & a_9 \\ a_{10} & a_{11} & a_{12} \\ a_{13} & a_{14} & a_{15} \\ a_{16} & a_{17} & a_{18} \end{bmatrix} \middle| a_i \in C(\langle Z_{26} \cup I \rangle) \right.$$

= {a + bi_F + cI + i_FdI | a, b, c, d ∈ $Z_{26}$}, $i_F^2$ = 25, $I^2$ = I and $(i_FI)^2$ = 25I}, 1≤ i ≤ 18} be a S-neutrosophic complex modulo integer vector space of over $Z_{26}$.

Consider

$$W_1 = \left\{ \begin{bmatrix} a_1 & a_2 & a_3 \\ 0 & 0 & 0 \\ \vdots & \vdots & \vdots \\ 0 & 0 & 0 \end{bmatrix} \middle| a_1, a_2, a_3 \in C(\langle Z_{22} \cup I \rangle) \right.$$

= {a + bi_F + cI + i_FdI | a, b, c, d ∈ $Z_{22}$}, $i_F^2$ = 21, $I^2$ = I and $(i_FI)^2$ = 21I}} ⊆ V,

$$W_2 = \left\{ \begin{bmatrix} 0 & 0 & 0 \\ a_1 & a_2 & a_3 \\ 0 & 0 & 0 \\ \vdots & \vdots & \vdots \\ 0 & 0 & 0 \end{bmatrix} \middle| a_1, a_2, a_3 \in C(\langle Z_{22} \cup I \rangle) \right\} \subseteq V,$$

$$W_3 = \left\{ \begin{bmatrix} 0 & 0 & 0 \\ 0 & 0 & 0 \\ a_1 & a_2 & a_3 \\ 0 & 0 & 0 \\ \vdots & \vdots & \vdots \\ 0 & 0 & 0 \end{bmatrix} \middle| a_1, a_2, a_3 \in C(\langle Z_{22} \cup I \rangle) \right\} \subseteq V,$$

$$W_4 = \left\{ \begin{bmatrix} 0 & 0 & 0 \\ 0 & 0 & 0 \\ 0 & 0 & 0 \\ a_1 & a_2 & a_3 \\ 0 & 0 & 0 \\ 0 & 0 & 0 \end{bmatrix} \middle| a_1, a_2, a_3 \in C(\langle Z_{22} \cup I \rangle) \right\} \subseteq V,$$

$$W_5 = \left\{ \begin{bmatrix} 0 & 0 & 0 \\ 0 & 0 & 0 \\ 0 & 0 & 0 \\ 0 & 0 & 0 \\ a_1 & a_2 & a_3 \\ 0 & 0 & 0 \end{bmatrix} \middle| a_1, a_2, a_3 \in C\left(\langle Z_{22} \cup I \rangle\right) \right\} \subseteq V$$

and

$$W_6 = \left\{ \begin{bmatrix} 0 & 0 & 0 \\ 0 & 0 & 0 \\ 0 & 0 & 0 \\ 0 & 0 & 0 \\ 0 & 0 & 0 \\ a_1 & a_2 & a_3 \end{bmatrix} \middle| a_1, a_2, a_3 \in C\left(\langle Z_{22} \cup I \rangle\right) \right\} \subseteq V$$

be S-neutrosophic complex modulo number vector subspaces of V. $V = \bigcup_{i=1}^{6} W_i$; $W_i \cap W_j = (0)$ if $i \neq j$; $1 \leq i, j \leq 6$. Thus V is the direct sum of S-neutrosophic complex modulo integer vector subspaces of V over $Z_{22}$.

*Example 3.42:* Let

$$V = \left\{ \begin{bmatrix} a_1 & a_2 & a_3 \\ a_4 & a_5 & a_6 \\ a_7 & a_8 & a_9 \end{bmatrix} \middle| a_i \in C\left(\langle Z_{34} \cup I \rangle\right) \right\}$$

$= \{a + bi_F + cI + i_F dI \mid a, b, c, d \in Z_{34}\}$, $i_F^2 = 33$, $I^2 = I$ and $(i_F I)^2 = 33I\}$, $1 \leq i \leq 9\}$ be a Smarandache neutrosophic complex modulo integer vector space over the S-ring $Z_{34}$.

Take

$$W_1 = \left\{ \begin{bmatrix} a_1 & 0 & a_2 \\ 0 & a_3 & 0 \\ 0 & 0 & 0 \end{bmatrix} \middle| a_1, a_2, a_3 \in C\left(\langle Z_{34} \cup I \rangle\right) \right\} \subseteq V,$$

$$W_2 = \left\{ \begin{bmatrix} a_2 & 0 & 0 \\ a_3 & a_1 & 0 \\ 0 & 0 & 0 \end{bmatrix} \middle| a_1, a_2, a_3 \in C(\langle Z_{34} \cup I \rangle) \right\} \subseteq V,$$

$$W_3 = \left\{ \begin{bmatrix} a_1 & 0 & 0 \\ 0 & a_2 & a_3 \\ 0 & 0 & 0 \end{bmatrix} \middle| a_1, a_2, a_3 \in C(\langle Z_{34} \cup I \rangle) \right\} \subseteq V,$$

$$W_4 = \left\{ \begin{bmatrix} a_1 & 0 & 0 \\ 0 & 0 & a_3 \\ a_2 & 0 & 0 \end{bmatrix} \middle| a_1, a_2, a_3 \in C(\langle Z_{34} \cup I \rangle) \right\} \subseteq V,$$

and

$$W_5 = \left\{ \begin{bmatrix} a_1 & 0 & 0 \\ 0 & 0 & 0 \\ 0 & a_2 & a_3 \end{bmatrix} \middle| a_1, a_2, a_3 \in C(\langle Z_{34} \cup I \rangle) \right\} \subseteq V$$

be S-neutrosophic complex modulo integer vector subspaces of V over the S-ring $Z_{34}$. Clearly

$$\bigcup_{i=1}^{5} W_i \; ; \; W_i \cap W_j \neq \begin{bmatrix} 0 & 0 & 0 \\ 0 & 0 & 0 \\ 0 & 0 & 0 \end{bmatrix}$$

if $i \neq j$; $1 \leq j$, $j \leq 5$. So V is not a direct sum of $W_1, W_2, \ldots, W_5$, but only a pseudo direct sum of $W_1, W_2, \ldots, W_5$ of V.

We can also define neutrosophic complex modulo integer vector space / linear algebra over $\langle Z_n \cup I \rangle$.

We will give only examples of them.

*Example 3.43:* Let

$$V = \left\{ \begin{bmatrix} a_1 & a_2 & a_3 & a_4 \\ a_5 & a_6 & a_7 & a_8 \\ a_9 & a_{10} & a_{11} & a_{12} \\ a_{13} & a_{14} & a_{15} & a_{16} \end{bmatrix} \middle| a_i \in C(\langle Z_{29} \cup I \rangle) \right\}$$

= {a + bi$_F$+ cI + i$_F$dI | a, b, c, d ∈ Z$_{29}$}, i$_F^2$ = 28, I$^2$ = I and (i$_F$I)$^2$ = 28I}, 1≤ i ≤ 16} be a neutrosophic complex neutrosophic modulo integer linear algebra over F = C (⟨Z$_{29}$ ∪ I⟩).

*Example 3.44:* Let

$$V = \left\{ \begin{bmatrix} a_1 & a_2 & a_3 \\ a_4 & a_5 & a_6 \\ \vdots & \vdots & \vdots \\ a_{31} & a_{32} & a_{33} \end{bmatrix} \middle| a_i \in C(\langle Z_{40} \cup I \rangle) \right\}$$

= {a + bi$_F$ + cI + i$_F$dI | a, b, c, d ∈ Z$_{40}$}, i$_F^2$ = 39, I$^2$ = I and (i$_F$I)$^2$ = 39I}, 1≤ i ≤ 33} be a neutrosophic complex neutrosophic modulo integer vector space over (⟨Z$_{40}$ ∪ I⟩).

*Example 3.45:* Let

$$V = \left\{ \begin{bmatrix} a_1 & a_2 & a_3 & a_4 \\ a_5 & a_6 & a_7 & a_8 \\ \vdots & \vdots & \vdots & \vdots \\ a_{37} & a_{38} & a_{39} & a_{40} \end{bmatrix} \middle| a_i \in C(\langle Z_{29} \cup I \rangle) \right\}$$

= {a + bi$_F$+ cI + i$_F$dI | a, b, c, d ∈ Z$_{29}$}, i$_F^2$ = 28, I$^2$ = I and (i$_F$I)$^2$ = 28I}, 1≤ i ≤ 40} be a complex neutrosophic complex modulo integer vector space over the complex ring C(Z$_{29}$) = {a + bi$_F$ | a, b ∈ Z$_{29}$}, i$_F^2$ = 28}.

*Example 3.46:* Let

$$V = \left\{ \begin{bmatrix} a_1 & a_2 & a_3 & a_4 \\ a_5 & a_6 & a_7 & a_8 \\ a_9 & a_{10} & a_{11} & a_{12} \\ a_{13} & a_{14} & a_{15} & a_{16} \end{bmatrix} \middle| a_i \in C(\langle Z_{12} \cup I \rangle) \right\}$$

= {a + bi$_F$+ cI + i$_F$dI | a, b, c, d ∈ Z$_{12}$}, i$_F^2$ = 11, I$^2$ = I and (i$_F$I)$^2$ = 11I}, 1≤ i ≤ 16} be a complex neutrosophic complex modulo integer linear algebra over the complex ring C (Z$_{12}$) = {a + bi$_F$ | a, b ∈ Z$_{12}$}, i$_F^2$ = 11}.

*Example 3.47:* Let

$$V = \left\{ \begin{bmatrix} a_1 & a_2 & a_3 & a_4 \\ a_5 & a_6 & a_7 & a_8 \\ a_9 & a_{10} & a_{11} & a_{12} \\ a_{13} & a_{14} & a_{15} & a_{16} \\ a_{17} & a_{18} & a_{19} & a_{20} \end{bmatrix} \middle| a_i \in C(\langle Z_{16} \cup I \rangle) \right.$$

$= \{a + bi_F + cI + i_F dI \mid a, b, c, d \in Z_{16}\}, i_F^2 = 15, I^2 = I$ and $(i_F I)^2 = 15I\}$, $1 \le i \le 20\}$ be a complex neutrosophic complex modulo integer vector space over the complex ring $C(Z_{16}) = \{a + bi_F \mid a, b \in Z_{16}\}, i_F^2 = 15\}$.

*Example 3.48:* Let

$$V = \left\{ \begin{bmatrix} a_1 & a_2 & a_3 & a_4 \\ a_5 & a_6 & a_7 & a_8 \\ a_9 & a_{10} & a_{11} & a_{12} \end{bmatrix} \middle| a_i \in C(\langle Z_{14} \cup I \rangle) \right.$$

$= \{a + bi_F + cI + i_F dI \mid a, b, c, d \in Z_{14}\}, i_F^2 = 13, I^2 = I$ and $(i_F I)^2 = 13I\}$, $1 \le i \le 12\}$ be a strong neutrosophic complex modulo integer vector space over the complex neutrosophic modulo integer ring $C(\langle Z_{14} \cup I \rangle)$. Clearly dimension of V over S is 12.

*Example 3.49:* Let

$$V = \left\{ \begin{bmatrix} a_1 & a_2 \\ a_3 & a_4 \\ \vdots & \vdots \\ a_{19} & a_{220} \end{bmatrix} \middle| a_i \in C(\langle Z_{13} \cup I \rangle) \right.$$

$= \{a + bi_F + cI + i_F dI \mid a, b, c, d \in Z_{13}\}, i_F^2 = 12, I^2 = I$ and $(i_F I)^2 = 12I\}$, $1 \le i \le 20\}$ be a strong neutrosophic complex modulo integer vector space over the modulo integer ring C $(\langle Z_{13} \cup I \rangle)$. V is finite dimensional and dimension of V over S is 20.

*Example 3.50:* Let

$$V = \left\{ \begin{bmatrix} a_1 & a_2 & a_3 \\ a_4 & a_5 & a_6 \\ a_7 & a_8 & a_9 \end{bmatrix} \middle| a_i \in C\left(\langle Z_{43} \cup I \rangle\right) \right.$$

= $\{a + bi_F + cI + i_F dI \mid a, b, c, d \in Z_{43}\}$, $i_F^2 = 42$, $I^2 = I$ and $(i_F I)^2 = 42I\}$, $1 \le i \le 9\}$ be a special complex neutrosophic modulo integer vector space over the neutrosophic complex modulo integer ring $S = C(\langle Z_{43} \cup I \rangle)$.

Consider

$$B = \left\{ \begin{bmatrix} 1 & 0 & 0 \\ 0 & 0 & 0 \\ 0 & 0 & 0 \end{bmatrix}, \begin{bmatrix} 0 & 1 & 0 \\ 0 & 0 & 0 \\ 0 & 0 & 0 \end{bmatrix}, \begin{bmatrix} 0 & 0 & 1 \\ 0 & 0 & 0 \\ 0 & 0 & 0 \end{bmatrix}, \begin{bmatrix} 0 & 0 & 0 \\ 1 & 0 & 0 \\ 0 & 0 & 0 \end{bmatrix}, \begin{bmatrix} 0 & 0 & 0 \\ 0 & 1 & 0 \\ 0 & 0 & 0 \end{bmatrix}, \right.$$

$$\left. \begin{bmatrix} 0 & 0 & 0 \\ 0 & 0 & 1 \\ 0 & 0 & 0 \end{bmatrix}, \begin{bmatrix} 0 & 0 & 0 \\ 0 & 0 & 0 \\ 1 & 0 & 0 \end{bmatrix}, \begin{bmatrix} 0 & 0 & 0 \\ 0 & 0 & 0 \\ 0 & 1 & 0 \end{bmatrix}, \begin{bmatrix} 0 & 0 & 0 \\ 0 & 0 & 0 \\ 0 & 0 & 1 \end{bmatrix} \right\} \subseteq V,$$

B is a basis of V over S. Clearly dimension of V over S is nine. However if S is replaced by $\langle Z_{43} \cup I \rangle$ or $Z_{43}$ or $C(Z_{43})$ the dimension is different from nine.

*Example 3.51:* Let

$$V = \left\{ \begin{bmatrix} a_1 & a_2 & a_3 & a_4 \\ a_5 & a_6 & a_7 & a_8 \\ a_9 & a_{10} & a_{11} & a_{12} \end{bmatrix} \middle| a_i \in C\left(\langle Z_{37} \cup I \rangle\right) \right.$$

= {a + bi$_F$+ cI + i$_F$dI | a, b, c, d ∈ Z$_{37}$}, i$_F^2$ = 36, I$^2$ = I and (i$_F$I)$^2$ = 36I}} be a strong neutrosophic complex modulo integer vector space over the ring S = C (⟨Z$_{37}$ ∪ I⟩).

Take

$$W = \left\{ \begin{bmatrix} a_1 & a_2 & a_3 & a_4 \\ 0 & 0 & 0 & 0 \\ a_5 & a_6 & a_7 & a_8 \end{bmatrix} \middle| a_i \in C(\langle Z_{37} \cup I \rangle); 1 \le i \le 8 \right\}$$

⊆ V, W is a strong neutrosophic complex modulo integer vector space.

Consider

$$P = \left\{ \begin{bmatrix} 0 & 0 & 0 & 0 \\ a_1 & a_2 & a_3 & a_4 \\ 0 & 0 & 0 & 0 \end{bmatrix} \middle| a_i \in C(Z_{37}) \right.$$

= {a + bi$_F$ | a, b ∈ C(Z$_{37}$)}, i$_F^2$ = 36, ⊆ C (⟨Z$_{37}$ ∪ I⟩)}; 1≤ i ≤ 4} ⊆ V, P is only a pseudo strong complex neutrosophic complex modulo integer vector subspace of V over C (Z$_{37}$). Clearly P is not a strong vector subspace of V over C (⟨Z$_{37}$ ∪ I⟩).

Consider

$$B = \left\{ \begin{bmatrix} 0 & a_1 & 0 & a_2 \\ 0 & 0 & a_3 & 0 \\ 0 & 0 & a_6 & a_5 \end{bmatrix} \middle| a_i \in C(\langle Z_{37} \cup I \rangle) \right.$$

= {a + bi$_F$ | a, b ∈ Z$_{37}$}, I$^2$ = I,; 1≤ i ≤ 6} ⊆ V; B is only a pseudo neutrosophic strong neutrosophic complex modulo integer vector subspace of V over the neutrosophic field ⟨Z$_{37}$ ∪ I⟩ ⊆ C(⟨Z$_{37}$ ∪ I⟩). Clearly B is not a strong vector subspace of V over C (⟨Z$_{37}$ ∪ I⟩).

Consider

$$C = \left\{ \begin{bmatrix} a_1 & a_2 & a_3 & a_4 \\ a_5 & a_6 & a_7 & a_8 \\ a_9 & a_{10} & a_{11} & a_{12} \end{bmatrix} \middle| a_i \in Z_{37}; 1 \le i \le 12 \right\}$$

$\subseteq$ V, C is only a usual modulo integer vector space over $Z_{43}$. Thus C is a pseudo strong usual modulo integer neutrosophic complex vector subspace of V over $Z_{37} \subseteq C(\langle Z_{37} \cup I \rangle)$. Clearly C is not a strong vector space over $C(\langle Z_{37} \cup I \rangle)$. Thus we have the following theorem.

**THEOREM 3.3:** *Let V be a strong neutrosophic complex modulo integer vector space / linear algebra over C ($\langle Z_p \cup I \rangle$).*
  *Then*
i. *V has pseudo complex strong neutrosophic complex modulo integer vector subspace / linear subalgebra over C $(Z_p) = \{a + bi_F \mid a, b \in Z_p, i_F^2 = p – 1 \subseteq C (\langle Z_p \cup I \rangle)$.*
ii. *V has pseudo neutrosophic strong neutrosophic complex modulo integer vector subspace / linear subalgebra over $\langle Z_p \cup I \rangle = \{a + bI \mid a, b \in Z_p, I^2 = I\} \subseteq C (\langle Z_p \cup I \rangle)$.*
iii. *V has pseudo real modulo integer strong neutrosophic complex vector subspace / linear subalgebra over $Z_p \subseteq C (\langle Z_p \cup I \rangle)$.*

Proof is simple and hence is left as an exercise to the reader.

***Example 3.52:*** Let

$$V = \left\{ \begin{bmatrix} a_1 & a_2 \\ a_3 & a_4 \end{bmatrix} \middle| a_i \in C (\langle Z_2 \cup I \rangle) \right\}$$

= $\{a + bi_F + cI + i_FdI \mid i_F^2 = 1, I^2 = I$ and $(i_FI)^2 = I\}$, $1 \le i \le 4\}$ be a strong neutrosophic complex modulo integer linear algebra over the neutrosophic complex modulo integer ring.

$$C(\langle Z_2 \cup I \rangle) \; M = \left\{ \begin{bmatrix} a_1 & a_2 \\ a_3 & a_4 \end{bmatrix} \middle| a_i \in C (Z_2); 1 \le i \le 4 \right\} \subseteq V$$

be a pseudo real strong neutrosophic complex modulo integer linear subalgebra of V over $Z_2 \subseteq C(\langle Z_2 \cup I \rangle)$.

 Clearly M is not a strong neutrosophic complex modulo integer linear subalgebra of V over $C(\langle Z_2 \cup I \rangle)$.

Consider

$$P = \left\{ \begin{bmatrix} a_1 & a_2 \\ a_3 & a_4 \end{bmatrix} \middle| a_i \in C(Z_2); \right.$$

$\{a + bi_F \mid i_F^2 = 1, a, b \in Z_2\}; 1 \le i \le 4\} \subseteq V$; P is a pseudo complex strong neutrosophic complex modulo integer linear subalgebra of V over $C(Z_2) \subseteq C(\langle Z_2 \cup I \rangle)$.

Clearly P is not a strong complex neutrosophic modulo integer linear subalgebra of V over $C(\langle Z_2 \cup I \rangle)$.

$$B = \left\{ \begin{bmatrix} a_1 & a_2 \\ a_3 & a_4 \end{bmatrix} \middle| a_i \in C(\langle Z_2 \cup I \rangle) \right.$$

= $\{a + bI \mid a, b \in Z_2, I^2 = I\} \subseteq C(\langle Z_2 \cup I \rangle); 1 \le i \le 4\} \subseteq V$, is a pseudo neutrosophic strong complex neutrosophic complex modulo integer linear subalgebra of V over $\langle Z_2 \cup I \rangle$. Clearly B is not a strong complex neutrosophic modulo integer linear subalgebra of V over $C(\langle Z_2 \cup I \rangle)$.

Now we proceed onto define the notion of set neutrosophic complex modulo integer vector space over a set.

**DEFINITION 3.3:** *Let V be a set of elements from $C(\langle Z_n \cup I \rangle)$ (the elements can be matrices with entries from $C(\langle Z_n \cup I \rangle)$ or polynomial with coefficients from $C(\langle Z_n \cup I \rangle)$). Suppose $S \subseteq Z_n$ be a subset of $Z_n$. If for all $v \in V$ and $s \in S$, $sv = vs \in V$ then we define V to be a set neutrosophic complex modulo integer vector space over the set $S \subseteq Z_n$.*

We will give examples of this situation.

*Example 3.53:* Let

$$V = \left\{ \sum_{i=0}^{10} a_i x_i, (x_1, x_2, x_3, x_4, x_5), \begin{bmatrix} x_1 \\ x_2 \\ \vdots \\ x_7 \end{bmatrix}, \begin{bmatrix} a_1 & a_2 & a_3 \\ a_4 & a_5 & a_6 \end{bmatrix} \right.$$

$a_i, x_j \in C(\langle Z_7 \cup I \rangle); 0 \le i \le 10, 1 \le j \le 7\}$ be a set complex neutrosophic vector space over the set $\{0, 1, 5\} \subseteq Z_7$.

*Example 3.54:* Let

$$V = \left\{ \begin{pmatrix} x_1 & x_2 & x_3 \\ x_4 & x_5 & x_6 \end{pmatrix}, \begin{bmatrix} x_1 & x_2 & x_3 \\ x_4 & x_5 & x_6 \\ \vdots & \vdots & \vdots \\ x_{28} & x_{29} & x_{30} \end{bmatrix}, \begin{bmatrix} x_1 & x_2 & \ldots & x_{10} \\ x_{11} & x_{12} & \ldots & x_{20} \\ x_{21} & x_{22} & \ldots & x_{30} \end{bmatrix} \right.$$

$x_i \in C(\langle Z_{19} \cup I \rangle); 1 \le i \le 30\}$ be a set neutrosophic complex modulo integer vector space over the set $S = \{0, 1, 11, 12, 4, 7\} \subseteq Z_{19}$.

*Example 3.55:* Let

$$V = \left\{ \begin{bmatrix} a_1 & a_2 & a_3 \\ a_4 & a_5 & a_6 \\ \vdots & \vdots & \vdots \\ a_{58} & a_{59} & a_{60} \end{bmatrix}, \begin{bmatrix} a_1 & a_2 & a_3 & a_4 & a_5 \\ a_6 & a_7 & a_8 & a_9 & a_{10} \\ a_{11} & a_{12} & a_{13} & a_{14} & a_{15} \end{bmatrix}, \begin{bmatrix} a_1 & a_2 & a_3 \\ a_4 & a_5 & a_6 \\ a_7 & a_8 & a_9 \end{bmatrix} \right.$$

$a_i \in C(\langle Z_{47} \cup I \rangle); \{a + bi_F + cI + i_F dI \mid a, b, c, d \in Z_{47}\}, i_F^2 = 46, I^2 = I$ and $(i_F I)^2 = 46I\}, 1 \le i \le 9\}$ be a set neutrosophic complex modulo integer vector space over the set $S = \{0, 1, 7, 16, 19, 42, 43\} \subseteq Z_{47}$.

*Example 3.56:* Let

$$V = \left\{ \begin{bmatrix} a_1 & a_2 & a_3 \\ a_4 & a_5 & a_6 \\ \vdots & \vdots & \vdots \\ a_{28} & a_{29} & a_{30} \end{bmatrix}, \begin{bmatrix} a_1 & a_2 & \ldots & a_{41} \\ a_{42} & a_{43} & \ldots & a_{82} \end{bmatrix}, (a_1, a_2, \ldots, a_{19}) \right.$$

$a_i \in C(\langle Z_{19} \cup I \rangle); \{a + bi_F + cI + i_F dI \mid a, b, c, d \in Z_{19}\}, i_F^2 = 18, I^2 = I$ and $(i_F I)^2 = 18 i_F\}, 1 \le i \le 82\}$ be a set complex neutrosophic modulo integer vector space over the set $S = \{0, 1, 2, 14, 10, 5, 16\} \subseteq Z_{19}$.

Let us now give examples of set complex neutrosophic vector subspaces and subset neutrosophic complex modular integer subspaces of a vector space. The definition is left as an exercise.

*Example 3.57:* Let

$$V = \left\{ \sum_{i=0}^{20} a_i x_i, \begin{bmatrix} a_1 \\ a_2 \\ \vdots \\ a_{17} \end{bmatrix}, \begin{bmatrix} a_1 & a_2 & \cdots & a_{10} \\ a_{11} & a_{12} & \cdots & a_{20} \end{bmatrix}, \right.$$

$$\left. \begin{bmatrix} a_1 & a_2 & a_3 & a_4 \\ a_5 & a_6 & a_7 & a_8 \\ a_9 & a_{10} & a_{11} & a_{12} \\ a_{13} & a_{14} & a_{15} & a_{16} \end{bmatrix}, \begin{bmatrix} a_1 & a_2 \\ a_3 & a_4 \end{bmatrix} \right|$$

$a_i \in C(\langle Z_{19} \cup I \rangle); \{a + bi_F + cI + i_F dI \mid a, b, c, d \in Z_{19}\}, i_F^2 = 18, I^2 = I$ and $(i_F I)^2 = 18I\}, 1 \leq i \leq 20\}$ be a set neutrosophic complex modulo integer vector space over the set $S = \{I, 0, 1, 2, 2 + 3I, 18, 15I, 6I, 8 + 5I\} \subseteq \langle Z_{19} \cup I \rangle$.

Consider

$$P = \left\{ \sum_{i=0}^{10} a_i x_i, \begin{bmatrix} a_1 \\ a_2 \\ a_3 \\ 0 \\ 0 \\ 0 \\ \vdots \\ 0 \\ a_4 \\ a_5 \end{bmatrix}, \begin{bmatrix} a_1 & 0 & 0 & \cdots & 0 & a_2 & a_3 \\ a_4 & 0 & 0 & \cdots & 0 & a_5 & a_6 \end{bmatrix}, \right.$$

$$\left[\begin{array}{cc} a_1 & 0 \\ a_2 & 0 \end{array}\right], \left[\begin{array}{cccc} a_1 & a_2 & 0 & 0 \\ 0 & 0 & a_3 & a_4 \\ a_6 & 0 & 0 & a_5 \\ 0 & a_7 & a_8 & 0 \end{array}\right]\right\}$$

$a_i \in C(\langle Z_{19} \cup I \rangle); 1 \leq i \leq 10\} \subseteq V$, P is a set complex neutrosophic modulo integer vector subspace of V over the set S.

Consider

$$B = \left\{ \sum_{i=0}^{8} a_i x_i, \begin{bmatrix} 0 \\ 0 \\ a_1 \\ a_2 \\ a_3 \\ 0 \\ 0 \\ 0 \\ \vdots \\ 0 \\ 0 \end{bmatrix}, \begin{bmatrix} a_1 & 0 & a_3 & 0 & a_5 & 0 & a_7 & 0 & a_9 & 0 \\ 0 & a_2 & 0 & a_4 & 0 & a_6 & 0 & a_8 & 0 & 0 \end{bmatrix}, \right.$$

$$\left. \begin{bmatrix} 0 & 0 \\ a_1 & a_2 \end{bmatrix}, \begin{bmatrix} 0 & 0 & 0 & a_1 \\ 0 & 0 & a_2 & 0 \\ 0 & a_3 & 0 & 0 \\ a_4 & 0 & 0 & 0 \end{bmatrix} \right\}$$

$a_i \in C(\langle Z_{19} \cup I \rangle); 1 \leq i \leq 9\} \subseteq V$, B is a subset complex neutrosophic vector subspace of V over the subset $S = \{0, 1, 2, 18, 6I, 8 + 5I\} \subseteq \langle Z_{19} \cup I \rangle$.

*Example 3.58:* Let

$$V = \left\{ \begin{bmatrix} a_1 & a_2 & a_3 \\ a_4 & a_5 & a_6 \\ \vdots & \vdots & \vdots \\ a_{28} & a_{29} & a_{30} \end{bmatrix}, \begin{bmatrix} a_1 & a_2 \\ a_3 & a_4 \end{bmatrix}, \sum_{i=0}^{20} a_i x^i \,\middle|\, a_i \in C(\langle Z_{43} \cup I \rangle) \right.$$

$= \{a + bi_F + cI + i_F dI \mid a, b, c, d \in Z_{43}\}$, $i_F^2 = 42$, $I^2 = I$ and $(i_F I)^2 = 42I\}$, $0 \leq i \leq 30\}$ be a set neutrosophic complex modulo integer vector space over the set $S = \{0, 2, 22, 19, I\} \subseteq C(\langle Z_{43} \cup I \rangle)$. Take

$$P = \left\{ \begin{bmatrix} a_1 & a_2 \\ a_3 & a_4 \end{bmatrix}, \sum_{i=0}^{12} a_i x^i \,\middle|\, a, b, c, d, a_i \in C(\langle Z_{43} \cup I \rangle); 0 \leq i \leq 12 \right\}$$

$\subseteq V$; P is a set neutrosophic complex modulo integer vector subspace of V over S.

*Example 3.59:* Let

$$V = \left\{ \begin{bmatrix} a_1 & a_2 & a_3 \\ a_4 & a_5 & a_6 \\ \vdots & \vdots & \vdots \\ a_{28} & a_{29} & a_{30} \end{bmatrix}, \begin{bmatrix} a_1 & a_2 & \ldots & a_{19} \\ a_{20} & a_{21} & \ldots & a_{38} \end{bmatrix} \begin{bmatrix} a_1 & a_2 & a_5 \\ a_3 & a_4 & a_6 \end{bmatrix} \right.$$

$a_i \in C(\langle Z_{23} \cup I \rangle) = \{a + bi_F + cI + i_F dI \mid a, b, c, d \in Z_{23}\}$, $i_F^2 = 22$, $I^2 = I$ and $(i_F I)^2 = 22I\}$, $1 \leq i \leq 30\}$ be a set neutrosophic complex modulo integer vector space over the set $S = \{0, 1, 5, 20, 18, 7\} \subseteq Z_{23}$.

Take

$$M = \left\{ \begin{bmatrix} 0 & 0 & 0 \\ a_1 & a_2 & a_3 \\ a_4 & a_5 & a_6 \\ 0 & 0 & 0 \\ \vdots & \vdots & \vdots \\ 0 & 0 & 0 \end{bmatrix}, \begin{bmatrix} 0 & a_1 & 0 & \ldots & 0 \\ 0 & a_2 & 0 & \ldots & 0 \end{bmatrix}, \begin{bmatrix} a_1 & 0 & a_2 \\ 0 & a_3 & 0 \end{bmatrix} \right.$$

$a_i \in C(\langle Z_{23} \cup I \rangle); 1 \leq i \leq 6\} \subseteq V$; M is a set complex neutrosophic modulo integer vector subspace of V over S.

Now we provide examples of set neutrosophic complex modulo integer linear algebra over the set S.

*Example 3.60:* Let

$$V = \left\{ \begin{bmatrix} a_1 & a_2 & a_5 \\ a_3 & a_4 & a_6 \end{bmatrix}, \begin{bmatrix} a_1 & a_2 & a_3 & a_4 \\ a_5 & a_6 & a_7 & a_8 \\ \vdots & \vdots & \vdots & \vdots \\ a_{69} & a_{70} & a_{71} & a_{72} \end{bmatrix}, \sum_{i=0}^{29} a_i x^i \right\}$$

$a_i \in C (\langle Z_{11} \cup I \rangle) = \{a + bi_F + cI + i_F dI \mid a, b, c, d \in Z_{11}\}$, $i_F^2 = 10$, $I^2 = I$ and $(i_F I)^2 = 10I\}$, $0 \le i \le 72\}$ be a set neutrosophic complex modulo integer vector space over the set $S = \{1, 2, 0, 1+3I, I, 2I+1\}$. Consider

$$B = \left\{ \begin{bmatrix} 0 & a_1 & 0 \\ a_2 & 0 & a_3 \end{bmatrix}, \sum_{i=0}^{10} a_i x^i, \begin{bmatrix} 0 & 0 & 0 & 0 \\ a_1 & a_2 & a_3 & a_4 \\ 0 & 0 & 0 & 0 \\ \vdots & \vdots & \vdots & \vdots \\ 0 & 0 & 0 & 0 \end{bmatrix} \right\}$$

$a_i \in C (\langle Z_{11} \cup I \rangle); 0 \le i \le 10\} \subseteq V$, B is a set neutrosophic complex modulo integer vector subspace over the set S.

Now we see in example 3.60 we cannot define addition on V. So V is not a linear algebra.

*Example 3.61:* Let $V = \left\{ \sum_{i=0}^{25} a_i x^i \,\middle|\, a_i \in C(\langle Z_{17} \cup I \rangle) = \{a + bi_F + cI + i_F dI \mid a, b, c, d \in Z_{17}\}, i_F^2 = 16, I^2 = I \text{ and } (i_F I)^2 = 16I\}, 0 \le i \le 25 \right\}$ be a set neutrosophic complex modulo integer linear algebra over the set $S = \{0, 1, 4, 10, 12\}$.

*Example 3.62:* Let $V = \left\{ \sum_{i=0}^{25} a_i x^i \,\middle|\, a_i \in C (\langle Z_{17} \cup I \rangle) = \{a + bi_F + cI + i_F dI \mid a, b, c, d \in Z_{17}\}, i_F^2 = 16, I^2 = I \text{ and } (i_F I)^2 = 16I\}, 0 \le i \le 25 \right\}$ be a set neutrosophic complex modulo integer linear algebra over the set $S = \{0, 1, 4, 10, 12\}$.

*Example 3.63:* Let

$$B = \left\{ \begin{bmatrix} a_1 & a_2 & a_3 \\ a_4 & a_5 & a_6 \\ \vdots & \vdots & \vdots \\ a_{28} & a_{29} & a_{30} \end{bmatrix} \middle| a_i \in C(\langle Z_{37} \cup I \rangle) \right.$$

$= \{a + bi_F + cI + i_FdI \mid a, b, c, d \in Z_{37}\}$, $i_F^2 = 36$, $I^2 = I$ and $(i_FI)^2 = 36I\}, 1 \le i \le 30\}$ be a set neutrosophic complex modulo integer linear algebra over the set $S = \{0, 1, 2, 3, 10, 18, 23, 31\}$.

*Example 3.64:* Let $V = \{(a_1, a_2, \ldots, a_{13}) \mid a_i \in C(\langle Z_{31} \cup I \rangle) = \{a + bi_F + cI + i_FdI \mid a, b, c, d \in Z_{31}\}$, $i_F^2 = 30$, $I^2 = I$ and $(i_FI)^2 = 30I\}$, $1 \le i \le 31\}$ be a set neutrosophic complex modulo integer linear algebra over the set $S = \{0, 1, 30\}$.

We give examples of linear subalgebras.

*Example 3.65:* Let $V = \left\{ \begin{bmatrix} a_1 & a_2 & a_3 \\ a_4 & a_5 & a_6 \\ a_7 & a_8 & a_9 \end{bmatrix} \middle| a_i \in C(\langle Z_{43} \cup I \rangle) = \{a \right.$

$+ bi_F + cI + i_FdI \mid a, b, c, d \in Z_{42}\}$, $i_F^2 = 42$, $I^2 = I$ and $(i_FI)^2 = 42I\}$, $1 \le i \le 9\}$ be a set complex neutrosophic modulo integer linear algebra over the set $S = \{0, 4, 8, 24, 9, 41, 39\}$. We see on V we can define yet another operation product; so V becomes a set neutrosophic complex modulo integer strong linear algebra over S. We have the following interesting observations related with such the set algebras.

**THEOREM 3.4:** *Let V be a set complex neutrosophic modulo integer strong linear algebra over the set S.*
 (i)  *V is a set complex neutrosophic modulo integer linear algebra over the set S.*
 (ii) *V is a set complex neutrosophic modulo integer vector space over the set S.*
*However the converse of both (i) and (ii) are not true in general.*

*Example 3.66:* Let

$$V = \left\{ \begin{bmatrix} a_1 & a_2 & a_3 \\ a_4 & a_5 & a_6 \\ a_7 & a_8 & a_9 \\ a_{10} & a_{11} & a_{12} \end{bmatrix} \middle| a_i \in C(\langle Z_{23} \cup I \rangle) \right.$$

$= \{a + bi_F + cI + i_F dI \mid a, b, c, d \in Z_{23}\}$, $i_F^2 = 22$, $I^2 = I$ and $(i_F I)^2 = 22I\}$, $1 \leq i \leq 12\}$ be a set complex neutrosophic modulo integer linear algebra over the set $S = \{20I, 0, 1 + 4I, 5\}$.

Clearly V is not a set complex neutrosophic modulo integer strong linear algebra over the set S.

However V is a set complex neutrosophic modulo integer vector space over the set S. This proves the converse part for the claim (1) of theorem 3.

*Example 3.67:* Let

$$V = \left\{ \begin{bmatrix} a_1 & a_2 \\ a_3 & a_4 \end{bmatrix}, (a_1, a_2, a_3, a_4, a_5, a_6, a_7, a_8, a_9), \begin{bmatrix} a_1 & a_2 & a_3 \\ a_4 & a_5 & a_6 \\ \vdots & \vdots & \vdots \\ a_{28} & a_{29} & a_{30} \end{bmatrix} \right.$$

$a_i \in C(\langle Z_{11} \cup I \rangle) = \{a + bi_F + cI + i_F dI \mid a, b, c, d \in Z_{30}\}$, $i_F^2 = 10$, $I^2 = I$ and $(i_F I)^2 = 11I\}$, $1 \leq i \leq 30\}$ be a set complex neutrosophic modulo integer vector space over the set $S = \{0, I, 1+I, 5+6I, 10I + 3\}$. Clearly V is not set complex neutrosophic modulo integer linear algebra over the set S. Further it is not a strong linear algebra as it is not even a linear algebra. Hence converse of (ii) of theorem is verified.

*Example 3.68:* Let $V = \{(a_1, a_2, \ldots, a_{11});$ where $a_i \in C(\langle Z_{67} \cup I \rangle) = \{a + bi_F + cI + i_F dI \mid a, b, c, d \in Z_{67}\}$, $i_F^2 = 66$, $I^2 = I$ and $(i_F I)^2 = 66I\}$, $1 \leq i \leq 11\}$ be a set complex neutrosophic modulo integer strong linear algebra over the set $S = \{0, 1, I, 20 + 33I, 40 + 4I, 5+17I, 20I + 41\}$. Clearly V is a set complex neutrosophic modulo integer linear algebra over the set S. V is also a set complex neutrosophic modulo integer vector space over S.

Now we will give examples of these structures.

*Example 3.69:* Let

$$V = \left\{ \begin{bmatrix} a_1 & a_2 \\ a_3 & a_4 \end{bmatrix} \middle| a_i \in C(\langle Z_{13} \cup I \rangle) \right.$$

$= \{a + bi_F + cI + i_FdI \mid a, b, c, d \in Z_{13}\}, i_F^2 = 12, I^2 = I$ and $(i_FI)^2 = 12I\}\}$ be a set complex neutrosophic modulo integer strong linear algebra over the set $S = \{0, 1, 2, 5, 7, 10, 0\}$.

$$W = \left\{ \begin{bmatrix} a_1 & a_2 \\ 0 & a_3 \end{bmatrix} \middle| a_i \in C(\langle Z_{13} \cup I \rangle); 1 \le i \le 3 \right\} \subseteq V$$

is a set complex neutrosophic modulo integer strong linear subalgebra of V over S.

$$P = \left\{ \begin{bmatrix} a_1 & a_2 \\ 0 & 0 \end{bmatrix} \middle| a_i \in C(\langle Z_{13} \cup I \rangle); 1 \le i \le 2 \right\} \subseteq V$$

is also set complex neutrosophic modulo integer strong linear subalgebra of V over S.

Take

$$M = \left\{ \begin{bmatrix} 0 & a_1 \\ a_2 & 0 \end{bmatrix} \middle| a_i \in C(\langle Z_{13} \cup I \rangle); 1 \le i \le 2 \right\} \subseteq V,$$

M is only a pseudo set complex neutrosophic strong linear subalgebra of V over S as M is not a strong set linear subalgebra of V over S as product is not defined on M.

If we take

$$P = \left\{ \begin{bmatrix} a_1 & 0 \\ 0 & 0 \end{bmatrix}, \begin{bmatrix} 0 & 0 \\ a_2 & 0 \end{bmatrix}, \begin{bmatrix} 0 & a_3 \\ 0 & 0 \end{bmatrix} \middle| a_1, a_2, a_3 \in C(\langle Z_{13} \cup I \rangle) \right\}$$

$\subseteq$ M, P is only a pseudo set complex neutrosophic modulo integer subspace of V over S as in P we see the sum of two elements.

$$x = \begin{bmatrix} a_1 & 0 \\ 0 & 0 \end{bmatrix} \text{ and } y \begin{bmatrix} 0 & a_3 \\ 0 & 0 \end{bmatrix}$$

is $\quad x + y = \begin{bmatrix} a_1 & 0 \\ 0 & 0 \end{bmatrix} + \begin{bmatrix} 0 & a_3 \\ 0 & 0 \end{bmatrix} = \begin{bmatrix} a_1 & a_3 \\ 0 & 0 \end{bmatrix} \notin P.$

Hence the claim.

***Example 3.70:*** Let

$$V = \left\{ \begin{bmatrix} a_1 & a_2 & a_3 \\ a_4 & a_5 & a_6 \\ a_7 & a_8 & a_9 \end{bmatrix} \middle| a_1, a_2, a_3 \in C(\langle Z_{19} \cup I \rangle) \right.$$

$= \{a + bi_F + cI + i_F dI \mid a, b, c, d \in Z_{19}\}, i_F^2 = 18, I^2 = I$ and $(i_F I)^2 = 18I\}\}$ be a set complex neutrosophic modulo integer strong linear algebra over the set $S = \{0, 1, I+3, 9+8I, 3, 2I\}$.

Consider

$$P_1 = \left\{ \begin{bmatrix} a & 0 & 0 \\ 0 & 0 & 0 \\ 0 & 0 & 0 \end{bmatrix} \middle| a \in (\langle Z_{19} \cup I \rangle) \right\} \subseteq V$$

and

$$P_2 = \left\{ \begin{bmatrix} 0 & 0 & 0 \\ 0 & 0 & 0 \\ 0 & 0 & b \end{bmatrix} \middle| b \in (\langle Z_{19} \cup I \rangle) \right\} \subseteq V$$

be a set complex neutrosophic modulo integer strong linear subalgebra over S of V.

However we may not be in a position to write V as a direct sum of strong linear subalgebras.

***Example 3.71:*** Let $V = \{(a_1, a_2, a_3, a_4, a_5, a_6, a_7, a_8) \mid a_i \in C (\langle Z_5 \cup I \rangle) = \{a + bi_F + cI + i_F dI \mid a, b, c, d \in Z_{59}\}, i_F^2 = 4, I^2 = I$ and $(i_F I)^2 = 4I\}\}$ be a set neutrosophic complex modulo integer strong linear algebra over the set $S = \{0, 1, I + 3, 2I + 4\}$.

Consider $P_1 = \{(a_1, a_2, 0\ 0\ 0\ 0\ 0\ 0) \mid a_1\ a_2 \in C (\langle Z_5 \cup I \rangle)\} \subseteq V$, $P_2 = \{(0, 0, a_1, a_2, 0\ 0\ 0\ 0) \mid a_1, a_2 \in C (\langle Z_5 \cup I \rangle)\} \subseteq V$; $P_3 = \{(0, 0, 0, 0, a_1, a_2, 0\ 0\ ) \mid a_1, a_2 \in C (\langle Z_5 \cup I \rangle)\} \subseteq V$ and $P_4 = \{(0\ 0\ 0\ 0\ 0, a_1, a_2) \mid a_1, a_2 \in C (\langle Z_5 \cup I \rangle)\} \subseteq V$; $P_1, P_2, P_3$ and $P_4$ are set complex neutrosophic modulo integer strong linear subalgebras of V over the set S.

Clearly V = Clearly $\bigcup_{i=1}^{4} P_i$ ; $P_i \cap P_j = (0\ 0\ 0\ 0\ 0\ 0\ 0\ 0)$ if $i \neq j$; $1 \leq j, j \leq 4$. V is a direct sum of strong linear subalgebras over the set S. Let us consider the following examples.

*Example 3.72:* Let

$$V = \left\{ \begin{bmatrix} a_1 & a_2 & a_3 \\ a_4 & a_5 & a_6 \\ a_7 & a_8 & a_9 \end{bmatrix} \middle| a_i \in C(\langle Z_3 \cup I \rangle) \right.$$

$= \{a + bi_F + cI + i_F dI \mid a, b, c, d \in Z_3\}$, $i_F^2 = 2$, $I^2 = I$ and $(i_F I)^2 = 2I\}$; $1 \leq i \leq 9\}$ be a set complex neutrosophic modulo integer strong linear algebra over the set $S = \{0, 1, 1+I, 2+2I, 2\}$.

Consider

$$P_1 = \left\{ \begin{bmatrix} 0 & 0 & a_1 \\ 0 & a_2 & 0 \\ a_3 & 0 & 0 \end{bmatrix} \middle| a_i \in C(\langle Z_3 \cup I \rangle); 1 \leq i \leq 3 \right\} \subseteq V,$$

$$P_2 = \left\{ \begin{bmatrix} a_1 & a_2 & 0 \\ a_3 & 0 & 0 \\ 0 & 0 & 0 \end{bmatrix} \middle| a_i \in C(\langle Z_3 \cup I \rangle); 1 \leq i \leq 3 \right\} \subseteq V \text{ and}$$

$$P_3 = \left\{ \begin{bmatrix} 0 & 0 & 0 \\ 0 & 0 & a_1 \\ 0 & a_2 & a_3 \end{bmatrix} \middle| a_i \in C(\langle Z_3 \cup I \rangle); 1 \leq i \leq 3 \right\} \subseteq V$$

are only pseudo set complex neutrosophic modulo integer strong linear subalgebras of V. They are pseudo strong as they are not closed under product, that is product is not closed on $P_1$, $P_2$ and $P_3$.

But we see $V = \bigcup_{i=1}^{3} P_i$ ; $P_i \cap P_j = (0)$ if $i \neq j$; $1 \leq j, j \leq 3$. Thus V is only a pseudo direct sum and not a direct sum of strong linear subalgebras. Now having seen properties one can define linear transformations provided they are defined on the same set.

*Example 3.73:* Let

$$V = \left\{ \begin{bmatrix} a_1 & a_2 & a_3 \\ a_4 & a_5 & a_6 \\ a_7 & a_8 & a_9 \end{bmatrix} \middle| a_i \in C(\langle Z_{13} \cup I \rangle) \right.$$

= {a + bi$_F$ + cI + i$_F$dI | a, b, c, d ∈ Z$_3$}, i$_F^2$ = 12, I$^2$ = I and (i$_F$I)$^2$ = 12I}; 1 ≤ i ≤ 9} be a set neutrosophic complex modulo integer strong linear algebra defined over the set S = {0, I, 4I, 2 + 7I, 9I, 8, 1 + 3I}.

Consider

$$M = \left\{ \begin{bmatrix} a_1 & a_2 & a_3 & a_4 \\ a_5 & a_6 & a_7 & a_8 \\ a_9 & a_{10} & a_{11} & a_{12} \\ a_{13} & a_{14} & a_{15} & a_{16} \end{bmatrix} \middle| a_i \in C(\langle Z_{13} \cup I \rangle) \right.$$

= {a + bi$_F$ + cI + i$_F$dI | a, b, c, d ∈ Z$_{13}$}, i$_F^2$ = 12, I$^2$ = I and (i$_F$I)$^2$ = 12I}; 1 ≤ i ≤ 16} be a set neutrosophic complex modulo integer strong linear algebra over the same set S.

We define T : V → M by

$$T\left(\begin{bmatrix} a_1 & a_2 & a_3 \\ a_4 & a_5 & a_6 \\ a_7 & a_8 & a_9 \end{bmatrix}\right) = \begin{bmatrix} a_1 & a_2 & a_3 & 0 \\ a_4 & a_5 & a_6 & 0 \\ a_7 & a_8 & a_9 & 0 \\ 0 & 0 & 0 & 0 \end{bmatrix}$$

T is a strong linear transformation from V to M.

Suppose we define P : V → M by

$$P\left(\begin{bmatrix} a_1 & a_2 & a_3 \\ a_4 & a_5 & a_6 \\ a_7 & a_8 & a_9 \end{bmatrix}\right) = \begin{bmatrix} a_1 & a_2 & a_5 & 0 \\ a_3 & a_6 & 0 & 0 \\ a_4 & 0 & 0 & 0 \\ 0 & 0 & 0 & 0 \end{bmatrix},$$

P is not a strong linear transformation from V to M only a linear transformation from V to M.

We can derive all related properties with appropriate modifications. This task is left as an exercise to the reader.

Now we proceed onto define the notion of semigroup neutrosophic complex modulo integer vector spaces / linear algebras over a semigroup.

**DEFINITION 3.4:** *Let V be a set neutrosophic complex modulo integer vector space over the set S. If S is an additive semigroup then we define V to be a semigroup complex neutrosophic modulo integer vector space over the semigroup S.*

We will give examples of this situation.

*Example 3.74:* Let

$$V = \left\{ \sum_{i=0}^{20} a_i x^i, \begin{bmatrix} a_1 & a_2 & \ldots & a_{10} \\ a_{11} & a_{12} & \ldots & a_{20} \\ a_{21} & a_{22} & \ldots & a_{30} \end{bmatrix}, \begin{bmatrix} a_1 \\ a_2 \\ \vdots \\ a_{15} \end{bmatrix} \middle| a_i \in C(\langle Z_{19} \cup I \rangle) \right.$$

$= \{a + bi_F + cI + i_F dI \mid a, b, c, d \in Z_{19}\}, i_F^2 = 18, I^2 = I$ and $(i_F I)^2 = 18I\}; 0 \leq i \leq 30\}$ be semigroup complex neutrosophic modulo integer vector space over the semigroup $Z_{19}$.

*Example 3.75:* Let

$$V = \left\{ \begin{bmatrix} a_1 & a_2 & a_3 & a_4 \\ a_5 & a_6 & a_7 & a_8 \\ a_9 & a_{10} & a_{11} & a_{12} \end{bmatrix}, (a_1, a_2, \ldots, a_{20}), \sum_{i=0}^{9} a_i x^i, \begin{bmatrix} a_1 & a_2 & a_3 \\ a_4 & a_5 & a_6 \\ \vdots & \vdots & \vdots \\ a_{31} & a_{32} & a_{33} \end{bmatrix} \right.$$

$a_i \in C(\langle Z_{61} \cup I \rangle) = \{a + bi_F + cI + i_F dI \mid a, b, c, d \in Z_{19}\}, i_F^2 = 60, I^2 = I$ and $(i_F I)^2 = 60I\}; 0 \leq i \leq 33\}$ be a semigroup complex neutrosophic modulo integer vector space over the semigroup $C(Z_{61}) = \{a + bi_F \mid a, b \in Z_{61}\}, i_F^2 = 60\}$.

*Example 3.76:* Let

$$V = \left\{ \begin{bmatrix} a_1 & a_2 & a_3 \\ a_4 & a_5 & a_6 \\ \vdots & \vdots & \vdots \\ a_{37} & a_{38} & a_{39} \end{bmatrix}, \begin{bmatrix} a_1 & a_2 & \ldots & a_{14} \\ a_{15} & a_{16} & \ldots & a_{28} \\ a_{29} & a_{30} & \ldots & a_{42} \end{bmatrix}, \begin{bmatrix} a_1 & a_2 & a_3 \\ a_4 & a_5 & a_6 \\ a_7 & a_8 & a_9 \end{bmatrix}, (a_1, a_2, \ldots, a_{25}) \right.$$

$a_i \in C (\langle Z_{23} \cup I \rangle) = \{a + bi_F + cI + i_F dI \mid i_F^2 = 22, I^2 = I$ and $(i_F I)^2 = 22I\}$; $1 \le i \le 42\}$ be a semigroup neutrosophic complex modulo integer vector space over the semigroup $S = (Z_{23})$.

We can define semigroup complex neutrosophic modulo integer linear algebra over a semigroup S under addition. We give only examples of them.

*Example 3.77:* Let

$$V = \left\{ \begin{bmatrix} a_1 & a_2 & a_3 & a_4 \\ a_5 & a_6 & a_7 & a_8 \\ \vdots & \vdots & \vdots & \vdots \\ a_{37} & a_{38} & a_{39} & a_{40} \end{bmatrix} \middle| a_i \in C (\langle Z_{29} \cup I \rangle) \right.$$

$= \{a + bi_F + cI + i_F dI \mid a, b, c, d \in Z_{29}\}$, $i_F^2 = 28$, $I^2 = I$ and $(i_F I)^2 = 28I\}$; $1 \le i \le 40\}$ be a semigroup neutrosophic complex modulo integer linear algebra over the semigroup $Z_{29} = S$.

*Example 3.78:* Let

$$V = \left\{ \begin{bmatrix} a_1 & a_2 & \ldots & a_{12} \\ a_{13} & a_{14} & \ldots & a_{24} \\ a_{23} & a_{24} & \ldots & a_{36} \\ a_{37} & a_{38} & \ldots & a_{48} \end{bmatrix} \middle| a_i \in C (\langle Z_{17} \cup I \rangle) \right.$$

$= \{a + bi_F + cI + i_F dI \mid a, b, c, d \in Z_{16}\}$, $i_F^2 = 16$, $I^2 = I$ and $(i_F I)^2 = 16I\}$; $1 \le i \le 48\}$ be a semigroup neutrosophic complex modulo integer linear algebra over the semigroup $S = Z_{17}$.

*Example 3.79:* Let

$$V = \left\{ \begin{bmatrix} a_1 & a_2 & a_3 & a_4 \\ a_5 & a_6 & a_7 & a_8 \\ a_9 & a_{10} & a_{11} & a_{12} \\ a_{13} & a_{14} & a_{15} & a_{16} \end{bmatrix} \middle| a_i \in C (\langle Z_7 \cup I \rangle) \right.$$

= {a + bi$_F$+ cI + i$_F$dI | a, b, c, d ∈ Z$_7$}, i$_F^2$ = 6, I$^2$ = I and (i$_F$I)$^2$ = 6I}; 1 ≤ i ≤ 16} be a semigroup neutrosophic complex modulo integer linear algebra over the semigroup S = Z$_7$.

***Example 3.80:*** Let V = {(a$_1$, a$_2$, …, a$_{25}$) | a$_i$ ∈ C (⟨Z$_{13}$ ∪ I⟩)= {a + bi$_F$+ cI + i$_F$dI | a, b, c, d ∈ Z$_{13}$}, i$_F^2$ = 12, I$^2$ = I and (i$_F$I)$^2$ = 12I} 1 ≤ i ≤ 25} be a semigroup complex neutrosophic modulo integer strong linear algebra over the semigroup S = Z$_{13}$.

We have three types of substructures associated with semigroup neutrosophic complex modulo integer strong linear algebra.

***Example 3.81:*** Let

$$V = \left\{ \begin{bmatrix} a_1 & a_2 & a_3 \\ a_4 & a_5 & a_6 \\ a_7 & a_8 & a_9 \end{bmatrix} \middle| a_i \in C(\langle Z_{19} \cup I \rangle) \right.$$

= {a + bi$_F$+ cI + i$_F$dI | a, b, c, d ∈ Z$_{19}$}, i$_F^2$ = 18, I$^2$ = I and (i$_F$I)$^2$ = 18I}; 1 ≤ i ≤ 9} be a semigroup complex neutrosophic modulo integer strong linear algebra over the semigroup S = ⟨Z$_{19}$ ∪ I⟩.

$$W = \left\{ \begin{bmatrix} a_1 & a_2 & a_3 \\ 0 & a_4 & a_5 \\ 0 & 0 & a_6 \end{bmatrix} \middle| a_i \in C(\langle Z_{19} \cup I \rangle)'; 1 \le i \le 6 \right\}$$

be semigroup complex neutrosophic modulo integer strong linear subalgebra of V over the semigroup S.

Take

$$M = \left\{ \begin{bmatrix} a_1 & 0 & 0 \\ 0 & a_2 & 0 \\ 0 & 0 & a_3 \end{bmatrix} \middle| a_i \in C(\langle Z_{19} \cup I \rangle)'; 1 \le i \le 3 \right\},$$

a semigroup complex neutrosophic modulo integer strong linear subalgebra of V over the subsemigroup T = Z$_{19}$ ⊆ S.

$$B = \left\{ \begin{bmatrix} a_1 & a_2 & 0 \\ a_3 & 0 & 0 \\ 0 & 0 & 0 \end{bmatrix} \middle| a_i \in C(\langle Z_{19} \cup I \rangle); 1 \leq i \leq 3 \right\} \subseteq V$$

is only a pseudo semigroup complex neutrosophic modulo integer strong linear subalgebra of V over S. Clearly B is not closed with respect to product as

$$\begin{bmatrix} a_1 & a_2 & 0 \\ a_3 & 0 & 0 \\ 0 & 0 & 0 \end{bmatrix} \begin{bmatrix} x & y & 0 \\ z & 0 & 0 \\ 0 & 0 & 0 \end{bmatrix} = \begin{bmatrix} a_1x + a_2z & a_1y & 0 \\ a_3x & a_3y & 0 \\ 0 & 0 & 0 \end{bmatrix} \notin B$$

where as

$$\begin{bmatrix} a_1 & a_2 & 0 \\ a_3 & 0 & 0 \\ 0 & 0 & 0 \end{bmatrix} \text{ and } \begin{bmatrix} x & y & 0 \\ z & 0 & 0 \\ 0 & 0 & 0 \end{bmatrix}$$

are in B. Hence B is only a semigroup complex neutrosophic modulo integer linear subalgebra over the semigroup which is not a strong linear subalgebra over the semigroup S. Take

$$P = \left\{ \begin{bmatrix} a_1 & 0 & 0 \\ 0 & 0 & 0 \\ 0 & 0 & 0 \end{bmatrix}, \begin{bmatrix} 0 & 0 & 0 \\ 0 & 0 & 0 \\ 0 & 0 & a_2 \end{bmatrix}, \begin{bmatrix} 0 & 0 & a_3 \\ 0 & 0 & 0 \\ 0 & 0 & 0 \end{bmatrix} \right\}$$

$a_1, a_2, a_3 \in C(\langle Z_{19} \cup I \rangle)\} \subseteq V$, P is only a pseudo semigroup neutrosophic complex modulo integer pseudo vector subspace of V over S. Clearly V is not closed with respect to addition or multiplication. Hence we see we can define several types of substructures in case of strong linear algebras defined over semigroups.

*Example 3.82:* Let

$$V = \left\{ \begin{bmatrix} a_1 & a_2 & a_3 \\ a_4 & a_5 & a_6 \\ a_7 & a_8 & a_9 \\ a_{10} & a_{11} & a_{12} \end{bmatrix}, (a_1, a_2, ..., a_{20}), \begin{bmatrix} a_1 & a_{13} \\ a_2 & a_{14} \\ \vdots & \vdots \\ a_{12} & a_{24} \end{bmatrix} \right.$$

$a_i \in C(\langle Z_{47} \cup I \rangle) = \{a + bi_F + cI + i_F dI \mid a, b, c, d \in Z_{47}\}$, $i_F^2 = 46$, $I^2 = I$ and $(i_F I)^2 = 46I\}\}$ be a semigroup neutrosophic complex modulo integer vector space over the semigroup $S = Z_{47}$.

Clearly V is not a semigroup complex neutrosophic modulo integer linear subalgebra or V is not a semigroup complex neutrosophic modulo integer strong linear subalgebra.

*Example 3.83:* Let

$$V = \left\{ \begin{bmatrix} a_1 & a_2 & a_3 \\ a_4 & a_5 & a_6 \\ a_7 & a_8 & a_9 \\ a_{10} & a_{11} & a_{12} \end{bmatrix} \middle| a_i \in C(\langle Z_{17} \cup I \rangle) \right\}$$

= $\{a + bi_F + cI + i_F dI \mid a, b, c, d \in Z_{17}\}$, $i_F^2 = 16$, $I^2 = I$ and $(i_F I)^2 = 16I\}$, $1 \le i \le 12\}$ be a semigroup neutrosophic complex modulo integer linear algebra over $Z_{17}$ the semigroup.

Clearly V is not a semigroup neutrosophic complex modulo integer strong linear algebra over the semigroup $Z_{17}$.

But consider

$$M = \left\{ \begin{bmatrix} a_1 & a_2 & a_3 \\ 0 & 0 & 0 \\ 0 & 0 & 0 \\ a_4 & a_5 & a_6 \end{bmatrix} \middle| a_i \in C(\langle Z_{17} \cup I \rangle) \right\}$$

= $\{a + bi_F + cI + i_F dI \mid a, b, c, d \in Z_{17}\}$, $i_F^2 = 16$, $I^2 = I$ and $(i_F I)^2 = 16I\}$, $1 \le i \le 6\} \subseteq V$ is only a semigroup complex neutrosophic modulo integer linear subalgebra of V over the semigroup $S = Z_{17}$.

We can define linear transformation of two semigroup neutrosophic complex linear algebra only they are defined over the same semigroup S.

We give examples of them.

*Example 3.84:* Let

$$V = \left\{ \begin{bmatrix} a_1 & a_2 & a_3 \\ a_4 & a_5 & a_6 \\ a_7 & a_8 & a_9 \end{bmatrix} \middle| a_i \in C(\langle Z_{29} \cup I \rangle) \right.$$

$= \{a + bi_F + cI + i_F dI \mid a, b, c, d \in Z_{29}\}$, $i_F^2 = 28$, $I^2 = I$ and $(i_F I)^2 = 28I\}\}$ be a semigroup neutrosophic complex modulo integer strong linear algebra over $Z_{29}$. $W = \{(a_1, a_2, a_3, a_4, a_5, a_6, a_7, a_8, a_9) \mid a_i \in C(\langle Z_{29} \cup I \rangle) = \{a + bi_F + cI + i_F dI \mid a, b, c, d \in Z_{29}\}$, $i_F^2 = 28$, $I^2 = I$ and $(i_F I)^2 = 28I\} 1 \le i \le 9\}$ be a semigroup complex neutrosophic modulo integer strong linear algebra over $S = Z_{29}$.

Define $T : V \to W$ where

$$T\left( \begin{bmatrix} a_1 & a_2 & a_3 \\ a_4 & a_5 & a_6 \\ a_7 & a_8 & a_9 \end{bmatrix} \right) = (a_1, a_2, a_3, a_4, a_5, a_6, a_7, a_8, a_9);$$

T is a linear transformation from V to W.

Let $\eta : V \to V$ be defined by

$$\eta \left( \begin{bmatrix} a_1 & a_2 & a_3 \\ a_4 & a_5 & a_6 \\ a_7 & a_8 & a_9 \end{bmatrix} \right) = \begin{bmatrix} a_1 & 0 & 0 \\ 0 & a_2 & a_5 \\ 0 & a_3 & a_4 \end{bmatrix}.$$

$\eta$ is a linear operator on V.

Now we can derive several related properties with some simple appropriate changes. Now we proceed onto define the notion of group neutrosophic complex modulo integer vector space /linear algebra / strong linear algebra over a group G.

If V is a set with zero of complex neutrosophic modulo integers and G to be a group of integers addition. We call V a group neutrosophic complex modulo integer vector space over the group G if

i) for every $v \in V$ and $g \in G$ and gv and vg are in V.
ii) $0.v = 0$ for every $v \in V$ and 0 the additive identity of G.

We give examples of them.

*Example 3.85:* Let

$$V = \left\{ \begin{bmatrix} a_1 & a_2 & a_3 \\ a_4 & a_5 & a_6 \end{bmatrix}, \begin{bmatrix} a_1 \\ a_2 \\ \vdots \\ a_{10} \end{bmatrix}, (a_1, a_2, ..., a_{24}), \begin{bmatrix} a_1 & a_2 \\ a_3 & a_4 \end{bmatrix} \right.$$

$a_i \in C(\langle Z_{47} \cup I \rangle) = \{a + bi_F + cI + i_F dI \mid a, b, c, d \in Z_{47}\}, i_F^2 = 46$, $I^2 = I$ and $(i_F I)^2 = 46I\}$; $1 \le i \le 24\}$ be a group complex neutrosophic modulo integer vector space over the group $G = Z_{47}$ under addition.

*Example 3.86:* Let

$$M = \left\{ \begin{bmatrix} a_1 & a_2 & ... & a_{10} \\ a_{11} & a_{12} & ... & a_{20} \\ a_{21} & a_{22} & ... & a_{30} \end{bmatrix} \middle| a_i \in C(\langle Z_{29} \cup I \rangle) \right.$$

$= \{a + bi_F + cI + i_F dI \mid a, b, c, d \in Z_{29}\}, i_F^2 = 28, I^2 = I$ and $(i_F I)^2 = 28I\}$; $1 \le i \le 30\}$ be a group neutrosophic complex modulo integer linear algebra over the semigroup $G = Z_{29}$.

Clearly M is also a group complex neutrosophic modulo integer vector space over the group. However V in example 3.78 is not a linear algebra only a vector space.

*Example 3.87:* Let

$$V = \left\{ \begin{bmatrix} a_1 & a_2 & a_3 & a_4 \\ a_5 & a_6 & a_7 & a_8 \\ a_9 & a_{10} & a_{11} & a_{12} \\ a_{13} & a_{14} & a_{15} & a_{16} \end{bmatrix} \middle| a_i \in C(\langle Z_{23} \cup I \rangle) \right.$$

= {a + bi$_F$+ cI + i$_F$dI | a, b, c, d ∈ Z$_{23}$}, i$_F^2$ = 22, I$^2$ = I and (i$_F$I)$^2$ = 22I}; 1 ≤ i ≤ 16} be a group neutrosophic complex modulo integer strong linear algebra over the group G = Z$_{23}$.

*Example 3.88:* Let

$$V = \left\{ \begin{bmatrix} a_1 & a_2 \\ a_3 & a_4 \end{bmatrix}, (a_1, a_2, ..., a_{12}), \begin{bmatrix} a_1 & a_2 & a_3 \\ a_4 & a_5 & a_6 \\ \vdots & \vdots & \vdots \\ a_{28} & a_{29} & a_{39} \end{bmatrix} \right\}$$

a$_i$ ∈ C (⟨Z$_{11}$ ∪ I⟩) = {a + bi$_F$+ cI + i$_F$dI | a, b, c, d ∈ Z$_{11}$}, i$_F^2$ = 10, I$^2$ = I and (i$_F$I)$^2$ = 10I}; 1 ≤ i ≤ 30} be a group neutrosophic complex modulo integer vector space over the group G = Z$_{11}$. V is not a group neutrosophic complex modulo integer linear algebra over G = Z$_{11}$.

Clearly V is not a group neutrosophic complex modulo integer strong linear algebra over the group G = Z$_{11}$.

*Example 3.89:* Let

$$V = \left\{ \begin{bmatrix} a_1 & a_2 & a_3 \\ a_4 & a_5 & a_6 \\ \vdots & \vdots & \vdots \\ a_{13} & a_{14} & a_{15} \end{bmatrix} \middle| a_i \in C(\langle Z_{23} \cup I \rangle) \right\}$$

= {a + bi$_F$+ cI + i$_F$dI | a, b, c, d ∈ Z$_{23}$}, i$_F^2$ = 22, I$^2$ = I and (i$_F$I)$^2$ = 22I}; 1 ≤ i ≤ 15} be a group neutrosophic complex modulo integer linear algebra over the group G = Z$_{23}$.

Clearly T is not a group complex neutrosophic modulo strong linear algebra over G.

*Example 3.90:* Let

$$V = \left\{ \begin{bmatrix} a_1 & a_6 \\ a_2 & a_7 \\ a_3 & a_8 \\ a_4 & a_9 \\ a_5 & a_{10} \end{bmatrix} \middle| a_i \in C(\langle Z_7 \cup I \rangle) \right\}$$

$= \{a + bi_F + cI + i_F dI \mid a, b, c, d \in Z_7\}$, $i_F^2 = 6$, $I^2 = I$ and $(i_F I)^2 = 6I\}$; $1 \le i \le 10\}$ and

$$P = \left\{ \begin{bmatrix} a_1 & a_2 & a_3 & a_4 \\ a_5 & a_6 & a_7 & a_8 \end{bmatrix} \middle| a_i \in C(\langle Z_7 \cup I \rangle) \right\}$$

$= \{a + bi_F + cI + Idi_F \mid a, b, c, d \in Z_7$, $i_F^2 = 6$, $(i_F I)^2 = 6I$, $I^2 = I\}$; $1 \le i \le 8\}$ be a group complex neutrosophic modulo integer linear algebra over the group $G = (Z_7, +)$.

Define $\eta : M \to P$

$$\eta \left( \begin{bmatrix} a_1 & a_6 \\ a_2 & a_7 \\ a_3 & a_8 \\ a_4 & a_9 \\ a_5 & a_{10} \end{bmatrix} \right) = \begin{bmatrix} a_1 & a_2 & a_3 & a_4 \\ a_5 & a_6 & a_7 & a_8 \end{bmatrix}.$$

$\eta$ is a linear transformation from M to P.
If $T : M \to M$ such that

$$T \left( \begin{bmatrix} a_1 & a_6 \\ a_2 & a_7 \\ a_3 & a_8 \\ a_4 & a_9 \\ a_5 & a_{10} \end{bmatrix} \right) = \begin{bmatrix} a_1 & a_2 \\ 0 & 0 \\ a_3 & a_4 \\ 0_4 & 0 \\ a_5 & a_6 \end{bmatrix}$$

then T is a linear operator on M.

We can derive almost all properties for group neutrosophic complex modulo integer (strong) linear algebra defined over the group G.

We can also find $\text{Hom}_G(V, W)$ and $\text{Hom}_G(V, V)$ study the algebraic structure enjoyed by them.

Also study the substructure and writing them as direct sum and pseudo direct sum can be taken as a routine exercise by the interested reader.

**Chapter Four**

# APPLICATIONS OF COMPLEX NEUTROSOPHIC NUMBERS AND ALGEBRAIC STRUCTURES USING THEM

In this chapter we study the probable applications of complex neutrosophic reals and complex neutrosophic modulo integers.

It is pertinent to keep on record such study is very dormant; we can say solving equations and finding solutions in $C(\langle R \cup I \rangle) = \{a + bi + cI + Idi \mid a, b, c, d \in R; i^2 = -1\}$ has not been carried out.

Clearly this set $C(\langle R \cup I \rangle)$ contain R and the algebraically closed field namely the field of complex number $C = \{a + ib$ where $a, b \in R\}$ as proper subsets.

Thus we can say this extended like field will also give the roots when the roots are indeterminates. We can denote the complex neutrosophic number by the 4 - tuple (a, b, c, d) the first coordinate represents the real value, the second coordinate the complex coefficients, the third coordinate the neutrosophic coefficient and forth coordinate the complex neutrosophic coefficient and (a, b, c, d) = a + bi + cI + idI.

Study of the eigen values as complex or indeterminate stands as one of the applications.

Further the notion of complex modulo integers $C(Z_n)$ is itself very new. Here the complex number $i_F$ is defined as $i_F^2 = n - 1$ when $Z_n$ is taken into account. As n varies the value of the finite square of the complex number also vary.

Hence $C(Z_n) = \{a + bi_F \mid a, b \in Z_n \text{ and } i_F^2 = n - 1\}$. Thus $C(Z_3) = \{a + bi_F \mid a, b \in Z_3 \text{ and } i_F^2 = 2\}$. Now in due course of time these new structures will find very many applications.

Finally the notion of complex neutrosophic modulo integers is defined. This is also represented as a 4-tuple with a special value for the finite complex number $i_F$, $i_F$ is defined as $i_F^2 = n - 1$ and $(Ii_F)^2 = (n - 1)I$.

Thus $C(\langle Z_n \cup I \rangle) = \{a + bi_F + cI + i_F dI \mid a, b, c, d \in Z_n; i_F^2 = n - 1, I^2 = I \text{ and } (i_F I)^2 = (n - 1)I\}$.

These new structures are given algebraic structures like groups semigroups, rings, vector spaces and linear algebras and they will find application in due course of time once this research becomes popular.

**Chapter Five**

# SUGGESTED PROBLEMS

In this chapter we suggest over 150 problems some at research level and some just routine exercises.

1. Obtain some interesting properties enjoyed by
   (a) neutrosophic complex reals.
   (b) neutrosophic complex modulo integers.
   (c) neutrosophic complex rationals.

2. Can any geometrical interpretation be given to the field of neutrosophic complex numbers $C(\langle Q \cup I \rangle)$?

3. Can $C(\langle Z \cup I \rangle)$ be a Smarandache ring?

4. Is $\{(a_1, a_2) \mid a_1, a_2 \in C(\langle Z \cup I \rangle)$ under product $\times\}$ a Smarandache semigroup?

5. Let $V = \left\{ \begin{bmatrix} a_1 \\ a_2 \\ a_3 \\ a_4 \end{bmatrix} \middle| a_i \in C(\langle Q \cup I \rangle), i = 1, 2, 3, 4, + \right\}$ be a group.
   i) Define a automorphism $\eta : V \to V$ so that ker $\eta$ is a nontrivial subgroup.
   ii) Is $V \cong C(\langle Q \cup I \rangle) \times C(\langle Q \cup I \rangle) \times C(\langle Q \cup I \rangle) \times C(\langle Q \cup I \rangle)$?

6. Let $M = \left\{ \begin{pmatrix} a_1 & a_2 \\ a_3 & a_4 \end{pmatrix} \middle| a_i \in C(\langle Q \cup I \rangle), 1 \le i \le 4 \right\}$ be a semigroup under multiplication.
   i) Prove M is a S-semigroup.
   ii) Is M commutative?
   iii) Find at least three zero divisors in M.
   iv) Does M have ideals?
   v) Give subsemigroups in M which are not ideals.

7. Let $S = \left\{ \begin{bmatrix} a_1 & a_2 & a_3 \\ a_4 & a_5 & a_6 \\ a_7 & a_8 & a_9 \\ a_{10} & a_{11} & a_{12} \end{bmatrix} \middle| a_i \in C(\langle Q \cup I \rangle), 1 \le i \le 12 \right\}$ be a neutrosophic complex rational semigroup under '+'.
   i) Find subsemigroups of S.
   ii) Can S have ideals?
   iii) Can S have idempotents?
   iv) Can S have zero divisors? Justify your claim.

8. Let $V = \left\{ \begin{pmatrix} a_1 & a_2 \\ a_3 & a_4 \end{pmatrix} \middle| a_i \in C(\langle Q \cup I \rangle), 1 \le i \le 4 \right\}$ be a semigroup of neutrosophic complex rationals under product.

      i) Is V commutative?
      ii) Can V have idempotents?
      iii) Does V have a subsemigroup which is not an ideal?
      iv) Give an ideal in V.
      v) Can V have zero divisors?
      vi) Is V a Smarandache semigroup?
      vii) Is V a Smarandache commutative semigroup?

9. Let $V = \left\{ \begin{pmatrix} a_1 & ... & a_{10} \\ a_2 & ... & a_{20} \end{pmatrix} \middle| a_i \in C(\langle Q \cup I \rangle), 1 \leq i \leq 20 \right\}$ be a complex neutrosophic group under addition.
    a) Is V commutative?
    b) Find subgroups of V.
    c) Find a subgroup H and describe V/H.
    d) Define $\eta : V \to V$ so that
$$\ker \eta \neq \begin{pmatrix} 0 & 0 & 0 & 0 & ... & 0 \\ 0 & 0 & 0 & 0 & ... & 0 \end{pmatrix}.$$
    e) Find $\eta : V \to V$ so that $\eta^{-1}$ exists.

10. Prove $C(\langle Z \cup I \rangle)$ is a group under addition and only a semigroup under multiplication. Is $C(\langle Z \cup I \rangle)$ an integral domain? Justify your claim.

11. Is $C(\langle Q \cup I \rangle) = \{a + bi + cI + idI \mid a, b, c, d \in Q\}$ a field? Is $C(\langle Q \cup I \rangle)$ a prime field?

12. Can one say for all polynomials with complex neutrosophic coefficients $C(\langle R \cup I \rangle)$ is the algebraically closed field?

13. Can $C(\langle R \cup I \rangle)[x]$ have irreducible polynomials?

14. Find irreducible polynomials in $C(\langle Q \cup I \rangle)$?

15. Find irreducible polynomials in $C(\langle Z \cup I \rangle)$? Is every ideal in $C(\langle Z \cup I \rangle)$ principal? Justify your claim.

16. Study the ring $P = \left\{ \begin{bmatrix} a & b \\ c & d \end{bmatrix} \middle| a, b, c, d, \in C(\langle Z \cup I \rangle) \right\}$.

17. Is $G = \left\{ \begin{bmatrix} a & b \\ c & d \end{bmatrix} \middle| ad - bc \neq 0, a, b, c, d, \in C(\langle Z \cup I \rangle) \right\}$ a group? Is G simple?

18. Find zero divisors in P mentioned in problem (16)

19. What are the advantages of using the algebraic structure $C(\langle R \cup I \rangle)$?

20. Give some uses of this complete algebraic structure $C(\langle R \cup I \rangle)$.

21. What is the cardinality of the semigroup $S = \left\{ \begin{bmatrix} a & b \\ c & d \end{bmatrix} \middle| a, b, c, d, \in C(\langle Z_2 \cup I \rangle), \times \right\}$?

22. Find the number of elements in $T = \left\{ \begin{bmatrix} a & b & c \\ d & e & f \end{bmatrix} \middle| a, b, c, d, e, f \in C(\langle Z_3 \cup I \rangle) \right\}$, the semigroup under addition.

23. Is $M = \left\{ \begin{bmatrix} a & b \\ c & d \end{bmatrix} \middle| a, b, c, d, \in C(\langle Z_5 \cup I \rangle) \right\}$ a ring?.
    a) Find ideals if any in M.
    b) Is M a S-ring?
    c) What is the order of M?
    d) Can M have S-zero divisors?
    e) Does M contain idempotents?

24. Is $C(\langle Z_{19} \cup I \rangle)$ a field?

25. Prove $C(\langle Z_{25} \cup I \rangle)$ can only be a ring.

26. Can $C(\langle Z_{12} \cup I\rangle)$ be a S-ring?  Justify.

27. Can $C(\langle Z_5 \cup I\rangle)$ have S-ideals?

28. Find ideals in $C(\langle Z_6 \cup I\rangle)$.

29. Find maximum ideals of $C(\langle Z_{18} \cup I\rangle)$.

30. Find S-zero divisors and S-units if any in $C(\langle Z_{24} \cup I\rangle)$.

31. Is $C(\langle Z_{11} \cup I\rangle)$ a Smarandache semigroup under product?

32. Does $S = C(\langle Z_{15} \cup I\rangle)$ have S-ideals where S is a semigroup under $\times$?

33. Let $R = C(\langle Z_{24} \cup I\rangle)$ be a ring.
    i) Find S-subrings of R.
    ii) Can R have S-ideals?
    iii) Does R contain S-subrings which are not ideals?
    iv) Find zero divisors in R.
    v) Is every zero divisor in R a S-zero divisor?
    vi) Find S-idempotent if any in R.
    vii) Determine the number of elements in R.
    viii) Is R a S-ring?
    ix) Find an ideal I in R so that R / I is a field.
    x) Does there exist an ideal J in R so that R/J is a S-ring?

34. Let $V = \{(x_1, x_2, x_3, x_4, x_5) \mid x_i \in C(\langle Z_{11} \cup I\rangle); 1 \le i \le 5\}$ be a vector space over $Z_{11}$.
    i) Find dimension of V.
    ii) Is V finite dimensional?
    iii) Find a basis for V.
    iv) Find subspaces of V.
    v) Find $\text{Hom}_{Z_{11}}(V, V)$.
    vi) Find a linear operator T on V so that $T^{-1}$ does not exist.
    vii) Write V as a direct sum.

viii) Write V as a pseudo direct sum.

ix) Define a T : V → V and verify Null T + Range T = dim V.

35. Let $V = \left\{ \begin{bmatrix} a_1 & a_2 & a_3 & a_4 & a_5 \\ a_6 & a_7 & a_8 & a_9 & a_{10} \end{bmatrix} \middle| a_i \in C(\langle Z \cup I \rangle); 1 \leq i \leq 10 \right\}$ be a set complex neutrosophic linear algebra over the set $S = 3Z \cup 5Z$.

i) Find a basis of V.
ii) Is V finite dimensional?
iii) Write V as a direct sum.
iv) Write V as a pseudo direct sum.
v) Obtain conditions on a linear operator T of V so that $T^{-1}$ exists.
vi) Can V have subset complex neutrosophic linear subalgebras?
vii) Find the algebraic structure enjoyed $\text{Hom}_S(V,V)$.
viii) Can V be a double linear algebra?

36. Let $V = \left\{ \begin{bmatrix} a_1 & a_2 \\ a_3 & a_4 \\ a_5 & a_6 \end{bmatrix}, (a_1, a_2, a_3), \begin{bmatrix} a_1 & a_2 & \ldots & a_{10} \\ a_{11} & a_{12} & \ldots & a_{20} \end{bmatrix}, \sum_{i=0}^{5} a_i x^i \right\}$ be a set neutrosophic complex vector space with $a_i \in C(\langle Q \cup I \rangle)$, $0 \leq i \leq 20$ over the set $S = 7Z \cup 3Z$.

i) Find a basis for V.
ii) Is V finite dimensional?
iii) Write V as a direct sum of subspaces.
iv) Find a linear operator on T so that $T^{-1}$ does not exist.
v) Does V contain subset vector subspace?
vi) Can V be written as a pseudo direct sum?
vii) Let $V = W_1 \oplus \ldots \oplus W_4$ find projection E on V and describe their properties.

37. Let $V = \left\{ \begin{bmatrix} a_1 & a_2 & \dots & a_6 \\ a_7 & a_8 & \dots & a_{12} \\ a_{13} & a_{14} & \dots & a_{18} \end{bmatrix}, \sum_{i=0}^{21} a_i x^i, \begin{bmatrix} a_1 & a_2 \\ a_3 & a_4 \end{bmatrix}, \begin{bmatrix} a_1 \\ a_2 \\ a_3 \end{bmatrix} \middle| a_i \right.$

$\in C(\langle Q \cup I \rangle), 0 \leq i \leq 21\}$ be a special set neutrosophic complex vector space over the set $S = ZI \cup 3Z \cup 5Z \cup C(Z)$.
   i) Find a basis.
   ii) What is the dimension of V over S?
   iii) Does there exists a subset vector subspace W of dimension less than five over a subset T of S?
   iv) Write V as a direct sum of set subspaces.
   v) Find a non invertible linear operator on V.
   vi) Write V as a pseudo direct sum and define projections. Is that possible?

38. Let $V = \left\{ \begin{pmatrix} a_1 & a_2 \\ a_3 & a_4 \end{pmatrix} \middle| a_i \in C(\langle Q \cup I \rangle), 1 \leq i \leq 4 \right\}$ be a

semigroup neutrosophic complex double linear algebra over $S = Z$.
   i) Find a basis of V over Z
   ii) What is the dimension of V over Z?
   iii) If Z is replaced by $C(Z)$ what is dimension of V?
   iv) If Z is replaced by $\langle Z \cup I \rangle$ what is dimension of V?
   v) If Z is replaced by $C(\langle Z \cup I \rangle)$ what is dimension of V?
   vi) Is Z is replaced by $C(Q)$ what is dimension of V?
   vii) If Z is replaced by Q what is dimension of V over Q?
   viii) If Z is replaced by $\langle Q \cup I \rangle$ what is dimension of V over $\langle Q \cup I \rangle$?
   ix) If Z is replaced by $C(\langle Q \cup I \rangle)$ what is dimension of V over $C(\langle Q \cup I \rangle)$?

Compare the dimension in (ii) to (viii) and derive conclusions based on it.

39. Let $V = \left\{ \begin{bmatrix} a_1 & a_2 & a_3 \\ a_4 & a_5 & a_6 \\ a_7 & a_8 & a_9 \end{bmatrix} \middle| a_i \in C (\langle Z \cup I \rangle), 1 \leq i \leq 9 \right\}$ be a group complex neutrosophic double linear algebra over the group Z.
   i) Find a basis of V.
   ii) Does V have complex neutrosophic pseudo double linear subalgebras?

40. Let $V = C (\langle Z_3 \cup I \rangle)$ be a semigroup under multiplication.
   i) Is V a S-semigroup?
   ii) Find order of V.
   iii) Find zero divisors if any in V.
   iv) Can V have S-zero divisors?
   v) Can V be a group?
   vi) Can V have S-units?
   vii) Can V have S-subsemigroups?

41. Let $M = \{a + bi + cI + idI \mid a, b, c, d \in Z_6\}$ be a semigroup under multiplication.
   Study questions (i) to (vii) suggested in problem 40 for V.

42. Prove $G = C (\langle Z_8 \cup I \rangle)$ can be a group under addition and only a semigroup under multiplication. Find order of G.

43. Let $G = C (\langle Z_7 \cup I \rangle)$;
   i) Is G a group?
   ii) Can G be a ring?
   iii) Can G be a field?
   iv) What is the strongest algebraic structure enjoyed by G?

44. Let $R = C (\langle Z_{10} \cup I \rangle)$ be a complex neutrosophic ring.
   i) Find the order of R.
   ii) Is R a S-ring?
   iii) Can R have S-subrings?
   iv) Can R have subrings which are not S-ideals?
   v) Can R have S-idempotents?

vi) Does R contain zero divisors which are not S-zero divisors?
vii) Define $\eta: R \to R$ so that $I = \ker \eta$ is a proper subring of R; so that R/I is a field.
viii) Can R have $Z_{10}$ as its subring?
ix) Is $C(Z_{10})$ a subring of R?
x) Can R have S-units?
xi) Does there exist an ideal I in R so that R/I is a field?

45. Let $V = C(\langle Z_3 \cup I \rangle) = \{a + bi + cI + idI \mid a, b, c, d \in Z_3\}$ be a complex neutrosophic vector space over the field $Z_3$.
    i) Find a basis for V.
    ii) Write V as a direct sum.
    iii) Write V as a pseudo direct sum.
    iv) Find subspaces of V so that V is neither a direct sum nor a pseudo direct sum.
    v) Find the algebraic structure enjoyed by $\text{Hom}_{Z_3}(V,V)$.
    vi) What is the algebraic structure of $L(V, Z_3)$?

46. Find some interesting properties of set complex neutrosophic vector space defined over the set S.

47. Study the special properties associated with the neutrosophic complex modulo integer ring $C(\langle Z_{90} \cup I \rangle)$.

48. What are the distinct properties of the complex modulo integer $C(Z_n)$; n not a prime?

49. Enumerate the properties of the complex modulo integer $C(Z_p)$; p a prime.

50. Find the zero divisors and units of $C(Z_{24})$.

51. Find an ideal I in $C(Z_{128})$ so that $C(Z_{128})/I$ is a field.

52. Does there exist an ideal I in $C(Z_{49})$ so that $C(Z_{49})/I$ is a field?

53. Find a subring S in $C(Z_n)$ so that S is not an ideal.

54. Is every complex modulo integer ring a S-ring?

55. Find a necessary and sufficient condition for a complex modulo integer ring $S = C(Z_n)$ to have ideals I such that the quotient ring is never a field.

56. Does every $C(Z_n)$ contain a zero divisor? Justify your claim.

57. Can every $C(Z_n)$ be a field?

58. Is every element in $C(Z_7)$ invertible?

59. Let $G = \left\{ \begin{bmatrix} a_1 & a_2 & a_3 & a_4 \\ a_5 & a_6 & a_7 & a_8 \\ a_9 & a_{10} & a_{11} & a_{12} \end{bmatrix} \middle| a_i \in C(Z_2) = \{a + bi_F \mid a, b \in Z_2, i_F^2 = 1\}; 1 \le i \le 12 \right\}$ be a complex modulo integer group under addition.
   i) Find the order of G?
   ii) Find three subgroups and verify Lagranges theorem for G.
   iii) Does G have p-Sylow subgroup?
   iv) Does G satisfy Cauchy theorem for finite abelian groups.
   v) What is the order of $A = \left\{ \begin{bmatrix} a & a & a & a \\ a & a & a & a \\ a & a & a & a \end{bmatrix} \text{ where } a \in C(Z_2) \right\}$?

60. Let $G = \left\{ \sum_{i=0}^{\infty} a_i x^i \middle| a_i \in C(Z_3) \right\}$ be a set.
   i) Can G be a group under addition?
   ii) Will G be a group under multiplication?
   iii) Will G be a semigroup of complex modulo integers under multiplication?
   iv) Can G have normal subgroup under +?

v) Is G a torsion free subgroup of torsion group?
vi) Can G have subgroups of finite order under +?

61. Let $S = \left\{ \begin{bmatrix} a_1 & a_2 & a_3 \\ a_4 & a_5 & a_6 \\ a_7 & a_8 & a_9 \end{bmatrix} \middle| a_i \in C(Z_{31}) = \{a + bi_F \mid a, b \in Z_{31}, i_F^2 = 30\} \right\}$ be a complex modulo integer semigroup under matrix multiplication.
   i) Find order of S.
   ii) Is S commutative?
   iii) Find ideals in S.
   iv) Is S a S-semigroup?
   v) Can S have subsemigroups which are not ideals?
   vi) Can S have zero divisors?
   vii) Can S have idempotents?
   viii) Give S-ideals in any in S.
   ix) Can S have nilpotent elements?

62. Let $S = \left\{ \begin{bmatrix} a_1 & a_2 & ... & a_{10} \\ a_{11} & a_{12} & ... & a_{20} \\ a_{21} & a_{22} & ... & a_{30} \end{bmatrix} \middle| a_i \in C(Z_{47}) = \{a + bi_F \mid a, \right.$
    $b \in Z_{47}, i_F^2 = 46\}; 1 \leq i \leq 30\}$ be a complex modulo integer vector space over the field $Z_{47}$.
    i) Find dimension of V over $Z_{47}$.
    ii) Find the order of V.
    iii) Write V as a direct sum.
    iv) Write V as a pseudo direct sum over $Z_{47}$.
    v) Define a linear operator T on V such that

    $W = \left\{ \begin{bmatrix} a_1 & a_2 & 0 & 0 & ... & 0 & a_{10} \\ a_3 & a_4 & 0 & 0 & ... & 0 & a_{20} \\ a_5 & a_6 & 0 & 0 & ... & 0 & a_{30} \end{bmatrix} \middle| a_i, a_j \in C(Z_{47}); \right.$
    $1 \leq i \leq 6, j=10, 20, 30\} \subseteq V$ so that $T(W) \subseteq W$.
    vi) Find S on V so that $S(W) \not\subseteq W$.

63. Let $V = \left\{ \begin{bmatrix} a_1 & a_2 & a_3 \\ a_4 & a_5 & a_6 \\ a_7 & a_8 & a_9 \end{bmatrix} \middle| a_i \in C(Z_{20}) = \{a + bi_F \mid a, b \in Z_{20}, i_F^2 = 19\}; 1 \le i \le 9 \right\}$ be a complex modulo integer group under addition.
   i) Find order of V.
   ii) Verify Lagrange's theorem.
   iii) What are the orders of subgroups of V?
   iv) Find quotient groups.
   v) Find order of $x = \begin{bmatrix} 3+i_F & 17 & 12+i_F \\ 7i_F & 13i_F+1 & 5+2i_F \\ 0 & 3+10i_F & 7 \end{bmatrix} \in V$.
   vi) Is Cauchy theorem true for x in V?

64. Let $V = \left\{ \begin{bmatrix} a_1 & a_2 & a_3 & a_4 \\ a_5 & a_6 & a_7 & a_8 \\ a_9 & a_{10} & a_{11} & a_{12} \\ a_{13} & a_{14} & a_{15} & a_{16} \end{bmatrix} \middle| a_i \in C(Z_{30}) = \{a + bi_F \mid a, b \in Z_{30}, i_F^2 = 29\}; 1 \le i \le 16 \right\}$ be a Smarandache complex modulo integer linear algebra over the S-ring $Z_{30}$.
   i) Find a basis of V.
   ii) Is V finite dimensional?
   iii) Is order of V finite?
   iv) Write V as a direct sum.
   v) Write V as a pseudo direct sum.
   vi) Find a sublinear algebra W of V such that for a linear operator T, $T(W) \subseteq W$.
   vii) Find a pseudo S-subvector space of V over $Z_{30}$.

65. Let $S = \left\{ \begin{bmatrix} a_1 & a_2 \\ a_3 & a_4 \\ \vdots & \vdots \\ a_{21} & a_{22} \end{bmatrix} \middle| a_i \in C(Z_3) = \{a + bi_F \mid a, b \in Z_3, i_F^2 = 2\}; 1 \le i \le 22 \right\}$ be a complex modulo integer group under addition.

   i) Find subgroups of S.
   ii) Can S be written as a interval direct sum of complex modulo integer subgroups?
   iii) Is every element in S is of finite order?
   iv) What is the order of S?
   v) Can S have p-Sylow subgroups?
   vi) Find all the p-Sylow subgroups of G.

66. Let $G = \left\{ \begin{bmatrix} a_1 & a_2 \\ a_3 & a_4 \end{bmatrix} \middle| a_i \in C(Z_{40}) = \{a + bi_F \mid a, b \in Z_{40}, i_F^2 = 39\}; 1 \le i \le 4 \right\}$ be a semigroup of complex modulo integers under matrix multiplication.

   i) Is G commutative? (Prove your claim)
   ii) What is the order of G?
   iii) Can G have S-subsemigroups?
   iv) Is G a S-semigroup?
   v) Can G have ideals?
   vi) Can G have S-idempotents?
   vii) Let $H = \left\{ \begin{bmatrix} a_1 & a_2 \\ a_3 & a_4 \end{bmatrix} \middle| a_1 a_4 - a_2 a_3 \ne 0, a_i \in C(Z_{40}); 1 \le i \le 4 \right\} \subseteq G$. Is H a subgroup of G.
   viii) What is the order of H?
   ix) Is H a commutative structure with respect to product?
   x) Can G have zero divisors which are not S-zero divisors?
   xi) Can G be written as a direct sum of subsemigroups?

67. Let $S = \left\{ \begin{bmatrix} a_1 & a_2 \\ a_3 & a_4 \end{bmatrix} \middle| a_i \in C(Z_{20}) = \{a + bi_F \mid a, b \in Z_{20}, i_F^2 = 19\}; 1 \le i \le 3 \right\}$ be a complex modulo integer ring of characteristic 20.
   i) Is S a commutative ring?
   ii) Is S a S-ring.?
   iii) Can S have S-ideals?
   iv) Can S have ideals which are not S-ideals?
   v) What is the order of S?
   vi) Can S have subrings which are not S-subrings?
   vii) Does S contain a maximal ideal?
   viii) Can S have zero divisors?
   ix) Can S have S-units?
   x) Does S contain subrings which are not S-ideals?
   xi) Find a homomorphism $\eta : S \to S$ so that
   $\ker \eta \ne \begin{pmatrix} 0 & 0 \\ 0 & 0 \end{pmatrix}$. Find S/ker $\eta$.

68. Let $M = \left\{ \sum_{i=0}^{\infty} a_i x^i \middle| a_i \in C(Z_5) = \{a + bi_F \mid a, b \in Z_5, i_F^2 = 4\} \right\}$ be a complex modulo integer polynomial ring.
   i) Is M a field ?
   ii) Can M be a S-ring?
   iii) Is M an integral domain?
   iv) Is M a principal ideal domain?
   v) Can M have irreducible polynomials?
   vi) Is the polynomial $p(x) = (3 + 2i_F) + (1 + 4i_F)x^2 + (3 + i_F)x^3 + 4i_F x^4$ in M a reducible polynomial order $C(Z_5)$?

69. Let $N = \left\{ \sum_{i=0}^{\infty} a_i x^i \middle| a_i \in C(Z_{12}) = \{a + bi_F \mid a, b \in Z_{12}, i_F^2 = 1\} \right\}$ be a complex modulo integer polynomial ring.
   i) Is N a field or an integral domain?
   ii) Is N a S-ring?
   iii) Can N have zero divisors?

iv) Is N a principal ideal domain?
v) Can N have S-subrings?
vi) Does N contain S-units?
vii) Can N have subrings which are not ideals?
viii) Can N have S-idempotents?
ix) Can N have ideals which are minimal?

70. Let $V = \left\{ \sum_{i=0}^{\infty} a_i x^i \,\middle|\, a_i \in C(Z_{16}) = \{a + bi_F \mid a, b \in Z_{16}, i_F^2 = 15\} \right\}$ be a complex modulo integer polynomial ring. Answer the problem (1) to (ix) of problem (69) for V.

71. Let $W = (C(Z_{12}) \times C(Z_{10}) \times C(Z_{43})) = \{a + i_F b, c + di_F, e + mi_F) \mid a, b \in Z_{12}, i_F^2 = 11, c, d \in Z_{10}, i_F^2 = 9$ and $e, m \in Z_{43}; i_F^2 = 42\}$.
   i) Is W a group of complex modulo integers under addition?
   ii) What is the order of W?
   iii) Is W a semigroup under product?
   iv) Is W as a semigroup of S-semigroup?
   v) Find zero divisors in W?
   vi) Find all p-Sylow subgroups of W treated as a group of complex modulo integers? (Is it possible?)
   vii) Can W be given a complex modulo integers ring structure?

72. Let $R = (Z_3 \times C(Z_5) \times C(Z_7) \times Z_{11}) = \{(a, b, c, d) \mid a \in Z_3, b \in C(Z_5) = \{x + i_F y \mid x, y \in Z_5, i_F^2 = 4\} c \in C(Z_7) = \{m + ni_F \mid m, n \in Z_7, i_F^2 = 6\}, d \in Z_{11}\}$ be a ring of modulo complex integers.
   i) Find order of R.
   ii) Can R have ideals?
   iii) Is R a S-ring?
   iv) Can R have S-units?
   v) Is R a principal ideal domain?
   vi) Can R have S-zero divisors?
   vii) Find ideals in R which are not S-ideal.

viii) Find an ideal I in R and determine R /I.
ix) Is it possible to have R/I to be a field?
x) Can R ever be an integral domain?

73. Let P = {$(Z_7 \times C(Z_{40}) \times C(Z_{11}))$ = (a, b, c) where a ∈ $Z_7$, b ∈ $C(Z_{40})$ = {d + $gi_F$ | d, g ∈ $Z_{40}$, $i_F^2$ = 39}, c ∈ $C(Z_{11})$ = {m + $ni_F$ | m, n ∈ $Z_{11}$; $i_F^2$ = 10}} be a group under addition.
   i) Mention some special properties enjoyed by P.
   ii) What is the order of P?
   iii) Is P abelian?
   iv) Find p-Sylow subgroups of P.
   v) Find subgroups of P which are only modulo integer subgroups and not complex modulo integer subgroups.

74. Let S = {$(Z_{20} \times C(Z_{12}) \times C(Z_{15}) \times Z_6)$ = (a, b, c, d) | a ∈ $Z_{20}$, b ∈ $C(Z_{12})$ = {e + $gi_F$ | e, g ∈ $Z_{12}$, $i_F^2$ = 11}, b ∈ $C(Z_{12})$ = {e + $gi_F$ | e, g ∈ $Z_{12}$; $i_F^2$ = 11}, c ∈ $C(Z_{15})$ = {m + $ni_F$ | m, n ∈ $Z_{15}$; $i_F^2$ = 14}, d ∈ $Z_6$} be a complex modulo integer semigroup under product.
   i) Find order of S.
   ii) Is S a S-semigroup?
   iii) Can S have zero divisors?
   iv) Can S have S-units?
   v) Find ideals in S.
   vi) Find S-subsemigroups which are not ideals in S.
   vii) Does S contain idempotents?
   viii) Can S have S-Lagrange subgroups?

75. Let T = {$(C(Z_3) \times C(Z_2) \times C(Z_6))$ = (a, b, c) | a ∈ $C(Z_3)$ = {x + $yi_F$ | x, y ∈ $Z_3$, $i_F^2$ = 2}, b ∈ $C(Z_2)$ = {m + $ni_F$ | m, n ∈ $Z_2$; $i_F^2$ = 1}, c ∈ $C(Z_6)$ = {t + $ui_F$ | t, u ∈ $Z_6$; $i_F^2$ = 5} be the ring of complex modulo integers.
   i) What is the order of T?
   ii) Is T a S-ring?
   iii) Find S-ideals if any in T.

iv) Find a subring in T which is not an ideal.
v) Is T a principal ideal domain?
vi) Can T be a unique factorization domain?
vii) Does T contain S-zero divisors?

76. Let M = {(C($Z_3$) × C($Z_4$) × C($Z_5$) × C($Z_6$) × C($Z_7$)) = (a, b, c, d, e) | a ∈ C($Z_3$) = {g + h$i_F$ | g, h ∈ $Z_3$, $i_F^2$ = 2}, b ∈ C($Z_4$) = {m + n$i_F$ | m, n ∈ $Z_4$; $i_F^2$ = 3}, c ∈ C ($Z_5$) = {t + u$i_F$ | t, u ∈ $Z_5$; $i_F^2$ = 4}, d ∈ C($Z_6$) = {r + s$i_F$ | r, s ∈ $Z_6$; $i_F^2$ = 5}, e ∈ C ($Z_7$) = {p + q$i_F$ | p, q ∈ $Z_7$; $i_F^2$ = 6}} be a group of complex modulo integers under addition.
    i) Find order of M.
    ii) Find subgroups of M.
    iii) Find all p-Sylow subgroups of M.
    iv) Prove every element in M is of finite order.
    v) Find an automorphism on M.

77. Can we have a unique factorization domain of complex modulo integers? Justify your claim.

78. Give an example of a complex modulo integers which is a principal ideal domain (Does it exist!).

79. What is the algebraic structure enjoyed by the complex modulo ring C($Z_{19}$)?

80. Obtain some interesting properties about complex modulo integer rings.

81. Obtain some unique properties enjoyed by rings built using the complex modulo integers C($Z_n$).

82. Is C(⟨Z ∪ I⟩) a unique factorization domain?

83. Can C(⟨R ∪ I⟩) be a principal ideal domain?

84. If $A = \left\{ \begin{bmatrix} a & b \\ c & d \end{bmatrix} \middle| a, b, c, d \in C(\langle Q \cup I \rangle); ad - bc \neq 0 \right\}$ be a

 ring of neutrosophic complex modulo integers.
 i) Is A commutative?
 ii) Is A a S-ring?
 iii) Can A have S-ideals?
 iv) Can A have zero divisors?
 v) Can A have S- subrings which are not ideals?
 vi) Can the concept of a.c.c. or d.c.c. to imposed on A?

85. Let $M = \left\{ \sum_{i=0}^{\infty} a_i x^i \middle| a_i \in C(\langle Z \cup I \rangle) \right\}$ be the neutrosophic

 complex polynomial ring.
 i) Find subrings in M which are not ideals?
 ii) Is M a S-ring?
 iii) Study the concept of irreducible and reducible polynomials in M.
 iv) Find a linearly reducible polynomial in M.
 v) Find a irreducible polynomial of degree three in M.
 vi) Will an irreducible polynomial in M generate a maximal ideal?
 vii) Is M a principal ideal domain?
 viii) Is M a unique factorization domain?

86. Let $R = \left\{ \sum_{i=0}^{\infty} a_i x^i \middle| a_i \in C(\langle Q \cup I \rangle) \right\}$ be a neutrosophic

 complex polynomial ring. Answer the questions (i) to (viii) mentioned in problem (85).

87. Let $F = \left\{ \sum_{i=0}^{\infty} a_i x^i \middle| a_i \in C \right\}$ be a ring. Find the differences

 between M in problem 82, R in problem (83) and F in this problem.

88. Suppose $P = \left\{ \sum_{i=0}^{\infty} a_i x^i \mid a_i \in C(\langle R \cup I \rangle) \right\}$ be the real neutrosophc ring of polynomials.
    i)  Study / answer questions (i) to (viii).
    ii) Compare P with R, F and M in problems (82) to (84).

89. Mention the properties distinctly associated with complex modulo integer vector spaces defined over $Z_p$, p a prime.

90. Find interesting properties enjoyed by Smarandache complex modulo integer linear algebras defined over the S-ring $Z_n$, n not a prime.

91. Let $P = \left\{ \begin{bmatrix} a_1 & a_2 & \ldots & a_{12} \\ a_{13} & a_{14} & \ldots & a_{24} \\ a_{25} & a_{26} & \ldots & a_{36} \\ a_{37} & a_{38} & \ldots & a_{48} \end{bmatrix} \middle| a_i \in C(Z_{24}) = \{a + bi_F \mid a, \right.$

    $\left. b \in Z_{24}, i_F^2 = 23\}, 1 \le i \le 48 \right\}$ be a Smarandache complex integer vector space over the S-ring $Z_{24}$.
    i)    Find a basis for P over $Z_{24}$.
    ii)   What is the dimension of P over $Z_{24}$?
    iii)  Find subspaces of P.
    iv)   Write P as a direct sum.
    v)    Write P as a pseudo direct sum.
    vi)   Find a linear operator T on P so that $T^{-1}$ exists.
    vii)  Does P contain pseudo S-complex modulo integer vector subspaces over a field $F \subseteq Z_{24}$?

92. Let $V = \left\{ \sum_{i=0}^{\infty} a_i x^i \middle| a_i \in C(\langle Z_7 \cup I \rangle) = \{a + bi_F + cI + i_F I d \mid a, \right.$

    $\left. b, c, d \in Z_7, i_F^2 = 6, I^2 = I, (i_F I)^2 = 6I \} \right\}$ be the neutrosophic complex polynomial ring.
    i)  Is V a S-ring?
    ii) Can V have reducible polynomials?

iii) Give examples of irreducible polynomials?
iv) Is V an integral domain?
v) Can V have zero divisors?
vi) Is V a principal ideal domain?
vii) Can V have subrings which are not ideals?

93. Let $W = \left\{ \sum_{i=0}^{\infty} a_i x^i \;\middle|\; a_i \in C(\langle Z_6 \cup I \rangle) = \{a + bi_F \mid a, b \in Z_6,$

   $i_F^2 = 5\}\}$ be a neutrosophic complex modulo integer ring.
   i) Find zero divisors in W.
   ii) Find two relatively prime polynomials in W.
   iii) Find ideals in W.
   iv) Is W a S-ring?
   v) Find subrings in W which are not ideals?
   vi) Can W have S-units?
   vii) Can W have S-zero divisors?
   viii) Is W a principal ideal ring?
   ix) Can W have maximal ideals?
   x) Does W contain minimal ideals?
   xi) Will a.c.c. condition on ideals be true in W?
   xii) Does W contain pure neutrosophic ideals?
   xiii) Can W have complex modulo integer ideals?

94. Let $G = \left\{ \sum_{i=0}^{\infty} a_i x^i \;\middle|\; a_i \in C(\langle Z_2 \cup I \rangle) = \{a + bi_F + cI + i_F dI \mid$

   $a, b, c, d \in Z_2,\; i_F^2 = 1\}\}$ be a group under addition.
   i) Find order of G.
   ii) Can G have Sylow subgroups?
   iii) Can G be written as a direct product of subgroups?
   iv) Let $H = \left\{ \sum_{i=0}^{10} a_i x^i \;\middle|\; a_i \in C(\langle Z_2 \cup I \rangle) \right\} \subseteq G$. Find G / H.
   v) What is the order of G / H?
   vi) Find $\sum_{i=0}^{19} a_i x^i H$.

95. Let $F = \left\{ \begin{bmatrix} a_1 & a_2 \\ a_3 & a_4 \end{bmatrix} \middle| a_i \in C(\langle Z_{12} \cup I \rangle) = \{a + bi_F \mid a, b \in Z_{12}, i_F^2 = 1\}, 1 \leq i \leq 4 \right\}$ be a ring.
   i) Is F a S-ring?
   ii) Is F is of finite order?
   iii) Find subrings in F which are not ideals.
   iv) Is F commutative?
   v) Find a left ideal in F which is not a right ideal in F.
   vi) Find S-zero divisors if any in F.
   vii) Can F be a ring with S-units?
   viii) What is the difference between F and
   $$H = \left\{ \begin{bmatrix} a_1 & a_2 \\ a_3 & a_4 \end{bmatrix} \middle| a_i \in Z_{12} \right\}?$$
   ix) Find some special properties enjoyed by F and not by
   $$T = \left\{ \begin{bmatrix} a & b \\ c & d \end{bmatrix} \middle| a, b, c, d \in C(Z_{12}) = \{a + bi_F \mid a, b \in Z_{12}, i_F^2 = 11\} \right\}.$$

96. Let $P = \left\{ \begin{bmatrix} a & b \\ c & d \end{bmatrix} \middle| a, b, c, d \in C(\langle Z_9 \cup I \rangle) = \{a + bi_F + cI + i_F dI \mid a, b, c, d \in Z_9, i_F^2 = 8, I^2 = I, (i_F I)^2 = 8I\} \right\}$ be a complex neutrosophic modulo integer semigroup under multiplication.
   i) Find order of P.
   ii) Is P commutative?
   iii) Find zero divisors if any in P.
   iv) Is P a S-semigroup?
   v) Can P have S-subsemigroups?
   vi) Find right ideals which are not left ideals in P.
   vii) Find S-units if any in P.

97. Let $M = \left\{ \begin{bmatrix} a & b \\ c & d \end{bmatrix} \middle| \, a, b, c, d \in C(\langle Z_{11} \cup I \rangle) = \{m + ni_F + tI + i_F sI \mid m, n, t, s \in Z_{11}, i_F^2 = 10, I^2 = I, (i_F I)^2 = 10I\}; ad - bc \neq 0 \right\}$ be a group under multiplication.
   i) Find order of M.
   ii) Prove M is non commutative.
   iii) Find normal subgroups in M.
   iv) Find subgroups which are not normal in M.
   v) Is Cauchy theorem true in M?
   vi) Show M can be embedded in a symmetric group $S_n$, give that n? (Cauley's theorem).

98. Find some interesting properties associated with neutrosophic complex modulo integer semigroups.

99. Let $S = \{(a_1, a_2, a_3, a_4, a_5) \mid a_i \in C(\langle Z_{24} \cup I \rangle) = \{a + bi_F + cI + i_F dI \mid a, b, c, d \in Z_{24}, i_F^2 = 23, I^2 = I, (i_F I)^2 = 23I\}; 1 \leq i \leq 5\}$ be a complex neutrosophic modulo integer semigroup under multiplication.
   i) Find order of S.
   ii) Find zero divisors in S.
   iii) Find ideals in S.
   iv) Is S a S-ring?
   v) Find S-subsemirings in S.

100. Suppose S in problem (98) is taken as addition what are the relevant differences you can find?

101. Let $V = \left\{ \begin{bmatrix} x_1 \\ x_2 \\ \vdots \\ x_{10} \end{bmatrix} \middle| \, x_i \in C(\langle Z_4 \cup I \rangle) = \{a + bi_F \mid a, b \in Z_4, i_F^2 = 3, I^2 = I, (i_F I)^2 = 3I\}; 1 \leq i \leq 10 \right\}$ be a group of complex neutrosophic modulo integers under addition.
   i) Find subgroups.

ii) Find order of V.
iii) Find p-Sylow subgroups of V.
iv) Write V as a direct sum.

102. Let P = {(C ($\langle Z_7 \cup I \rangle$)) × C ($\langle Z_{11} \cup I \rangle$) × C ($\langle Z_3 \cup I \rangle$)) = (a, b, c) | a ∈ C ($\langle Z_7 \cup I \rangle$) = {d + gi$_F$ + hI + i$_F$ kI | d, g, h, k ∈ Z$_7$, i$_F^2$ = 6, I$^2$ = I, (i$_F$I)$^2$ = 6I}} b ∈ C ($\langle Z_{11} \cup I \rangle$) = {m + ni$_F$ + tI + i$_F$ sI | m, n, t, s ∈ Z$_{11}$, i$_F^2$ = 10, I$^2$ = I, (i$_F$I)$^2$ = 10I} and C ∈ C ($\langle Z_3 \cup I \rangle$) = {(a + bi$_F$ + eI + i$_F$ dI) | i$_F^2$ = 2, I$^2$ = I, (i$_F$I)$^2$ = 2I}} be a neutrosophic complex modulo integer semigroup under product.
   i) Find order of P.
   ii) Find ideal of P.
   iii) Is P a S-semigroup?
   iv) Is T = T = {Z$_7$ × Z$_{11}$ × Z$_3$} ⊆ P a pseudo modulo integer subsemigroup of P? Can T be an ideal of P?
   v) Can W = {C(Z$_7$) × C (Z$_{11}$) × C (Z$_3$)} ⊆ P be a pseudo complex ideal of P?
   vi) Can B = {$\langle Z_7 \cup I \rangle$ × C ($\langle Z_{11} \cup I \rangle$) × C ($\langle Z_3 \cup I \rangle$)} ⊆ P be a pseudo neutrosophic ideal of P?

103. Suppose P in problem 101 is taken as a neutrosophic complex modulo integer group under addition then study the basic properties associated with P. Compare P as a group under + and semigroup under multiplication.

104. Let M = {Z$_7$ × C (Z$_5$) × C ($\langle Z_{12} \cup I \rangle$) × ($\langle Z_{10} \cup I \rangle$) = {(a, b, c, d) | a ∈ Z$_7$, b ∈ C (Z$_5$) = a + pi$_F$; a, p ∈ Z$_5$, i$_F^2$ = 4}; c ∈ C ($\langle Z_{12} \cup I \rangle$) = {a + bi$_F$ + dI + ei$_F$I | a, b, d, e ∈ Z$_{12}$, i$_F^2$ = 11, I$^2$ = I, (i$_F$I)$^2$ = 11I}, d ∈ $\langle Z_{10} \cup I \rangle$ = {a + bI | a, b ∈ Z$_{10}$, I$^2$ = I}} be a semigroup under multiplication.
   i) Find order of M.
   ii) Find zero divisors in M.
   iii) Is M a S-semigroup?
   iv) Find S-ideals if any in M.
   v) Can M have S-units?

105. Let $L = \left\{ \begin{bmatrix} a_1 & a_2 & \ldots & a_9 \\ a_{10} & a_{11} & \ldots & a_{18} \\ a_{19} & a_{20} & \ldots & a_{27} \end{bmatrix} \middle| a_i \in C(Z_{40}) = \{a + bi_F \mid a, b \in Z_{40}, i_F^2 = 39\}; 1 \le i \le 27 \right\}$ be a complex modulo integer group under addition.

i) Find all p-Sylow subgroups of L.
ii) Can L be written as a direct union sum of subgroups?
iii) If $H = \left\{ \begin{bmatrix} a_1 & 0 & \ldots & 0 & a_4 \\ a_2 & 0 & \ldots & 0 & a_5 \\ a_3 & 0 & \ldots & 0 & a_6 \end{bmatrix} \middle| a_i \in C(Z_{40}); 1 \le i \le 6 \right\}$
$\subseteq L$. Is H a subgroup? If so find the coset of H in L.

106. Let $M = \left\{ \begin{bmatrix} a_1 & a_2 & a_3 \\ a_4 & a_5 & a_6 \\ a_7 & a_8 & a_9 \end{bmatrix} \middle| a_i \in C(\langle Z_{12} \cup I \rangle) = \{a + bi_F + cI + i_F dI \mid a, b, c, d \in Z_{12}, i_F^2 = 11\}, I^2 = I; (Ii_F)^2 = 11I\}; 1 \le i \le 9 \right\}$ be a neutrosophic complex modulo integer ring.

i) Is M a S-ring?
ii) Write S-units in M.
iii) Find ideals in M.
iv) Find S-ideals if any in M.
v) Find subrings which are not ideals.
vi) What is the order of M?
vii) Does M contain S-subrings which are not S-ideals?
viii) Find a right ideal which is not a left ideal of M.

107. Let $T = \left\{ \begin{bmatrix} a_1 & a_2 & a_3 & a_4 \\ a_5 & a_6 & a_7 & a_8 \\ a_9 & a_{10} & a_{11} & a_{12} \\ a_{13} & a_{14} & a_{15} & a_{16} \end{bmatrix} \middle| a_i \in C(\langle Z_{22} \cup I \rangle) = \{a + bi_F + cI + i_F dI \mid a, b, c, d \in Z_{22}, i_F^2 = 21, I^2 = I; (i_F I)^2 = 21I\}; 1 \le i \le 16 \right\}$ be a neutrosophic complex modulo integer semigroup.

i)   Find order of T.
   ii)  Find ideals in T.
   iii) Find subrings in T which are not ideals of T.
   iv)  Find a right ideal of T and left ideal of T.
   v)   Find S-zero divisors if any in T.
   vi)  Is T a S-semigroup?
   vii) Can T have idempotents?

108. Let B = $\{(x_1, x_2, x_3, x_4, x_5, x_6) \mid x_i \in C(\langle Z_{12} \cup I \rangle) = \{a + bi_F + cI + i_F dI \mid a, b, c, d \in Z_{12}, i_F^2 = 11, I^2 = I; (i_F I)^2 = 11I\}; 1 \leq i \leq 6\}$ be a ring of neutrosophic complex modulo integers.
   i)    Find order of B.
   ii)   Find ideals of B.
   iii)  Can B have subrings which are not ideals?
   iv)   Can B have S-ideals?
   v)    Is B a S-ring?
   vi)   Can B have S-subrings?
   vii)  Can I = $\{(x_1\ 0\ x_2\ 0\ x_3\ 0) \mid x_i \in C(\langle Z_{12} \cup I \rangle); 1 \leq i \leq 3\} \subseteq B$ be an ideal?
   viii) Find B/I.

109. Let M = $\left\{ \begin{bmatrix} a & b \\ c & d \end{bmatrix} \middle| a, b, c, d \in C(\langle Z_{14} \cup I \rangle) = \{a + bi_F + cI + i_F dI \mid i_F^2 = 13, I^2 = I; (i_F I)^2 = 13I\} \right\}$ be a ring of neutrosophic complex modulo integer.
   i)    Is M a S-ring?
   ii)   Find order of M.
   iii)  Find S-ideals if any in M.
   iv)   Can M have S-zero divisors?
   v)    Find idempotents if any in M.
   vi)   Find a left ideal of M which is not a right ideal of M and vice versa.
   vii)  Does M have a zero divisor which is not a S-zero divisor?
   viii) Find a S-subring which is not an ideal of M.

110. Let $T = \left\{ \sum_{i=0}^{10} a_i x^i \mid a_i \in C(\langle Z_3 \cup I \rangle) = \{a + bi_F + cI + i_F dI \mid a, b, c, d \in Z_3; i_F^2 = 2, I^2 = I; (i_F I)^2 = I\} \right\}$ be a semigroup (or group under addition).
   i) Find order of T.
   ii) What are the stricking differences when T is considered as a group and as a semigroup?
   iii) Find all p-Sylow subgroups of T.
   iv) Find $S_n$ (exact n) so that T is embeddable in $S_n$.
   v) What are the distinct features as a group T enjoys?
   vi) Can T have ideals, (T as a semigroup)?

111. Let M = {all 10 × 10 upper triangular matrices with entries from $C(\langle Z_{25} \cup I \rangle) = \{a + bi_F + cI + i_F dI \mid a, b, c, d \in Z_{25}; i_F^2 = 24, I^2 = I; (i_F I)^2 = 24I\}$} be a neutrosophic complex modulo integer ring.
   i) Is M commutative?
   ii) Is M a S-ring?
   iii) Is M finite?
   iv) Find subrings in M which are not ideals?
   v) Find S-ideals if any.
   vi) Does M have zero divisors?
   vii) Can M have S-units?
   viii) Find any special property enjoyed by M.

112. Let $V = \{(x_1, x_2, \ldots, x_{10}) \mid x_i \in C(\langle Z_{11} \cup I \rangle) = \{a + bi_F + cI + i_F dI \mid a, b, c, d \in Z_{11}, i_F^2 = 10, I^2 = I; (i_F I)^2 = 10I\}\}$ be a neutrosophic complex modulo integer linear algebra over the field $Z_{11}$.
   i) Find subspaces of V.
   ii) Write V as a direct sum of subspaces (linear subalgebras).
   iii) Find a basis of V.
   iv) If V is defined over $\langle Z_{11} \cup I \rangle$. find the basis.
   v) If V is defined over $C(Z_{11})$ what is the dimension of V over $C(Z_{11})$?

vi) What is dimension of V if V is defined over $C(\langle Z_{11} \cup I \rangle)$?

113. Let $M = \left\{ \begin{bmatrix} a_1 & a_2 & \ldots & a_8 \\ a_9 & a_{10} & \ldots & a_{16} \end{bmatrix} \middle| a_i \in C(\langle Z_{23} \cup I \rangle) = \{a + bi_F + cI + i_F dI \mid a, b, c, d \in Z_{23}, i_F^2 = 22, I^2 = I; (i_F I)^2 = 22I\}; ; 1 \leq i \leq 16 \right\}$ be a neutrosophic complex modulo integer vector space over $Z_{23}$.

   i) Find a basis of M over $Z_{23}$.
   ii) What is dimension of M over $Z_{23}$?
   iii) Write M as a direct sum.
   iv) Write M as a pseudo direct sum.
   v) Let $W = \left\{ \begin{bmatrix} a_1 & 0 & a_2 & 0 & a_3 & 0 & a_4 & 0 \\ 0 & a_5 & 0 & a_6 & 0 & a_7 & 0 & a_8 \end{bmatrix} \middle| a_i \in C(\langle Z_{23} \cup I \rangle) \right\} 1 \leq i \leq 8\} \subseteq M$ be a subspace of M. Find a linear operator on V so that $T(W) \subseteq W$.
   vi) Study (1) and (ii) if $Z_{23}$ if replaced by $C(\langle Z_{23} \cup I \rangle)$.

114. Let $P = \left\{ \begin{pmatrix} a_1 & a_2 \\ a_3 & a_4 \end{pmatrix} \middle| a_i \in C(\langle Z_{12} \cup I \rangle) = \{a + bi_F + cI + i_F dI \mid a, b, c, d \in Z_{12}, i_F^2 = 11, I^2 = I; (i_F I)^2 = 11I\}; ; 1 \leq i \leq 4 \right\}$ be a S-neutrosophic complex modular integer linear algebra over the S-ring $Z_{12}$.

   i) Find a basis of P.
   ii) What is the dimension of P over $Z_{12}$?
   iii) Write P as a direct sum.
   iv) Write P as a pseudo direct sum.
   v) Find the algebraic structure enjoyed by $\text{Hom}_{Z_{12}}(V,V)$.
   vi) Give a linear operator T such that $T^{-1}$ exists and keeps no subspace of P invariant under it.
   vii) Does such a linear operator exist?

115. Let $V = \left\{ \begin{bmatrix} a_1 & a_2 \\ a_3 & a_4 \end{bmatrix}, (a_1 \ a_2 \ a_3), \sum_{i=0}^{5} a_i x^i \;\middle|\; a_i \in C \right.$

$(\langle Z_{20} \cup I \rangle) = \{a + bi_F + cI + i_F dI \mid i_F^2 = 19, I^2 = I; (i_F I)^2 = 19I\}; \; 0 \le i \le 5\}$ be a set neutrosophic complex integer modulo vector space over the set $S = \{0, 2, 5, 10, 11\} \subseteq Z_{20}$.

i) Find the number of elements in V.
ii) Find a basis of V over S.
iii) Does V have subspaces?
iv) Can V be written as a direct sum? If so do it.
v) Write V as a pseudo direct sum.

116. Let P =

$\left\{ \sum_{i=0}^{5} a_i x^i, \begin{bmatrix} a_1 \\ a_2 \\ a_3 \\ a_4 \\ a_5 \\ a_6 \\ a_7 \end{bmatrix}, \begin{bmatrix} a_1 & a_3 & a_5 & a_7 & a_9 \\ a_2 & a_4 & a_6 & a_8 & a_{10} \end{bmatrix}, \begin{bmatrix} a_1 & a_2 & a_3 \\ a_4 & a_5 & a_6 \\ a_7 & a_8 & a_9 \end{bmatrix}, \right.$

$\left. \begin{bmatrix} a_1 & a_2 & a_3 & a_4 \\ a_5 & a_6 & a_7 & a_8 \\ a_9 & a_{10} & a_{11} & a_{12} \\ a_{13} & a_{14} & a_{15} & a_{16} \\ a_{17} & a_{18} & a_{19} & a_{20} \\ a_{21} & a_{22} & a_{23} & a_{24} \end{bmatrix} \;\middle|\; a_i \in C \;(\langle Z_{11} \cup I \rangle) = \{a + bi_F + cI + i_F dI \mid a, b, c, d \in Z_{11}, i_F^2 = 10, I^2 = I; (i_F I)^2 = 10I\}; \; 0 \le i \le 24 \right\}$ be a set neutrosophic complex vector space of modulo integers over the set $S = \{0, 3, 5, 7\} \subseteq Z_{11}$.

i) Find dimension of P over S.
ii) Find a basis of P over S.
iii) Write P as a direct sum.

iv) Find T : P → P such that T keep atleast two subspaces invariant.
v) Prove P can have pure neutrosophic subspaces over S.
vi) Find complex modulo integer set vector subspaces of P over S.
vii) Find just modulo integer vector subspaces of P over S.
viii) Find the structure enjoyed by $Hom_S$ (P, P).
ix) Can a linear operator on P be such that it keeps every subspace invariant?

117. Let M = $\left\{ \begin{bmatrix} a_1 & a_2 & a_3 & a_4 & a_5 & a_6 \\ a_7 & a_8 & a_9 & a_{10} & a_{11} & a_{12} \\ a_{13} & a_{14} & a_{15} & a_{16} & a_{17} & a_{18} \end{bmatrix} \middle| a_i \in C \right.$

$(\langle Z_{43} \cup I \rangle) = \{a + bi_F + cI + i_F dI \mid a, b, c, d \in Z_{43}, i_F^2 = 42,$ $I^2 = I; (i_F I)^2 = 42I\}; \ 1 \le i \le 18\}$ be a set neutrosophic complex linear algebra over the set S = {0, 1, 4, 41, 9, 42} ⊆ $Z_{43}$.
i) Find a basis for M over S.
ii) Write M as a direct sum of sublinear algebras over S.
iii) Find dimension of M over S.
iv) Find two disjoint sublinear algebras of M so that the intersection is the 3 × 6 zero matrix.
v) Find a linear operator T on M so that $T^{-1}$ does not exist.

118. Let V = $\left\{ \begin{bmatrix} a_1 & a_2 & a_3 \\ a_4 & a_5 & a_6 \\ a_7 & a_8 & a_9 \\ a_{10} & a_{11} & a_{12} \\ a_{13} & a_{14} & a_{15} \end{bmatrix} \middle| a_i \in C \ (\langle Z_{23} \cup I \rangle) = \{a + bi_F + cI \right.$

$+ i_F dI \mid a, b, c, d \in Z_{23}, i_F^2 = 22, I^2 = I; (i_F I)^2 = 22I\}; \ 1 \le i \le 15\}$ be a set neutrosophic complex modulo integer linear algebra over the set S = {0, 1, 20, 19, 21, 3, 5, 7} ⊆ $Z_{23}$.

i) Find a basis for V.
ii) What is the dimension of V over S?
iii) Write V as a direct sum of set sublinear algebras.
iv) Study $Hom_S (V, V)$.
v) Study $L (V, S)$.
vi) Introduce and study any other properties related with V.

119. Let $G = \left\{ \begin{bmatrix} a_1 & a_3 & a_5 & a_7 & a_9 \\ a_2 & a_4 & a_6 & a_8 & a_{10} \end{bmatrix} \middle| a_i \in C (\langle Z_5 \cup I \rangle) = \{a + bi_F + cI + i_F dI \mid a, b, c, d \in Z_5, i_F^2 = 4, I^2 = I; (i_F I)^2 = 4I\}; 1 \le i \le 10 \right\}$ be a group of neutrosophic complex modulo integers under addition.
i) What is the order of G?
ii) Find all p-Sylow subgroups of G.
iii) Let $H = \left\{ \begin{bmatrix} a_1 & a_2 & 0 & 0 & 0 \\ 0 & 0 & a_3 & a_4 & 0 \end{bmatrix} \middle| a_i \in C(<Z_5 \cup I>); 1 \le i \le 4 \right\}$
$\subseteq$ be a subgroup of G. find G / H. What is o (G/H)?

120. Let $P = \left\{ \begin{bmatrix} a & a \\ a & a \end{bmatrix} \middle| a \in C (\langle Z_{25} \cup I \rangle) = \{a + bi_F + cI + i_F dI \mid a, b, c, d \in Z_{25}, i_F^2 = 24, I^2 = I; (i_F I)^2 = 24I\} \right\}$ be a group under addition of neutrosophic complex modulo integers.
i) Find order of P.
ii) Obtain all p-Sylow subgroups of P.
iii) Give a subgroup of P which is not a p-Sylow subgroup.
iv) Can P be a group under product?
v) Is it possible to write P as a direct sum of subgroups?
vi) Can $\eta : P \to P$ be such that ker $\eta$ is nontrivial? ($\eta$ - a group homomorphism).

121. Give any nice and interesting property enjoyed by neutrosophic complex modulo integer groups.

122. Obtain some classical theorems of finite groups in case of complex modulo neutrosophic integer groups.

123. Construct a class of neutrosophic complex modulo integer semigroups which are not S-semigroups.

124. Give an example of a complex neutrosophic modulo integer semigroup which is a S-semigroup.

125. Does there exist rings built using neutrosophic complex modulo integers which is not a S-ring?

126. Study the algebraic structure of $\text{Hom}_{Z_p}(V, V)$ where V is a neutrosophic complex modulo integer vector space over $Z_p$.

127. Study the algebraic structure of $L(V, Z_p)$ where V is the neutrosophic complex modulo integer vector space over $Z_p$.

128. Let $V = \left\{ \begin{bmatrix} a_1 & a_2 & ... & a_{10} \\ a_{11} & a_{12} & ... & a_{20} \\ a_{21} & a_{22} & ... & a_{30} \end{bmatrix} \middle| a_i \in C(\langle Z_{43} \cup I \rangle) = \{a + bi_F + cI + i_F dI \mid a, b, c, d \in Z_{43}, i_F^2 = 42, I^2 = I \text{ and } (i_F I)^2 = 42I\}; 1 \leq i \leq 30 \right\}$ be a neutrosophic complex modulo integer vector space over the field $Z_{43}$. Find $L(V, Z_{43})$.

129. Let $V = \left\{ \begin{bmatrix} a_1 & a_2 & ... & a_6 \\ a_7 & a_8 & ... & a_{12} \end{bmatrix}, (a_1, a_2, ..., a_{10}), \begin{bmatrix} a_1 \\ a_2 \\ \vdots \\ a_7 \end{bmatrix} \middle| a_i \right.$ $\in C(\langle Z_{47} \cup I \rangle) = \{a + bi_F + cI + di_F I \mid a, b, c, d \in Z_{47}, i_F^2 = 46, I^2 = I \text{ and } (i_F I)^2 = 46I\}; 1 \leq i \leq 12 \right\}$ be a set complex

neutrosophic modulo integer vector space over the set S = {0, 42I, 3+7I, 27I, 5+40I, 1}.
i) Find a basis.
ii) Is V finite or infinite dimensional over S?
iii) Write V as a direct sum.
iv) Write V as a pseudo direct sum.
v) Find $Hom_S$ (V,V).
vi) Find atleast one linear operator T on V so that $T^{-1}$ exists.

130. Find some nice applications of set complex neutrosophic modulo integer vector spaces defined over a set S.

131. What is the advantage of using set complex neutrosophic modulo integer strong linear algebras over the set S?

132. Let V be a set neutrosophic complex modulo integer strong linear algebra defined over a set S.
Find $Hom_S$ (V,V).

133. Let M = $\left\{ \begin{pmatrix} a_1 & a_2 \\ a_3 & a_4 \end{pmatrix} \middle| a_i \in C (\langle Z_3 \cup I \rangle), 1 \leq i \leq 4 \right\}$ where C $(\langle Z_3 \cup I \rangle) = \{a + bi_F + cI + i_F dI \mid a, b, c, d \in Z_3, i_F^2 = 2, I^2 = I$ and $(i_F I)^2 = 2I\}$ be a group complex neutrosophic modulo integer strong linear algebra over the group G = $Z_3$.
i) Find a basis of M over G.
ii) Show the number of elements in M is finite (Find o (M)).
iii) Write M as a pseudo direct sum.
iv) Is it possible to write M as a direct sum of strong linear subalgebras?

134. Find some interesting properties associated with group modulo integer neutrosophic complex vector spaces defined over a group G.

135. Compare group neutrosophic complex modulo integer vector spaces with group neutrosophic complex modulo integer strong linear algebras.

136. Find the difference between the semigroup complex neutrosophic modulo integer vector spaces and semigroup complex neutrosophic modulo integer strong linear algebras.

137. Let $V = \{(a_1, a_2, \ldots, a_9) \mid a_i \in C(\langle Z_5 \cup I \rangle) = \{a + bi_F + cI + i_F dI \mid a, b, c, d \in Z_5, i_F^2 = 4, I^2 = I \text{ and } (i_F I)^2 = 4I\}; 1 \leq i \leq 9\}$ be a semigroup neutrosophic complex modulo integer strong linear algebra over the semigroup $S = C(Z_5) = \{a + i_F b \mid a, b \in Z_5\}$ under addition.
   i) Find the number of elements in V.
   ii) Find a basis of V.
   iii) What is the dimension of V over S?
   iv) Write V as a direct sum of sublinear algebras.
   v) Write V as a pseudo direct sum.

138. Find the algebraic structure enjoyed by $\text{Hom}_S(V,V)$; V a set complex neutrosophic modulo integers over the set S.

139. If V is a neutrosophic complex modulo integer linear algebra over $Z_p$ find $\text{Hom}_{Z_p}(V, V)$.

140. Let $M = \left\{ \begin{bmatrix} a_1 & a_2 & a_3 & a_4 & a_5 \\ a_6 & a_7 & a_8 & a_9 & a_{10} \\ a_{11} & a_{12} & a_{13} & a_{14} & a_{15} \\ a_{16} & a_{17} & a_{18} & a_{19} & a_{20} \end{bmatrix} \right.$ $a_i \in C(\langle Z_2 \cup I \rangle) = \{a + bi_F + cI + i_F dI \mid a, b, c, d \in Z_2, i_F^2 = 1, I^2 = I \text{ and } (i_F I)^2 = I\}; 1 \leq i \leq 20\}$ be a complex neutrosophic modulo integer vector space over the field $Z_2$.
   i) Find dimension of M over $Z_2$.
   ii) What is the order of M?
   iii) Find a basis of M over $Z_2$.
   iv) Is M a linear algebra?

v) Can M be a strong linear algebra?
vi) Write M as a direct sum of subspaces.
vii) Find $\text{Hom}_{Z_2}$ (M, M).
viii) Find L (M, $Z_2$).
ix) Write M as a pseudo direct sum of subspaces.
x) If $Z_2$ is replaced by C ($\langle Z_2 \cup I \rangle$) what is the structure of M?

141. Let $M = \left\{ \begin{bmatrix} a_1 & a_2 & \ldots & a_{10} \\ a_{11} & a_{12} & \ldots & a_{20} \\ a_{21} & a_{22} & \ldots & a_{30} \end{bmatrix} \middle| a_i \in C(\langle Z_{11} \cup I \rangle) = \{a + bi_F + cI + i_F dI \mid a, b, c, d \in Z_{11}, i_F^2 = 10, I^2 = I \text{ and } (i_F I)^2 = 10I\}; 1 \leq i \leq 30 \right\}$ be a set complex neutrosophic modulo integer linear algebra over the set $S = \{0, 1\}$.

i) Find a basis of M.
ii) Find the number of elements in M.
iii) Write M as a direct sum of subspaces.
iv) Find $\text{Hom}_S$ (M, M).
v) Find $T \in \text{Hom}_S$ (M, M) such that $T^{-1}$ does not exist.

142. If S in problem 140 is replaced by $Z_{11}$, answer questions (1) to (v).

143. If S in problem 140 is replaced by C ($Z_{11}$) then study the questions (1) to (v).

144. If S in problem 140 is replaced by C ($\langle Z_{11} \cup I \rangle$), will S be a vector space?

145. Let $V = \left\{ \begin{bmatrix} a_1 \\ a_2 \\ a_3 \\ a_4 \\ a_5 \\ a_6 \end{bmatrix} \middle| a_i \in C(\langle Z_{20} \cup I \rangle) = \{a + bi_F + cI + i_F dI \mid a, \right.$
b, c, d $\in Z_{20}$, $i_F^2 = 19$, $I^2 = I$ and $(i_F I)^2 = 19I\}$; $1 \leq i \leq 6\}$

be a Smarandache neutrosophic complex modulo integer vector space over the S-ring $Z_{20}$.
i) Find a basis of V.
ii) Find number of elements of V.
iii) Write V as a direct sum of subspaces.
iv) Find $\text{Hom}_{Z_{20}}$ (V, V).
v) Find a projection operator on V.
vi) Does there exists a linear operator on V which keeps every subspace of V invariant?

146. Let $M = \left\{ \begin{bmatrix} a & b \\ c & d \end{bmatrix} \middle| a, b, c, d \in C(\langle Z_{10} \cup I \rangle) = \{m + ni_F + rI + i_F sI \mid m, n, r, s \in Z_{10}, i_F^2 = 9, I^2 = I \text{ and } (i_F I)^2 = 9I\}\right\}$ be a Smarandache complex neutrosophic modulo integer linear algebra over the S-ring $S = Z_{10}$.
i) Find a basis of M.
ii) Find dimension of M over S - ring.
iii) If M is treated only as a vector space will dimension of M over S be different? Justify your claim.
iv) Can M be written as a direct sum of linear subalgebras?
v) Find the algebraic structure enjoyed by $\text{Hom}_S$ (M, M).
vi) Write M as a pseudo direct sum of linear subalgebras.
vii) Is $P = \left\{ \begin{bmatrix} 0 & a \\ b & 0 \end{bmatrix} \middle| a, b \in C(\langle Z_{10} \cup I \rangle) \right\} \subseteq M$ be a linear subalgebra?

147. Let $V = \left\{ \begin{bmatrix} a & b & c & d \\ e & f & g & h \end{bmatrix} \middle| a, b, c, d, e, f, g, h \in C(Z_{46}) \right\}$ be the Smarandache complex modulo integer vector space over the S-ring $Z_{46}$.
i) Find a basis of V over $Z_{46}$.
ii) Is V finite dimensional over $Z_{46}$?
iii) Find $\text{Hom}_{Z_{46}}$ (V, V).
iv) If $Z_{46}$ is replaced by $C(Z_{46})$ will the dimension of V

over C ($Z_{46}$) vary?

148. Let $T = \left\{ \begin{bmatrix} a_1 & a_2 & a_3 \\ a_4 & a_5 & a_6 \\ a_7 & a_8 & a_9 \end{bmatrix} \middle| a_i \in C(Z_7); 1 \leq i \leq 9 \right\}$ be a complex modulo integer linear algebra over the field $F = Z_7$.
   i) Find a basis of T over $Z_7$.
   ii) Find L (T, $Z_7$).
   iii) Find $\text{Hom}_{Z_7}$ (T, T).
   iv) Find the number of elements in T.
   v) Write T as a pseudo direct sum of sublinear algebras.
   vi) Can T be written as a direct sum of sublinear algebras?

149. Enumerate some interesting properties enjoyed by complex modulo integer vector space.

150. What is the distinct features enjoyed by Smarandache complex modulo integer vector spaces defined over a S-ring ($Z_n$) and a complex modulo integer vector spaces defined over a field $Z_p$.

151. Let $V = \left\{ \begin{bmatrix} a_1 & a_2 & a_3 & a_4 \\ a_5 & a_6 & a_7 & a_8 \\ a_9 & a_{10} & a_{11} & a_{12} \\ a_{13} & a_{14} & a_{15} & a_{16} \end{bmatrix} \middle| a_i \in C(\langle Z_{23} \cup I \rangle) = \{a + bi_F + cI + i_F dI \mid a, b, c, d \in Z_{23}, i_F^2 = 22, I^2 = I \text{ and } (i_F I)^2 = 22I\}; 1 \leq i \leq 16 \right\}$ and W = {($a_1, a_2, \ldots, a_{12}$) | $a_i \in C(\langle Z_{23} \cup I \rangle)$; $1 \leq i \leq 12$} be two complex neutrosophic modulo integer linear algebra defined over the field $Z_{23}$.
   i) Find a basis of V and a basis of W.
   ii) What is the dimension of V?
   iii) Find $\text{Hom}_{Z_{23}}$ (V, W) = S.
   iv) Find $\text{Hom}_{Z_{23}}$ (W, V) = R.
   v) Is R ≅ S?

152. Obtain some interesting applications of complex modulo integer linear algebras.

153. Determine some nice applications of complex neutrosophic modulo integer vector space.

154. Find any property enjoyed by S - complex neutrosophic modulo integer vector spaces defined over a S-ring.

155. Let $V = \left\{ \begin{bmatrix} a_1 & \ldots & a_{20} \\ a_{21} & \ldots & a_{40} \end{bmatrix} \middle| a_i \in C(\langle Z_{29} \cup I \rangle) = \{a + bi_F + cI + i_F dI \mid a, b, c, d \in Z_{29}, i_F^2 = 28, I^2 = I \text{ and } (i_F I)^2 = 28I\}; 1 \leq i \leq 40 \right\}$ be a set complex neutrosophic linear algebra of modulo integer defined over the set $S = \{0, 1, I\}$.
   i) Find a basis of V over S.
   ii) Will the dimension of V change if V is replaced by $2Z_{29}$.
   iii) Find $\text{Hom}_S(V, V)$.

156. Let $P = \{(a_1, \ldots, a_{25}) \mid a_i \in C(\langle Z_{13} \cup I \rangle) = \{a + bi_F + cI + i_F dI \mid a, b, c, d \in Z_{13}, i_F^2 = 12, I^2 = I \text{ and } (i_F I)^2 = 12I\}; 1 \leq i \leq 25\}$ be a set complex neutrosophic modulo integer strong linear algebra over $S = \{0, I\}$.
   i) Find a basis of P over S.
   ii) Find dimension of P over S.
   iii) Write P as a direct sum of strong linear subalgebras.
   iv) If S is replaced by $T = \{0, 1\}$ will dimension of P over T different?

157. Find applications of complex neutrosophic real linear algebras.

158. Let S = $\left\{ \begin{bmatrix} a & b \\ c & d \end{bmatrix} \middle| a, b, c, d \in C(\langle R \cup I \rangle) = \{m + ni + sI + Iri \mid i^2 = -1, I^2 = I \text{ and } (iI)^2 = -I\} \right\}$ be a complex neutrosophic linear algebra over R.
   i) Find a basis of S over R.
   ii) Is S finite dimensional over R?
   iii) Find $\text{Hom}_R(S,S)$.
   iv) Find L (S, R).
   v) If R is replaced by $\langle R \cup I \rangle$ study results (1) to (v).
   vi) If R is replaced by C (R) = {a + bi | a, b ∈ R, $i^2 = -1$} = C study questions (1) to (iv).
   vii) If R is replaced by C ($\langle R \cup I \rangle$) What is dimension of S? Study (1) to (iv) questions.
   viii) Compare and distinguish between the spaces given in questions (v) (vi) and (vii).

159. Obtain some interesting results about complex neutrosophic semivecor spaces defined over Z.

160. Let P = $\left\{ \begin{bmatrix} a & b \\ c & d \end{bmatrix} \middle| a, b, c, d \in C(\langle Z \cup I \rangle) = \{m + ni + rI + isI \mid i^2 = -1, I^2 = I \text{ and } (iI)^2 = -I\}, m, n, r, s, \text{ are in } Z \right\}$. Is P a semifield? Justify your answer.

# FURTHER READINGS

# INDEX

# ABOUT THE AUTHORS

**Dr.W.B.Vasantha Kandasamy** is an Associate Professor in the Department of Mathematics, Indian Institute of Technology Madras, Chennai. In the past decade she has guided 13 Ph.D. scholars in the different fields of non-associative algebras, algebraic coding theory, transportation theory, fuzzy groups, and applications of fuzzy theory of the problems faced in chemical industries and cement industries. She has to her credit 646 research papers. She has guided over 68 M.Sc. and M.Tech. projects. She has worked in collaboration projects with the Indian Space Research Organization and with the Tamil Nadu State AIDS Control Society. She is presently working on a research project funded by the Board of Research in Nuclear Sciences, Government of India. This is her $58^{th}$ book.

On India's 60th Independence Day, Dr.Vasantha was conferred the Kalpana Chawla Award for Courage and Daring Enterprise by the State Government of Tamil Nadu in recognition of her sustained fight for social justice in the Indian Institute of Technology (IIT) Madras and for her contribution to mathematics. The award, instituted in the memory of Indian-American astronaut Kalpana Chawla who died aboard Space Shuttle Columbia, carried a cash prize of five lakh rupees (the highest prize-money for any Indian award) and a gold medal.
She can be contacted at vasanthakandasamy@gmail.com
Web Site: http://mat.iitm.ac.in/home/wbv/public_html/
or http://www.vasantha.in

---

**Dr. Florentin Smarandache** is a Professor of Mathematics at the University of New Mexico in USA. He published over 75 books and 200 articles and notes in mathematics, physics, philosophy, psychology, rebus, literature.

In mathematics his research is in number theory, non-Euclidean geometry, synthetic geometry, algebraic structures, statistics, neutrosophic logic and set (generalizations of fuzzy logic and set respectively), neutrosophic probability (generalization of classical and imprecise probability). Also, small contributions to nuclear and particle physics, information fusion, neutrosophy (a generalization of dialectics), law of sensations and stimuli, etc. He can be contacted at smarand@unm.edu